\documentclass[12pt]{article}

\usepackage{amsmath,amstext,amsfonts,amsthm,amssymb,eucal,bm,diagrams}

\newcommand{\CFT}{\mbox{CohFT}}
\newcommand{\CFTs}{\mbox{CohFTs}}
\newcommand{\dL}{\mathbb{L}}
\newcommand{\compl}{\widehat{\phantom{x}}}
\newcommand{\petlya}{\operatorname{loop}}

\newcommand{\maf}{\operatorname{mf}}

\newcommand{\nat}{\operatorname{nat}}
\newcommand{\be}{\begin{equation}}
\newcommand{\ee}{\end{equation}}

\newcommand{\bP}{{\mathbf P}}

\newcommand{\vareps}{\varepsilon}

\newcommand{\Id}{\operatorname{Id}}
\newcommand{\del}{\partial}

\newcommand{\ch}{\operatorname{ch}}

\renewcommand{\mod}{\operatorname{mod}}
\newcommand{\Com}{\operatorname{Com}}
\newcommand{\Cone}{\operatorname{Cone}}
\newcommand{\Tr}{\operatorname{Tr}}

\newcommand{\OO}{{\cal O}}
\newcommand{\Sym}{\operatorname{Sym}}

\newcommand{\pr}{\operatorname{pr}}

\newcommand{\coker}{\operatorname{coker}}
\newcommand{\QMF}{\operatorname{QMF}}
\newcommand{\Si}{\Sigma}

\newcommand{\DD}{{\cal D}}

\newcommand{\KK}{{\cal K}}
\newcommand{\JJ}{{\cal J}}

\newcommand{\tot}{\operatorname{tot}}

\newcommand{\G}{{\mathbb G}}

\newcommand{\mg}{{\frak m}}
\newcommand{\hra}{\hookrightarrow}

\newcommand{\lan}{\langle}
\newcommand{\ran}{\rangle}
\newcommand{\Coh}{\operatorname{Coh}}
\newcommand{\GG}{{\cal G}}
\newcommand{\CC}{{\cal C}}
\newcommand{\YY}{{\cal Y}}

\newcommand{\Spec}{\operatorname{Spec}}

\renewcommand{\P}{{\mathbb P}}

\newcommand{\si}{\sigma}

\newcommand{\Pic}{\operatorname{Pic}}
\newcommand{\Ind}{\operatorname{Ind}}
\newcommand{\ga}{\gamma}
\newcommand{\de}{\delta}
\newcommand{\eps}{\epsilon}

\renewcommand{\ker}{\operatorname{ker}}
\newcommand{\im}{\operatorname{im}}

\numberwithin{equation}{section}

\newtheorem{theor}{Theorem}[subsection]
\newtheorem{thm}[theor]{Theorem}
\newtheorem{lem}[theor]{Lemma}

\newtheorem{prop}[theor]{Proposition}

\newtheorem{cor}[theor]{Corollary}
{  \theoremstyle{definition}
           \newtheorem{defi}[theor]{Definition}
           \newtheorem{rem}[theor]{Remark}
           \newtheorem{rems}[theor]{Remarks}
          \newtheorem{ex}[theor]{Example}
          
}

\newcommand{\Irr}{\operatorname{Irr}}
\newcommand{\Pf}{\noindent {\it Proof}}
\newcommand{\id}{\operatorname{id}}

\newcommand{\DQMF}{\operatorname{DQMF}}
\newcommand{\DCMF}{\operatorname{DCMF}}
\newcommand{\ov}{\overline}
\newcommand{\we}{\wedge}

\newcommand{\rk}{\operatorname{rk}}
\newcommand{\ra}{\rightarrow}
\renewcommand{\AA}{{\cal A}}

\newcommand{\FF}{{\cal F}}

\newcommand{\MM}{{\cal M}}
\newcommand{\TT}{{\cal T}}
\newcommand{\XX}{{\cal X}}
\newcommand{\HH}{{\cal H}}

\newcommand{\VV}{{\cal V}}
\renewcommand{\SS}{{\cal S}}

\newcommand{\LL}{{\cal L}}
\newcommand{\WW}{{\cal W}}
\newcommand{\st}{\operatorname{st}}

\newcommand{\Om}{\Omega}

\newcommand{\Hom}{\operatorname{Hom}}

\newcommand{\Ext}{\operatorname{Ext}}
\newcommand{\End}{\operatorname{End}}
\newcommand{\Res}{\operatorname{Res}}

\newcommand{\Barc}{\operatorname{Bar}}
\newcommand{\PPic}{{\cal P}ic}

\renewcommand{\a}{\alpha}
\renewcommand{\b}{\beta}
\newcommand{\bfe}{{\bf e}}
\newcommand{\bfm}{{\bf m}}
\newcommand{\bfd}{{\bf d}}
\newcommand{\om}{\omega}
\newcommand{\De}{\Delta}
\newcommand{\la}{\lambda}
\renewcommand{\th}{\theta}
\newcommand{\C}{{\mathbb C}}

\newcommand{\A}{{\mathbb A}}

\newcommand{\Z}{{\mathbb Z}}
\newcommand{\Q}{{\mathbb Q}}
\newcommand{\La}{\Lambda}
\newcommand{\Ga}{\Gamma}
\newcommand{\Td}{\operatorname{Td}}

\newcommand{\wt}{\widetilde}
\newcommand{\ot}{\otimes}

\newcommand{\sub}{\subset}
\newcommand{\com}{\operatorname{com}}
\newcommand{\ed}{\qed\vspace{3mm}}

\newcommand{\Qcoh}{\operatorname{Qcoh}}

\newcommand{\tr}{\operatorname{tr}}

\newcommand{\pa}{\partial}

\newcommand{\unit}{{\mathbf 1}}
\newcommand{\Sg}{\operatorname{Sg}}
\newcommand{\tree}{\text{tree}}
\newcommand{\ev}{\operatorname{ev}}
\newcommand{\rig}{\operatorname{rig}}
\newcommand{\MF}{{\operatorname{MF}}}
\newcommand{\HMF}{\operatorname{HMF}}
\newcommand{\DMF}{\operatorname{DMF}}
\newcommand{\mf}{matrix factorization}
\newcommand{\w}{{\boldsymbol w}} 

\newcommand{\Homb}{\mathcal{H}om} 
\newcommand{\Per}{\operatorname{Per}}
\newcommand{\fC}{{\mathfrak C}}
\newcommand{\LHZ}{\operatorname{LHZ}}
\newcommand{\bdeg}{{\bf deg}}
\newcommand{\bth}{{\bm \theta}}
\newcommand{\bq}{{\bf q}}

\newcommand{\Bun}{\operatorname{Bun}}
\newcommand{\str}{\operatorname{str}}

\title{Matrix factorizations and Cohomological Field Theories}
\author{Alexander Polishchuk \and Arkady Vaintrob}
\date{}

\begin{document}
\maketitle
\begin{abstract}
We give a purely algebraic construction of a cohomological field theory associated
with a quasihomogeneous isolated hypersurface singularity $W$ and a subgroup $G$ of the
diagonal group of symmetries of $W$. This theory can be viewed as an analogue of
the Gromov-Witten theory for an orbifoldized Landau-Ginzburg model for $W/G$.
The main geometric ingredient for our construction is provided by the moduli of
curves with $W$-structures introduced by Fan, Jarvis and Ruan.
We construct certain matrix factorizations on the products of these moduli stacks with affine spaces
which play a role similar to that of the virtual fundamental classes in the Gromov-Witten theory.
These matrix factorizations are used to produce functors from the categories of
equivariant matrix factorizations to the derived categories of coherent sheaves on the
Deligne-Mumford moduli stacks of stable curves. The structure maps of our cohomological
field theory are then obtained by passing to the induced maps on Hochschild homology.
We prove that for simple singularities a specialization of our theory gives
the cohomological field theory constructed by Fan, Jarvis and Ruan using analytic tools.
\end{abstract}

\tableofcontents

\section*{Introduction}

The notion of a cohomological field theory was introduced in \cite{KonM} (see also \cite{Manin}) 
to axiomatize properties of quantum cohomology and Gromov-Witten invariants
and to provide a basis for formulating some mathematical aspects
of mirror symmetry. 
Recall that a (complete) cohomological field theory (\CFT) is an 
algebraic structure on a vector space $\HH$ (called the state space of the theory)
with a collection of operations indexed by homology classes of 
the Deligne-Mumford moduli spaces $\ov{\MM}_{g,n}$ of stable curves with marked points.
The main example of a \CFT\ is provided by the Gromov-Witten theory 
associated with a smooth projective variety $X$ (or more generally, with a compact symplectic
manifold). In this context mirror symmetry can be viewed as 
an isomorphism of \CFTs\ originating from two different geometric inputs (the so called A- and
B-models, where the GW-theory corresponds to the A-model and the B-model is related
to deformations of complex structures of $X$).
The part of the \CFT\ data corresponding to genus
zero curves can be described as a formal Frobenius manifold structure, which leads in the case
of GW-theory to the notion of quantum cohomology, a certain deformation of the cohomology ring of
$X$.
Through Frobenius manifolds \CFTs\  
are related to integrable hierarchies 
of systems of partial differential equations (see \cite{DZ}).

The first example of \CFTs\ besides GW-theory was provided by the theory of $r$-spin
curves constructed in \cite{Witten-Nm}, \cite{Witten-2dimgr},
\cite{JKV}, \cite{PV} and \cite{P-class} (see also \cite{Moch}).
The corresponding Frobenius manifolds are isomorphic to the ones constructed by Saito
for simple singularities of type $A_{r-1}$, and the corresponding integrable hierarchies are
the Gelfand-Dickey hierarchies (see \cite{JKV}, \cite{FSZ}).

Starting from the work of Givental \cite{Giv-highergenus} 
it was realized (see \cite{FSZ}, \cite{Teleman}) that 
an arbitrary generically semisimple Frobenius manifold extends to a unique \CFT\ .
This construction 
can be applied to Saito's Frobenius manifold of any quasihomogeneous isolated singularity $\w$
to give a \CFT\ for such a singularity, which corresponds to the B-side of the Landau-Ginzburg
model related to $\w$ (see \cite{LVW}). 
A \CFT\ corresponding to the A-side for these LG-models was recently constructed by
Fan, Jarvis and Ruan in \cite{FJR}, \cite{FJR2}. We will refer to this \CFT\ as 
{\it FJR-theory}.\footnote{Sometimes this theory is referred to as FJRW-theory since Fan, Jarvis and Ruan attribute some key ideas of their model to Witten.}
Their construction is based on the study of a certain PDE 
over coverings of moduli spaces of stable curves (that generalize the moduli
spaces of $r$-spin curves). More precisely, they construct a virtual fundamental cycle
on the moduli space of solutions of this PDE corresponding to a linear perturbation of the potential $\w$.
The dependence of this virtual class on the perturbation is governed by the state space of the theory
which is given by the orbifoldized Milnor ring of $(\w, G)$ (see \cite{Wall1}, \cite{Kauf}), where 
$G$ is a finite group of diagonal symmetries of $\w$.
One of the main results of \cite{FJR} is that 
in the case of ADE simple singularities this \CFT\ is isomorphic to
Givental's \CFT\ associated with Saito's Frobenius manifold of the dual singularity 
(which is the same singularity for the series $E$ and for series $D$ with non-maximal
symmetry group). Using the work of Frenkel-Givental-Milanov (see \cite{GM}, \cite{FGM})
it is shown in \cite{FJR} that the associated total potential function is a $\tau$-function 
of the corresponding Kac-Wakimoto hierarchy.
Another recent result on the mirror symmetry for the Landau-Ginzburg models is due to
Krawitz \cite{Krawitz} who established the isomorphism between Frobenius algebras associated
with A- and B-side of the dual invertible quasi-homogeneous potentials. 

In this paper we present a purely algebraic (and perhaps more general) version of the FJR-theory. 
The main role in
our construction is played by \mf s. These are generalizations of complexes obtained by
replacing the condition $d^2=0$ with $d^2=\w$
(see \cite{Eisen}, \cite{Bu}).
Matrix factorizations appeared in physics in connection with open-closed topological string theories
(see \cite{KL1}, \cite{KL2}, \cite{KR}).
They are also related to an important invariant of the potential $\w$, the {\it singularity category}
(see \cite{Orlov}, \cite{Posic-mf}, \cite{LP}). The category of \mf s for an isolated singularity fits
naturally into the framework of noncommutative geometry developed from the point of view
of dg-categories or $A_{\infty}$-algebras (see \cite{Dyck}, \cite{KKP}).

Our main construction gives a \CFT\ whose state space is 
built from the Hochschild homology of the dg-categories
of equivariant \mf s associated with a quasihomogeneous isolated singularity $\w$ 
and a finite group of symmetries $G$. 
This \CFT\ carries a priori more information than the FJR-theory: it is a \CFT\ {\it with coefficients }
in the representation ring $R$ of $G$. Conjecturally, the reduced theory
obtained by the specialization $R\to\C$ is equivalent to the FJR-theory. We show that this is
true for all simple singularities.

Let $\w(x_1,\ldots,x_n)$ be a quasihomogeneous polynomial
with an isolated singularity at the origin. We fix the degrees of quasihomogeneity
$d_j=\deg(x_j)$ and we set $q_j=d_j/d$, where $d=\deg(\w)$.
We denote the Milnor ring of $\w$ by $\AA_\w=\C[x_1,\ldots,x_n]/(\partial_1\w,\ldots,\pa_n\w)$
and set
$$H(\w)=\AA_\w\ot (dx_1\we\ldots\we dx_n).$$
This space is canonically isomorphic to the Hochschild homology of the category of \mf s of $\w$.
Let $G_\w\sub(\C^*)^n$ be the group of all diagonal symmetries of $\w$.
For each $\ga\in G_\w$ consider the subspace of invariants $(\A^n)^{\ga}\sub\A^n$ and
set $\w_{\ga}=\w|_{(\A^n)^{\ga}}$. Then $\w_{\ga}$ still has an isolated singularity at zero.
Let $G\sub G_\w$ be a finite subgroup containing 
the exponential grading element
$$J=(\exp(2\pi i q_1),\ldots,\exp(2\pi i q_n))\in G_\w,$$
and let $R=\C[\widehat{G}]$ be the representation ring of $G$.
Our construction associates with such $G$ a \CFT\ with coefficients in $R$
on the state space
$$\HH(\w,G)=\bigoplus_{\ga,\ga'\in G}H(\w_{\ga,\ga'})^{G},$$
where $\w_{\ga,\ga'}$ is the restriction of $\w$ to the subspace of invariants
$(\A^n)^{\{\ga,\ga'\}}\sub\A^n$. We view this space as an $R$-module via the
$G$-grading given by $\ga'$.
The component corresponding to $\ga'=1$ admits a certain twist by a Todd class, and conjecturally
this reduced \CFT\ is isomorphic to the FJR-theory.
In \cite{ChiR} the FJR-theory for the quintic threefold is related to a certain specialization of
the corresponding Gromov-Witten theory. One can speculate that our \CFTs\ associated
with other elements $\ga'\in G$ are also related to some specializations of the Gromov-Witten
invariants in the Calabi-Yau case.

As in \cite{FJR},
the main geometric ingredient of our construction is the moduli stacks of the so-called $\w$-curves
which are orbicurves $\CC$ with marked (orbi-)points together with a collection of line bundles
$\LL_1,\ldots,\LL_n$ satisfying certain constraints coming from the monomials of $\w$.
However, we observe that these coverings should really be viewed as an attribute
of the group of symmetries $G$ rather than of the potential $\w$. More precisely, we reformulate
the notion of a $\w$-structure using principal $\Ga$-bundles, where $\Ga$ is an extension of
$\G_m$ by $G$. This leads to the moduli stacks of $\Ga$-spin curves that replace $\w$-curves
considered in \cite{FJR}. This technical device allows us to remove the assumption that $G_\w$ is
finite imposed in \cite{FJR} and to show that every finite subgroup of $G_\w$ containing $J$
is admissible in the sense of \cite{FJR}.

Our results give also a categorification of these \CFTs\ in a certain weak sense.
Namely, we construct a collection of functors 
inducing the \CFT\ maps after passing to Hochschild homology
(up to rescaling). These functors are given by kernels which are certain \mf s, 
called {\it fundamental \mf s}, on the
product of the moduli spaces of $\Ga$-spin curves with affine spaces.
In some sense these fundamental \mf s play a role similar to that of the virtual fundamental
class in the GW-theory. 
The factorization axiom holds on the categorified level after passing to appropriate
finite covers of the relevant moduli spaces. It seems plausible that there should also be a version
of quantum K-theory in our setup (see \cite{Lee}).

Note that representations of functors between categories of \mf s by kernels
are discussed extensively in the work of Ballard, Favero and Katzarkov \cite{BFK} (in
the context of graded algebras), where the authors give interesting applications to Homological
Mirror Symmetry and Hodge theory.

One of the obstacles that prevents us from comparing the specialization of our \CFT\
at $\ga'=1$ with the FJR-theory is
that we were able to prove the analog of the Dimension Axiom of \cite[Sec.\ 4]{FJR} only in
some special cases. In particular, we verified it for
all simple singularities, which together with the computation of the 
corresponding Frobenius algebras enables us to show that in this case
our reduced \CFT\ is isomorphic to that of Fan-Jarvis-Ruan (see Section \ref{simple-sing-sec}). 

The paper is organized as follows. In Section \ref{mf-sec} we review basics of
the theory of \mf s. In particular, we discuss the relation between the derived category of \mf s and
the singularity category, prove auxiliary results about push-forwards of \mf s and establish some properties of Koszul \mf s, crucial for our main construction. 
In Section \ref{quasihom-mf-sec} we specialize
to the case of a quasihomogeneous isolated singularity. We consider the dg-category of
equivariant \mf s and compute its Hochschild homology space with the canonical
pairing. We also discuss functors between categories of \mf s given by kernels.
In Section \ref{simple-constr-sec} we recast the notion of a $\w$-structure from \cite{FJR} in
terms of torsors over some commutative algebraic groups and consider the corresponding moduli
stacks of $\Ga$-spin curves. Section \ref{main-sec} is the technical core of the paper.
Here we construct the fundamental \mf\ over the product of the moduli space of $\Ga$-spin curves
with an affine space. 
The construction shares some of the features with 
the construction of the Witten's top Chern class in \cite{PV}, however, it uses also some new ingredients, notably the push-forward of \mf s. Also, the proof of independence 
of the fundamental \mf\ of choices of resolutions is based on a different idea (we use properties
of regular Koszul \mf s). In Section \ref{axioms-sec} we define the \CFT s associated with
a pair $(\w,G)$ using the functor associated with the fundamental \mf\ and passing to the induced
map on Hochschild homology.
We also prove for our theory analogs of all properties established in \cite[Sec.\ 4]{FJR}
for the FJR-theory except the Dimension Axiom, which we prove in some particular cases (see \ref{dim-sec}). 
In Section \ref{calc-sec} we give a recipe for calculating genus zero three-point correlators
for our theory which are responsible for the Frobenius algebra structure on the state space. 
In Section \ref{simple-sing-sec} we compute all such correlators in the case
when $\w$ is a simple singularity of type $A$, $D$, $E_6$, $E_7$ or $E_8$. We also prove
in Section \ref{comparison-sec} that our reduced \CFT\ is isomorphic to the FJR-theory in this case.
In Appendix we discuss the constructions of functoriality for Hochschild homology. In particular,
we prove that the construction used in \cite{PV-mf} is compatible with the standard one
(see Theorem \ref{hoch-funct-thm}).

\bigskip

\noindent
{\it Notation and conventions}.
We work with schemes and stacks over $\C$ and all our (dg-) categories are $\C$-linear. 
For a finite abelian group $G$ we denote by $\widehat{G}$ the dual abelian group. For a
commutative algebraic group $\Ga$ we denote by $X(\Ga)$ the group of characters of $\Ga$.
We say that a triangulated category $\DD$ 
is generated by a set of objects $(E_i)$ if the minimal full triangulated subcategory containing
$E_i$ and closed under taking direct summands is the entire $\DD$.
For an additive category $\CC$ we denote by $\Com(\CC)$ the category of complexes
over $\CC$.

We always assume that our algebraic stacks are Noetherian and semi-separated.
For such an algebraic stack $X$ we denote by $\Coh(X)$ (resp., $\Qcoh(X)$; resp., $D^b(X)$) the category of coherent sheaves (resp., quasicoherent sheaves; resp., bounded derived
category of coherent sheaves) on $X$. By \cite[Cor. 2.11]{AB}, $D^b(X)$ 
is equivalent to the full subcategory
of the bounded derived category of $\OO_X$-modules consisting of complexes with coherent
cohomology.
By a vector bundle on a stack we mean a locally free sheaf of $\OO$-modules of finite rank.

\vspace{2mm}

\noindent
{\it Acknowledgments}.
Both authors would like to thank the IHES, where part of this work was done,
for hospitality and stimulating atmosphere. 
At various stages of this project we benefited from conversations with Alessandro Chiodo,
Alexander Givental, Bernhard Keller, Maxim Kontsevich, Dmitry Orlov, Tony Pantev
and Leonid Positselski. 
We are also grateful to Matthew Ballard, David Favero and Ludmil Katzarkov for sharing
with us a preliminary version of their work \cite{BFK}.
The first author was partially supported by the NSF grant DMS-1001364.

\section{Matrix factorizations on stacks}\label{mf-sec}

In this section we review the theory of \mf s on stacks from \cite{PV-stacks}. We also
establish some technical results on push-forwards of \mf s and on Koszul \mf s.

\subsection{Categories of \mf s}\label{mf-cat-sec}

Let us recall some basic definitions from \cite{PV-stacks}.

\begin{defi}
Let $X$ be an algebraic stack, 
$L$ a line bundle on $X$, and $W\in H^0(X,L)$ a section (called a {\it potential}). 
A {\it \mf\ } $\bar{E}=(E_\bullet,\de_\bullet)$ of $W$ on $X$ consists of a pair of vector bundles 
(i.e., locally free sheaves of finite rank) $E_0$, $E_1$
on $X$ together with homomorphisms 
$$\de_1:E_1\to E_0 \ \text{ and } \ \de_0:E_0\to E_1\ot L,$$ 
such that
$\de_0\de_1=W\cdot \id$ and $\de_1\de_0=W\cdot\id$.

Sometimes we will assume that the potential $W$ is not a zero divisor, i.e., the morphism
$W:\OO_X\to L$ is injective. 
\end{defi}

In the case $W=0$ we have $\de_0\de_1=\de_1\de_0=0$, so we can define cohomology
of a \mf\ $\bar{E}$ by
\begin{equation}\label{mf-coh-eq}
H^0(\bar{E})=\ker(\de_0)/\de_1(E_1), \ \ H^1(\bar{E})=\ker(\de_1)/\de_0(E_0\ot L^{-1}).
\end{equation}

\begin{defi} We define the dg-category $\MF(X,W)$ of  \mf s of $W$ as follows.
For a pair of \mf s $\bar{E}$ and $\bar{F}$ we define a $\Z$-graded complex of morphisms
$\Homb_\MF(\bar{E},\bar{F})$ by setting
$$\Homb_\MF(\bar{E},\bar{F})^{2n}=\Hom(E_0,F_0\ot L^n)\oplus\Hom(E_1,F_1\ot L^n),$$
$$\Homb_\MF(\bar{E},\bar{F})^{2n+1}=\Hom(E_0,F_1\ot L^{n+1})\oplus\Hom(E_1,F_0\ot L^n).$$
The differential on $\Homb_\MF(\bar{E},\bar{F})$, raising the degree by $1$, is given by
\begin{equation}\label{diff-mf-eq}
d f = \de_F\circ f - (-1)^{|f|} f \circ \de_E~.
\end{equation}
\end{defi}

We denote by $\HMF(X,W)=H^0\MF(X,W)$ the corresponding homotopy category. We will usually
omit $X$ from the notation.
As in the standard case considered in \cite{Orlov} the category $\HMF(W)$ has a triangulated structure
(see \cite[Def.\ 1.3]{PV-stacks}).

One can conveniently write the complexes $\Homb_\MF(\bar{E},\bar{F})$ using
the half-twist notation (see \cite[Def.\ 1.2]{PV-stacks}):
$$\Homb_\MF(\bar{E},\bar{F})^i=\Hom_{i \mod 2}(E(L^{1/2}), F(L^{1/2})\ot L^{i/2}),$$
where $E(L^{1/2})=E_0\oplus (E_1\ot L^{1/2})$ and $\Hom_{i \mod 2}$ denotes morphisms
of $\Z/2$-graded bundles, homogeneous of degree $i \mod 2$.

We also consider the dg-category $\MF^\infty(X,W)$ and the corresponding homotopy category
$\HMF^\infty(X,W)$  of {\it quasi-\mf s} defined using locally free sheaves
of not necessarily finite rank (see \cite[Def.\ 1.4]{PV-stacks}). 
An even larger dg-category $\QMF(X,W)$ of {\it quasicoherent \mf s}
is obtained if we allow $E_0$ and $E_1$ to be arbitrary quasicoherent sheaves.
This category 
is featured prominently in more recent treatments of \mf s 
in the non-affine case (see \cite{Posic-mf} and \cite{LP}).
We will mostly use it in the case when $W=0$.
Note that in this case we have a natural class of quasi-isomorphisms in the corresponding homotopy
category (defined using the cohomology \eqref{mf-coh-eq}). 
Localizing with respect to quasi-isomorphisms we get the derived category
$\DQMF(X,0)$. We denote by $\DQMF_c(X,0)\sub\DQMF(X,0)$ 
the full subcategory of quasicoherent \mf s of $0$ with coherent cohomology.
Similarly, replacing quasicoherent sheaves with coherent sheaves one can talk about {\it coherent
\mf s} of $0$ and define the derived category of coherent \mf s $\DCMF(X,0)$.

Let us introduce some natural operations on \mf s.

For a bounded complex of vector bundles on $X$, $(C^{\bullet},\de_C)$, and 
a quasicoherent \mf\  $\bar{E}=(E_\bullet,\de_\bullet)$ of $W\in H^0(X,L)$ we define the \mf\ 
$C^{\bullet}\ot\bar{E}$ of $W$ by setting
\begin{equation}\label{complex-tensor-eq}
(C^{\bullet}\ot E)(L^{1/2})=C(L^{-1/2})\ot E(L^{1/2}),
\end{equation}
with the differential $\de_C\ot\id+\id\ot\de$,
where $C(L^{-1/2})=\oplus_{n\in\Z} C^n\ot L^{-n/2}$ with the $\Z/2$-grading induced by the
$\Z$-grading and the induced differential $\de_C:C(L^{-1/2})\to C(L^{-1/2})\ot L^{1/2}$.
Explicitly,
$$(C^{\bullet}\ot E)_i=\oplus_{n\in\Z}C^n\ot E(L^{1/2})_{n+i}\ot L^{-(n+i)/2}, \ \text{ for }
i=0,1.$$
Note that we have natural isomorphisms
$$(C^{\bullet}[1])\ot\bar{E}\simeq C^{\bullet}\ot (\bar{E}[1])\simeq (C^{\bullet}\ot\bar{E})[1].$$

Let $X_0\sub X$ be the zero locus of $W$.
With
a quasicoherent \mf\ $\bar{E}=(E_\bullet,\de_\bullet)$ we associate a $\Z$-graded complex
 of vector bundles on $X_0$
\begin{equation}\label{mf-Z-gr-com-eq}
\com(\bar{E}): \ \ \ldots\to (E_0\ot L^{-1})|_{X_0} \rTo{\ov{\de}_0} E_1|_{X_0}\rTo{\ov{\de}_1} E_0|_{X_0}\rTo{\ov{\de}_0}
(E_1\ot L)|_{X_0}\to\ldots
\end{equation}
where $\ov{\de}_i$ is induced by $\de_i$, and $E_0|_{X_0}$ is placed in degree $0$ (for quasi-\mf s
this construction was considered in \cite[Sec.\ 1]{PV-stacks}).
By \cite[Lem.\ 1.5]{PV-stacks}, this complex is exact provided $W$ is not a zero divisor and
$\bar{E}$ is a {\it quasi-\mf}.
This construction extends to a dg-functor
$$\com:\QMF(W)\to\Com(\Qcoh(X_0))$$
that induces an exact functor between the corresponding homotopy categories.
It is easy to see that for a bounded complex $C^\bullet$ of vector bundles on $X$ one has
a natural isomorphism of complexes on $X_0$
\begin{equation}\label{tensor-com-isom}
\com(C^\bullet\ot \bar{E})\simeq C^\bullet|_{X_0}\ot_{\OO_{X_0}} \com(\bar{E}).
\end{equation}
In the case $W=0$ the complex $\com(\bar{E})$ satisfies 
$$H^{2n}\com(\bar{E})=H^0(\bar{E})\ot L^n, \ \ H^{2n-1}\com(\bar{E})=H^1(\bar{E})\ot L^n.$$
In particular, 
a closed morphism $q$ of quasicoherent \mf s of $0$ is a quasi-isomorphism if and only if
$\com(q)$ is a quasi-isomorphism.

\begin{defi}
For a pair of potentials $W,W'\in H^0(X,L)$
the tensor product dg-functor 
$$\MF(W)\otimes \MF(W')\to \MF(W+W')$$
is defined as follows.
For $\bar{E}=(E,\de_E)$ and $\bar{F}=(F,\de_F)$ we set
\begin{equation}\label{usual-tensor-product}
(\bar{E}\ot\bar{F})_0=E_0\ot F_0\oplus E_1\ot F_1\ot L \ \text{ and }E_0\ot F_1\oplus E_1\ot F_0
\end{equation}
with the differential $\de_E\ot\id_F+\id_E\ot \de_F$.
Note that
$$(E\ot F)(L^{1/2})=E(L^{1/2})\ot F(L^{1/2}).$$
\end{defi}

\begin{defi}
For a bounded complex $(C^\bullet,\de_C)$ of vector bundles and a line bundle $L$
on $X$ let us define $(\maf(C^\bullet),\de)$, a \mf\ of $0\in H^0(X,L)$, by setting
\begin{equation}\label{maf-eq}
\maf(C^\bullet)_0=\bigoplus_n C^{2n}\ot L^{-n}, \ \ 
\maf(C^\bullet)_1=\bigoplus_n C^{2n-1}\ot L^{-n}
\end{equation}
with the differential $\de$ induced by $\de_C$.
\end{defi}

A straightforward check shows that the tensor product operations \eqref{complex-tensor-eq} and
\eqref{usual-tensor-product} are consistent.

\begin{lem}\label{complex-tensor-lem} 
One has a natural isomorphism 
$$C^\bullet\ot\bar{E}\simeq\maf(C^\bullet)\ot\bar{E}$$
in $\MF(W)$,
on the left we use the operation \eqref{complex-tensor-eq}.
Hence, by \eqref{tensor-com-isom} we have
$$\com(\maf(C^\bullet)\ot\bar{E})\simeq C^\bullet|_{X_0}\ot\com(\bar{E}).$$
\end{lem}
\ed

We have the duality dg-functor
\begin{equation}\label{duality-eq}
\MF(W)^{op}\to \MF(-W)
\end{equation}
sending $\bar{E}=(E,\de_E)$ to $\bar{E}^*=(E(L^{1/2})^{\vee}(L^{-1/2}),\de^*)$. 
In other words, the even part of $\bar{E}^*$ is $E_0^{\vee}$
and its odd part is $E_1^{\vee}\ot L^{-1}$.

\begin{lem}\label{Hom-mf-lem}
For a pair of \mf s $\bar{E}$ and $\bar{F}$ in $\MF(X,W)$ we have an isomorphism of complexes
$$\Homb_\MF(\bar{E},\bar{F})\simeq H^0(X,\com(\bar{E}^*\ot\bar{F})),$$
where $\bar{E}^*\ot\bar{F}\in\MF(X,0)$.
\end{lem}
\ed

For a morphism of stacks $f:X'\to X$, a line bundle $L$ over $X$
and a section $W\in H^0(X,L)$ we have natural pull-back functors on \mf s: 
a dg-functor
$$f^*:\MF(X,W)\to\MF(X',f^*W),$$
where $f^*W$ is the induced section of $f^*L$ on $X'$, and the induced exact functor
$$f^*:\HMF(X,W)\to\HMF(X,f^*W).$$

\begin{defi}\label{ext-tens-def}
The {\it external product} of $(X,W)$ and $(X',W')$ (where $W\in H^0(X,L)$
and $W'\in H^0(X',L'))$ is defined as a pair $(U,W\wt{\oplus}W')$, where
$U\to X\times X'$ is the $\G_m$-torsor 
associated with the line bundle $L'\ot L^{-1}$.
Let $p_1:U\to X$ and $p_2:U\to X'$ be the natural
projections. Then we have an isomorphism $L_U=p_1^*L\simeq p_2^*L'$, and
we define
$$W\wt{\oplus}W'=p_1^*W+p_2^*W'\in H^0(U,L_U).$$
\end{defi}

Combining the pull-back and tensor product functors defined above we obtain the 
{\it external tensor product} functor
\begin{equation}\label{ext-tensor-eq}
\MF(X,W)\otimes \MF(X',W')\to MF(U,p_1^*W+p_2^*W).
\end{equation}

We can define Koszul \mf s $\{\a,\b\}$ in our setting (see \cite{KhR}, \cite{PV-mf}).

\begin{defi} Let $X$, $L$ and $W$ be as above, and let
$V$ be a vector bundle $V$ on $X$.
For global sections 
$$\a\in H^0(X,V\ot L),\ \b\in H^0(X,V^{\vee}) \ \text{ such that }\lan\a,\b\ran=W$$ 
we define define the {\it Koszul matrix factorization} 
$\{\a,\b\}$ of $W$ by
\be\label{Koszul-mf-def}
\{\a,\b\}=\left({\bigwedge}^\bullet(V\ot L^{1/2})(L^{-1/2}), \de_{\a,\b}\right),
\end{equation}
with the $\Z/2$-grading on ${\bigwedge}^\bullet(V\ot L^{1/2})$ induced by the 
$\Z$-grading. The differential is given by
$$\de_{\a,\b}=\a\we?+\iota(\b),$$
where $\iota(\b)$ is the contraction by $\b$.
\end{defi}

Explicitly, 
$$\{\a,\b\}_0=\OO_X\oplus (\we^2 V\ot L)\oplus (\we^4V\ot L^2)\oplus\ldots,$$
$$\{\a,\b\}_1=V\oplus (\we^3 V\ot L)\oplus (\we^5V\ot L^2)\oplus\ldots.$$

\begin{lem}\label{Koszul-complex-lem}
For $\a\in H^0(X,V\ot L)$ let 
\begin{equation}\label{usual-Koszul-complex-eq}
K^\bullet(\a)=({\bigwedge}^\bullet(V\ot L),\a\we?)
\end{equation}
be the Koszul complex. Then one has an isomorphism of \mf s in $\MF(0)$
$$\{\a,0\}\simeq\maf(K^\bullet(\a)).$$
\end{lem}

The proof is straightforward.

\subsection{Equivariant \mf s}

Let $X$ be a stack and $\Ga$ an affine algebraic
group acting on $X$. Let also $W$ be a regular function on $X$, semi-invariant with respect to $\Ga$.
Thus, we have a character $\chi:\Ga\to\G_m$ such that 
$$W(\ga\cdot x)=\chi(\ga)W(x)$$
for $\ga\in \Ga$, $x\in X$. 
Recall that {\it $\Ga$-equivariant \mf s of $W$ with respect to the character } $\chi$
are defined as pairs of $\Ga$-equivariant vector bundles $E_0, E_1$ on $X$ together with
$\Ga$-invariant homomorphisms
$$\de_1:E_1\to E_0\ \text{ and } \de_0:E_0\to E_1\ot\chi,$$
such that $\de_0\de_1=W\cdot \id$ and
$\de_1\de_0=W\cdot \id$. 

The dg-category $\MF_{\Ga,\chi}(X,W)$ of $\Ga$-equivariant \mf s is naturally equivalent to
the category $\MF(X/\Ga,\ov{W})$, 
where $\ov{W}$ is the section induced by $W$ of the line bundle over $X/\Ga$ associated with $\chi$
(see \cite[Prop.\ 2.2]{PV-stacks}).

\begin{ex}
Let $V$ be a $\Ga$-equivariant vector bundle on $X$.
For $\Ga$-invariant sections $\a\in H^0(X,V\ot\chi)^G$ and
$\b\in H^0(U,V^{\vee})^\Ga$ 
the Koszul \mf\ $\{\a,\b\}$ (see \eqref{Koszul-mf-def}) is $\Ga$-equivariant.
\end{ex}

We will often consider the following special situation.
Let $\Ga$ be a commutative algebraic group with a
surjective homomorphism $\chi:\Ga\to \G_m$ such that $G:=\ker(\chi)$ is finite, and let $X$ be
a stack with the trivial action of $\Ga$.

Let $\chi_1,\ldots,\chi_d\in \widehat{\Ga}$ be a system of representatives 
for the cosets of the subgroup $\lan\chi\ran\sub\widehat{\Ga}$ of characters of $\Ga$.
A $\Ga$-equivariant quasicoherent \mf\ $\bar{E}$ consists of a pair of $\widehat{\Ga}$-graded
quasicoherent sheaves
$$E_0=\oplus_{\xi\in \widehat{\Ga}} E_{0,\xi}, \ \ E_1=\oplus_{\xi\in \widehat{\Ga}} E_{1,\xi}$$ 
and a differential $\de$ on $E_0\oplus E_1$, such that 
$$\de(E_{1,\xi})\sub E_{0,\xi} \ \text{ and } \de(E_{0,\xi})\sub E_{1,\xi\chi^{-1}}.$$
Now we associate with $\bar{E}$ the complex of $G$-equivariant quasicoherent sheaves on $X$
\begin{equation}\label{com-G-eq}
\com_G(\bar{E}):=\bigoplus_{i=1}^d\com(\bar{E})_{\chi_i}, \text{ where}
\end{equation}
$$\com(\bar{E})_{\chi_i}:\ldots E_{0,\chi_i\chi}\to E_{1,\chi_i}\to E_{0,\chi_i}\to E_{1,\chi_i\chi^{-1}}\to E_{0,\chi_i\chi^{-1}}\to\ldots$$
where the action of $G$ is given by restricting the action of $\Ga$.
Note that we have an isomorphism of $\Z/2$-graded complexes 
$$\com_G(\bar{E})\simeq E_\bullet.$$
This implies the following result.

\begin{prop}\label{W0-prop}
In the above situation the functor 
$$\com_G:\QMF_{\Ga,\chi}(X,0)\to\Com(\Qcoh_G(X))$$
is an equivalence of dg-categories, 
which restricts to an equivalence
$$\MF_{\Ga,\chi}(X,0)\to \Com^b(\Bun_G(X)),$$
where $\Bun_G(X)$ is the category of $G$-equivariant
vector bundles on $X$ and $\Com^b(?)$ denotes the category of bounded complexes. 
The functor $\com_G$ also induces an equivalence 
between the derived category of quasicoherent
\mf s with coherent cohomology 
$\DQMF_c(X/\Ga,0)$ and the derived category $D^b_c(\Qcoh_G(X))$
of complexes with bounded coherent cohomology.
\end{prop}
\ed

\begin{rem}
If we choose $\chi_1$ to be the trivial character 
then
the $G$-invariant part of $\com_G(\bar{E})$ is
isomorphic to $\Ga$-invariants of the complex $\com(\bar{E})$ (see \eqref{mf-Z-gr-com-eq}). 
\end{rem}

\begin{rem}\label{com-G-mf-rem}
In the above situation let $C^\bullet$ be a bounded complex of $\Ga$-equivariant vector bundles on $X$.
Recall that we can associate with it a \mf\ 
$\maf(C^\bullet)\in\MF_{\Ga,\chi}(X,0)$ (see \eqref{maf-eq}). Then we have natural isomorphisms of
$G$-equivariant sheaves
\begin{equation}\label{coh-maf-eq}
\begin{array}{l}
H^{even}\com_G(\maf(C^\bullet))\simeq H^0(\maf(C^\bullet))\simeq H^{even}(C^\bullet),\\
H^{odd}\com_G(\maf(C^\bullet))\simeq H^1(\maf(C^\bullet))\simeq H^{odd}(C^\bullet).
\end{array}
\end{equation}
\end{rem}

\subsection{Connection with categories of singularities}\label{sing-sec}

As was proved by Orlov \cite{Orlov},
the homotopy category of \mf s of a nonzero 
function $W$ on a smooth affine variety $X$
is equivalent to a certain triangulated category $D_{\Sg}(X_0)$ that ``measures the singularity"
of the hypersurface $X_0=(W=0)$. The {\it singularity category} 
$D_{\Sg}(X_0)$ is the quotient of
the derived category $D^b(X_0)$ of coherent sheaves on $X_0$ by the subcategory of
perfect complexes. This definition also makes sense for stacks.
In \cite{PV-stacks} we proved an extension of Orlov's result to smooth stacks (satisfying certain
technical assumptions) by replacing the homotopy category of \mf s with the appropriate
derived category (see below). The stacks that we allow are called FCDRP-stacks 
(see \cite[Def.\ 3.1]{PV-stacks}), where FCD stands for ``finite cohomological dimension" and
RP for ``resolution property".

Let $X$ be a stack, and let
$W\in H^0(X,L)$ be a potential, where $L$ is a line bundle on $X$.
Assume that $W$ is not a zero divisor, and let $X_0=W^{-1}(0)$ be the zero locus of $W$. As in \cite{Orlov}, we consider the natural functor 
\begin{equation}\label{coker-functor-eq}
\fC:\HMF(X,W)\to D_{\Sg}(X_0)
\end{equation}
that associates with a \mf\ $(E_{\bullet},\de)$ the cokernel of $\de_1:E_1\to E_0$. 
This functor is exact (see Lemma 3.12 of \cite{PV-stacks}).

In the case when $X$ is a smooth affine scheme and $L$ is trivial, the functor $\fC$ is an equivalence
by \cite[Thm.\ 3.9]{Orlov}.
In the non-affine case we need to localize the category $\HMF(X,W)$. Namely, we consider the 
full subcategory
$$\LHZ(X,W)\sub\HMF(X,W)$$
consisting of 
\mf s $\bar{E}$ that are locally contractible
(i.e., there exists an open covering $U_i$ of $X$ in smooth topology such that
$\bar{E}|_{U_i}=0$ in $\HMF(U_i,W|_{U_i})$).
Then we define the {\it derived category
of \mf s} as the quotient 
$$\DMF(X,W)=\HMF(X,W)/\LHZ(X,W).$$
We proved in \cite[Thm.\ 3.14]{PV-stacks} that in the case when $X$ is a smooth FCDRP-stack
the functor $\fC$ induces an exact equivalence 
\begin{equation}\label{main-equiv-eq}
\ov{\fC}:\DMF(W)\to D_{\Sg}(X_0).
\end{equation}

We will need the following property of the functor $\fC$.

\begin{prop}\label{sing-tensor-lem} Assume that $W$ is not a zero divisor.
For any
$\bar{E}=(E,\de)\in\HMF(X,W)$ and a bounded complex $C^\bullet$ of vector bundles on $X$
there is an isomorphism
$$\fC(C^\bullet\ot\bar{E})\simeq C^\bullet|_{X_0}\ot\fC(\bar{E})$$
in $D_{\Sg}(X_0)$, where we use the operation of tensor
product of a \mf\ with a complex of vector bundles on $X$ (see \eqref{complex-tensor-eq}).
\end{prop}

\Pf . 
Since $\fC$ commutes with the translation functors, it is enough to consider the case
when $C^\bullet$ is concentrated in non-positive degrees. 
Recall that
$$(C^\bullet\ot E)_i=\bigoplus_{n\ge 0}C^{-n}\ot E(L^{1/2})_{n-i}\ot L^{(n-i)/2}, \text{ for } i=0,1,$$
and the (injective) differential $\wt{\de}_1:(C^\bullet\ot E)_1\to (C^\bullet\ot E)_0$ is induced
by $\de$ and the differential $d_C$ on $C^\bullet$.
Let us consider the coherent sheaves on $X_0$,
$$\FF:=\fC(\bar{E})=\coker(\de_1:E_1\to E_0) \ \text{ and }
\wt{\FF}=\fC(C^\bullet\ot\bar{E})=\coker(\wt{\de}_1).$$
We are going to construct an exact triple of bounded complexes of coherent sheaves on $X_0$,
concentrated in non-negative degrees, of the form
$$0\to S^\bullet\to\GG^\bullet\rTo{p} C^\bullet|_{X_0}\ot\FF\to 0,$$ 
such that the terms of $S^\bullet$ are locally free and $\GG^\bullet$ is a resolution for 
$\wt{\FF}$. This will imply that $\GG^\bullet$ is isomorphic
to $C^\bullet|_{X_0}\ot \FF$ in $D_{\Sg}(X_0)$, and the assertion will follow.
We define the complex $\GG^\bullet$ by
$$\GG^{-i}=C^{-i}\ot\FF\oplus\bigoplus_{n>i} C^{-n}\ot E(L^{1/2})_{n-i}\ot L^{(n-i)/2}|_{X_0}
\ \text{  for } i\ge 0,$$
with the differential induced by $\de$, $d_C$, and the embedding $\FF\to E_1\ot L|_{X_0}$,
The map 
$$p:\GG^\bullet\to C^\bullet\ot \FF=C^\bullet|_{X_0}\ot \FF$$ 
is defined as the natural projection.
It remains to show that $\GG^\bullet$ is a resolution of $\wt{\FF}$.
It is easy to check that the canonical map
$(C^\bullet\ot E)_0\to\wt{\FF}$ factors through a map $\GG^0\to\wt{\FF}$ that extends
to a morphism of complexes $\GG^\bullet\to\wt{\FF}$. To see that it is a quasi-isomorphism,
we use the increasing filtrations on both sides induced by
the stupid filtration $[C^\bullet]_{\ge -n}$ on $C^\bullet$.
The associated quotients of the induced filtration on $\wt{\FF}$ are
\begin{equation}\label{nth-quotient-sh}
C^{-n}\ot\coker\bigl((E(L^{1/2})_{n-1}\ot L^{(n-1)/2}\to E(L^{1/2})_n\ot L^{n/2}\bigr),
\end{equation}
where $n\ge 0$.
On the other hand, the associated quotients
of the induced filtration of $\GG^\bullet$ are
the complexes
\begin{equation}\label{nth-quotient-complex}
C^{-n}\ot \FF_0\to C^{-n}\ot E_1\ot L|_{X_0}\to\ldots\to C^{-n}\ot E(L^{1/2})_n\ot L^{n/2}
\end{equation}
concentrated in degrees $[-n,0]$ with $n\ge 0$.
It remains to observe that by \cite[Lem.\ 1.5]{PV-stacks}, the complex \eqref{nth-quotient-complex}
is a resolution of the sheaf \eqref{nth-quotient-sh}.
\ed

\begin{rem}\label{W0-der-rem} 
In the case $W=0$ the definition of the derived category of \mf s still makes sense.
For instance, as in Proposition \ref{W0-prop}, let us consider the case
$X=Y/\Ga$ with $\Ga$ acting trivially on $Y$ and equipped with
a surjective character $\chi:\Ga\to\G_m$ (which defines the line bundle $L$).
Assume also that $G=\ker(\chi)$ is finite.
Then the category $\DMF_{\Ga,\chi}(Y,0)=\DMF(X,0)$ is equivalent to the category of
$\widehat{G}$-graded objects in the usual derived category of bounded ($\Z$-graded)
complexes of vector bundles on $Y$ (since a bounded acyclic complex of
projective modules is contractible). 
\end{rem}

As we have shown in \cite[Sec.\ 4]{PV-stacks}, the functor $\fC$ extends naturally to quasi-\mf s. More precisely,
we have a functor
$$\fC^{\infty}:\DMF^\infty(X,W)\to D'_{\Sg}(X_0),$$
where $\DMF^\infty(X,W)$ is the derived category of quasi-\mf s defined as the quotient
of the homotopy category $\HMF^\infty(X,W)$ by the subcategory of objects that are locally
homotopic to zero, and $D'_{\Sg}(X_0)$ is the quotient of $D^b(\Qcoh(X_0))$ by the subcategory of
bounded complexes of locally free sheaves.

If $f:X\to Y$ is a smooth affine morphism with integral fibers,
where $X$ and $Y$ are stacks, and $W\in H^0(Y,L)$ is a potential,
then we have a natural push-forward functor that takes a quasi-\mf\ of $f^*W$ on $X$ to
a quasi-\mf\ of $W$ on $Y$ (see \cite[Def.\ 4.8]{PV-stacks}). Furthermore, if $W$ is not a zero divisor and
$Y$ is smooth with the resolution property, then there is an induced functor of derived categories
$f_*:\DMF^\infty(X,f^*W)\to\DMF^\infty(Y,W)$, so that we have
a commutative diagram
\begin{diagram}
\DMF^\infty(X,f^*W)&\rTo{\fC^\infty}&D'_{\Sg}(X_0)\\
\dTo{f_*}&&\dTo{g_*}\\
\DMF^\infty(Y,W)&\rTo{\fC^\infty}&D'_{\Sg}(Y_0)
\end{diagram}
where $Y_0$ (resp., $X_0$) is the zero locus
of $W$ (resp., $f^*W$), $g:X_0\to Y_0$ is the morphism induced by $f$, and
the right vertical arrow is induced by the push-forward functor 
$g_*: D^b(\Qcoh(X_0))\to D^b(\Qcoh(Y_0))$.

\subsection{Supports}\label{supports-sec}

Let $X$ be a stack and let $W\in H^0(X,L)$ be a section.
We denote by $\ov{\DMF}(X,W)$ the idempotent closure of $\DMF(X,W)$.
For a closed substack $Z\subset X$ the full subcategory 
$\ov{\DMF}(X,Z;W)\sub \ov{\DMF}(X,W)$ of \mf s with support on $Z$
is defined as the common kernel of the restriction functors
$$\bar{E}\mapsto i_x^*\com(\bar{E})$$ 
for closed points $x\in X_0\setminus Z$,
where $X_0=W^{-1}(0)\sub X$ (see \cite[(5.1)]{PV-stacks}).
In the case when $W$ is a non-zero-divisor 
and $X$ is a smooth FCDRP-stack
we proved that
a \mf\ $\bar{E}=(E,\de)$ belongs to $\ov{\DMF}(X,Z;W)$ 
if and only if it restricts to zero in 
$\ov{\DMF}(X\setminus Z,W|_{X\setminus Z})$
(see the proof of \cite[Prop.\ 5.6]{PV-stacks}).

\begin{lem}\label{Koszul-mf-support-lem} 
Let $\{\a,\b\}$ be the Koszul \mf\ of $W\in H^0(X,L)$
associated with sections $\a\in V\ot L$ and $\b\in V^{\vee}$
such that $\lan\a,\b\ran=W$ (see \eqref{Koszul-mf-def}). Then
$\{\a,\b\}$ is supported on the zero locus of the section $(\a,\b)\in (V\ot L)\oplus V^{\vee}$. 
\end{lem}

\Pf . By definition of support of a \mf,
we have to check that if $V$ is a vector space and $\a\in V$ and $\b\in V^*$ are 
such that $\lan \a,\b\ran=0$ and $(\a,\b)\neq (0,0)$, then the complex
$$\{\a,\b\}=({\bigwedge}^\bullet V,\a\we ?+\iota(\b))$$
is acyclic. Since $\lan\a,\b\ran=0$, we can find
 a direct sum decomposition $V=V_1\oplus V_2$ such that
$\a\in V_1$ and $\b\in V_2^*$. Then $\{\a,\b\}$ becomes the total complex of the tensor product
of the complexes $({\bigwedge}^\bullet V_1,\a\we ?)$ and 
$({\bigwedge}^\bullet V_2, \iota(\b))$, at least one of which is acyclic because $(\a,\b)\neq (0,0)$.
\ed

Recall that for a pair of potentials $W,W'\in H^0(X,L)$ 
we have the tensor product bifunctor (see \eqref{usual-tensor-product})
\be\label{HMF-tensor-functor}
\HMF(X,W)\times\HMF(X,W')\to\HMF(X,W+W').
\end{equation}

\begin{prop}\label{tensor-prod-lem} 
Let $X$ be a stack 
and $W,W'\in H^0(X,L)$ two sections. 

\noindent
(i) For a a pair of closed substacks $Z\sub X$ and $Z'\sub X$ 
the bifunctor \eqref{HMF-tensor-functor}
induces an exact bifunctor
$$\ov{\DMF}(X,Z;W)\times \ov{\DMF}(X,Z';W')\to \ov{\DMF}(X,Z\cap Z';W+W').$$

\noindent
(ii) Assume in addition that $X$ is a smooth FCDRP stack and $W$ and $W'$ are not
zero divisors. The the tensor product induces a bifunctor
$$\ov{\DMF}(X,W)\times \ov{\DMF}(X,W')\to \ov{\DMF}(X,\Si\cap \Si';W+W').$$
where $\Sigma$ (resp., $\Sigma'$)
is the singularity locus of the hypersurface $W=0$ (resp., $W'=0$).
\end{prop}

\Pf . (i) The tensor product functor \eqref{HMF-tensor-functor} is compatible with pull-backs, hence,
if at least one of the \mf s
$\bar{E}\in\HMF(X,W)$ or $\bar{F}\in\HMF(X,W')$ is locally contractible
then so is $\bar{E}\ot\bar{F}$. Therefore, this functor descends to derived categories.
Similarly, for $\bar{E}\in\ov{\DMF}(X,Z;W)$ and $\bar{F}\in\ov{\DMF}(X,Z';W')$ 
the tensor product $\bar{E}\ot\bar{F}$ is contractible in a neighborhood of $x\not\in Z\cap Z'$ 
(since either $\bar{E}$ of $\bar{F}$ is contractible near $x$). 

\noindent
(ii) This follows from part (i) and \cite[Cor.\ 5.3]{PV-stacks}.
\ed

\subsection{Push-forwards}\label{push-forward-sec}

In \cite[Sec.\ 6]{PV-stacks} we defined the push-forwards for \mf s with relatively proper support. 
Let $f:X\to Y$ be a representable morphism of
smooth FCDRP-stacks and $W\in H^0(Y,L)$ a potential
such that $W$ and $f^*W$ are not zero divisors.
Let $Z\sub X_0$ be a closed substack of the zero locus of $f^*W$, such that
the induced morphism $f:Z\to Y$ is proper. Let $Y_0\sub Y$ denote the zero locus of $W$, and
let $f_0:X_0\to Y_0$ be the map induced by $f$. 
The derived push-forward functor
$$Rf_{0*}:D^b(X_0,Z)\to D^b(Y_0,f(Z))$$
induces a functor
$$\ov{D_{\Sg}(X_0,Z)}\to\ov{D_{\Sg}(Y_0,f(Z))},$$
and hence, by \cite[Prop.\ 5.6]{PV-stacks}, a functor
\begin{equation}\label{mf-push-forward}
Rf_*:\ov{\DMF}(X,Z;f^*W)\to \ov{\DMF}(Y,f(Z);W).
\end{equation}
In the case when $f$ is proper
this functor preserves the usual derived categories inside their idempotent completions
(see \cite[Rem.\ 6.2]{PV-stacks}).
We proved in \cite[Prop.\ 6.3]{PV-stacks} that this functor is compatible with the natural push-forwards
for quasi-\mf s with respect to smooth affine morphisms with integral fibers
(see \cite[Def.\ 4.8]{PV-stacks}).
More generally, if $f:X\to Y$ is an arbitrary affine morphism, then
for any $W\in H^0(Y,L)$ we have the naive
push-forward functor
\begin{equation}\label{naive-push-forward}
f_*:\QMF(X,f^*W)\to\QMF(Y,W): \bar{E}=(E,\de)\mapsto (f_*E,f_*\de)
\end{equation}
for quasicoherent \mf s (see Section \ref{mf-cat-sec}).
Note that
in the case when $W=0$ this functor respects quasi-isomorphisms.

\begin{rem}\label{coh-mf-rem} 
In the case when $f$ is a finite morphism, the functor \eqref{mf-push-forward} is compatible
with the naive push-forward \eqref{naive-push-forward} in the following sense.
For a \mf\ $\bar{E}=(E,\de)$ the naive push-forward 
$f_*\bar{E}$ is a coherent \mf\ with an additional property that 
the multiplication by $W$ is an injective endomorphism.
The cokernel functor $\fC:\MF(Y,W)\to D_{\Sg}(Y_0)$
extends naturally to the category of such coherent \mf s, so we can view $f_*\bar{E}$
as an object of the derived category $\DMF(Y,W)$ 
(cf. \cite[Def.\ 3.21, Rem.\ 6.2]{PV-stacks}).
Since $\fC(f_*\bar{E})\simeq f_{0*}(\bar{E})$, we obtain an isomorphism in $\DMF(Y,W)$
of $f_*\bar{E}$ with the push-forward of $\bar{E}$ given by \ref{mf-push-forward}.
\end{rem}

\begin{rem} As it was pointed out to us by Leonid Positselski,
an alternative construction of push-forwards can be given using exotic
derived categories of \cite{Posic} and the results of \cite{Posic-mf}.
One can start with the natural push-forward functors between the coderived categories
of quasicoherent \mf s which are defined using injective resolutions (see \cite[Sec.\ 3.7]{Posic}).
Then one can use the fact that for a regular scheme of finite Krull dimension
the coderived category of quasi-\mf s is equivalent 
to the coderived category of quasicoherent \mf s (see \cite[Thm.\ 1(a)]{Posic-mf}). Note that
the absolute derived category of \mf s, which is a full subcategory in the above coderived category,
is equivalent to the corresponding hypersurface singularity category 
(see \cite[Thm.\ 2]{Posic}, \cite[Prop.\ 2.13]{LP} or \cite{O-nonaff}) 
and hence to our category $\DMF(X,W)$.
\end{rem}

The following lemma gives an important relation between the push-forward with respect to the embedding of the zero locus of $W$ and the tensor
product of \mf s of $W$ and of $-W$.

\begin{lem}\label{complex-quasi-isom-lem} 
Let $X$ be a smooth FCDRP-stack, $W\in H^0(X,L)$ a potential and $\bar{E}$ a \mf\ of $W$.

\noindent
(i) Let us define the \mf\ $C_\bullet(\bar{E})$ of $0$ on $X$ as follows:
$$C_0(\bar{E})=E_0\oplus E_1\otimes L,\ \ C_1(\bar{E})=E_1\oplus E_0,$$
with the differential $\de(x,y)=(\de_E(x)+y,-W(x)-\de_E(y))$.
Then $C_\bullet(\bar{E})$ is contractible (i.e., homotopy equivalent to $0$).

\noindent
(ii) Let $\bar{F}$ be a \mf\ of $-W$, and let $i:X_0\hra X$
be the inclusion of the zero locus of $W$.
Assume that $W$ is not a zero divisor.
Then the map
\begin{equation}\label{DCMF-0-fun}
\FF\mapsto i_*(\FF\ot_{\OO_{X_0}}i^*\bar{F})
\end{equation}
gives a well-defined functor
$D_{\Sg}(X_0)\to \DCMF(X,0)$,
where $\DCMF(X,0)$ is the derived category of coherent \mf s (see Section \ref{mf-cat-sec}).
Also, the natural map of coherent \mf s of $0$ on $X$
$$q:\bar{E}\ot\bar{F}\to i_*\left(\fC(\bar{E})\ot_{\OO_{X_0}}i^*\bar{F})\right),$$
induced by the projection $E_0\to\fC(\bar{E})=\coker(E_1\to E_0)$,
is a quasi-isomorphism. 
\end{lem}

\Pf . (i) The contracting homotopy sends $(x,y)$ to $(0,y)$.

\noindent
(ii) By \cite[Lem.\ 1.5]{PV-stacks},
 the complex $\com(i^*\bar{F})=\com(\bar{F})$ of bundles on $X_0$ is exact. Therefore,
from the isomorphism~\eqref{tensor-com-isom} it follows
that  for a
perfect complex $C^\bullet$ on $X_0$ the coherent \mf\ of zero
$C^\bullet\ot i^*\bar{F}$ is acyclic. 
Hence, \eqref{DCMF-0-fun} gives a well-defined functor.
We have to check that the map $\com(q)$ is a quasi-isomorphism (see \eqref{mf-Z-gr-com-eq}).
But this map fits into an exact sequence of complexes of sheaves on $X$
$$0\to K^\bullet\to\com(\bar{E}\ot\bar{F})\rTo{\com(q)} i_*\left(\fC(\bar{E})\ot\com(\bar{F})\right)\to 0,$$
where $K^\bullet=E_1\ot\com(C_\bullet(\bar{F}))$. Hence, by part (i), $K^\bullet$ is acyclic. 
\ed

We will need the push-forward functors in the following situation. 

\begin{ex}\label{G-push-forward-constr} 
Let $\pi:E\to X$ be a smooth affine morphism with integral fibers, 
where $X$ is a smooth FCDRP-stack, and let $W\in H^0(X,\OO_X)$ be a potential. 
Suppose that we have a commutative algebraic group 
$\Ga$ acting on both $E$ and $X$ compatibly, so that $W$ is semi-invariant with respect to $\Ga$
and a character $\chi:\Ga\to\G_m$, where $\chi$ is surjective with finite kernel $G$. 
Let $Z\sub E$ be a $\Ga$-invariant closed substack 
such that $\pi^*W|_Z=0$ and $Z$ is proper over $X$.
We will define the push-forward functor in the following two cases.

\noindent
{\bf Case 1.} Assume that $W$ is not a zero divisor.
Then by \cite[Prop.\ 6.1]{PV-stacks} (applied to the morphism $E/\Ga\to X/\Ga$), 
we have the push-forward functor
$$\pi_*:\ov{\DMF}_\Ga(E,Z;\pi^*W)\to\ov{\DMF}_{\Ga}(X,W).$$
If in addition, a subgroup $I\sub G$ acts trivially on $X$ 
then we can combine
the above functor with taking $I$-invariants to get a functor
\begin{equation}\label{G-mf-push-forward}
\pi_*^I:\ov{\DMF}_{\Ga}(E,Z;\pi^*W)\to\ov{\DMF}_{\Ga/I}(X,W).
\end{equation}

\noindent
{\bf Case 2.}
Assume that $W=0$ and the action of $\Ga$ on $X$ (but not necessarily on $E$) is trivial.
Then we can still define the push-forward functor as follows. Given a $\Ga$-equivariant \mf\ 
$\bar{P}$ of $0$ on $E$, supported on $Z$, we can consider the push-forward
$\pi_*\bar{P}$ (see \cite[Def.\ 4.8]{PV-stacks}) 
which will be a $\Ga$-equivariant 
quasi-\mf\ of $0$ on $X$.
Let us consider the corresponding complex $\com_G(\pi_*\bar{P})$
of $G$-equivariant quasicoherent sheaves defined by \eqref{com-G-eq} 
(using representatives for the $\lan\chi\ran$-cosets in
the character group of $\Ga$). The assumption on the support of $\bar{P}$ implies that
this complex has bounded coherent cohomology. Since the category of such complexes 
is equivalent to the bounded derived category of coherent
sheaves, we obtain a functor
\begin{equation}\label{G-mf-0-push-forward}
\pi_*:\ov{\DMF}_\Ga(E,Z;0)\to D^b_G(X)\simeq \DMF_{\Ga}(X,0).
\end{equation}
The map on the Grothendieck groups induced by this functor
sends the class of a \mf\  $\bar{P}$ to $[\pi_*H^0(P,\de)]-[\pi_*H^1(P,\de)]$.

In both cases the diagram of functors
\begin{equation}\label{push-forward-coh-quasicoh-diagram}
\begin{diagram}
\ov{\DMF}_\Ga(E,Z;\pi^*W)&\rTo{\pi_*}& \ov{\DMF}_{\Ga}(X,W)\\
\dTo{}&&\dTo{}\\
\DMF^\infty_\Ga(E,f^*W) &\rTo{\pi_*}& \DMF^\infty_\Ga(X,W)
\end{diagram}
\end{equation}
is commutative. Indeed, for $W=0$ (and $\Ga$ acting
trivially on $X$) this is clear from the definition, and for $W\neq 0$ this follows from
\cite[Prop.\ 6.3]{PV-stacks}.

Note that if we have
a $\Ga$-equivariant closed substack $i:E'\to E$ such that $Z\sub E'$, such that
$\pi'=\pi|_{E'}$ is still smooth with integral fibers,
then we can consider the push-forward functors
\eqref{mf-push-forward} associated with the projections $\pi:E\to X$ and
$\pi':E'\to X$ and also the functor $i_*:\ov{\DMF}(E',f^*W|_{E'})\to \ov{\DMF}(E,f^*W)$. 
In this situation one has an isomorphism of functors 
$$\pi_*\circ i_*\simeq \pi'_*$$
from $\ov{\DMF}(E',Z;(\pi')^*W|_{E'})$ to $\ov{\DMF}(X,W)$.
\end{ex}

We have the following analog of the projection formula.

\begin{prop}\label{proj-formula-prop} 
Let $\pi:E\to X$, $Z\sub E$, $\Ga$ and $W$ be as
in one of the two cases of Example \ref{G-push-forward-constr}, where $\Ga$ acts
trivially on $X$. Then for $\bar{P}\in\ov{\DMF}_\Ga(E,Z;\pi^*W)$ and
$\bar{Q}\in\ov{\DMF}_\Ga(X,-W)$ one has a functorial isomorphism
$$\pi_*(\bar{P}\ot \pi^*\bar{Q})\simeq \pi_*(\bar{P})\ot \bar{Q}.$$
\end{prop}

\Pf . 
The case $W=0$ immediately reduces to the usual projection formula, so we will assume that
$W$ is not a zero divisor.
By Lemma \ref{complex-quasi-isom-lem}(ii), we have natural quasi-isomorphisms
$$\pi_*(\bar{P}\ot \pi^*\bar{Q})\simeq 
i_*\pi_{0*}\left(\fC(\bar{P})\ot_{\OO_{E_0}}\pi_0^*i^*\bar{Q})\right) 
\text{ and}$$
$$\pi_*(\bar{P})\ot\bar{Q}\simeq i_*\left(\fC(\pi_*\bar{P})\ot_{\OO_{X_0}}i^*\bar{Q}\right),$$
where $E_0=\pi^{-1}(X_0)$, $\pi_0:E_0\to X_0$ is the restriction of $\pi$, and $i:X_0\to X$
is the natural embedding. It remains to 
use the isomorphism
$$\pi_{0*}\fC(\bar{P})\simeq \fC(\pi_*\bar{P})$$ 
and the usual projection formula for $\pi_0$.
\ed

Our push-forward functors also have the following base change property. 

\begin{prop}\label{base-change-prop}
(i) Suppose we have a cartesian diagram of smooth FCDRP-stacks
\begin{diagram}
X' &\rTo{v}& X\\
\dTo{f'}&&\dTo{f}\\
Y' &\rTo{u}& Y
\end{diagram}
where $f$ is a flat representable morphism, 
and let $W\in H^0(Y,L)$ be a potential such that its pull-backs 
$W_{X}$, $W_{X'}$ and $W_{Y'}$ to
$X$, $X'$ and $Y'$, respectively, are not zero-divisors. Let $X_0$, $Y_0$, $X'_0$ and $Y'_0$ be
the zero loci of these potentials, $Z\sub X_0$ a closed substack, proper over $Y$,
and $Z'\sub X'_0$ the induced closed substack. Then the diagram of functors
\begin{diagram}
 \ov{\DMF}(X,Z,W_X)&\rTo{v^*}&\ov{\DMF}(X',Z',W_{X'})\\
 \dTo{Rf_*}&&\dTo{Rf'_*}\\
 \ov{\DMF}(Y,W)&\rTo{u^*}&\ov{\DMF}(Y',W_{Y'})
\end{diagram}
is commutative.

\noindent
(ii) Let $(\pi:E\to X, W, \Ga, \chi, G, Z)$ be as in Example \ref{G-push-forward-constr}, where
either $W$ is not a zero-divisor or $W=0$ and $\Ga$ acts trivially on $X$.
Suppose we have a cartesian diagram
 \begin{diagram}
E' &\rTo{v}& E\\
\dTo{\pi'}&&\dTo{\pi}\\
X' &\rTo{u}& X
\end{diagram}
of stacks with $\Ga$-action, where all the maps are $\Ga$-equivariant.
Assume that the action of $\Ga$ on $X'$ is trivial and $u^*W=0$. Then for any \mf\ 
$\bar{P}$ of $W$ on $E$, supported on $Z$, there is an isomorphism
$$u^*\pi_*\bar{P}\simeq \pi'_*v^*\bar{P}$$
in $D^b_{G}(X')\simeq\DMF_\Ga(X',0)$.
\end{prop}

\Pf . (i) This follows easily from the base change formula for stacks (see \cite[Prop.\ 3.10]{BFN}).

\noindent
(ii) It is enough to check the corresponding isomorphism in $\DMF^\infty_\Ga(X',0)$.
Therefore, the statement follows from the commutativity of the diagrams 
\eqref{push-forward-coh-quasicoh-diagram} for $\pi$ and $\pi'$ together
with the usual base change formula.
\ed

As in the case of usual sheaves, the base change formula leads to a relative 
K\"unneth isomorphism for push-forwards of \mf s.

\begin{prop}\label{rel-Kun-lem} Let $f_1:X_1\to Y_1$ and $f_2:X_2\to Y_2$ be smooth
morphisms with connected fibers of FCDRP-stacks over a smooth FCDRP-stack $S$.
Let $W_1\in H^0(Y_1,L_1)$ and $W_2\in H^0(Y_2,L_2)$ be potentials,
and let $W=p_{Y_1}^*W_1+p_{Y_2}^*W_2$ be the corresponding potential on the $\G_m$-torsor $Y\to Y_1\times_S Y_2$
associated with the line bundle $L_1\ot L_2^{-1}$,
where $p_{Y_i}:Y\to Y_i$ are the projections.
Let $Z_1\sub X_1$ (resp., $Z_2\sub X_2$) be a closed substack, proper over $Y_1$ (resp., $Y_2$).
Consider the map $f:X\to Y$ fitting into the Cartesian square
\begin{diagram}
X&\rTo{f}&Y\\
\dTo{}&&\dTo{}\\
X_1\times_S X_2&\rTo{f_1\times_S f_2}& Y_1\times_S Y_2
\end{diagram}
Assume that $W_1$, $W_2$ and $W$ are non-zero-divisors. 
Then for $\bar{E}_1\in\ov{DMF}(X_1,Z_1,f_1^*W_1)$ and
$\bar{E}_2\in\ov{\DMF}(X_2,Z_2,f_2^*W_2)$
there is a functorial isomorphism
$$Rf_*\bigl(p_{X_1}^*(\bar{E}_1)\ot p_{X_2}^*(\bar{E}_2)\bigr)\simeq 
p_{Y_1}^*(Rf_{1*}(\bar{E}_1))\otimes p_{Y_2}^*(Rf_{2*}(\bar{E}_2))$$
in $\ov{\DMF}(Y,W)$, where $p_{X_i}:X\to X_i$ are the projections.

The same assertion holds if one (or both) of the morphisms
$f_i$ is of the form $E/\Ga\to X/\Ga$, where $E\to X$ is a smooth affine morphism,
$\Ga$ is a commutative group acting trivially on $Z$ and $W_i=0$ (see Case 2 of Example 
\ref{G-push-forward-constr}).
\end{prop}

\Pf . Let us denote $\bP=Rf_*\bigl(p_{X_1}^*(\bar{E}_1)\ot p_{X_2}^*(\bar{E}_2)\bigr)$ and
$\bP_i=Rf_{i*}(\bar{E}_i)$ for $i=1,2$.
Consider the commutative diagram with a cartesian square
\begin{equation}\label{sum-sing-diag}
\begin{diagram}
X &\rTo{\wt{f}_2}&X_1\times_{Y_1} Y&\rTo{\wt{f}_1}& Y\\
&\rdTo{p_{X_1}}&\dTo{p_1}&&\dTo{p_{Y_1}}\\
&&X_1&\rTo{f_1}&Y_1
\end{diagram}
\end{equation}
where $\wt{f}_2$ and $\wt{f}_1$ are obtained from $\id_{X_1}\times f_2$ and $f_1\times\id_{Y_2}$ by
the base change.
Note that the composition of the arrows in the first row is equal to $f$. Thus, we have
$$\bP
\simeq R(\wt{f}_1)_*R(\wt{f}_2)_*(p_{X_1}^*(\bar{E}_1)\ot p_{X_2}^*(\bar{E}_2))
\simeq R(\wt{f}_1)_*\bigl(p_1^*(\bar{E}_1)\ot R(\wt{f}_2)_*p_{X_2}^*(\bar{E}_2)\bigr).
$$
Applying the base change formula in the cartesian square
\begin{diagram}
X &\rTo{p_{X_2}}& X_2\\
\dTo{\wt{f}_2}&&\dTo{f_2}\\
X_1\times_{Y_1} Y&\rTo{p_{Y_2}}& Y_2
\end{diagram}
we get
$$R(\wt{f}_2)_*p_{X_2}^*(\bar{E}_2)\simeq p_{Y_2}^*\bP_2.$$
Hence,
$$\bP\simeq R(\wt{f}_1)_*(p_1^*(\bar{E}_1))\ot p_{Y_2}^*(\bP_2),$$
Finally, the base change formula in the cartesian
square of \eqref{sum-sing-diag} shows that
$$R(\wt{f}_1)_*(p_1^*(\bar{E}_1))\simeq p_{Y_1}^*(\bP_1).$$
\ed

\subsection{Regular Koszul \mf s}

Here we study Koszul \mf s $\{\a,\b\}$ in the case when $\b$ is a regular section.
Such \mf s should be viewed as deformations of the Koszul complex $\{0,\b\}$. 

In this section we assume that $X$ is a smooth FCDRP-stack.
We fix a potential $W\in H^0(X,L)$, a vector bundle $V$ on $X$, and sections
$\a\in H^0(X,V\ot L)$ and $\b\in H^0(X,V^{\vee})$, such that
$\lan \a,\b\ran=W$.

\begin{defi} The Koszul \mf\ (see \eqref{Koszul-mf-def})
$$\{\a,\b\}=\left({\bigwedge}^\bullet(V\ot L^{1/2})(L^{-1/2}), \de_{\a,\b}\right),$$
is called {\it regular} if $\b$ is a regular section of $V^\vee$.
\end{defi}

In this section we always assume that $\{\a,\b\}$ is regular 
(with the exception of  Proposition \ref{koszul-deformed-prop}) and denote by
$i:X'\hra X$ the embedding of the zero locus of $\b$. 
Note that $X'$ is contained in $X_0$, the zero locus of $W$.

\begin{lem}\label{koszul-main-lem} 
(i) If $W$ is not a zero divisor then
$$\fC(\{\a,\b\})\simeq \OO_{X'}\ \text{ in } D_{\Sg}(X_0),$$
where $\fC$ is the cokernel functor \eqref{coker-functor-eq}, and
$\OO_{X'}$ is viewed as a coherent sheaf on $X_0$. 

\noindent
(ii) If $W=0$ then 
$$H^0(\{\a,\b\})\simeq i_*\OO_{X'}\ \text{ and } H^1(\{\a,\b\})=0,$$ 
where $H^i(\bar{P}):=H^i(\com(\bar{P}))$ (see \eqref{mf-Z-gr-com-eq}). 
In other words, the natural morphism of quasicoherent \mf s
$$\{\a,\b\}\to i_*\maf(\OO_{X'})$$
is a quasi-isomorphism (where $\maf(\OO_{X'})_0=\OO_{X'}$, $\maf(\OO_{X'})_1=0$).
\end{lem}

\Pf .
(i) We have a natural map $\fC(\{\a,\b\})\to i_*\OO_{X'}$ induced by the projection 
${\bigwedge}^\bullet\to {\bigwedge}^0$
and by the map $\OO_X\to \OO_{X'}=\coker(\b^\vee:V\to\OO_X)$.
It is enough to prove that this map is an isomorphism locally, 
so the statement reduces to the affine case proved in \cite[Prop. 2.3.1]{PV-mf}. 

\noindent
(ii) Since $\lan\a,\b\ran=W=0$,
the complex $\com(\{\a,\b\})$ 
can be identified with the total complex of the bicomplex 
$$K=\bigoplus_{i,j} K^{i,j}, \ \text{ where } 
K^{i,j}={\bigwedge}^{i-j}V\ot L^i,$$
with differentials given by $\iota(\b)$ and $\a\we?$. 
The regularity of $\b$ implies that the cohomology of the differential
$\iota(\b)$ is concentrated along the diagonal $\bigoplus_i K^{i,i}$.
From the spectral sequence we
immediately see that 
$$H^{2n+1}(\com(\{\a,\b\}))=0\ \text{ and }H^{2n}(\com\{\a,\b\})\simeq i_*\OO_{X'}\ot L^n.$$
\ed

\begin{prop}\label{koszul-reg-prop}
(i) Assume that $X'$ (the zero locus of $\b$) is smooth. Let 
$W_1\in H^0(X,L)$ be another potential such that $W+W_1$ and $W_1|_{X'}$ are not zero divisors. 
Then for every \mf\ $\bar{P}=(P,\de)\in\MF(X,W_1)$ we have a functorial isomorphism
in $\DMF(X,W+W_1)$
\begin{equation}\label{koszul-reg-isom}
q: \bar{P}\ot\{\a,\b\}\rTo{\sim} i_*i^*\bar{P},
\end{equation}
where on the right-hand side $i_*$ is the push-forward functor
$$i_*:\DMF(X',W_1|_{X'})\to\DMF(X,W+W_1).$$

\noindent
(ii) Assume that $W$ is not a zero divisor. Then for every \mf\ $\bar{P}\in\MF(X,-W)$ we have
a quasi-isomorphism in $\QMF(X,0)$
$$q: \bar{P}\ot\{\a,\b\}\to i_*i^*\bar{P}.$$
\end{prop}

\Pf . (i) The projection $\bigwedge^{\bullet}(V\ot L^{1/2})(L^{-1/2})\to\OO_X$
induces a natural morphism 
$$q: P\ot{\bigwedge}^{\bullet}(V\ot L^{1/2})(L^{-1/2})\to i_*i^*(P\ot{\bigwedge}^{\bullet}(V\ot L^{1/2})(L^{-1/2}))\to i_*i^*P
$$ 
of coherent \mf s of $W+W_1$, where we use the fact that the naive push-forward $i_*i^*P$ is compatible with the push-forward functor \eqref{mf-push-forward}
(see Remark \ref{coh-mf-rem}).
To show that $q$ induces an isomorphism
in $\DMF(X,W+W_1)$, we can argue locally.
Thus, we can assume that $L$ and $V$ are trivial bundles. 
We will use induction in the rank $r$ of $V$. In the case when $r=1$, i.e., $V=\OO$, we have 
$\b=f$, $\a=g$, where $f$ and $g$
are functions on $X$ such that $W=fg$ and $i:X'=Z(f)\hra X$ is a divisor.
By definition,
$$\fC(\bar{P}\ot\{g,f\})=\coker(D:P_1\oplus P_0\to P_0\oplus P_1),$$
where 
$$D(p_1,p_0)=(\de(p_1)+f\cdot p_0,\de(p_0)-g\cdot p_1).
$$
Let us consider the exact triple of two-term complexes
$$0\to [0\to P_1]\to [P_1\oplus P_0\rTo{D}P_0\oplus P_1]\to [P_1\oplus P_0\rTo{\ov{D}}P_0]\to 0,$$
where $\ov{D}(p_1,p_0)=\de(p_1)+f\cdot p_0$.
Since $D$ is a part of the differential of the \mf\ $\bar{P}\ot\{g,f\}$ of $W+W_1$, it is injective. Thus, from
the above exact triple we get an exact sequence of sheaves
\begin{equation}\label{exact-seq-coker}
0\to\ker(\ov{D})\rTo{\ga} P_1\to\coker(D)\to\coker(\ov{D})\to 0,
\end{equation}
where $\ga(p_1,p_0)=\de(p_0)-g\cdot p_1$.
Now 
$$\coker(\ov{D})\simeq P_0/(\de(P_1)+fP_0)\simeq i_*\fC(i^*\bar{P}).$$
Since $W_1|_{X'}\neq 0$, the morphism $\de|_{X'}:P_1/fP_1\to P_0/fP_0$ is
injective. In other words, $\de^{-1}(fP_0)\cap P_1=fP_1$. 
This implies that the map $p_1\mapsto (fp_1,-\de(p_1))$
gives an isomorphism
$$\la: P_1\wt{\ra}\ker(\ov{D}).$$ 
The composition of $\ga$ with $\la$ sends
$p_1$ to $-\de^2(p_1)-fg\cdot p_1=-(W+W_1)\cdot p_1$. Thus, the exact sequence
\eqref{exact-seq-coker} is isomorphic to
$$0\to P_1\rTo{-W-W_1} P_1\to \fC(\bar{P}\ot \{g,f\})\to i_*\fC(i^*\bar{P})\to 0.$$
Since $\fC(i_*i^*\bar{P})\simeq i_*\fC(i^*\bar{P})$, we get an exact sequence
$$0\to P_1/(W+W_1)P_1\to \fC(\bar{P}\ot \{g,f\})\rTo{\fC(q)} \fC(i_*i^*\bar{P})\to 0$$
that implies that the map $\fC(q)$ is an isomorphism in $D_{\Sg}(X_0)$, and so
by \cite[Thm.\ 3.14]{PV-stacks}, $q$ is an isomorphism. This gives the base of induction.

When $r=\rk V>1$, decompose $V$ as $V=\OO_X\oplus V'$,
and let $\{\a,\b\}=\{g,f\}\ot\{\a',\b'\}$ be the corresponding decomposition,
where $\b'\in (V')^{\vee}$, $\a'\in V'$, $f,g\in H^0(\OO_X)$. 
Let $j:X''\hra X$ be the zero locus of $\b'$. 
Note that $X'\sub X''$ is the zero locus of the function $j^*f$.
Since $\b'$ is a regular section of $V'$, we can apply the
induction hypothesis to the \mf\ $\{\a',\b'\}$ on $X$ and conclude that
the natural map
$$\bar{P}\ot\{g,f\}\ot\{\a',\b'\}\to j_* (j^*\bar{P}\ot\{j^*g,j^*f\})$$
is an isomorphism.
Applying the case $r=1$ to the \mf\ $\{j^*g,j^*f\}$, we see
that 
$$j^*\bar{P}\ot \{j^*g,j^*f\}\simeq k_*k^*j^*\bar{P}\simeq k_*i^*\bar{P},$$
where $k:X'\hra X''$ is the natural embedding. This establishes the induction step.
Now the assertion follows.

\noindent
(ii) 
Since $\fC(\{\a,\b\})\simeq i_*\OO_{X'}$ in $D_{\Sg}(X_0)$ (see Lemma 
\ref{koszul-main-lem}(i)), this follows from Lemma \ref{complex-quasi-isom-lem}(ii). 
\ed

The following result deals with the situation when the section $\b\in H^0(X,V^\vee)$ is not regular but
is the image of a regular section of a subbundle of $V^\vee$.

\begin{prop}\label{koszul-deformed-prop}
Let $U\sub V$ be a subbundle.
Assume that we have a regular section $\b'\in H^0(X,(V/U)^\vee)$ such that
$\b=\iota(\b')$,
where $\iota:(V/U)^\vee\to V^\vee$ is the natural inclusion.

\noindent
(i) Assume that $W$ is not a zero divisor. Then
we have the following equality in the Grothendieck group of $D_{\Sg}(X_0)$:
$$[\fC(\{\a,\b\})]=[i'_*i^*{\bigwedge}^\bullet(U\ot L^{1/2})(L^{-1/2})],$$
where $i':X'\hra X_0$ is the inclusion.

\noindent
(ii) Assume that $W=0$. Then
$$[\{\a,\b\}]=[i_*\maf i^*{\bigwedge}^\bullet(U\ot L^{1/2})(L^{-1/2})],$$
in the Grothendieck group of the derived category of coherent \mf s (see Section \ref{mf-cat-sec}).
\end{prop}

\Pf . Note that the differential $\de_{\a,\b}$ is compatible with the filtration of the exterior algebra
${\bigwedge}^\bullet(V\ot L^{1/2})(L^{-1/2})$ by powers of the ideal generated by
$U$. Hence, in the Grothendieck group we can replace $\{\a,\b\}$ by the associated quotient,
which is isomorphic to the tensor product $C^\bullet\ot\{\ov{\a},\b'\}$, where
$$C^\bullet={\bigwedge}^\bullet(U\ot L^{1/2})(L^{-1/2})$$
and $\ov{\a}$ is the section of $(V/U)\ot L$ induced by $\a$.
Since $\b'$ is regular, in case (i) we have $\fC(\{\ov{\a},\b'\})\simeq i'_*\OO_{X'}$ by Lemma
\ref{koszul-main-lem}(i). By Proposition \ref{sing-tensor-lem}, this implies that
$$\fC(C^\bullet\ot\{\ov{\a},\b'\}\simeq i'_*i^*C^\bullet.$$
In case (ii) we deduce from Lemma \ref{koszul-main-lem}(ii) combined with
Lemma \ref{complex-tensor-lem}
 the following quasi-isomorphisms of coherent \mf s:
$$C^\bullet\ot\{\ov{\a},\b'\}\simeq C^\bullet\ot i_*\maf(\OO_{X'})\simeq
i_*(i^*C^\bullet\ot\maf(\OO_{X'}))\simeq i_*\maf(i^*C^\bullet).$$
\ed

\section{Matrix factorizations of a quasihomogeneous isolated singularity}
\label{quasihom-mf-sec}

Throughout this section we fix a quasihomogeneous potential $\w$ on $\A^n$
with an isolated singularity at $0$. Recall that the latter condition means that the quotient
$\C[x_1,\ldots,x_n]/(\pa_1\w,\ldots,\pa_n\w)$ is finite-dimensional (where $\pa_i$ denotes
the partial derivative with respect to $x_i$).
We fix the set of coprime positive degrees $(d_1,\ldots,d_n)$, such that $\w$ is homogeneous
of degree $d$ with respect to the grading $\deg(x_i)=d_i$.

From now on we denote by $\w$ such a potential. 

In this section we calculate the Hochschild homology of the dg-category
of $\Ga$-equivariant
\mf s of $\w$, where $\Ga$ is a certain one-dimensional subgroup of $\G_m^n$,
and the canonical bilinear form on this Hochschild
homology. 
This is a $\Z$-graded analog of the computations of the Hochschild homology
of the $\Z/2$-dg-category of \mf s of $\w$ and of the canonical bilinear form on it
performed in \cite{Dyck} and \cite{PV-mf}.
We also discuss functors between categories of equivariant
\mf s given by kernels.

\subsection{Symmetry groups}\label{sym-sec}

With $\w$ we associate certain natural groups as follows. Let us write
$$\w=\sum_{k=1}^N c_k M_k,$$
where $M_k(x_1,\ldots,x_k)$ are monomials and $c_k\in\C^*$.
We have a homomorphism
$$\rho:\G_m^n\to \G_m^N:(\la_\bullet)\mapsto (M_k(\la_\bullet)).$$
Let $\Ga_\w\sub\G_m^n$ be the preimage of the diagonal $\G_m\subset\G_m^N$
under $\rho$. 
In other words, $\Ga_\w$ is the
maximal subgroup of diagonal transformations of $\A^n$
under which $\w$ is semi-invariant. Let $\chi_w:\Ga_\w\to\G_m$ be the natural character.
It is easy to see that $G_\w=\ker(\rho)=\ker(\chi_w)$, so  
we have a canonical extension of commutative algebraic groups
$$1\to G_\w\to\Ga_\w\rTo{\chi_\w}\G_m\to 1.$$
The choice of coprime degrees $\bfd=(d_1,\ldots,d_n)$ defines
an injective homomorphism
\begin{equation}\label{quasihom-hom-eq}
i_{\bfd}:\G_m\to\Ga_\w:\la\to (\la^{d_1},\ldots,\la^{d_n})
\end{equation} 
such that $\chi_\w\circ i_{\bfd}(\la)=\la^d$.
Thus, the intersection of $G_\w$ with $i_{\bfd}(\G_m)$ is the cyclic subgroup of order $d$ generated by
the {\it exponential grading element}
\begin{equation}\label{J-eq}
J=(\exp(2\pi i q_1),\ldots,\exp(2\pi i q_n))\in G_\w, \text{ where } q_j=d_j/d.
\end{equation}
Note that we have a short exact sequence
$$1\to \Z/d\rTo{i} G_\w\times\G_m\rTo{(\iota,i_{\bfd})}\Ga_\w\to 1$$ 
where $\iota:G_\w\to\Ga_\w$ is the natural embedding and $i(1)=(J,\exp(-2\pi i/d))$.

We will often use the following correspondence between certain subgroups of $\Ga_\w$ and
subgroups of $G_\w$.

\begin{lem}\label{J-lem}
There is a natural bijection between the set of algebraic
subgroups $\Ga\sub\Ga_\w$ containing $i_{\bfd}(\G_m)$
and the set of algebraic subgroups $G\sub G_\w$ containing the element $J$ that
associates with $\Ga$ the intersection $G=\Ga\cap G_\w$ and with $G$ the image
$\Ga=(\iota,i_{\bfd})(G\times\G_m)$.
\end{lem}

\Pf . Let $G\sub G_\w$ be a subgroup containing $J$. Then intersection of 
$\Ga=(\iota,i_{\bfd})(G\times\G_m)$ with $G_\w=\ker(\chi_w)$ consists of the elements
$\iota(g)i_{\bfd}(\la)$, where $g\in G$, $\la\in\G_m$, such that
$\chi_\w(\iota(g)i_{\bfd}(\la))=\la^d=1$. Since in this case $i_{\bfd}(\la)$ belongs to the subgroup
generated by $J$, it follows that $\Ga\cap G_\w=G$.
Now let $\Ga\sub\Ga_\w$ be any subgroup containing $i_{\bfd}(\G_m)$.
Consider the subgroup 
$$\wt{\Ga}=(\iota,i_{\bfd})^{-1}(\Ga)\sub G_\w\times\G_m.$$
Then $\wt{\Ga}$ contains $1\times \G_m$, so we have $\wt{\Ga}=G\times\G_m$.
Furthermore, $\wt{\Ga}$ contains the element $i(1)\in\ker(\iota,i_{\bfd})$, hence $G$ contains $J$.
\ed

Let us fix a commutative algebraic group $\Ga$ equipped with a homomorphism
$\xi:\Ga\to\Ga_\w$ 
such that the composition 
$$\chi=\chi_\w\circ\xi:\Ga\to\G_m$$ is surjective and
the kernel of $\chi$ is finite.
We set
$G=\ker(\chi)\sub\Ga$, so that we have an exact sequence
\begin{equation}\label{G-Ga-ex-seq}
1\to G\to \Ga\rTo{\chi} \G_m\to 1.
\end{equation}

We will work with the dg-category of $\Ga$-equivariant \mf s of $\w$
$$\MF_{\Ga}(\w)=\MF_{\Ga,\chi}(\A^n,\w).$$ 
The corresponding homotopy category $\HMF_{\Ga}(\w)=\HMF(\A^n/\Ga,\w)$ is equivalent to
the equivariant category of singularity of the hypersurface $\w=0$ 
(see \cite[Prop.\ 3.19]{PV-stacks}), i.e., in this case we have
$\DMF(\A^n/\Ga,\w)=\HMF(\A^n/\Ga,\w)$.

Note that if we view $\MF_{\Ga}(\A^n,\w)$ as a $\Z/2$-graded dg-category then
it can be identified with the full subcategory in $\MF_G(\A^n,\w)$.

\subsection{Functors defined by kernels}\label{kernel-sec}

By analogy with Fourier-Mukai transforms on derived categories of coherent sheaves we
can use \mf s of the external product of potentials as kernels representing functors
between the categories of \mf s.

Let $W$ be a function on a smooth FCDRP-stack $X$, which is not a zero divisor
and is semi-invariant with respect to a group 
$H$ acting on $X$ and a surjective character $\chi_W:H\to\G_m$ with finite kernel.
Also let $K$ be a subgroup of $H\times \Ga$ such that
the restrictions of the characters $\chi_W\times 1$ and $1\times\chi_\w$ to $K$ are equal.
Consider the potential
$$(W\oplus \w):=p_1^*W+p_2^*\w$$ 
on $X\times\A^n$, where $p_1:X\times\A^n\to X$ and $p_2:X\times\A^n\to\A^n$ are
the projections. 
Then $W\oplus\w$ is semi-invariant with respect to $K$.
Note that if we take $K$ to be maximal then 
$((X\times\A^n)/K, W\oplus \w)$ can be identified the external product of $(X/H,W)$
with $(\A^n/\Ga,\w)$ (see Definition \ref{ext-tens-def}).
Since the singularity locus of $W\oplus \w$ is a subset of $X_0\times\{0\}$,
Proposition \ref{tensor-prod-lem} implies that for any
$P\in\ov{\DMF}_K(W\oplus \w)$ and $Q\in\ov{\DMF}_{\Ga}(\A^n,-\w)$ the tensor product
$P\ot Q\in\ov{\DMF}_K(p_1^*W)$ belongs to the subcategory
$$\ov{\DMF}_{X_0\times\{0\}}(X\times\A^n/K ,p_1^*W)\sub \ov{\DMF}(X\times\A^n/K,p_1^*W).$$
Assume moreover that the projection $K\to H$ is surjective with the kernel $K_0$.
Then we can apply the functor $p_{1*}^{K_0}$ to $P\ot Q$ (see \eqref{G-mf-push-forward}),
so we obtain a functor
\begin{equation}\label{kernel-fun-eq}
\Phi_P:\ov{\DMF}_{\Ga}(\A^n,-\w)\to\ov{\DMF}_{H}(X,W): Q\mapsto p_{1*}(P\ot Q)^{K_0}.
\end{equation}

The same construction 
works in the case $W=0$ if we assume that $H$ acts trivially on $X$ (see Example  
\ref{G-push-forward-constr}).

\begin{rems}\label{dg-fun-rem}
1. The functor $\Phi_P$ has a natural dg-version. Namely, the category $\ov{\DMF}_H(X,W)$ can be identified with
the full subcategory of the derived category of quasi-matrix factorizations $\DMF^\infty_H(X,W)$ (see
Cor.\ 4.5 and Lem.\ 4.6 of \cite{PV-stacks}). The latter category has a natural dg-version 
$\DMF^{\infty,dg}_H(X,W)$ obtained using the
construction of dg-quotient (see \cite[Thm.\ 4.8]{Keller-dg}). 
Thus, we get a dg-version of $\ov{\DMF}_H(X,W)$ by
taking the corresponding full dg-subcategory of $\DMF^{\infty,dg}_H(X,W)$. Now to obtain a dg-functor
$\Phi_P^{dg}$ inducing $\Phi_P$ we can use the natural push-forward dg-functor for quasi-matrix
factorizations (see \cite[Def.\ 4.8]{PV-stacks}). Namely, $\Phi_P^{dg}$ is defined by the same formula
\eqref{kernel-fun-eq}, where $Q$ is a quasi-matrix factorization of $-w$ and $p_{1*}$ is the naive
push-forward. By the results of \cite{PV-stacks} this construction is compatible with the functor between
derived categories, hence, it factors through the dg-functor defined on $\DMF^{\infty,dg}_H(X,W)$.
Furthermore, all the isomorphisms of functors of the form $\Phi_P$
discussed below will be induced by morphisms between dg-functors.

\noindent
2. In the case $W=0$ we can view the functor $\Phi_P$ as taking values in $D_{H_0}(X)$, where
$H_0=\ker(\chi_W)$.
By Corollary \ref{Groth-cor}, the map on Hochschild homology induced by the dg-functor
$\Phi_P^{dg}$ depends only on the class of $P$ in the Grothendieck group.
In particular, since 
a different choice of representatives for $\widehat{H}/\lan\chi_W\ran$ would change $\Phi_P$, by an
even power of the translation functor on each component of the decomposition
$$D_{H_0}(X)=\bigoplus_{\eta\in \widehat{H_0}}D(X)\ot\eta,$$
the induced map on Hochschild homology does not depend on these choices.
\end{rems}

\begin{lem}\label{stable-functor-lem}
Let  $U\subset \A^n$ be a $\Ga$-invariant linear subspace such that $U\sub S:=\w^{-1}(0)$, 
and let $X$ be a smooth FCDRP-stack $X$ with a trivial $\Ga$-action.
Denote by $k:S \to \A^n$, $k':U\to S$, $\pr:X\times\A^n\to \A^n$,
$p_1:X\times U\to X$ and $p_2:X\times U \to U$ the corresponding
embeddings and projections.
If  $\bar{P}\in \MF_{\Ga}(X\times \A^n,\pr^*\w)$ is such that
$\fC(\bar{P})\simeq (\id_X\times k')_*(C^\bullet)$
in $D_{\Sg}(X\times S/\Ga)$, where $C^\bullet$ is a bounded
complex of vector bundles on $X\times U$
(and $\fC$ is the cokernel functor \eqref{coker-functor-eq}),
then for $\bar{E}\in \ov{\DMF}_{\Ga}(\A^n,-\w)$ we have
a functorial isomorphism
$$\Phi_{\bar{P}}(\bar{E})\simeq \com_G(p_{1*}(C^\bullet\ot p_2^*k^*\bar{E}))\in D_G^b(X).$$
In particular, this holds if $\bar{P}=\wt{C}^\bullet\ot\{\a,\b\}$,
where $\b$ is a regular section with the zero locus $X\times U$ and
$\wt{C}^\bullet$ is a bounded complex of
vector bundles on $X\times\A^n$.
\end{lem}

\Pf . By Lemma \ref{complex-quasi-isom-lem}(ii), we have a canonical isomorphism
$$\bar{P}\ot \pr^*\bar{E}\simeq \left(\id_X\times i)_*(\fC(\bar{P})\ot_{\OO_{X\times S}}
(\id_X\times i)^*\pr^*\bar{E}\right),$$
where $i:S\hra \A^n$ is the natural embedding. Now using the isomorphism
$\fC(\bar{P})\simeq (\id_X\times k')_*(C^\bullet)$ and the projection formula we obtain an isomorphism
$$\bar{P}\ot \pr^*\bar{E}\simeq \left(\id_X\times k)_*(C^\bullet\ot_{\OO_{X\times U}}
p_2^*k^*\bar{E}\right),$$
which gives the required isomorphism after the push-forward to $X$. The fact that the assumptions
are satisfied for $\bar{P}=\wt{C}^\bullet\ot\{\a,\b\}$ with $\b$ regular follows from Lemma \ref{koszul-main-lem}(i).
\ed

\subsection{Generators of categories of \mf s}\label{gener-sec}

Let us define an object $\C^{\st}\in\MF_{\Ga}(\A^n,\w)$, the stabilized residue field, 
similarly to the $\Z/2$-graded considered in \cite[Sec.\ 2.5]{PV-mf}.
Let $T=(\mg/\mg^2)^*$ be the tangent space to $\A^n$ at the origin. 
The projection $\mg\to T^*$
admits a $\Ga$-equivariant splitting $s:T^*\to\mg$, 
which defines an element $\psi\in (T\ot \C[x])^{\Ga}$.
If we choose generators $x_1,\ldots,x_n$ of $\mg\sub\C[x]$ 
as $x_i=s(e_i)$, where $e_1,\ldots,e_n$ is a basis of $T^*$, then
$\psi=\sum_i e_i^*\ot x_i$, where $(e_i^*)$ is the dual basis of $T$.
On the other hand, since the $\Ga$-equivariant map
$$\lan ?,\psi\ran: T^*\ot \C[x]\to\mg$$
is surjective, we can find a $\Ga$-equivariant element 
$\phi\in T^*\ot \C[x]\ot\chi$, such that $\lan\phi,\psi\ran=\w$.
The pair $\phi$, $\psi$ defines a $\Ga$-equivariant Koszul \mf\ 
\begin{equation}\label{C-st-eq}
\C^{\st}=\{\phi,\psi\}
\end{equation}
of $\w$ on $\A^n$ with the respect to the character $\chi$.

Let $\lan\chi\ran\sub \widehat{\Ga}$ denote the subgroup (isomorphic to $\Z$) generated by $\chi$
in the group of characters of $\Ga$.
Let $\chi_1,\ldots,\chi_r$ be a set of representatives for $\lan\chi\ran$-cosets in $\widehat{\Ga}$.

\begin{prop}\label{gen-prop}
The \mf\  $\bigoplus_{i=1}^r \C^{\st}\ot \chi_i$ 
is a generator of the triangulated category $\ov{\DMF}_{\Ga}(\A^n,\w)=\ov{\HMF}_{\Ga}(\A^n,\w)$. 
\end{prop}

\Pf . Under the equivalence of $\DMF_{\Ga}(\A^n,\w)$
with the equivariant singularity category 
of the hypersurface $S=(\w=0)$ (see \cite[Thm.\ 3.14]{PV-stacks})
the stabilized residue field $\C^{\st}$ corresponds to the skyscraper sheaf at the origin. Since the singularity locus of $S$ is the origin, by
\cite[Cor.\ 5.3]{PV-stacks}, the category $\ov{D_{\Sg}(X/\Ga)}$ is equivalent to
$\ov{D_{\Sg}(X/\Ga, 0/\Ga)}$. The latter category is generated by the skyscraper sheaf at the origin
twisted by characters of $\Ga$. Since the twisting by $\chi$ is isomorphic to the square of the
translation functor, it is enough to consider representatives of $\widehat{\Ga}/\lan\chi\ran$
(cf. \cite[sec. 12]{Seidel} for a similar reasoning). 
\ed

\begin{cor}\label{ext-tens-prop}
Let $\w'$ be another quasi-homogeneous potential on $\A^m$ with an isolated singularity
at $0$, semi-invariant with respect to $(\Ga',\chi')$. Let $\Pi\sub\Ga'\times\Ga$
denote the preimage of the diagonal under the homomorphism 
$\chi'\times\chi:\Ga'\times\Ga\to\G_m\times\G_m$. Then
the external tensor product dg-functor
$$\MF_{\Ga'}(\A^m,\w')\otimes\MF_{\Ga}(\A^n,\w)\to 
\MF_{\Pi}(\A^m\times\A^n,\w'\oplus \w)$$
induces an equivalence of perfect derived categories (i.e., of the categories
of compact objects in the corresponding derived categories).
\end{cor}

\Pf . Use Proposition \ref{gen-prop} and argue as in \cite[Sec.\ 6.1]{Dyck}.
\ed

\begin{rem} 
A more general version of the above Corollary is proved independently in the
work of Ballard, Favero and Katzarkov \cite[Prop.\ 6.7]{BFK}.
\end{rem}

\subsection{The diagonal matrix factorization}
\label{diag-sec}

Here we construct a $\Ga_\w$-equivariant version of the 
diagonal matrix factorization representing the identity functor (see \cite{Dyck}, \cite{PV-mf}).
We keep the notation of the previous section.
Consider the $\Ga_\w$-invariant elements 
$\a\in V\ot k[x,y]\ot\chi$ and $\b\in V^*\ot k[x,y]$, given by
$$\a=\sum_j e_j\ot \w_j, \ \ \b=\sum_j e_j^*\ot (y_j-x_j),$$
where $\w_j(x,y)$ are polynomials such that
$$\w(y)-\w(x)=\sum_{j=1}^n(y_j-x_j)\w_j(x,y)$$
(such polynomials exist because $\Ga_\w$ is reductive).
We define a $\Ga_\w$-equivariant \mf\ of $\wt{\w}=-\w(x)+\w(y)$ on $\A^n\times\A^n$ 
(with respect to the diagonal action of $\Ga_\w$ on $\A^n\times\A^n$) by 
$$\De^{\st}_\w=\{\a,\b\}.$$

Now let $X$ be a smooth FCDRP-stack, and let $W$ be a function on $X$, 
as in section \ref{kernel-sec}, i.e., $W$ is semi-invariant with respect to $(H,\chi_W)$ and
we have a subgroup $K\sub H\times \Ga$ such that $\chi_W\times\id|_K=\id\times \chi|_K$ and
the projection $K\to H$ is surjective.
We also assume that either $W$ is not a zero divisor, or $W=0$ and $H$ acts trivially on $X$.

\begin{prop}\label{diag-mf-lem}
(i) The functor 
$$\Phi_{\De^{\st}_\w}:\DMF_{\Ga}(\A^n,\w)\to\DMF_{\Ga}(\A^n,\w)$$
associated with the kernel 
$\De^{\st}_\w\in\DMF_{\Ga}(\A^n\times\A^n,\wt{\w})$  by \eqref{kernel-fun-eq},
is isomorphic to the identity functor.

\noindent
(ii) Let $\pi:\A^n\to pt$ be the projection and 
$$\wt{\De}^{\st}_\w=p_{23}^*\De^{\st}_\w\in\DMF_K(X\times \A^n\times\A^n, 0\oplus(-\w)\oplus \w)$$
be the pull-back of $\De^{\st}_\w$ under the projection 
$p_{23}:X\times \A^n\times\A^n\to\A^n\times\A^n$.
Then the following diagram of functors is commutative up to an isomorphism:
\begin{diagram}
\ov{\DMF}_K(X\times \A^n\times\A^n, W\oplus\w\oplus(-\w))&\rTo{?\ot \wt{\De}^{\st}_\w}&
\ov{\DMF}_{K,X\times\{(0,0)\}}(X\times \A^n\times\A^n, W\oplus 0\oplus 0)\\
\dTo{(\id_X\times\De)^*}&&\dTo{(\id_X\times\pi\times\pi)_*^{K_0}}\\
\ov{\DMF}_{K,X\times\{0\}}(X\times\A^n,W\oplus 0)&\rTo{(\id_X\times\pi)_*^{K_0}}&
\ov{\DMF}_{H}(X,W)
\end{diagram}
where $\De:\A^n\to\A^n\times \A^n$ is the diagonal embedding, $K_0=\ker(K\to H)$,
and we use the push-forward functors combined with taking $K_0$-invariants
as in \eqref{G-mf-push-forward}.
\end{prop}

\Pf . (i) This can be derived from Proposition \ref{gen-prop} similarly to 
the non-equivariant case considered in \cite{Dyck} (see also \cite{PV-mf}).

\noindent
(ii) The proof is based on the fact that 
$\wt{\De}^{\st}_\w$ is a regular Koszul \mf\ with the zero locus 
$X\times\De(\A^n)\sub X\times\A^n\times\A^n$.
Assume first that $W$ is not a zero divisor. Then by Proposition \ref{koszul-reg-prop}(i), we have
$$\bar{P}\ot\wt{\De}^{\st}_\w\simeq (\id_X\times\De)_*(\id_X\times\De)^*\bar{P},$$
which implies the result since 
$$(\id_X\times\pi\times\pi)_*\circ (\id_X\times\De)_*=(\id_X\times\pi)_*.$$
In the case $W=0$ the proof is similar, but we use Proposition \ref{koszul-reg-prop}(ii) instead.
\ed

\begin{cor} The dg-category $\MF_{\Ga}(\A^n,\w)$ is dg-Morita equivalent to
a smooth proper dg-algebra.
\end{cor}

\Pf . This follows from the existence of a compact generator (see Proposition \ref{gen-prop}), 
from Corollary \ref{ext-tens-prop} 
and from the fact that the diagonal bimodule is represented by a \mf\
(cf. \cite[Sec.\ 7]{Dyck}).
 
\subsection{Hochschild homology and the Chern character for dg-categories}\label{dg-Hoch-sec}

Below we will use the formalism of \cite[Sec.\ 1]{PV-mf} (see also 
\cite{Keller-derived} and \cite{Toen} and Appendix to this paper).

Let $\CC$ be a dg-category over a field $k$. 
We denote by $D(\CC)$ the derived category of right $\CC$-modules,
by $\Per(\CC)\sub D(\CC)$ the perfect derived category, and by $\Per_{dg}(\CC)$ the dg-category
of homotopically finitely presented right $\CC$-modules (see \cite[Sec.\ 7]{Toen}).
The Hochschild homology of $\CC$ is given by
$$HH_*(\CC)=\Tr_{\CC}(\De_\CC),$$
where $\De_{\CC}$ is the diagonal $\CC-\CC$-bimodule $E\ot F^{\vee}\mapsto\Hom_{\CC}(F,E)$
and 
\begin{equation}\label{Tr-functor}
\Tr_{\CC}:D(\CC^{op}\ot\CC)\to D(k)
\end{equation}
is the trace functor given by the derived tensor product with $\De_{\CC}$.

As in \cite[Sec.\ 1.2]{PV-mf} we consider only dg-categories $\CC$
such that the $\CC-\CC$-bimodule $\De_\CC$ is perfect, the complexes $\Hom_\CC(A,B)$
for $A,B\in\CC$ have finite dimensional cohomology, and the derived category $D(\CC)$
has a compact generator. 
Such dg-categories are dg Morita equivalent to homologically smooth
and proper dg-algebras and can be characterized by the condition that $\Per_{dg}(\CC)$
is {\it saturated}, i.e., proper, smooth and triangulated (see \cite[Sec.\ 2.2]{TV}). 
 
Any dg-functor $F:\Per_{dg}(\CC)\to\Per_{dg}(\DD)$ between dg-categories of the above type
comes from a kernel in $\Per_{dg}(\CC^{op}\ot\DD)$ (see \cite[Sec.\ 2.2.]{TV}, \cite[Sec.\ 5.4]{Toen-dg}) and
induces a map $F_*:HH_*(\CC)\to HH_*(\DD)$. This map can be defined in several
equivalent ways. We will use a slight modification of the construction given in \cite[Sec.\ 1.2]{PV-mf}
(see Appendix for the details).

First, consider the induced functor
$$F^{(2)}:\Per(\CC^{op}\ot\CC)\to\Per(\DD^{op}\ot\DD)$$
that sends the representable module $h_{C_1^\vee\ot C_2}$
to $h_{F(C_1)^\vee\ot F(C_2)}$.

There is a canonical morphism of functors 
\begin{equation}\label{tr-phi2-abs-eq}
\Tr_{\CC}\to\Tr_{\DD}\circ F^{(2)}
\end{equation}
and a canonical morphism
\be\label{phi2-Delta-abs-eq}
F^{(2)}(\De_{\CC})\to\De_{\DD}
\end{equation}
in $\Per(\DD^{op}\ot\DD)$.

Now the map $F_*$ is defined as the composition
$$\Tr_{\CC}(\De_{\CC})\to \Tr_{\DD}F^{(2)}(\De_{\CC})\to \Tr_{\DD}(\De_{\DD}),$$
where the first arrow is induced by \eqref{tr-phi2-abs-eq} and the second is induced by
\eqref{phi2-Delta-abs-eq}. 
The map $F_*$ depends only on the quasi-isomorphism class of the $\CC-\DD$-bimodule
producing $F$ (in fact, only on its class in the Grothendieck group).

For an object $E\in\Per_{dg}(\CC)$ we define its Chern character 
using the functor $\unit_E: \Per_{dg}(k)\to \Per_{dg}(\CC)$ sending $k$ to $E$ as follows
$$\ch(E)=(\unit_E)_*(1)\in\HH_0(\CC).$$

The maps $F_*$ are compatible with the composition (see \cite[Lem.\ 1.2.1]{PV-mf}), which implies
the functoriality of the Chern character
$$\ch(F(E))=F_*(\ch(E))$$
for $E\in\Per(\CC)$.

Another application of the functoriality is the construction of the canonical perfect bilinear pairing
\begin{equation}\label{abs-can-pairing-eq}
(\cdot,\cdot)_{\CC}:HH_*(\CC^{op})\ot HH_*(\CC)\to k.
\end{equation}
Namely, it is induced by the dg-version of the trace functor \eqref{Tr-functor} 
restricted to perfect bimodules (see \cite[Sec.\ 1.2]{Shk} and \cite[Sec.\ 1.2]{PV-mf}).
Using this pairing one can express the map on the Hochschild homology induced by a functor
in terms of the Chern character of the corresponding kernel.

\begin{lem}\label{kernel-pairing-lem}
Let $F:\Per_{dg}(\CC)\to \Per_{dg}(\DD)$ be the dg-functor associated with a kernel
$K\in\Per_{dg}(\CC^{op}\ot\DD)$. Then the induced map $F_*:HH_*(\CC)\to HH_*(\DD)$ is given by
$$a\mapsto \tr_{12}(a \ot \ch(K)),$$
where 
$$\tr_{12}: HH_*(\CC)\ot HH_*(\CC^{op})\ot HH_*(\DD)\to HH_*(\DD)$$
sends $a\ot b\ot c$ to $(a,b)_{\CC^{op}} c$, where 
$(\cdot,\cdot)_{\CC^{op}}$ 
is the canonical pairing for $\CC^{op}$.
\end{lem}

\Pf . For $M\in\Per_{dg}(\CC)$ we have
$$F(M)=M\ot_{\CC} K\simeq \De_{\CC^{op}}\ot_{\CC\ot\CC^{op}} (M\boxtimes K)
=(\Tr_{\CC^{op}}\ot\Id)(M\boxtimes K).$$
Hence, $F_*$ is the composition of the map
$$HH_*(\CC)\to HH_*(\CC)\ot HH_*(\CC^{op}\ot\DD)\simeq HH_*(\CC)\ot HH_*(\CC^{op})\ot HH_*(\DD)$$
induced by $M\mapsto M\boxtimes K$, with $\tr_{12}$. The assertion follows immediately from this.
\ed

\subsection{Hochschild homology of $\MF_{\Ga}(\w)$}\label{Hoch-sec}

Let  $\AA_\w=\C[x]/(\partial_1\w,\ldots,\pa_n\w)$
be the Milnor ring of the isolated singularity $\w$, and consider the space
$$H(\w)=\AA_\w\ot dx,$$
where $dx=dx_1\we\ldots\we dx_n$, equipped with the action of $\Ga$, induced
by its action on $\A^n$. Note that since the partial derivatives of $\w$ form a regular sequence,
$H(\w)$ is identified with the only nonzero cohomology (at the $n$th term) of the complex
$\Om^\bullet_{\C[x]/\C}$ with the differential $\a\mapsto \a\we d\w$.

Recall that 
the Hochschild homology of the $\Z/2$-dg-category $\MF(\w)$ is isomorphic to 
the space $H(\w)$ in degree $n\mod (2)$ (see \cite{Dyck}). 
The Hochschild homology of the $\Z/2$-dg-category $\MF_G(\w)$ is given by 
\begin{equation}\label{G-hoch-isom}
HH_*(\MF_G(\w))\simeq \bigoplus_{\ga\in G}H(\w_{\ga})^G
\end{equation}
where $\w_{\ga}$ is the restriction of $\w$ to the subspace $(\A^n)^{\ga}\sub\A^n$ (see \cite{PV-mf}).
Here we establish a version of this isomorphism for the $\Z$-graded dg-category 
$\MF_{\Ga}(\w)$.

Let $\Ga^{(2)}\sub \Ga\times \Ga$ denote the preimage of the diagonal under
the homomorphism $\chi\times\chi:\Ga\times \Ga\to\G_m\times\G_m$.
Let $\chi^{(2)}:\Ga^{(2)}\to\G_m$ be the character induced by $\chi$ and
by one of the projections $\Ga^{(2)}\to \Ga$.
Let us consider the dg-category $\MF_{\Ga^{(2)}}(\A^n\times\A^n,\wt{\w})$,
where $\wt{\w}(x,y)=-\w(x)+\w(y)$. Similar to the $\Z/2$-graded case (see \cite[sec. 6.1]{Dyck})
we can interpret the corresponding perfect derived category as the category
of dg-functors $\MF_{\Ga}(\A^n)\to\MF_{\Ga}(\A^n)$.
Namely, to a kernel $K\in\MF_{\Ga^{(2)}}(\A^n\times\A^n,\wt{\w})$ we associate
the dg-functor \eqref{kernel-fun-eq}
$$\Phi_K:\MF^{\infty}_{\Ga}(\A^n,\w)\to\MF_{\Ga}^{\infty}(\A^n,\w):\bar{E}\mapsto 
p_{2*}(p_1^*\bar{E}\ot K)^{G\times \{1\}}.$$
Note that here the invariants are taken
with respect to the action of the group $G$ on the
first factor of the product $\A^n\times\A^n$. 
Since $\Per_{dg}(\MF_{\Ga}(\A^n,\w))$ is saturated,
Corollary \ref{ext-tens-prop} implies that
every dg-functor from this category to itself is represented by a \mf\ of $\wt{\w}$.

Now we are ready to compute the Hochschild homology of $\MF_{\Ga}(\w)$.
Let $\widehat{G}$ be the dual group to $G$. The exact sequence \eqref{G-Ga-ex-seq} 
induces an exact sequence
$$0\to\Z\rTo{n\mapsto\chi^n}\widehat{\Ga}\to \widehat{G}\to 0.$$
Note that we have the natural action of $\widehat{\Ga}$ on the category $\DMF_{\Ga}(\w)$ 
given by tensor multiplication with $1$-dimensional representations of $\Ga$.
Furthermore, by definition of
the triangulated structure on $\DMF_{\Ga}(\w)$ we have
$$\bar{E}\ot\chi\simeq \bar{E}[2].$$
Hence, the induced action of $\widehat{\Ga}$ on $HH_*(\MF_{\Ga}(\w))$
factors through an action of $\widehat{G}$. 
In other words, $HH_*(\MF_{\Ga}(\w))$ has a natural structure of $R$-module for
$R=\C[\widehat{G}]$.

\begin{thm}\label{hoch-prop}
(i) The Hochschild homology of the dg-category $\MF_{\Ga}(\w)$ is given by 
\begin{equation}\label{graded-G-hoch-isom}
HH_*(\MF_{\Ga}(\w))\simeq\bigoplus_{\ga\in G} H(\w_\ga)^G,
\end{equation}
where $\w_\ga$ is the restriction of $\w$ to the subspace of $\ga$-invariants $(\A^n)^\ga$,
with the $\Z$-grading is given by
$$H(\w_\ga)^G=\bigoplus_{i\in\Z} H(\w_\ga)^{G}_{\chi^{-i}}[n_\ga-2i],$$
where
$n_\ga=\dim(\A^n)^\ga$. We have an isomorphism of $\Z/2$-graded spaces 
$$HH_*(\MF_{\Ga}(\w))\simeq HH_*(\MF_{G}(\w))$$ 
identifying the decompositions
\eqref{graded-G-hoch-isom} and \eqref{G-hoch-isom}.

\noindent
(ii) The decomposition of $HH_*(\MF_{\Ga}(\w))$ 
into $\ga$-isotypical subspaces (where $\ga\in G$ is viewed as a character of $\widehat{G}$) coincides
with the decompositions \eqref{graded-G-hoch-isom} and
\eqref{G-hoch-isom}.

\noindent
(iii) Let $\Ga'\sub \Ga$ be a subgroup such that the restriction of $\chi$ to $\Ga'$ is surjective,
and let $G'=\Ga'\cap G_\w$ be the corresponding subgroup of $G$.
Let
$$\Res^G_{G'}:HH_*(\MF_{\Ga}(\w))\to HH_*(\MF_{\Ga'}(\w))$$
be the map induced by the forgetful functor 
\begin{equation}\label{forget-functor-Ga-eq}
\Phi:\MF_{\Ga}(\w)\to\MF_{\Ga'}(\w).
\end{equation}
Then the restriction of $\Res^G_{G'}$ to the component of the decomposition 
\eqref{graded-G-hoch-isom} corresponding to an element $\ga\in G$ is equal to zero if
$\ga\not\in G'$ and 
to the canonical embedding $H(\w_{\ga})^G\to H(\w_{\ga})^{G'}$
if $\ga\in G'$.
\end{thm}

\Pf .
(i) First, let us check that the trace functor (see Section \ref{dg-Hoch-sec})
$$\Tr:\MF_{\Ga^{(2)}}(\A^n\times\A^n, \wt{\w})\to \Com_f(\C-\mod)$$
associates with a \mf\ $\bar{E}$ of $\wt{\w}$ the $\Ga$-invariants in the global sections of
the restriction of the complex $\com(\bar{E})$ 
(see \eqref{mf-Z-gr-com-eq})
to the diagonal $(y=x)$ in $\A^n\times\A^n$. Indeed, this follows from the isomorphism
$$H^0(\A^n,\com(\bar{E}^*\boxtimes\bar{E'})|_{y=x})^{\Ga}\simeq 
H^0(\A^n,\com(\bar{E}^*\ot\bar{E'}))^{\Ga}\simeq \Hom_{\HMF_{\Ga}}(\bar{E},\bar{E'})$$
for $\bar{E},\bar{E'}\in\MF_{\Ga}(\A^n,\w)$ (see Lemma \ref{Hom-mf-lem}).

Next, we observe that the identity functor on $\MF_{\Ga}(\A^n,\w)$ is represented by
\begin{equation}\label{tilde-De-eq}
\De^{\st}_G:=\bigoplus_{\ga\in G}(\id\times \ga)^*\De^{\st}\in \MF_{\Ga^{(2)}}(\wt{\w}),
\end{equation}
where $\De^{\st}=\De^{\st}_\w$,
with the $\Ga^{(2)}$-equivariant structure induced by the
$\Ga$-equivariant structure on $\De^{\st}$  via the  diagonal embedding of $\Ga$ into $\Ga^{(2)}$.
This follows immediately from Proposition \ref{diag-mf-lem}(i) because of the exact sequence
$1\to G\to \Ga^{(2)}\to \Ga\to 1$.

Finally, to compute the Hochschild homology we have to apply the functor $\Tr$ to $\De^{\st}_G$. Let 
us show that
\begin{equation}\label{Z-Hoch-com-eq}
\com(\De^{\st}_G)|_{y=x}\simeq \bigoplus_{\ga\in G}\bigoplus_{i\in\Z} (\KK_{\w_\ga}\ot\chi^{-i})[2i],
\end{equation}
where for any potential $W$ on $\A^m$, equivariant with respect to a character $\chi:\Ga\to\G_m$,
we denote by $\KK_W$ is the complex
\begin{equation}\label{K-w-eq}
\KK_W=[\Om_{\A^m}^0\to \Om_{\A^m}^1\ot\chi\to \ldots\to\Om_{\A^m}^m\ot\chi^m]
\end{equation}
placed in degrees $[0,m]$ with the differential given by $dW\we ?$. 
The summand of \eqref{Z-Hoch-com-eq}
corresponding to $\ga=1$ is
\begin{equation}\label{Z-w-Hoch-com-eq}
\com(\De^{\st})|_{y=x}\simeq \bigoplus_{i\in\Z} (\KK_{\w}\ot\chi^{-i})[2i].
\end{equation}
Using the identifications
$$(\De^{\st})_0|_{y=x}\simeq \bigoplus_{i\ge 0}\Om_{\A^n}^{2i}\cdot\chi^i,\ \text{ and }\ 
(\De^{\st})_1|_{y=x}\simeq \bigoplus_{i\ge 0}\Om_{\A^n}^{2i+1}\cdot\chi^i,$$
the complex $\com(\De^{\st})|_{y=x}$ can be presented as
\begin{diagram}
\Om^0\chi^{-1}&&&&\Om^0&&&&\Om^0\chi\\
&\rdTo{}&&&&\rdTo{}\\
&&\Om^1&&&&\Om^1\chi\\
&&&\rdTo{}&&&&\rdTo{}\\
\Om^2&&&&\Om^2\chi&&&&\Om^2\chi^2\\
&\rdTo{}&&&&\rdTo{}\\
\ldots&&&&\ldots&&&&\ldots
\end{diagram}
where all the arrows are given by $d\w\we ?$. This immediately gives the isomorphism
\eqref{Z-w-Hoch-com-eq}. The equality of the summands in \eqref{Z-Hoch-com-eq} 
corresponding to $g\neq 1$ is verified similarly 
by explicitly calculating $(\id\times \ga)^*\De^{\st}|_{y=x}\simeq \De^{\st}|_{y=\ga x}$
as in \cite[sec. 2.5]{PV-mf}.

Since each $\w_\ga$ is an isolated singularity (see \cite[Lem.\ 2.5.3(i)]{PV-mf}),
the cohomology of the complex $\KK_{\w_\ga}$ \eqref{K-w-eq}
is isomorphic to $H(\w_{\ga})\ot \chi^{n_{\ga}}$ concentrated in (cohomological) degree $n_{\ga}$.
Hence, using \eqref{Z-Hoch-com-eq} we see that the cohomology of
$\com(\De^{\st}_G)$
is isomorphic to 
$$\bigoplus_{\ga\in G}\bigoplus_{i\in\Z} (H(\w_\ga)\ot \chi^{n_\ga-i})[2i-n_\ga].$$ 
Passing to $\Ga$-invariants and 
substituting $i\mapsto n_\ga-i$ we get the result. The last assertion follows from 
the computation in \cite[Sec.\ 2.5]{PV-mf}.

\noindent (ii)
In general, if $\a:\CC\to\CC$ is an autoequivalence of a dg-category,
the induced automorphism $\a_*$ of $HH_*(\CC)$ is defined as follows (see section \ref{dg-Hoch-sec}).
We have an induced equivalence of $\a^{(2)}$ of $\CC^{op}\ot\CC$ and
natural isomorphisms 
$$\psi:\Tr_{\CC}\circ (\a^{(2)})\simeq\Tr_{\CC},$$
\begin{equation}\label{alpha-De-eq}
\phi:(\a^{(2)})(\De_{\CC})\simeq \De_{\CC}
\end{equation}
that induce an automorphism
$$\a_*:HH_*(\CC)=\Tr_{\CC}(\De_{\CC})\rTo{\phi}\Tr_{\CC}((\a^{(2)})(\De_{\CC}))\rTo{\psi}\Tr_{\CC}(\De_{\CC})=
HH_*(\CC).$$
Now let us specialize to the case of $\CC=\MF_{\Ga}(\w)$ and the autoequivalence $\a$ given by the tensoring with a character $\eta$ of $\Ga$. Under the identification of the perfect derived category
$\Per(\CC^{op}\ot\CC)$ with 
$\ov{\HMF}_{\Ga^{(2)}}(\A^n\times\A^n,\wt{\w})$ (see Corollary \ref{ext-tens-prop}) the functor $\a^{(2)}$ corresponds to tensoring with the
character $\eta^{-1}\times\eta|_{\Ga^{(2)}}$. Recall that by Proposition \ref{diag-mf-lem},
the kernel $\De_{\CC}$ representing the identity functor in
this case is $\De^{\st}_G=\bigoplus_{\ga\in G}(\id\times \ga)^*\De^{\st}$.
It is easy to check that the isomorphism $\phi$
in this case is given by the multiplication by $\eta(\ga)^{-1}$ on the component
$(\id\times \ga)^*\De^{\st}$. Since, the automorphism $\a_*$ is obtained by the restriction of
$\phi$ to the diagonal in $\A^n\times\A^n$, we obtain that $\a_*$ acts as $\eta(\ga)^{-1}$
on the component of $HH_*(\MF_{\Ga}(\w))$ coming from the term 
$$(\id\times \ga)^*\De^{\st}|_{y=x}\simeq\De^{\st}|_{y=\ga x}$$ 
of $\De^{\st}_G$, which is exactly 
the term corresponding to $\ga$ in the decomposition \eqref{graded-G-hoch-isom}.

\noindent
(iii) Let us apply the general construction of Section \ref{dg-Hoch-sec} (see also Appendix) to the forgetful functor
$\Phi:\CC\to \DD$, where
$\CC=\MF_{\Ga}(\w)$ and $\DD=\MF_{\Ga'}(\w)$.
By Corollary \ref{ext-tens-prop}, we have natural equivalences 
\be\label{square-equivalences}
\Per(\CC\ot\CC^{op})\simeq \ov{\HMF}_{\Ga^{(2)}}(\w\oplus(-\w)),\ \ \ 
\Per(\DD\ot\DD^{op})\simeq \ov{\HMF}_{(\Ga')^{(2)}}(\w\oplus(-\w)).
\end{equation}
Under the equivalences \eqref{square-equivalences} the induced functor
$\Phi^{(2)}:\Per(\CC\ot\CC^{op})\to\Per(\DD\ot\DD^{op})$ is identified with
the forgetful functor corresponding to the restriction from $\Ga^{(2)}$ to $(\Ga')^{(2)}$.
The map on Hochschild homology
$\Phi_*$ is given by the composition
$$\Tr_{\CC}(\De_{\CC})\to \Tr_{\DD}\Phi^{(2)}(\De_{\CC})\to \Tr_{\DD}(\De_{\DD})$$
via the natural transformation
\begin{equation}\label{tr-phi2-eq}
\Tr_{\CC}\to\Tr_{\DD}\circ\Phi^{(2)}
\end{equation}
and the natural morphism
\be\label{phi2-Delta-eq}
\Phi^{(2)}(\De_{\CC})\to\De_{\DD}.
\end{equation}
As we have seen in the beginning of the proof of (i), the functor
$\Tr_{\CC}$ (resp., $\Tr_{\DD}$) can be identified under the equivalences
\eqref{square-equivalences} with the restriction to the diagonal
followed by taking $\Ga$-invariants (resp., $\Ga'$-invariants). Since
the morphism \eqref{tr-phi2-eq} on objects of the form $E_1\ot E_2^\vee\in\CC\ot\CC^{op}$
corresponds to the natural embedding
$\Hom_{\CC}(E_2,E_1)\to\Hom_{\DD}(\Phi(E_2),\Phi(E_1))$,
we see that the morphism
\eqref{tr-phi2-eq} corresponds to the natural embedding of $\Ga$-invariants into
$\Ga'$-invariants. 
On the other hand, using equivalences \eqref{square-equivalences} and \eqref{tilde-De-eq}
we have
$$\De_{\CC}\simeq \bigoplus_{\ga\in G}(\id\times \ga)^*\De^{\st}_\w, \ \ 
\De_{\DD}\simeq \bigoplus_{\ga'\in G'}(\id\times \ga')^*\De^{\st}_\w.$$
We claim that the morphism \eqref{phi2-Delta-eq} corresponds
under these identifications to the natural projection (identity on summands corresponding
to elements $\ga\in G'$, and zero on all the other summands). This will immediately imply the desired
statement.
By choosing a $\Ga$-equivariant
generator of the category of (non-equivariant) \mf s of $\w$, we can reduce the task to a
similar question for a dg-algebra $A$ with an action of the group
$G$. The analog of a representation of the identity functor via the kernel \eqref{tilde-De-eq}
is the functorial isomorphism of $A[G]$-modules
\begin{equation}\label{A[G]-M-eq}
M\wt{\to} (A[G]\ot_A M)^G: m\mapsto\sum_{\ga\in G}\ga^{-1}\ot \ga m,
\end{equation}
where $M$ is any module over $A[G]$, the twisted group algebra of $G$.
The $G$-invariants on the right-hand side of \eqref{A[G]-M-eq} 
are taken with respect to the action of $G$ on $A[G]\ot_A M$ given
by $\ga\cdot (x\ot m)=x\ga^{-1}\ot \ga m$, while the 
$A[G]$-structure is induced by the the left action of $A[G]$ on itself. 
The morphism \eqref{phi2-Delta-eq} is obtained via the natural isomorphism
$$\Phi^{(2)}(\De_{\CC})\simeq \Phi\circ \Psi,$$
where $\Psi:A[G']-\mod\to A[G]-\mod$ is the right adjoint functor to 
the restriction functor $\Phi:A[G]-\mod\to A[G']-\mod$. For an $A[G']$-module $N$
we have a functorial isomorphism of $A[G]$-modules
$$\Psi(N)\simeq (A[G]\ot_A N)^{G'}$$
(with the same conventions as in \eqref{A[G]-M-eq}). Namely, for any $A[G]$-module $M$
the isomorphism
$$\Hom_{A[G]}(M, (A[G]\ot_A N)^{G'})\wt{\to} \Hom_{A[G']}(M,N)$$
associates with $f:M\to (A[G]\ot_A N)^{G'}$ its composition with the morphism of $A[G']$-modules
$$\pi:(A[G]\ot_A N)^{G'}\to N,$$ 
which itself is the composition of the embedding into $A[G]\ot_A N$ with the projection to $N$
sending $1\ot n$ to $n$ and $\ga\ot n$ to zero for $\ga\neq 0$.
Thus, the above morphism $\pi$ can be identified with the adjunction map $(\Phi\circ\Psi)(N)\to N$.
Now our claim follows from the commutativity of the triangle
\begin{diagram}
(A[G]\ot_A N)^{G'} &\rTo{}& (A[G']\ot_A N)^{G'}\\
&\rdTo{\pi}&\dTo{}\\
&& N
\end{diagram}
where the vertical arrow is the inverse of the isomorphism \eqref{A[G]-M-eq} for the group $G'$ and the $A[G']$-module $N$,
while the horizontal arrow
is induced by the projection $A[G]\to A[G']$ sending $[\ga]$ to zero for 
all $\ga\in G\setminus G'$.
\ed

Let us denote by
\begin{equation}\label{w-hoch-decomp-eq}
\HH(\w)=HH_*(\MF_{\Ga}(\w))=\bigoplus_{\ga\in G} H(\w_\ga)^G
\end{equation}
the Hochschild homology space computed in Theorem \ref{hoch-prop}(i).

\begin{cor}\label{Kunneth-cor}
Let $\w'$ be a quasi-homogeneous potential on $\A^m$ with an isolated singularity
at $0$, semi-invariant with respect to the same group $\Ga$ and the same character
$\chi:\Ga\to\G_m$. Then the tensor product functor induces an isomorphism
\begin{equation}\label{Kunneth-mf-eq}
\HH(\w')\ot_R \HH(\w)\to HH_*(\MF_{\Ga}(\A^m\times\A^n,\w'\oplus \w))^{G\times G},
\end{equation}
where $R=\C[\widehat{G}]$ and the $G\times G$-action on the Hochschild homology is induced by the action of
$G\times G$ on $\A^m\times\A^n$.
\end{cor}

\Pf .
By Corollary \ref{ext-tens-prop} together with the K\"unneth formula for Hochschild homology
(see \cite[Sec.\ 2.4]{Shk}, \cite[1.1.4]{PV-mf}), the tensor product induces an isomorphism
$$\HH(\w')\ot_\C \HH(\w)\to HH_*\left(\MF_{\Ga^{(2)}}(\A^m\times\A^n,\w'\oplus\w)\right).$$
Now by Theorem \ref{hoch-prop}(iii), the forgetful functor
$$\MF_{\Ga^{(2)}}(\A^m\times\A^n,\w'\oplus\w)\to \MF_{\Ga}(\A^m\times\A^n,\w'\oplus\w)$$
induces a map
$$HH_*(\MF_{\Ga^{(2)}}(\A^m\times\A^n,\w'\oplus\w))\to 
HH_*(\MF_{\Ga}(\A^m\times\A^n,\w'\oplus\w))^{G\times G}$$
given by the projection to the components associated with the image of the diagonal embedding
$G\to G\times G$ (where we use the natural
$G\times G$-grading on the $\Ga^{(2)}$-equivariant Hochschild
homology). Similarly, the map
$$\HH(\w')\ot_\C \HH(\w)\to \HH(\w')\ot_R\HH(\w)$$
can be identified with the projection to 
$\bigoplus_{\ga\in G}e_\ga\HH(\w')\ot e_\ga\HH(\w)$, so the assertion follows.
\ed

\begin{ex}\label{Hoch-ex}
Consider $\Ga=\G_m$ embedded naturally into
$\Ga_\w$ via $\la\mapsto (\la^{d_1},\ldots,\la^{d_n})$, so that the induced character
$\chi:\G_m\to\G_m$ is $\la\mapsto\la^d$, and $G=\Z/d$.
Then
$$HH_*(\MF_{\G_m}(\w))\simeq\bigoplus_{j\in\Z/d, i\in\Z} H(\w_j)_{di}[n_j-2i],
$$
where $n_j$ is the number of $s\in[1,n]$ such that $d|jd_s$ (we use the $\Z$-grading
of $H(\w_j)$ induced by the $\Z$-grading of the variables $x_i$).
\end{ex}

\begin{rem} In the case $n=0$ and $\w=0$ the category $\DMF_{\Ga}(0)$ is (noncanonically)
equivalent to $D(G-\mod)$. Similarly, $\DMF_{\Ga^{(2)}}(0)$ is equivalent
to $D(G^2-\mod)$. The diagonal object $\De^{\st}_G$ in this case can be identified with
$$\De^{\st}_G\simeq \Ind_G^{G^2}\unit_G\simeq \oplus_{\eta\in \widehat{G}}\eta\boxtimes\eta^{-1},$$
where $\unit_G$ is the trivial representation of $G$.
\end{rem}

\subsection{The Chern character and the canonical bilinear pairing}\label{bilinear-sec}

Recall that there is a canonical pairing on Hochschild homology
\begin{equation}\label{can-pairing-eq}
(\cdot,\cdot):HH_*(\MF_{\Ga}(\w)^{op})\ot HH_*(\MF_{\Ga}(\w))\to\C
\end{equation}
(see \eqref{abs-can-pairing-eq}).
Under the identification of $\Z/2$-graded spaces $HH_*(\MF_{\Ga}(\w))\simeq HH_*(\MF_G(\w))$
this pairing coincides with
the nondegenerate bilinear pairing on Hochschild homology of $\MF_G(\w)$
calculated in \cite{PV-mf} (see Theorem  \ref{hoch-prop}).
Note that the duality \eqref{duality-eq} 
gives a natural equivalence $\MF_{\Ga}(\w)^{op}\simeq\MF_{\Ga}(-\w)$,
so that the canonical pairing is induced by the dg-functor
$$\MF_{\Ga}(-\w)\ot \MF_{\Ga}(\w)\to\Com_f(\C-\mod):(E,F)\mapsto \com(E\ot F)^{\Ga},$$
where $\Com_f(\C-\mod)$ is the category of complexes of $\C$-vector spaces with finite-dimensional
total cohomology.

\begin{defi} Let $R=\C[\widehat{G}]$ be the group algebra of the dual group $\widehat{G}$. 
We define an $R$-bilinear version of the canonical pairing \eqref{can-pairing-eq}
\begin{equation}\label{R-bil-form-eq}
(\cdot,\cdot)^R:\HH(-\w)\ot_R \HH(\w)\to R
\end{equation}
as the map on Hochschild homology induced by the tensor product functor
$$\MF_{\Ga}(-\w)\ot\MF_{\Ga}(\w)\to\Com_f(G-\mod):(E,F)\mapsto \com_G(E\ot F)
$$
where $\com_G$ is given by \eqref{com-G-eq}, $\Com_f(G-\mod)$ is the category of complexes of
$G$-modules with finite-dimensional total cohomology.
\end{defi}

We have an isomorphism
$$\com(E\ot F)^{\Ga}\simeq\com_G(E\ot F)^G,$$
which implies that
\begin{equation}\label{tr-form-eq}
(\cdot,\cdot)=\tr\circ (\cdot,\cdot)^R,
\end{equation}
where $\tr:R\to\C$ is given by 
\begin{equation}\label{tr-R-eq}
\tr(\sum_{\eta\in \widehat{G}}c_\eta\cdot [\eta])=c_1.
\end{equation}

\begin{prop}\label{Casimir-prop}
The pairing $(\cdot,\cdot)^R$ is perfect and the corresponding Casimir element
$T_{\w}$ is given by
$$T_{\w}=\frac{1}{|G|}\cdot\sum_{g\in G}(\id\times g)^*\ch(\De^{\st}_\w)\in HH_*(\MF_{\Ga}(\A^n\times\A^n,\wt{\w}))^{G\times G}\simeq
\HH(-\w)\ot_R \HH(\w),$$
where the last isomorphism comes from Corollary \ref{Kunneth-cor}.
\end{prop}

\Pf . 
Since the functor 
$\Phi_{\De^{\st}_\w}$ given by the kernel $\De^{\st}_\w$
 is
isomorphic to the identity functor (see Proposition \ref{diag-mf-lem}(i)), 
the composition 
$$\MF_{\Ga}(\w)\rTo{(\ot p_{23}^*\De^{\st}_\w)\circ p_1^*}
\MF_{\Ga}(\w\oplus(-\w)\oplus \w)\rTo{\De_{12}^*} \MF_{\Ga}(\w)
$$
is isomorphic to the identity. Hence, the composition of the induced maps on Hochschild homology
$$\HH(\w)\rTo{\a} HH_*(\MF_{\Ga}(\w\oplus(-\w)\oplus \w))\rTo{\b}
\HH(\w)$$
is equal to the identity.
Since $\b$ is equivariant with respect to the $G$-action on 
$HH_*(\MF_{\Ga}(\w\oplus(-\w)\oplus \w))$ given by the embedding
of $1\times 1\times G\sub G\times G\times G$ 
and the trivial $G$-action on $\HH(\w)$, we 
have 
$\b\circ
(|G|^{-1}\cdot\sum_{g\in  G}(\id\times\id\times g)^*\circ\a) = \id.
$ 
But the element
$$T_\w=\frac{1}{|G|}\cdot\sum_{g\in G}(\id\times g)^*\ch(\De^{\st}_\w)\in HH_*(\MF_{\Ga}((-\w)\oplus \w))$$
is invariant under $G\times G$ (since $\De^{\st}_\w$ is equivariant with respect to the diagonal
action of $G$). Therefore, we have
$$\b(x\ot T_\w)=x$$
for any $x\in\HH(\w)$. It remains to observe that $x\ot T_\w$ belongs to the subspace
$$HH_*(\MF_{\Ga}(\w\oplus(-\w)\oplus \w))^{G\times G\times G}\simeq\HH(\w)\ot_R\HH(-\w)\ot_R\HH(\w),$$
and the restriction of $\b$ to this subspace is equal to $(\cdot,\cdot)^R\ot\id$.

By definition, the pairing $(\cdot,\cdot)^R$ is obtained as the restriction to the space of
$G\times G$-invariants of the map
$$HH_*(\MF_{\Ga}((-\w)\oplus\w))\to R$$
 induced by the functor of restricting to the diagonal.

\ed

\begin{ex}\label{k-st-ex}
Let $i:\{0\}\hra \A^n$ denote the natural embedding, and let
$$\kappa:HH_*(\MF_{\Ga}(\w))\to R$$
be the $R$-linear functional induced by the restriction functor
$$\MF_{\Ga}(\w)\to\Com_f(G-\mod):\bar{E}\mapsto \com_G(i^*\bar{E}).$$
Then 
$$\kappa(x)=\left(x,\ch(\C^{\st})\right)^R.$$
Indeed, this follows from Lemma \ref{stable-functor-lem} by spelling out the
definitions.
\end{ex}

Using \eqref{tr-form-eq} and the formula for the canonical pairing $(\cdot,\cdot)$ from
\cite[Thm.\ 4.2.1]{PV-mf}, we deduce the following explicit formula for the $R$-bilinear pairing $(\cdot,\cdot)^R$.

\begin{lem}\label{R-bil-form-lem}
Let us consider the standard residue pairing 
$$\lan f\ot dx, g\ot dx\ran_\w=(-1)^{{n\choose 2}}\Res_0(f\cdot g)$$
on the twisted Milnor ring $\AA_\w\ot dx=H(\w)$ (where $\Res_0$ is the
Grothendieck residue on $\AA_\w$).
For $h\in\HH(\w)$ let $h_\ga\in H(\w_\ga)$ be the component of $h$ with respect to
the decomposition \eqref{w-hoch-decomp-eq}.
Then the $R$-bilinear canonical pairing \eqref{R-bil-form-eq} is given by
\begin{equation}\label{R-pairing-formula}
(h,h')^R=\sum_{\ga\in G}c_\ga\cdot \lan h_\ga,h'_\ga\ran_{\w_{\ga}}\cdot e_\ga,
\end{equation}
where 
\begin{equation}\label{e-idemp}
e_\ga:=\frac{1}{|G|}\cdot \sum_{\eta\in \widehat{G}}\eta^{-1}(\ga)[\eta].
\end{equation}
$$c_\ga=\det[\id-\ga,T^*/(T^*)^\ga]^{-1},$$
and $T^*=\mg/\mg^2$.
\end{lem}

\Pf . By Theorem \ref{hoch-prop}, the decomposition \eqref{w-hoch-decomp-eq} 
coincides with the decomposition
$$\HH(\w)=\bigoplus_{\ga\in G} e_{\ga}\HH(\w)$$
induced by the $R$-module structure on $\HH(\w)$. Since both sides of
\eqref{R-pairing-formula} are $R$-bilinear, it is enough to check the equality after applying
$\tr$ (see \eqref{tr-R-eq}) which holds by \cite[Thm.\ 4.2.1]{PV-mf}.
\ed

\section{$\Ga$-spin curves and $\w$-structures}\label{simple-constr-sec}

\subsection{Abstract $\w$-structures}\label{abs-w-str-sec}

Let  
$$\w(x_1,\ldots,x_n)=\sum_{k=1}^N c_k M_k$$ 
be a Laurent polynomial in $x_1,\ldots,x_n$,
where $c_k\in\C^*$ and 
$$M_k=\prod_{i=1}^n x_i^{m_{ki}}$$
are monomials.
We denote by $\bfm_\w:\Z^n\to\Z^N$ the map given
by the matrix $(m_{ki})$ of exponents.

\begin{defi} Let $\bfd=(d_1,\ldots,d_n)\in\Z^n$ be a primitive vector and let $d$ be
an integer. A Laurent polynomial is {\it quasihomogeneous} of degree $d$ with respect to $\bfd$
if
$$\w(\la^{d_1}x_1,\ldots,\la^{d_n}x_n)=\la^d \w(x_1,\ldots,x_n).$$
\end{defi}

The above equation is equivalent to
\be\label{w-degrees-eq}
\bfm_\w(\bfd)=d\cdot \bfe,
\end{equation}
where $\bfe\in\Z^N$ the vector with all components equal to $1$.
Let us consider the dual homomorphisms $\bfm_\w^*:\Z^N\to\Z^n$,
$\bfe^*:\Z^N\to\Z$ and $\bfd^*:\Z^n\to\Z$. Then \eqref{w-degrees-eq} implies that
$$\bfd^*\circ\bfm_\w^*=d\cdot \bfe^*.$$
Let us consider the subgroup $P_\w=\im(\bfm_\w^*)\sub\Z^n$. Then we obtain that the restriction
$\bfd^*|_{P_\w}$ is divisible by $d$ and we can define the homomorphism
\begin{equation}\label{deg-P-Z-eq}
\deg=\frac{1}{d}\bfd^*:P_\w\to\Z,
\end{equation}
so that $\deg\circ\bfm_\w^*=\bfe^*$.

In the following definition
we view abelian monoids (such as $\Z^n$, $\Z^n_{\ge 0}$, etc.)
as symmetric monoidal categories with objects corresponding to elements, 
with only identity morphisms,
and the tensor operation given by the monoid structure.

\begin{defi}
Let $\CC$ be a symmetric monoidal
category with a unit object $\unit$, and let $E\in\CC$ be an invertible object.
Let also $\w$ be a Laurent polynomial in $x_1,\ldots,x_n$,  quasihomogeneous of degree $d$
with respect to $\bfd$.

(i) A {\it $(\w,\bfd)$-structure} in $\CC$ with respect to $E$ is a monoidal functor
$$\Phi:\Z^n\to\CC$$
together with an isomorphism of monoidal functors
$$\Phi|_{P_\w}\to\La\circ\deg,$$
where $\La:\Z\to\CC$ is the monoidal functor $i\mapsto E^i=E^{\ot i}$, and
$\deg$ is the homomorphism \eqref{deg-P-Z-eq}.

(ii) Assume in addition that $\w$ is a polynomial.
A {\it weak $(\w,\bfd)$-structure} in $\CC$ with respect to $E$ is a monoidal functor
$$\Phi:\Z^n\to\CC$$
together with a morphism of monoidal functors
\begin{equation}\label{phi-La-eq}
\phi:\Phi|_{P_{\w}^+}\to \La\circ\deg|_{P_{\w}^+},
\end{equation}
where $P_{\w}^+=P_\w\cap \Z_{\ge 0}^n$. 
\end{defi}

The following proposition gives a more down-to-earth interpretation of $\w$-structures.

\begin{prop}\label{w-str-prop} 
Let $\w$ be a Laurent polynomial, quasihomogeneous of degree $d$ with respect to $\bfd$,
and let $(\CC,E)$ be as above.

\noindent
(i) Let $v_1,\ldots,v_n$ be a basis of $\Z^n$ such that
$k_1v_1,\ldots,k_rv_r$ is a basis of $P_\w$ for  
positive integers $k_1,\ldots, k_r$, where $r=\rk P_\w\le n$. Then isomorphism
classes of
$\w$-structures with respect to $E$ correspond to isomorphism classes of
collections of invertible objects
$\Phi(v_1),\ldots,\Phi(v_n)$  in $\CC$ together with isomorphisms 
$$\phi_i: \Phi(v_i)^{\ot k_i}\to E^{\deg(k_iv_i)}, \ \ i=1,\ldots,r.$$

\noindent
(ii) Assume in addition that $\w$ is a polynomial.
For every weak $(\w,\bfd)$-structure $(\Phi,\phi)$ in $\CC$ with respect to $E\in\CC$,
such that $\phi$ is an isomorphism,
there exists a unique extension of $\phi$ to an isomorphism of monoidal functors
$$\wt{\phi}:\Phi|_{P_\w}\to\La\circ\deg,$$
i.e., to a $(\w,\bfd)$-structure.
\end{prop}

\Pf . (i) Any monoidal functor $\Phi:\Z^n\to\CC$ is determined up to an isomorphism by
the collection of invertible objects $\Phi(v_1),\ldots,\Phi(v_n)$. The same is true
for the group $P_\w\simeq\Z^r$ with the basis $k_1v_1,\ldots,k_rv_r$, which implies the result.

\noindent
(ii) For $p\in P_{\w}^+$ we have an isomorphism
$$\Phi(-p)\simeq\Phi(p)^{-1}\simeq E^{-\deg(p)}=E^{\deg(-p)},$$
where the first isomorphism comes from the monoidal structure on $\Phi$ and the second
is induced by $\phi$. Since 
$P_\w$ is generated by the vectors $d_\w^*(e_s)\in\Z_{\ge 0}^n$, $s=1,\ldots,N$,
every element of $P_\w$ can be represented 
as $p-p'$ with $p,p'\in P_{\w}^+$. The monoidal structure on $\Phi$ gives an isomorphism
$\Phi(p-p')\simeq \Phi(p)\ot\Phi(p')^{-1}$. Now 
$\phi(p)$ and $\phi(p')$ induce an isomorphism
$$\wt{\phi}(p-p'):\Phi(p-p')\to E^{\deg(p)}\ot E^{-\deg(p')}\simeq E^{\deg(p-p')},$$
which is the unique extension of $\phi$ to $P_\w$.
\ed

\begin{defi} The $(\w,\bfd)$-structure in $\CC$ with respect to $E=\unit$ given by 
$\Phi(?)=\unit$ with the identity
isomorphisms $\phi$ is called the {\it trivial $(\w,\bfd)$-structure}.
\end{defi} 

\begin{prop}\label{W-triv-cor}
Suppose that $\End(\unit)$ is an algebraically closed field,
$E\simeq\unit$, and a $(\w,\bfd)$-structure $(\Phi,\phi)$ satisfies
$\Phi(e_j)\simeq\unit$ for every $j$. Then this $(\w,\bfd)$-structure is isomorphic to
the trivial $(\w,\bfd)$-structure.
\end{prop}

\Pf . By Proposition \ref{w-str-prop}, such a $(\w,\bfd)$-structure corresponds to a collection of
isomorphisms 
$\phi_i:\Phi(v_i)^{\ot k_i}\to \unit$, $i=1,\ldots,r$.
But $\Phi(v_i)\simeq \unit$, so $\phi_i$ can be viewed as an element of $\End(\unit)$.
Choose $\xi_i\in\End(\unit)$ such that $\xi_i^{k_i}=\phi_i$. Then the morphisms
$$\Phi(v_i)\simeq \unit\rTo{\xi_i}\unit$$
induce an isomorphism with the trivial $(\w,\bfd)$-structure.
\ed

\begin{rem}
We will sometimes omit the vector of degrees $\bfd$ from notation and talk simply of $\w$-structures,
when this vector is fixed. In the case when $\bfm_\w$ is injective the vector $\bfd$ is uniquely determined
by $\w$ up to a sign.
\end{rem}

\subsection{$\Ga$-spin curves and their moduli} 
\label{w-moduli-sec}

Let us fix an algebraic subgroup $\Ga\sub\G_m^n$ with  a surjective character $\chi:\Ga\to\G_m$
such that $G=\ker(\chi)$ is finite. Thus, we have an exact sequence of commutative algebraic groups
$$1\to G\to \Ga\rTo{\chi} \G_m\to 1.$$
With these data we will associate a finite covering
$$\SS_{g,r,\Ga,\chi}\to \ov{\MM}_{g,r}$$
of the Deligne-Mumford moduli stacks of stable curves.
The stacks $\SS_{g,r,\Ga,\chi}$ are slight generalizations of the 
moduli spaces of $\w$-curves considered in \cite{FJR}.

Recall (see \cite[Sec.\ 4]{AV}, \cite[Sec.\ 2.1]{FJR}) that an {\it orbicurve with marked points} 
$(\CC,p_1,\ldots,p_r)$  is a proper
Deligne-Mumford stack $\CC$ whose coarse moduli
space is a (connected) nodal curve $C$, equipped with marked orbipoints $p_1,\ldots,p_r\sub\CC$,
such that the projection $\rho:\CC\to C$ is an isomorphism away from the marked points and from the
nodes. It is also required that each node is locally modeled by a quotient stack of the form
$\{xy=0\}/(\Z/n)$, where the action of $\Z/n$ is given by $(x,y)\mapsto 
(\exp(2\pi i/n)x,\exp(-2\pi i/n)y)$.
We say that an orbicurve $\CC$ is {\it smooth} if the curve $C$ is smooth.
We denote by 
$$\om^{\log}_{\CC}=\rho^*(\om_{C}(p_1+\ldots+p_r))$$
the log-canonical line bundle with respect to $p_1,\ldots,p_r$
(see \cite[Def. 2.1.2]{FJR}). The following definition is a coordinate-free version of the notion of
$\w$-curves introduced in \cite{FJR}. We use principal bundles of algebraic groups. For
a homomoprhism $f:G_1\to G_2$ of algebraic groups and a principal
$G_1$-bundle $P$ we denote by $f_*P$ the 
pushout
of $P$ with respect to $f$, which is a principal $G_2$-bundle.

\begin{defi}\label{spin-def} 
(i) A {\it $(\Ga,\chi)$-spin curve} (a {\it $\Ga$-spin curve} for short) is an 
orbicurve with marked points $(\CC,p_1,\ldots,p_r)$ together with
a principal $\Ga$-bundle $P$ over $\CC$ and an isomorphism of $\G_m$-bundles
\begin{equation}\label{Gamma-spin-eq}
\vareps:\chi_*P\to P(\om^{\log}_{\CC}),
\end{equation}
where $P(\om^{\log}_{\CC})$ is the principal $\G_m$-bundle associated with the line bundle $\om_{\CC}^{\log}$.
An {\it isomorphism} between two $\Ga$-spin curves is an isomorphism of curves
with marked points $f:(\CC,p_1,\ldots, p_r)\to(\CC',p'_1,\ldots,p'_r)$ 
and an isomorphism of $\Ga$-bundles
$t: P\to f^*P'$ compatible with isomorphisms \eqref{Gamma-spin-eq} for $P$ and $P'$.

(ii)  Let $(P,\eps)$ be a $\Ga$-spin structure on an orbicurve $(\CC,p_1,\ldots,p_r)$.
For each marked point $p_i$ let us consider the homomorphism
\begin{equation}\label{G-p-i-hom-eq}
G(p_i)\to\Ga\sub(\C^*)^n
\end{equation}
from the local automorphism group of $\CC$ at $p_i$
associated with its action on the fiber of $P$ at $p_i$.
We denote by 
\begin{equation}\label{gamma-i-eq}
\ga_i=\ga_i(P)=(\ga_{i1},\ldots,\ga_{in})\in (\C^*)^n
\end{equation}
the image of the canonical generator of the cyclic group $G(p_i)$ under \eqref{G-p-i-hom-eq}.
The collection $\ov{\ga}=(\ga_1,\ldots,\ga_r)$ is called the {\it type} of the $\Ga$-spin curve.

(iii) We say that a $\Ga$-spin curve is {\it stable} if $(\CC,p_1,\ldots,p_r)$ is stable and
all the homomorphisms \eqref{G-p-i-hom-eq} are injective.
\end{defi}

Since the restriction of $\om^{\log}_{\CC}$ to each marked point is trivial, the isomorphisms
\eqref{Gamma-spin-eq} imply that the elements $\ga_i\in \Ga$ belong to the subgroup
$\ker(\chi)=G$.
For a stable $\Ga$-spin curve we will identify $G(p_i)$ with
the subgroup $\lan\ga_i\ran\sub G$ generated by $\ga_i$.

It is useful to rewrite the definition of a $\Ga$-spin structure in terms of principal bundles that have
only algebraic tori as structure groups.
Namely, consider the algebraic torus $T=\G_m^n/G$.
The embedding $\Ga\hra\G_m^n$ induces an embedding
$\varphi:\G_m=\Ga/G\hra T$. Thus, we have a commutative diagram with exact rows
\begin{equation}\label{Ga-spin-diagram-eq}
\begin{diagram}
1\rTo{} & G &\rTo{}& \Ga &\rTo{\chi}&\G_m&\rTo{}1\\
&\dTo{\id} && \dTo{} &&\dTo{\varphi}\\
1\rTo{} & G &\rTo{}& \G_m^n &\rTo{\pi}& T&\rTo{}1
\end{diagram}
\end{equation}
such that $\pi$ induces an isomorphism $\G_m^n/\Ga\to T/\varphi(\G_m)$.

The principal $\Ga$-bundle $P$ in Definition \ref{spin-def} gives rise to a $\G_m^n$-bundle $P'$ via the embedding
$\Ga\to\G_m^n$. We are going to rewrite the definition of a $\Ga$-spin structure in terms of $P'$ and the
homomorphisms $\pi$ and $\varphi$ from diagram \eqref{Ga-spin-diagram-eq} (see Proposition
\ref{spin-prop}(i) below).
On the other hand, we can view $P'$ a collection of line bundles $(\LL_1,\ldots,\LL_n)$ on $\CC$.
We will show that the isomorphism \eqref{Gamma-spin-eq}
can be interpreted as a $(\w,\bfd)$-structure
in the category of line bundles on $\CC$ with respect to $\om_{\CC}^{\log}$ for
some quasihomogeneous Laurent polynomial $\w$ (see Section \ref{abs-w-str-sec})

Let $\Ga_0$ be the connected component of $1$ in $\Ga$. 
$\Ga_0$ is a one-dimensional torus, so we can choose an identification $\Ga_0=\G_m$. 
The embedding $\Ga_0=\G_m\hra\G_m^n$ takes form $\la\mapsto (\la^{d_1},\ldots,\la^{d_n})$
for some primitive vector $\bfd=(d_1,\ldots,d_n)\in\Z^n$. Furthermore, the restriction
of the character $\chi$ to $\Ga_0=\G_m$ is given by $\la\mapsto\la^d$ for some $d$.
Note that the degrees $(\bfd,d)\in\Z^{n+1}$ are determined by $\Ga$ uniquely up to a sign.

Consider the exact sequence of algebraic tori
$$1\to \G_m\rTo{\varphi} T\to T'\to 1.$$
Since we work over $\C$, we can find a splitting $T\simeq\G_m\times T'$. Consider the collection 
of characters of $T$ 
\begin{equation}\label{T-char-eq}
\eta_0=(\id,1), \eta_1=(\id,\eps_1), \ldots, \eta_{n-1}=(\id,\eps_{n-1}),
\end{equation}
where $\id$ is the identity character of $\G_m$ and $(\eps_1,\ldots,\eps_{n-1})$ is
a basis of the group of characters of the torus $T'$. 
Then we have 
$$\bigcap_{i=0}^{n-1}\ker(\eta_i)=1.$$ 
Using the projection $\pi:\G_m\to T$ we get characters of $\G_m^n$,
\begin{equation}\label{M-i-def-eq}
M_i=\eta_i\circ\pi:\G_m^n\to \G_m \text{ for }i=0,\ldots,n-1,
\end{equation}
satisfying $M_i|_{\Ga}=\eta_i\circ\varphi\circ\chi=\chi$ and such that
$$\bigcap_{i=0}^{n-1}\ker(M_i)=G.$$

\begin{prop}\label{spin-prop} 
(i) The category of $\Ga$-spin structures on $(\CC,p_1,\ldots,p_r)$ is equivalent to the category
of pairs $(P',\vareps')$, where $P'$ is a principal $\G_m^n$-bundle and $\vareps'$ is an isomorphism
\begin{equation}\label{Gamma-spin-bis-eq}
\vareps':\pi_*P'\to\varphi_*P(\om^{\log}_{\CC}).
\end{equation}

\noindent
(ii) There exists a Laurent polynomial $\w(x_1,\ldots,x_n)$, quasihomogeneous of degree $d$ with
respect to the grading of the variables given by $\bfd$, 
such that $G$ is equal to the group of diagonal symmetries of $\w$.

\noindent
(iii) For any $\w$ as in (ii) the category of $(\w,\bfd)$-structures in the category of line bundles on $\CC$ 
with respect to $\om_{\CC}^{\log}$ is equivalent to the category of $(\Ga,\chi)$-spin structures on 
$(\CC,p_1,\ldots,p_r)$.

\noindent
(iv) Let $\w$ be a quasihomogeneous polynomial of degree $d>0$ with respect to the grading given by
$\bfd$, and let
$G$ be a finite subgroup of the group $G_\w$ of diagonal symmetries of $\w$, such that
$G$ contains the exponential grading element $J$ (see \eqref{J-eq}).
Let also $\Ga\sub\G_m^n$ be the subgroup associated with $G\sub G_\w$ by Lemma \ref{J-lem}. 
Then to any $\Ga$-spin structure
on $(\CC,p_1,\ldots,p_r)$ there corresponds a natural $(\w,\bfd)$-structure in the category of line bundles on $\CC$
with respect to $\om_{\CC}^{\log}$.
\end{prop}

\Pf . (i) A $\Ga$-bundle $P$ can be viewed as a
$\G_m^n$-bundle $P'$ together with a trivialization of the induced $\G_m^n/\Ga$-bundle.
Commutativity of the right  square in \eqref{Ga-spin-diagram-eq} shows that
$$\pi_*P'\simeq\varphi_*\chi_*P.$$
Hence, the isomorphism \eqref{Gamma-spin-eq} gives rise to an isomorphism \eqref{Gamma-spin-bis-eq}.
Conversely, starting with a pair $(P',\vareps')$ we observe that the isomorphism
\eqref{Gamma-spin-bis-eq} induces a trivialization of the $T/\varphi(\G_m)$-bundle obtained by
the pushout from $\pi_*P'$,
or equivalently a trivialization of the $\G_m^n/\Ga$-bundle obtained by pushout from $P'$. 
Thus, we can reduce the structure to $\Gamma$ and obtain a $\Ga$-bundle $P$.
Now both parts of \eqref{Gamma-spin-bis-eq} become pushouts of the corresponding parts of
\eqref{Gamma-spin-eq} with respect to $\varphi$. Since $\vareps'$ is compatible with the trivializations
of the pushouts with respect to the projection $T\to T/\varphi(\G_m)$, it induces an isomorphism
\eqref{Gamma-spin-eq}.

\noindent
(ii) Recall that the restriction of $\chi$ to $\Ga_0=\G_m$ sends $\la$ to $\la^d$.
Therefore, each of the characters \eqref{M-i-def-eq} can be viewed as a Laurent monomial in $x_1,\ldots,x_n$ of
degree $d$ with respect to $\bfd$. Thus, we can take $\w=\sum_{i=0}^{n-1} M_i$.

\noindent
(iii) First, we observe that for an algebraic torus $T$ the category of 
principal $T$-bundles on $\CC$ is equivalent to the category of monoidal functors from
 with the character group $X(T)$ to the monoidal category $\PPic(\CC)$ of line bundles on $\CC$. 
Namely, with a $T$-bundle $P$ we associate the monoidal
functor $\Phi_P$ sending a character $\eta:T\to\G_m$ to the line bundle corresponding to the induced
$\G_m$-torsor $\eta_*P$. Indeed, a choice of a basis of $X(T)$ shows that
both structures are equivalent to collections of line bundles on $\CC$.
If $f:T_1\to T_2$ is a homomorphism of tori, then the monoidal functor $\Phi_{f_*P_1}:X(T_2)\to\PPic(\CC)$ 
associated with the pushout $f_*P_1$ of a principal $T_1$-bundle is isomorphic to the composition 
$\Phi_{P_1}\circ f^*$, where $f^*:X(T_2)\to X(T_1)$ is the induced homomorphism of the character
groups.

Given a Laurent polynomial $\w=\sum_{s=1}^N M_s$ as in (ii), we have
$$G=\ker((\bfm_\w)_*:\G_m^n\to\G_m^N),$$
where $\bfm_\w:\Z^n\to\Z^N$ is the linear map defined by the exponents of the monomials $M_s$ (see Section 
\ref{abs-w-str-sec}). 
Hence, the map $(\bfm_\w)_*$ factors through the projection $\pi:\G_m^n\to T$ followed by
an embedding of tori $T\hra\G_m^N$. The induced homomorphisms of character groups
$\Z^N\to X(T)$ and $\pi^*:X(T)\to \Z^n$ are surjective and injective, respectively. Therefore,
$\pi^*$ induces an isomorphism of $X(T)$ with $P_\w=\im(\bfm_\w)\sub\Z^n$ such that
the homomorphism $\deg:P_\w\to\Z$ gets identified with
the homomorphism on the character groups $\varphi^*:X(T)\to \Z$ induced by the embedding
of tori $\varphi:\G_m\to T$. Now the assertion follows from (i).

\noindent
(iv) Let us choose a Laurent polynomial $\ov{\w}(x_1,\ldots,x_n)$ as in (ii). Note that adding to
$\ov{\w}$ a linear combination of $G$-invariant monomials will still give us a Laurent polynomial with the diagonal group of symmetries equal to $G$. Thus, we can assume that 
$\ov{\w}$ contains all the monomials in $\w$. By (iii), a $\Ga$-spin structure induces a $(\ov{\w},\bfd)$-structure,
which in turn gives a $(\w,\bfd)$-structure.
\ed

\begin{rems}\label{quasihom-monom-rem}
1. 
Proposition \ref{spin-prop}(ii) implies that for a quasihomogeneous polynomial
$\w$, every finite subgroup $G$ of the group of diagonal symmetries $G_\w$ such that $G$ contains the
exponential grading element $J$, is admissible in the sense of \cite[Def.\ 2.3.2]{FJR}.

2. Let $(P,\vareps)$ be a $\Ga$-spin structure on $(\CC,p_1,\ldots,p_r)$.
Any Laurent monomial $M:\G_m^n\to\G_m$ of degree $d$
with respect to $\bfd$, such that $M|_G=1$, restricts to the character $\chi$ on $\Ga$.
Indeed, the condition that $M|_G=1$ implies that $M|_\Ga=\chi^a$ for some $a\in\Z$. 
The fact that $a=1$ follows from
the equality $M|_{\Ga_0}=\chi|_{\Ga_0}$ which in turn follows from the condition that $M$ has the
degree $d$ with respect to $\bfd$. 
Thus, the isomorphism 
\eqref{Gamma-spin-eq} gives rise to isomorphisms
\begin{equation}\label{w-str-isom-eq}
M(\LL_1,\ldots,\LL_n)\rTo{\sim}\om^{\log}_{\CC},
\end{equation}
for any Laurent monomial  $M$  of degree $d$ with respect to $\bfd$ such that $M|_G=1$. 
These are the kind of isomorphisms that appear in the original definition of a $\w$-curve in \cite{FJR}.
More precisely, if $\w(x_1,\ldots,x_n)$ is a Laurent polynomial,  
quasihomogeneous of degree $d$ with respect to $\bfd$, such that $G$ is equal to the group $G_\w$
of diagonal symmetries
of $\w$, then by Proposition \ref{spin-prop}(iii), the notion
of a $\Ga$-spin curve is equivalent to the notion of a $\w$-curve defined in \cite[Sec.\ 2.1]{FJR} (generalized to the case when $\w$ is a {\it Laurent} polynomial).
Indeed, 
the Smith normal forms appearing in \cite[Def. 2.1.10]{FJR} is just a way of recording an isomorphism of
monoidal functors from a free abelian group in coordinates (see Proposition \ref{w-str-prop}).
More generally, if $G$ is just a subgroup of $G_\w$ containing $J$, then by Proposition
\ref{spin-prop}(iv), any $\Ga$-spin curve
has a natural structure of a $\w$-curve.
\end{rems}

We will give now yet another way to describe $\Ga$-spin structures.

\begin{cor}\label{quasihom-monom-cor}
The category of $\Ga$-spin structures on $(\CC,p_1,\ldots,p_r)$
is equivalent to the category of collections of $n$ line bundles $\LL_1,\ldots,\LL_n$ on $\CC$
together with
isomorphisms \eqref{w-str-isom-eq} for 
the Laurent monomials $M=M_0,\ldots,M_{n-1}$ given by \eqref{M-i-def-eq}.
\end{cor}

\Pf . The line bundle $(\LL_1,\ldots,\LL_n)$ are obtained from a $\Ga$-spin structure
as discussed before Proposition
\ref{spin-prop}. The isomorphisms \eqref{w-str-isom-eq} for the monomials
$M_0,\ldots,M_{n-1}$
correspond to an isomorphism of $T$-bundles \eqref{Gamma-spin-bis-eq}
under the identification
$$(M_0, M_1, \ldots, M_{n-1}):T\rTo{\sim}\G_m^n.$$ 
Conversely, assume we have line bundles $\LL_1,\ldots,\LL_n$ equipped with isomorphisms 
\eqref{w-str-isom-eq} for $M_0,\ldots,M_{n-1}$. Then we obtain trivializations of the line bundles
$(M_i/M_0)(\LL_1,\ldots,\LL_n)$. Note that the homomorphism
$$(M_1/M_0,\ldots,M_n/M_0):\G_m^n\to \G_m^{n-1}$$
is equal to the composition 
$$\G_m^n\rTo{\pi}\to T\to T'\simeq \G_m^{n-1}.$$
Hence, the exact sequence of groups
$$1\to \Ga\to \G_m^n\to T'\to 1$$
shows that the $\G_m^n$-bundle $(\LL_1,\ldots,\LL_n)$ comes from a $\Ga$-bundle.
\ed

We are going to work with the moduli space of $\Ga$-spin curves.

\begin{defi}\label{fam-spin-def}
A genus-$g$, stable $\Ga$-spin curve with $k$ marked points over a base $T$ is a flat family
$\CC\to T$ of genus-$g$ orbicurves with gerbe markings $p_1,\ldots,p_k\sub \CC$ and sections
$\si_i:T\to p_i$ inducing isomorphism of $T$ with the coarse moduli of $p_i$, together with a relative
$\Ga$-spin structure $(P,\vareps)$ on $\CC$, such that all the fibers over closed points of $T$ are
stable $\Ga$-spin curves.
Here $P$ is a $\Ga$-bundle on $\CC$, and
$$\vareps:\chi_*P\simeq P(\om^{\log}_{\CC/T})$$
is an isomorphism of $\G_m$-bundles on $\CC$, where $\om^{\log}_{\CC/T}$ is the relative
log-canonical line bundle. 
These structures naturally form a stack $\SS_{g,r}=\SS_{g,r,\Ga,\chi}$,
which is the disjoint union of the open and closed substacks
$\SS_g(\ov{\ga})$ for $\ov{\ga}\in G^r$, parametrizing $\Ga$-spin curves of type $\ov{\ga}$.
\end{defi}

\begin{prop}\label{spin-stack-prop} 
The stack $\SS_{g,r}$ is a smooth proper DM-stack 
over $\C$ with projective coarse moduli, and the natural forgetful morphism 
$\SS_{g,r}\to\ov{\MM}_{g,r}$ is quasi-finite.
If $G$ is equal to the group $G_\w$ of diagonal symmetries of a Laurent polynomial
$\w$,  quasihomogeneous of degree $d$ with respect to $\bfd$, then
$\SS_{g,r}$ is naturally isomorphic to the stack $\WW_{g,r}(\w)$ of $\w$-curves constructed in \cite[Sec.\ 2.2]{FJR}.
\end{prop}

\Pf . First, we observe that although Fan-Jarvis-Ruan assume that $\w$ is a polynomial,
the definition of a $\w$-structure and the results of
\cite[Sec.\ 2.2]{FJR} are valid also in the case of a quasihomogeneous {\it Laurent} polynomial with finite group
of diagonal symmetries
(in fact, this extension is used in \cite[Sec.\ 2.3]{FJR} to define the moduli spaces associated with admissible
subgroups of $G_\w$). 

Applying Proposition \ref{spin-prop}(ii,iii), we see that our moduli stack $\SS_{g,r}$ is naturally isomorphic
to $\WW_{g,r}(\w)$ for some quasihomogeneous Laurent polynomial $\w$. Now we can use
\cite[Thm.\ 2.2.6]{FJR} to derive the required properties of the stack $\SS_{g,r}$.
(Alternatively, by modifying the arguments of \cite[Thm.\ 2.2.6]{FJR} one can work directly
with $\Ga$-spin structures.) 
\ed

Next, we consider rigidifications of $\Ga$-spin curves that were introduced in \cite[Sec.\ 2.2.3]{FJR}.

\begin{defi}
A {\it rigidification of a $\Ga$-spin structure $(P,\vareps)$ on an orbicurve $\CC$ at a marked point $p_i$}
is a trivialization of $P|_{p_i}$, i.e., an isomorphism
$$P|_{p_i}\simeq \Ga/\lan\ga_i\ran$$
compatible with the canonical trivialization
of $\om_{\CC}^{\log}|_{p_i}$ via the isomorphism \eqref{Gamma-spin-eq} (note that we can view the marked point $p_i$ as the gerbe $pt/\lan\ga_i\ran$).
A {\it rigidification of a $\Ga$-spin curve} $(\CC,p_1,\ldots,p_r; P,\vareps)$ consists
of a collection of rigidifications of $(P,\vareps)$ at every marked point $p_i$. 
\end{defi}

The group $\prod_{i=1}^r G/\lan\ga_i\ran$ acts simply transitively on the set of rigidifications
of a given $\Ga$-spin curve.
Thus, the moduli stack of rigidified $\Ga$-spin curves is a $\prod_{i=1}^r G/\lan\ga_i\ran$-torsor over
$\SS_g(\ov{\ga})$ that we denote by
\begin{equation}\label{rig-eq}
\SS^{\rig}_g(\ov{\ga})\to\SS_g(\ov{\ga})
\end{equation}
Let $(\LL_1,\ldots,\LL_n)$ be the line bundles associated with a $\Ga$-spin structure $P$.
For each $(p_i,j)$ such that the $j$th component of $\ga_i$ is trivial,
a rigidification structure induces a well-defined trivialization of $\LL_j|_{p_i}$.
Below we will define a different version of a rigidification structure that keeps track only of these
trivializations.

First, we define restrictions of $\Ga$-spin structures associated with coordinate projections
$\G_m^n\to\G_m^k$. 

\begin{defi}
Let $I=\{i_1,\ldots,i_k\}\sub\{1,\ldots,n\}$ be a nonempty
subset such that the vector $\bfd_I=(d_{i_1},\ldots,d_{i_k})$ is not zero,
and let 
$$p_I: \G_m^n\to\G_m^k$$ 
be the corresponding coordinate projection.
We say that $I$ is 
{\it $(\Ga,\chi)$-admissible} if the character $\chi$ factors through the projection
$\Ga\to \Ga_I:=p_I(\Ga)$. 
\end{defi}

The following statement is an immediate corollary of the definitions.

\begin{lem} Assume that $I$ is $(\Ga,\chi)$-admissible.
Let $\chi_I:\Ga_I\to\G_m$ be the character of $\Ga_I$ induced by $\chi$.
Then any $(\Ga,\chi)$-spin structure $(P,\vareps)$ on $(\CC,p_1,\ldots,p_r)$ naturally induces a 
$(\Ga_I,\chi_I)$-spin structure
$(P_I,\vareps_I)$ on $(\CC,p_1,\ldots,p_r)$ with $P_I=(p_I)_*P$.
Also, a rigidification of $(P,\vareps)$ induces a rigidification of $(P_I,\vareps_I)$.
\end{lem}

\ed

For every $\ga\in G\sub\G_m^n$ let us denote by $I(\ga)\sub\{1,\ldots,n\}$ the set of all $j$ such that
the $j$th component of $\ga$ is trivial.

\begin{defi}\label{restr-rig-def} Assume that the degree $d_i$ is nonzero for every $i=1,\ldots,n$.

(i) A collection $\ov{\ga}=(\ga_1,\ldots,\ga_r)\in G^r$ is called {\it $(\Ga,\chi)$-admissible} if 
for every $i=1,\ldots,r$ the subset $I(\ga_i)$ is $(\Ga,\chi)$-admissible, if nonempty.

(ii) If $\ov{\ga}$ is $(\Ga,\chi)$-admissible, we define
a {\it restricted rigidification} of a $\Ga$-spin structure $(P,\vareps)$ of type $\ov{\ga}$
as a collection of rigidifications
of the induced $\Ga_{I(\ga_i)}$-structure $(P_{I(\ga_i)},\vareps_{I(\ga_i)})$ at $p_i$ for $i=1,\ldots,r$
such that $I(\ga_i)$ is nonempty.
\end{defi}

Note that by definition, a restricted rigidification of a $\Ga$-spin
structure consists of trivializations of $\LL_j|_{p_i}$ for $j\in I(\ga_i)$ and $i=1,\ldots,r$,
satisfying certain compatibilities.

For a $(\Ga,\chi)$-admissible collection
$\ov{\ga}=(\ga_1,\ldots,\ga_r)\in G^r$ let us set 
$$G(\ov{\ga})=\prod_{i=1}^r G_{I(\ga_i)}, \text{ where }\ G_I=p_I(G).$$ 
Then the group $G(\ov{\ga})$
acts simply transitively on the set of restricted rigidifications of a given $\Ga$-spin curve of type $\ov{\ga}$.
We denote by 
$\SS_g^{\rig,0}(\ov{\ga})\to\SS_g(\ov{\ga})$
the $G(\ov{\ga})$-torsor of restricted rigidifications. 
We have a natural surjective morphism
\begin{equation}\label{sur-rig-mor-eq}
\SS_g^{\rig}(\ov{\ga})\to\SS_g^{\rig,0}(\ov{\ga})
\end{equation} 
compatible with the homomorphism
$$r_{\ov{\ga}}:\prod_{i=1}^r G/\lan\ga_i\ran\to G(\ov{\ga}).$$
Therefore, the $G(\ov{\ga})$-torsor $\SS_g^{\rig,0}(\ov{\ga})$ is isomorphic to the
pushout of \eqref{rig-eq} with respect to $r_{\ov{\ga}}$.

\begin{lem}\label{admissibility-lem} 
Let $\w(x_1,\ldots,x_n)$ be a polynomial with an isolated singularity at the origin, 
quasihomogeneous
of degree $d$ with respect to $\bfd=(d_1,\ldots,d_n)$, where $d_i\neq 0$ for $i=1,\ldots,n$.
Let also $G\sub G_\w$ be a finite subgroup containing the exponential grading element $J$, 
and let $\Ga\sub\G_m^n$ be the corresponding
extension of $\G_m$ by $G$ (see Lemma \ref{J-eq}). Then every $\ov{\ga}=(\ga_1,\ldots,\ga_r)\in G^r$ is
$(\Ga,\chi)$-admissible.
\end{lem}

\Pf . It is enough to check that for every $\ga\in G_\w$ the set 
$I(\ga)=\{i_1,\ldots,i_k\}\sub\{1,\ldots,n\}$ is $(\Ga,\chi)$-admissible,
provided it is nonempty.
To do this we use the fact that the restriction of $\w$ to the subspace of $\ga$-invariants still has an
isolated singularity at the origin (see \cite[Lem.\ 2.5.3(i)]{PV-mf}). 
In particular, this restriction is nonzero.
Let us take any monomial $M(x_{i_1},\ldots,x_{i_k})$
occurring in the restriction of $\w$ to the subspace of $\ga$-invariants. 
Since $M$ also occurs in $\w$, the action of $\Ga$ on $M$ rescales it by the character $\chi$.
But $M$ factors through $p_{I(\ga)}:\G_m^n\to\G_m^k$, hence, $\chi$ also factors through $p_{I(\ga)}$.
\ed

\subsection{Invariants of smooth $\Ga$-spin curves}\label{degrees-sec}

We keep the assumptions and notation of the beginning of Section \ref{w-moduli-sec}.
Let $(\CC,p_1,\ldots,p_r; P, \vareps)$ be a $\Ga$-spin curve of type 
$\ov{\ga}=(\ga_1,\ldots,\ga_r)\in G^r$ with a smooth orbicurve $\CC$, and let
$\LL_1,\ldots,\LL_n$ be the line bundles associated with the $\G_m^n$-bundle $P$.
 Consider the map $\rho:\CC\to C$, where $C$ is
 the smooth curve obtained by forgetting the
orbi-structure at the marked points, and the line bundles $L_j=\rho_*(\LL_j)$ on $C$.

For $\ga\in G$ we define
$$\bth_{\ga}=(\th_1,\ldots,\th_n)\in \Q^n$$
as the unique vector with $0\le \th_j<1$ for $j=1,\ldots,n$ such that
$$\ga=\exp(2\pi i\bth_{\ga})=(\exp(2\pi i\th_1),\ldots,\exp(2\pi i\th_n))\in (\C^*)^n.$$
Let $J=\exp(2\pi i\bq)\in (\C^*)^n$ be the exponential grading element \eqref{J-eq},
where $\bq=(q_1,\ldots,q_n)\in\Q^n$ with $q_j=d_j/d$. 
Note that by definition, $J\in\G_m=\Ga_0\sub \Ga$. Furthermore, $\chi(J)=1$, so
$J$ belongs to $G=\ker(\chi)\sub\Ga$.
The following result is essentially contained in \cite[Prop. 2.1.23, 2.2.8]{JKV}.

\begin{prop}\label{deg-prop} 
Let $g$ be the genus of $C$. 

\noindent
(i) One has the following identity in $G$:
\begin{equation}\label{gamma-J-eq}
\ga_1\cdot\ldots\cdot\ga_r=J^{2g-2+r}.
\end{equation}
Furthermore, consider the vector
$$\bdeg=(\deg L_1,\ldots,\deg L_n)\in\Z^n.$$ 
Then
\begin{equation}\label{deg-eq}
\bdeg=(2g-2+r)\bq-\bth_1-\ldots-\bth_r,
\end{equation}
where $\bth_s=\bth_{\ga_s}$.

\noindent
(ii) There exists a $\Ga$-spin structure on $(\CC,p_1,\ldots,p_r)$ of type $\ov{\ga}$ if and 
only if \eqref{gamma-J-eq} is satisfied.	

\noindent
(iii)
There exists a simple transitive action of the group $H^1(C,G)$ on the
set of isomorphism classes of $\Ga$-spin structures on $(\CC,p_1,\ldots,p_r)$ of type $\ov{\ga}$. 
In particular, if this set is nonempty, then it has $|G|^{2g}$ elements.
\end{prop}

\Pf . (i) Let $M=x_1^{k_1}\ldots x_n^{k_n}$ be a Laurent monomial of degree $d$
with respect to $\bfd$, such that $M|_G=1$, and
let $l_M:\Z^n\to\Z$ denote the corresponding linear form $\sum_j k_j e_j^*$.
Then it follows immediately from the definition \eqref{J-eq} that
\begin{equation}\label{l-M-J-eq}
l_M(\bq)=1.
\end{equation}
Recall that by \eqref{w-str-isom-eq}, the line bundle
$M(\LL_1,\ldots,\LL_n)$ is isomorphic to $\om^{\log}_{\CC}$.
Considering this isomorphism near each marked point we obtain that 
$$l_M(\bth_s)\in\Z \text{ for } s=1,\ldots,r.$$ 
We claim that the isomorphism \eqref{w-str-isom-eq} induces an isomorphism
\begin{equation}\label{M-L-C-main-isom}
M(L_1,\ldots,L_n)\simeq\om_C^{\log}(-l_M(\bth_1)p_1-\ldots-l_M(\bth_r)p_r).
\end{equation}
Indeed, let $p=p_s$ be one of the marked points, and let
$z$ be a local coordinate near $p$ on $\CC$, so that the generator $g_p$
of the local group $G(p)$ acts on $z$ by the multiplication with $\exp(-2\pi i/m)$.
Then we can view $z^m$ as a local coordinate near $p=p_s$ on $C$.
For each $j$ let $e_j(p)$ denote a generator of $\LL_j$ as an $\OO_{\CC}$-module near $p$.
For every $j=1,\ldots,n$, we have
$$g_p\cdot e_j(p)=\exp(2\pi i \th_{sj})\cdot e_j(p),$$
where $\bth_s=(\th_{s1},\ldots,\th_{sn})$. 
The line bundle $L_j=\rho_*\LL_j$ is generated near $p$
by $z^{m\th_{sj}}\cdot e_j(p)$.
The isomorphism \eqref{w-str-isom-eq}
implies that the action of $G(p)$ on $M(e_\bullet(p))$ is trivial. 
Hence, the line bundle
$\rho_*(M(\LL_\bullet))$ is generated near $p$ by $M(e_\bullet(p))$.
On the other hand, $L_j^{\ot k_j}$ is generated by
$z^{mk_j\th_{sj}}\cdot e_j(p)^{\ot k_j}$ near this point.
Thus, we have an isomorphism 
$$M(L_\bullet)\rTo{\sim} \rho_*(M(\LL_\bullet))(-\sum_{s=1}^r a_s p_s), \text{ where}$$
$$a_s=\sum_{j=1}^n k_j\th_{sj}.$$
Since $\rho_*(M(\LL_\bullet))\simeq \rho_*(\om_{\CC}^{\log})\simeq\om_C^{\log}$,
this gives \eqref{M-L-C-main-isom}.

Comparing the degrees in \eqref{M-L-C-main-isom}, we get a system of equations
$$l_M(\bdeg)=2g-2+r-l_M(\bth_1+\ldots+\bth_r),$$
where $M$ runs over Laurent monomials of degree $d$ with respect to $\bfd$ such that $M|_G=1$. 
Using \eqref{l-M-J-eq} we can rewrite this system as
$$l_M(\bdeg)=l_M\bigl((2g-2+r)\bq-\bth_1-\ldots-\bth_r)\bigr).$$
Now \eqref{deg-eq} follows from the fact that among $l_M$
there exist $n$ linearly independent forms:
it is enough to take the monomials $M$ corresponding to the characters 
\eqref{M-i-def-eq}. Since $\bdeg\in\Z^n$,
this also implies \eqref{gamma-J-eq}.

\noindent
(ii) Given $(\ga_1,\ldots,\ga_r)\in G^r$ satisfying \eqref{gamma-J-eq}, by 
Corollary \ref{quasihom-monom-cor},
we have to construct a collection of line bundles $(\LL_1,\ldots,\LL_n)$
on $\CC$, such that 
the action of the local group at every marked point $p_s$ on the fiber of $\bigoplus_j\LL_j$ at $p_s$ is
given by $\ga_s$, and isomorphisms
\begin{equation}\label{M-j-om-log-eq}
M_i(\LL_1,\ldots,\LL_n)\simeq\om_{\CC}^{\log} \text{ for }\ i=0,\ldots,n-1,
\end{equation}
where $M_0,\ldots,M_{n-1}$ are the Laurent monomials corresponding to the characters \eqref{M-i-def-eq}.
First, we claim that there exists a collection of line bundles $(L_1,\ldots,L_n)$ on $C$ together with isomorphisms
\eqref{M-L-C-main-isom} for $M=M_0,\ldots,M_{n-1}$. Indeed, since the characters \eqref{T-char-eq}
form a basis in the character lattice of $T$, the matrix of exponents of the monomials
$M_0,\ldots,M_{n-1}$ is nondegenerate. It is well known that the group $\Pic^0(C)$ of line bundles
of degree $0$ is divisible (recall that $C$ is smooth). Thus, it is enough to check that
the system of equations on degrees $\bdeg=(\deg(L_1),\ldots,\deg(L_n))$ imposed by
\eqref{M-j-om-log-eq} has a solution.
As we have seen in the proof of (i), 
if the condition \eqref{gamma-J-eq} is satisfied, then $\bdeg$ defined by \eqref{deg-eq}
gives a solution of this system. Finally, for $j=1,\ldots,n$, we define $\LL_j$ as the unique line bundle on $\CC$ with 
$\rho_*\LL_j\simeq L_j$ such that for every $s=1,\ldots,r$ the generator of the local group at $p_s$ acts
on the fiber of $\LL_j$ at $p_s$ by the $j$th component of $\ga_s$. As in part (i), this implies that
$$\rho_*M_i(\LL_1,\ldots,\LL_n)\simeq \om_C^{\log} \text{ for }\ i=0,\ldots,n-1,$$
and that the action of the local groups at the marked points on $M_i(\LL_1,\ldots,\LL_n)$ is trivial.
But this is equivalent to \eqref{M-j-om-log-eq}.

\noindent
(iii) Let $(\LL_1,\ldots,\LL_n)$ be a collection of line bundles on $\CC$ satisfying \eqref{M-j-om-log-eq},
such that the action of the local group on the fiber of $\bigoplus_j\LL_j$ at each marked point $p_s$ is given by
$\ga_s$. Any other such collection has the form $(\LL_1\ot\KK_1,\ldots,\LL_n\ot\KK_n)$,
where for each $j$ the local groups of the marked points act trivially on the fibers of $\KK_j$ and
$$M_i(\KK_1,\ldots,\KK_n)\simeq\OO_{\CC} \text{ for }\ i=0,\ldots,n-1.$$
Thus, different choices of $\LL_\bullet$ correspond to collections of line bundles $\KK_\bullet$ on $\CC$
such that $\KK_j=\rho^*\rho_*\KK_j$ for $j=1,\ldots,n$ and the line bundles $K_j=\rho_*\KK_j$ on $C$ satisfy
$$M_i(K_1,\ldots,K_n)\simeq\OO_C \text{ for }\ i=0,\ldots,n-1.$$
Since $G$ is the kernel of the homomorphism $\G_m^n\to\G_m^n$ given by $(M_0,\ldots,M_{n-1})$,
this is equivalent to a choice of a $G$-bundle.
Now the assertion follows from the fact that the group 
of isomorphism classes of $G$-bundles on $C$
is isomorphic $H^1(C,G)$.
\ed

As a consequence of the relation \eqref{deg-eq}, following \cite{FJR} we deduce the formula for the Euler characteristic of the bundle $\bigoplus_j L_j$. For $\ga\in G$ define the {\it degree shifting number} $\iota_{\ga}$ 
as the sum of coordinates of the vector $\bth_{\ga}-\bq$.

\begin{cor} One has
\begin{equation}\label{index-eq}
-\sum_{j=1}^n\chi(C,L_j)=D_g(\ga_1,\ldots,\ga_r)=(g-1)\hat{c}+\iota_{\ga_1}+\ldots+\iota_{\ga_r},
\end{equation}
where
$$\hat{c}=\sum_{j=1}^n (1-2q_j).$$
\end{cor}

We will also use the modified quantities
\begin{equation}\label{mod-index-eq}
\wt{D}_g(\ga_1,\ldots,\ga_r)=D_g(\ga_1,\ldots,\ga_r)+\frac{1}{2}\cdot \sum_{i=1}^r N_{\ga_i}
\end{equation}
with 
$$N_{\ga}=\dim (\A^n)^{\ga}.$$ 
They satisfy the following factorization properties.

\begin{lem}\label{mod-index-factor-lem} 
For any $\ov{\ga}=(\ga_1,\ldots,\ga_r)\in G^r$, $\ov{\ga'}=(\ga'_1,\ldots,\ga'_{r'})\in G^{r'}$ and $\ga\in G$ one has
\begin{equation}\label{mod-index-tree-eq}
\wt{D}_{g_1}(\ov{\ga},\ga)+\wt{D}_{g_2}(\ov{\ga'},\ga^{-1})=\wt{D}_{g_1+g_2}(\ov{\ga},\ov{\ga'}),
\end{equation}
\begin{equation}\label{mod-index-loop-eq}
\wt{D}_g(\ov{\ga},\ga,\ga^{-1})=\wt{D}_{g+1}(\ov{\ga}).
\end{equation}
\end{lem}

\Pf . This follows immediately from the simple relation 
$$\iota_{\ga}+\iota_{\ga^{-1}}=\hat{c}-N_{\ga}$$
established in \cite[Prop.\ 3.2.4]{FJR}.
\ed

\section{Matrix factorizations from $\w$-structures}\label{main-sec}

\subsection{$\w$-structures with respect to the canonical bundle}
\label{canonical-sec}

Let $\TT$ be a symmetric monoidal category with split
projectors over a field of characteristic zero.
To a monomial $M(x_1,\ldots,x_n)=x_1^{m_1}\ldots x_n^{m_n}$
we associate a polyfunctor on $\TT^n$ 
$$M_{\TT}:\TT^n\to \TT:
(A_1,\ldots,A_n)\mapsto M(A_\bullet):=S^{m_1}A_1\otimes\ldots\otimes S^{m_n}A_n,$$
where $S^m(?)$ denote the symmetric powers in $\TT$.
In particular, we are going to use this operation in the case when $\TT$ is the derived category
of coherent sheaves with the monoidal structure given by the derived tensor product.

Fix a quasihomogeneous polynomial 
\begin{equation}\label{monom-index-eq}
\w(x_1,\ldots,x_n)=\sum_{k=1}^N c_k M_k.
\end{equation}
Let
$\pi:C\to S$ be a family of nodal curves with a weak $\w$-structure $(\Phi,\phi)$
in the category of line bundles over $C$ with respect to the relative canonical bundle $\om_{C/S}$.
This structure can be specified by a collection of line bundles $L_j=\Phi(e_j)$, $j=1,\ldots,n$
and a collection of morphisms 
$$\phi_k:M_k(L_\bullet)\to\om_{C/S}$$
for every monomial $M_k$ appearing in $\w$.
Assume that the restriction of each $\phi_k$ to any fiber of $\pi$
is an isomorphism outside a finite number of points. 

For each $k=1,\ldots,N$, the morphism $\phi_k$
induces a morphism in the derived category $D(S)$
\begin{equation}\label{main-map}
M_k(R\pi_*(L_1),\ldots,R\pi_*(L_n))
\to R\pi_*(\om_{C/S})\to\OO_S[-1].
\end{equation}
Assume that each $R\pi_*(L_j)$ is represented by a complex 
\be\label{A-i-B-i-eq}
A_j\stackrel{\b_j}{\to} B_j
\end{equation}
of vector bundles on $S$, in such a way that the morphism \eqref{main-map} is realized on
the level of complexes. Recall that the $m$-th symmetric power of \eqref{A-i-B-i-eq} is the complex
$$S^{m}A_j\to S^{m-1}A_j\otimes B_j\to\ldots$$
concentrated in the degrees $0,1,\ldots,m$.
Therefore, the source of the map \eqref{main-map} is represented by the complex
$$M_k(A_\bullet)\stackrel{\de}{\to} \bigoplus_{j=1}^n B_j\ot\del_jM_k(A_\bullet)\to\ldots,$$
where $\del_jM_s$ are the ``partial derivatives" of the monomial 
$M_k=x_1^{m_{k1}}\ldots x_n^{m_{kn}}$:
$$\del_1M_k=x_1^{m_{k1}-1}x_2^{m_{k2}}\ldots x_n^{m_{kn}}, \ldots,
\del_nM_k=x_1^{m_{k1}}x_2^{m_{k2}}\ldots x_n^{m_{kn}-1}.
$$ 
The differential $\de$ is given by
$$\de(f_1\ot\ldots\ot f_n)=(\de_1(f_1)\ot f_2\ldots\ot f_n,\ldots,f_1\ot\ldots\ot \de_n(f_n)),$$
where $f_j\in S^{m_{kj}}A_j$, and $\de_j:S^{m_{kj}}A_j\to B_j\ot S^{m_{kj}-1}A_j$ is induced
by the Koszul differential $S^{m_{kj}}A_j\to A_j\ot S^{m_{kj}-1}A_j$ and by the map 
$\b_j:A_j\to B_j$, i.e.,
$$\de_j(a_j^{m_{kj}})=m_{kj} \b_j(a_j)\ot a_j^{m_{kj}-1}.$$
By our assumption, the map \eqref{main-map} is realized as a chain map of complexes, so
we have a map 
$$\a_k:\bigoplus_j B_j\ot\del_jM_k(A_\bullet)\to \OO_S$$
such that $\a_k\circ\de=0$.
The components of $\a_k$ can be viewed as morphisms
$$\a_{kj}:\del_j M_k(A_\bullet)\to B_j^{\vee},$$
and the condition $\a_k\circ\de=0$ can be expressed as follows:
$$\sum_{j=1}^n m_{kj} \lan \a_{kj}(a_1^{m_{k1}}\ot\ldots\ot a_j^{m_{kj}-1}\ot\ldots\ot 
a_n^{m_{kn}}), \b_j(a_j)\ran=0,
$$
where $\lan\cdot,\cdot\ran$ is the evaluation pairing.
Let $p:X\to S$ be the total space of the vector bundle $A_1\oplus\ldots\oplus A_n$ over $S$.
Then we can view the maps $\a_{kj}$ as sections of the induced bundles 
$p^*B_j^{\vee}$ on $X$ and the maps $\b_j$ as sections
of $p^*B_j$ on $X$. 
The above equation can be viewed as the following identity of functions on $X$:
\be\label{orthog-eq-1}
\sum_{j=1}^n m_{kj}\lan\a_{kj}, \b_j\ran=0.
\end{equation}
Now let us define sections $\a_j(\w)\in\Gamma(X,p^*B_j^{\vee})$ by setting
$$\a_j^{\w}=\sum_{k=1}^N c_k m_{kj}\a_{kj}.$$
Then from \eqref{orthog-eq-1} we obtain that
$$\sum_{j=1}^n \lan \a_j^{\w}, \b_j\ran=0,$$
so the sections
$$\a=(\a_1^{\w},\ldots,\a_n^{\w})\in \bigoplus_{j=1}^np^*B_j^{\vee} \ \text{ and }
\b=(\b_1,\ldots,\b_n)\in \bigoplus_{j=1}^np^*B_j$$
satisfy $\lan\a,\b\ran=0$.
Recall that $X$ is the total space of the vector bundle $\oplus_i A_i$ over $S$, so
it contains $S$ as the zero section.

\begin{prop}\label{zero-locus-prop}
Suppose that $\w$ has an isolated singularity at the origin. Then the common vanishing locus
of $\a$ and $\b$, $Z(\a,\b)$, coincides with the zero section in $X$. 
\end{prop}

\Pf . It is obvious that $Z(\a,\b)$ contains the zero section $S\sub X$.
It is enough to check the opposite inclusion in the case when $S$ is a point, i.e.,
$\CC=C$ is a single curve $C$ and $X$ is a vector space.
To compute $Z(\a,\b)$, first observe that $\ker(\b_j)$ is isomorphic to $H^0(C,L_j)$ (resp.,
$\coker(\b_j)\simeq H^1(C,L_j)$). For each monomial $M_k=x_1^{m_{k1}}\ldots x_n^{m_{kn}}$
the map \eqref{main-map} induces morphisms
$$R\Ga(C,L_j)\ot\del_jM_k(R\Ga(C,L_\bullet))\to \C[-1]$$
and in particular, well-defined morphisms
$$H^1(C,L_j)\ot\del_jM_k(H^0(C,L_\bullet))\to \C.$$
These maps are also induced by the maps $\a_{kj}$, so we deduce that
the restriction of $\a_{kj}$ to $\del_jM_k(H^0(L_\bullet))$ factors through
the embedding $\kappa_{kj}:H^1(L_j)^{\vee}\hookrightarrow B_j^{\vee}$. Note that by Serre duality
$H^1(L_j)^{\vee}\simeq\Hom(L_j,\om_C)$, and the map
$$\del_j M_k(H^0(L_\bullet))\to \Hom(L_j,\om_C)$$
corresponding to $\a_{kj}$ is the composition of the product map
$$\del_j M_k(H^0(L_\bullet))\to H^0(\del_jM_k(L_\bullet))$$
with the natural map 
$$\phi_{kj}:H^0(\del_jM_k(L_\bullet))\to\Hom(L_j,\om_C)$$
induced by $\phi_k$.
Since $\kappa_{kj}$ is an embedding, 
the condition that $\a_j^{\w}({\bf a})=0$ where ${\bf a}=(a_1,\ldots,a_n)$ and $a_j\in H^0(L_j)$, 
can be rewritten as
\begin{equation}\label{der-eq}
\sum_{k=1}^N c_k \phi_{kj}(\del_jM_k({\bf a}))=0,
\end{equation}
where the expression
$$M_k({\bf a})=a_1^{m_{k1}}\ot \ldots\ot a_n^{\ot m_{kn}}$$ 
is a section of $M_k(L_\bullet)=L_1^{m_{k1}}\ldots L_n^{m_{kn}}$. 

We know that all $\phi_k$ are isomorphisms at
the generic point $x\in C$ of any irreducible component of $C$.
Hence, there exist a collection of trivializations of $L_j$'s at $x$
and a trivialization of $\om_x$ such that each map $\phi_k$ becomes equal to the identity
under these trivializations. This follows easily from Proposition \ref{W-triv-cor}, since
the restriction of our weak $\w$-structure to $x$ becomes a $\w$-structure with respect to
the unit object.
Using these trivializations we can view ${\bf a}(x)=(a_1(x),\ldots,a_n(x))$ as 
a point in $\C^n$. Now the left-hand side of \eqref{der-eq} evaluated at $x$ 
becomes the $j$th partial derivative of $\w$ at ${\bf a}(x)$ (see \eqref{monom-index-eq}).
Thus, $\del_j\w({\bf a}(x))=0$ for all $k=1,\ldots,n$. This implies
that ${\bf a}(x)=0$, since $\w$ has an isolated singularity at $0$. Therefore, ${\bf a}=0$ on
a dense subset of $C$ and since $L_i$ have no torsion, this implies that ${\bf a}=0$.
So $Z(\a,\b)$ is contained in $S$. 
\ed

\subsection{Fundamental matrix factorizations}\label{constr-sec}

Let $\w$ be a quasihomogeneous polynomial with an isolated singularity, $G\sub G_\w$ a subgroup
containing $J$ and $\Ga\sub\Ga_\w$ the corresponding extension of $\G_m$ by $G$ 
(see Lemma \ref{J-lem}).
A key role in our construction of the \CFT\ associated with $\w$ and $G$ 
will be played by certain matrix factorizations of $-p^*\w_{\ov{\ga}}$
on $\SS_{g,\Ga}^{\rig}(\ov{\ga})\times \A^{\ov{\ga}}$, where 
$\A^{\ov{\ga}}=\prod_{i=1}^r (\A^n)^{\ga_i}$ and
$\w_{\ov{\ga}}=\sum_{i=1}^r \pr_i^*\w_{\ga_i}$.
Here $p:\SS_{g,\Ga}^{\rig}(\ga_1,\ldots,\ga_r)\times \A^{\ov{\ga}}\to\A^{\ov{\ga}}$ and
$\pr_i:\A^{\ov{\ga}}\to (\A^n)^{\ga_i}$
are the projections.
This collection of objects can be viewed as a categorified
analog of the fundamental class in Gromov-Witten theory.\footnote{Such a categorified version of
the fundamental class exists in GW-theory as well and can be constructed as in \cite{Lee}}

The construction roughly goes as follows. First, using the universal $\w$-structure
over $\SS_g^{\rig}(\ov{\ga})=\SS_{g,\Ga}^{\rig}(\ov{\ga})$ provided by Proposition \ref{spin-prop}(iv)
(together with the rigidification structure), we construct a Koszul \mf\ $\{\a,\b\}$ of the pull-back of 
$-p^*\w_{\ov{\ga}}$
on a certain affine bundle $X\to \SS_g^{\rig}(\ov{\ga})\times \A^{\ov{\ga}}$.
We prove that this \mf\ is supported on a section of this affine bundle over
$\SS_g^{\rig}(\ov{\ga})\times\{0\}$. This allows us to apply the construction of Section \ref{push-forward-sec}
to define the fundamental \mf\  
$\bP^{\rig}_g(\ov{\ga})\in \ov{\DMF}_{\Ga_\w}(\SS_g^{\rig}(\ov{\ga})\times \A^{\ov{\ga}}, -p^*\w_{\ov{\ga}})$ 
as the
push-forward of $\{\a,\b\}$. The space $X$ is the total space
of the vector bundle $\bigoplus_{j=1}^n A_j$ over $\SS_g^{\rig}(\ov{\ga})$
for an appropriate choice of
resolutions $[d_j:A_j\to B_j]$ of the derived push-forwards of the line bundles associated with
the $\Ga$-spin structure. The element $\b$ is the section of the pull-back of 
$\bigoplus B_j$ to $X$ induced by the differentials $d_j$, 
while the map $X\to \A^{\ov{\ga}}$ is induced by the rigidification
structure. The construction of the element $\a$, which is a section of $\bigoplus B_j^\vee$,
is more involved.
First, we use the isomorphisms \eqref{w-str-isom-eq} of the universal $\w$-structure to construct the
data needed for $\a$ on the level of derived categories. Then we use quasi-projectivity
of the coarse moduli space to lift these data to the level of complexes.

We start with a family of $\Ga$-spin curves over a base $S$ of type $\ov{\ga}=(\ga_1,\ldots,\ga_r)\in G^r$. 
Recall (see Definition \ref{fam-spin-def}) 
that this is a family $\pi:\CC\to S$ of nodal orbicurves with marked orbipoints
$p_1\hra \CC,\ldots, p_r\hra \CC$ and a principal $\Ga$-bundle $P$ together with an
isomorphism $\vareps:\chi_*P\simeq\varphi_*P(\om^{\log}_{\CC/S})$. Let $(\LL_1,\ldots,\LL_n)$ 
be a collection of line bundles on $\CC$ associated with $P$. Then $\vareps$
induces isomorphisms 
\begin{equation}\label{phi-M-isom-eq}
\phi_M:M(\LL_\bullet)=M(\LL_1,\ldots,\LL_n)\to\om_{\CC}^{\log}
\end{equation}
for every monomial $M$ occurring in $\w$ (see Remark \ref{quasihom-monom-rem}.2).
We also assume that our $\Ga$-spin structure $(P,\vareps)$ 
is equipped with a restricted rigidification (see Definition \ref{restr-rig-def}).
Note that this notion is well-defined in our 
situation by Lemma \ref{admissibility-lem}.

Let $(C\to S,p_1,\ldots,p_r)$ be the family of orbicurves with marked points,
obtained by forgetting the orbistructure at $p_i$'s 
(and keeping the orbistructure at the nodes),
and let $\rho:\CC\to C$ be the natural projection. The family
$(C,p_1,\ldots,p_r)$ can be constructed as follows. Present $\wt{C}$ as the 
union of two open substacks $\wt{C}^{reg}$
and $\wt{C}^0$, obtained by taking the complements of the nodes and of the marked points, respectively.
Then take the coarse moduli $C^{reg}$ of $\wt{C}^{reg}$ with the induced marked points (see
\cite{KM} or \cite{Conrad}). Note that by universality, $C^{reg}$ still maps to $S$. Our family $C$ is obtained by gluing $C^{reg}$ with $\wt{C}^0$.

For each $j=1,\ldots,n,$ the push-forward $L_j=\rho_*\LL_j$ 
is a line bundle on $C$. 
By abuse of notation we will denote the projection $C\to S$ also by $\pi$. Note that
we have natural isomorphisms $R\pi_*(\LL_j)\simeq R\pi_*(L_j)$. 

Let us introduce some more notation.
The $r$-tuple $\ov{\ga}=(\ga_1,\ldots,\ga_r)\in G^r$ determines a relation
\begin{equation}\label{Sigma-eq}
\Sigma=\Sigma(\ov{\ga})=\{(p_i,j)\ |\ \ga_{ij}=1\}\sub \{p_1,\ldots,p_r\}\times\{1,\ldots,n\}.
\end{equation}
For a $\Ga$-spin structure of type $\ov{\ga}$ one has $(p_i,j)\in\Sigma$ if and only if the action of
$G(p_i)$ on $\LL_j|_{p_i}$ is trivial.

Note that the cross-section of $\Sigma$ at $p_i$ coincides with the set $I(\ga_i)=\{j \ |\ \gamma_{ij}=1\}$ labeling the coordinates of the subspace of $\ga_i$-invariants $(\A^n)^{\ga_i}\subset\A^n$
Hence, 
the affine space
$$\A^{\ov{\ga}}=\prod_{i=1}^r(\A^n)^{\ga_i}$$
has coordinates $x_j(i)$ labeled by $(p_i,j)\in\Si$. 

Let us equip the space $\A^{\ov{\ga}}$ with the potential
$$\w_{\ov{\ga}}=\sum_{i=1}^r\pr_i^*\w_{\ga_i},$$ 
where $\w_{\ga_i}$ is the restriction of $\w$ to $(\A^n)^{\ga_i}$ and
$\pr_i:\A^{\ov{\ga}}\to(\A^n)^{\ga_i}$ is the projection.

Let 
\begin{equation}\label{Sigma-j-eq}
\Si_j=\{p_i\ |\ (p_i,j)\in\Si\}\sub\{p_1,\ldots,p_r\}
\end{equation}
be the 
cross-section of $\Si$ at $j$,  
and for each monomial $M$ 
in $\w$ let us set 
$$\Si_M=\bigcap_{\deg_{x_j}\!M>0}\Si_j$$
where $\deg_{x_j}$ is the degree in $x_j$. Thus, we take intersection over
those $j$ for which $x_j$ occurs in $M$.

\begin{lem}\label{W-Sigma-lem} We have
$$\w_{\ga_i}=\sum_{k: i\in\Si_{M_k}}c_k\cdot M_k \text{ and}$$
$$\w_{\ov{\ga}}=\sum_{k=1}^N c_k\cdot M_k^{\oplus \Si_{M_k}},$$
where $\displaystyle M^{\oplus \Si_{M}}=\sum_{i\in\Si_{M}}M(x_1(i),\ldots,x_n(i))$.
\end{lem}

\Pf . By definition, a monomial $M$ of $\w$ occurs in $\w_{\ga_i}$ if and only if
$\ga_i$ acts trivially on all variables $x_j$ such that $\deg_{x_j}M>0$, i.e., 
if and only if $i\in\Si_M$. 
\ed

The following 
important observation will eventually lead to a connection with the function $\w_{\ov{\ga}}$.

\begin{lem}\label{W-orbi-lem} 
For each monomial $M$ occurring in $\w$
the isomorphism \eqref{phi-M-isom-eq} induces 
an injective morphism
\begin{equation}\label{over-phi-eq}
\ov{\phi}_M:M(L_1,\ldots,L_n)\to \om_{C}(\Sigma_M)
\end{equation}
on $C$, which is an isomorphism in a neighborhood of $\Sigma_M$. In the above formula we view
$\Sigma_M\sub \{p_1,\ldots,p_r\}$ as an effective divisor on $C$.
\end{lem}

\Pf . To simplify notation we consider the case of a single curve (i.e., $S$ is a point).
As we have seen in the proof of Proposition \ref{deg-prop}(i),
the map 
$$\phi'_M:M(L_\bullet)\to \rho_*(M(\LL_\bullet))\simeq\om_C(p_1+\ldots+p_r),$$
induced by $\phi_M$, vanishes at $p_i$ to the order 
$$\frac{1}{m}\sum_{j=1}^n k_j\th_{ij}\ge 0,$$
where we use the notation of Proposition \ref{deg-prop}(i).
Furthermore, $\phi'_M$ is an isomorphism near $p_i$ if and only if
$\th_{ij}=0$ for all $j$ with $\deg_{x_j}M>0$. 
This immediately implies the assertion since $(p_i,j)\in\Sigma$ if and only if $\th_{ij}=0$.
\ed

The restricted rigidification structure on $(P,\vareps)$ induces trivializations
$$e_j(i):\OO_{p_i}\to \LL_j|_{p_i} \text{ for }\ (p_i,j)\in\Si,$$
such that for for every monomial $M$ occurring in $\w$ 
and every $p_i\in\Si_M$ the induced morphism
$$\OO_{p_i}\rTo{M(e_\bullet(i))} M(L_\bullet)\rTo{\ov{\phi}_M}\om_C(\Sigma_M)|_{p_i}\simeq\OO_{p_i}$$
is the identity.
For $j=1,\ldots,n$ let 
\begin{equation}\label{e-j-Si-j-eq}
\bfe_j:\OO_S^{\Si_j}\to \pi_*(L_j|_{\Si_j})
\end{equation}
denote the isomorphism induced by $(e_j(i))$.
From now on we will work exclusively with the data $(L_1,\ldots,L_n;\bfe_1,\ldots,\bfe_n)$
on the family of orbicurves $C/S$ that has trivial orbi-structure at the marked points $p_1,\ldots,p_r$,
together with the morphisms \eqref{over-phi-eq}.

Note that since $\Si_M\sub\Si_j$ 
for every $j=1,\ldots,n,$  for every $j$ with $\deg_{x_j}M>0$,
the map $\ov{\phi}_M$ induces a morphism
\begin{equation}\label{M-L-Si-om-eq}
M(L_\bullet(-\Si_\bullet))=M(L_1(-\Si_1),\ldots,L_n(-\Si_n))\to\om_{C/S},
\end{equation}
where we view $\Sigma_j$ as a subdivisor of $p_1+\ldots+p_r$. 
The collection $(L_j(-\Si_j))$ equipped with morphisms \eqref{M-L-Si-om-eq} is
a weak $\w$-structure with respect to $\om_{C/S}$.
Therefore, the construction of Section 
\ref{canonical-sec} gives a natural morphism
\begin{equation}\label{t-M-eq}
t_M: M(R\pi_*(L_\bullet(-\Si_\bullet)))\to\OO_S[-1]
\end{equation}
induced by \eqref{M-L-Si-om-eq} (see \eqref{main-map}).

For each $j$ we have an exact triangle in $D^b(S)$
$$R\pi_*(L_j(-\Sigma_j))\to R\pi_*(L_j)\rTo{r_j} \pi_*(L_j|_{\Sigma_j})\to\ldots
$$
The isomorphism \eqref{e-j-Si-j-eq} gives 
 a map
\begin{equation}\label{Z-j-eq}
Z_j:R\pi_*(L_j)\to \OO_S^{\Sigma_j},
\end{equation}
such that 
$$r_j=\bfe_j\circ Z_j.$$
Let us denote by 
$$Z_M:M(R\pi_*(L_\bullet))\rTo{M(Z_\bullet)} 
M(\OO_S^{\Sigma_\bullet})\to M(\OO_S^{\Si_M},\ldots,\OO_S^{\Si_M})\to \OO_S^{\Sigma_M}$$
the map induced by the composition
$R\pi_*(L_j)\stackrel{Z_j}{\to}\OO_S^{\Si_j}\to\OO_S^{\Si_M}$ and
by the algebra structure on $\OO_S^{\Si_M}$.

We split the construction of the fundamental \mf\ into four steps.
First, we extend the morphisms \eqref{t-M-eq} to morphisms
$\tau_M:E_M\to\OO_S[-1]$ where $E_M\in D^b(S)$ are certain modifications of 
$M(R\pi_*(L_\bullet))$. In Step 2 we realize these morphisms on the level of complexes
using appropriate resolutions $[A_j\to B_j]$ of $R\pi_*(L_j)$. At the same time we realize morphisms
$Z_j$ by maps of vector bundles $A_j\to \OO_S^{\Si_j}$. These maps combine into a morphism
$Z:X\to\A^{\ov{\ga}}$, where $X$ is the total space of the bundle $\bigoplus_j A_j$.
In Step 3 we construct a Koszul \mf\ $\{\a,\b\}$ of $-Z^*\w_{\ov{\ga}}$ and in Step 4 we show that
it can be pushed forward to a \mf\ of $-p^*\w_{\ov{\ga}}$ on $S\times\A^{\ov{\ga}}$, where $p$
is the projection $S\times\A^{\ov{\ga}}\to \A^{\ov{\ga}}$.

\

\noindent
{\bf Step 1}. For each monomial $M$ appearing in $\w$ we construct 
a canonical commutative diagram in $D^b(S)$ with an exact triangle as a middle row
\begin{equation}\label{main-diagram-eq}
\begin{diagram}
M(R\pi_*(L_\bullet(-\Si_\bullet)))\\
&\rdTo{\iota_M}\rdTo(4,2)^{\nat}\\
\dTo{t_M}&&E_M&\rTo{\eps}& M(R\pi_*(L_\bullet))&\rTo{Z_M}& \OO_S^{\Si_M}&\rTo{i}& E_M[1]\\
&\ldTo{\tau_M}&&&&&&\rdTo_{\Tr}&\dTo_{\tau_M[1]}\\
\OO_S[-1]&&&&&&&&\OO_S
\end{diagram}
\end{equation}
where 
$$\nat:M(R\pi_*(L_\bullet(-\Si_\bullet)))\to M(R\pi_*(L_\bullet))$$
is the natural map and
the map $\Tr:\OO_S^{\Si_M}\to\OO_S$ is given by the sum of components,
and $t_M$ is the morphism \eqref{t-M-eq}.

Let $M=x_1^{k_1}\ldots x_n^{k_n}$ and set $|M|=\deg M=k_1+\ldots+k_n$.
Consider the relative $|M|$th power $\pi^M:C^M:=C^{|M|}\to S$ of $C$ over $S$ and define
the line bundle $L_\bullet^M$ on $C^M$ by
$$L_\bullet^M:=L_1^{\boxtimes k_1}\boxtimes\ldots\boxtimes L_n^{\boxtimes k_n}.$$
Also, consider the product of symmetric groups 
\begin{equation}\label{Sym-M-eq}
\Sym(M)=S_{k_1}\times\ldots\times S_{k_n}.
\end{equation}
From the K\"unneth isomorphism we get an identification 
$$M(R\pi_*(L_\bullet))\simeq \left(R\pi^M_*(L_\bullet^M)\right)^{\Sym(M)},$$
where we take invariants with respect to the natural $\Sym(M)$-action.
Let 
$$\si_M:\Si_M\to C\rTo{\De_M} C^M$$
be the closed embedding,
where $\De_M$ is the diagonal map, and let $\JJ_M\sub\OO_{C^M}$ denote
the ideal sheaf of $\si_M(\Si_M)$.
Define a $\Sym(M)$-equivariant coherent sheaf $\FF_M$ on $C^M$ by
\begin{equation}\label{F-M-eq}
\FF_M=\JJ_M L_{\bullet}^M
\end{equation}
and set
$$E_M=\left(R\pi^M_*(\FF_M)\right)^{\Sym(M)}.$$
Note that we have an exact sequence
\begin{equation}\label{F-M-seq}
0\to\FF_M\to L_\bullet^M\to (\si_M)_*\si_M^*L_\bullet^M\to 0.
\end{equation}
The map $\ov{\phi}_M$ induces an isomorphism 
$$\si_M^*L_\bullet^M\simeq M(L_\bullet)|_{\Si_M}\simeq \OO_{\Si_M}$$ on
$\Si_M$, which coincides with the isomorphism obtained from the trivializations
of $L_j|_{\Si_j}$. Hence, taking the push-forward of the sequence \eqref{F-M-seq} to $S$, and
considering $\Sym(M)$-invariants
we get the horizontal exact triangle of the diagram \eqref{main-diagram-eq}.

We have an embedding of sheaves on $C^M$
\begin{equation}\label{L-M-FF-emb}
\iota_M^C:L_\bullet(-\Si_\bullet)^M\hra \FF_M.
\end{equation}
After taking the push-forward to $S$ and passing to
the $\Sym(M)$-invariants this embedding induces the map
$\iota_M:M(R\pi_*(L_\bullet(-\Si_\bullet)))\to E_M$.

The map $\tau_M$ is constructed similarly
using the morphism of sheaves
\begin{equation}\label{kappa-M-eq}
\kappa_M:\FF_M\to(\De_M)_*\om_{C/S}
\end{equation}
defined as follows. 
Let $\JJ_{\De}\sub\OO_{C^M}$
be the ideal sheaf of the diagonal $\De_M(C)\subset C^M$. 
Let $\psi_M:\De_{M}^*\FF_M\to\om_{C/S}$ be the composition of
the isomorphism
\begin{equation}\label{De-FF-M-eq}
\De_{M}^*\FF_M\simeq \JJ_M L_\bullet^M/\JJ_{\De} L_\bullet\simeq
M(L_\bullet)(-\Si_M)
\end{equation}
with the map $\ov{\phi}_M:M(L_\bullet)(-\Si_M)\to \om_{C/S}$.
Now
we define $\kappa_M$ as the morphism $\FF_M\to(\De_M)_*\om_{C/S}$
corresponding to $\psi_M$ by adjunction.

\begin{lem}\label{L-FF-M-kappa-lem} 
The composition $\kappa_M\circ \iota_M^C$ of
morphisms \eqref{L-M-FF-emb} and \eqref{kappa-M-eq} 
coincides with the map $L_\bullet(-\Si_\bullet)^M\to (\De_M)_*\om_{C/S}$ 
corresponding by adjunction to the map
\begin{equation}\label{De-L-M-om-eq}
\De_M^*L_\bullet(-\Si_\bullet)^M\simeq M(L_\bullet(-\Si_\bullet))\to\om_{C/S}.
\end{equation}
given by \eqref{M-L-Si-om-eq}.
\end{lem}

\Pf . This follows by adjunction from the fact that 
the map \eqref{De-L-M-om-eq} coincides with the composition of $\De_M^*\iota_M^C$
with the map 
$$\psi_M:\De_M^*\FF_M\simeq M(L_\bullet)(-\Si_M)\rTo{\ov{\phi}_M} \om_{C/S}.$$
\ed

We have a morphism of exact sequences
\begin{diagram}
0 & \rTo{} & \FF_M &\rTo{} &L_\bullet^M &\rTo{}& (\si_M)_*\OO_{\Si_M} &\rTo{} & 0\\
&&\dTo{\kappa_M} &&\dTo{}&&\dTo{\id}\\
0 & \rTo{} &(\De_M)_*\om_{C/S}&\rTo{}&(\De_M)_*\om_{C/S}(\Sigma_M)& \rTo{}& 
(\si_M)_*\OO_{\Si_M}&\rTo{} & 0
\end{diagram}
where the middle vertical arrow corresponds to $\ov{\phi}_M$ by adjunction.
It induces a morphism of exact triangles
\begin{diagram}
E_M&\rTo{\eps}&M(R\pi_*(L_\bullet))&\rTo{Z_M}&\OO_S^{\Sigma_M}&
\rTo{i}&E_M[1]\\
\dTo{\tau'_M}&&\dTo{}&&\dTo{\id}&&\dTo{\tau'_M}\\
R\pi_*(\om_{C/S})&\rTo{}&R\pi_*(\om_{C/S}(\Sigma_M))&\rTo{}&\OO_S^{\Sigma_M}&
\rTo{}& 
R\pi_*(\om_{C/S})[1]
\end{diagram}
On the other hand, we have a morphism of exact triangles
\begin{diagram}
R\pi_*(\om_{C/S})&\rTo{}&R\pi_*(\om_{C/S}(\Sigma_M))&\rTo{}&\OO_S^{\Sigma_M}&
\rTo{}& 
R\pi_*(\om_{C/S})[1]\\
\dTo{\Tr_{C/S}}&&\dTo{}&&\dTo{\id}&&\dTo{\Tr_{C/S}[1]}\\
\OO_S[-1]&\rTo{}&[\OO_S^{\Si_M}\to\OO_S]&\rTo{}&\OO_S^{\Si_M}&\rTo{\Tr}&\OO_S
\end{diagram}
where $\Tr_{C/S}:R\pi_*(\om_{C/S})\to\OO_S[-1]$
is the Grothendieck trace map (see \cite{Hart-RD}). 
Composing these two morphisms of exact triangles we get a diagram
\begin{equation}\label{morphism-triang-eq}
\begin{diagram}
E_M&\rTo{\eps}&M(R\pi_*(L_\bullet))&\rTo{Z_M}&\OO_S^{\Sigma_M}&
\rTo{i}&E_M[1]\\
\dTo{\tau_M}&&\dTo{}&&\dTo{\id}&&\dTo{\tau_M[1]}\\
\OO_S[-1]&\rTo{}&[\OO_S^{\Si_M}\to\OO_S]&\rTo{}&\OO_S^{\Si_M}&\rTo{\Tr}&\OO_S
\end{diagram}
\end{equation}
and in particular, the canonical map $\tau_M:E_M\to \OO_S[-1]$.
This finishes the construction of the diagram \eqref{main-diagram-eq}.
The equality 
\begin{equation}\label{tau-t-M-eq}
t_M=\tau_M\circ\iota_M
\end{equation}
(the commutativity of the leftmost triangle in
the diagram \eqref{main-diagram-eq}) follows from Lemma \ref{L-FF-M-kappa-lem}.

Note that when $\Sigma_M=\emptyset$ 
the map $\tau_M$ coincides with
the map \eqref{main-map} constructed from the data $(L_1,\ldots,L_n,\ov{\phi}_M)$.

\

\noindent
{\bf Step 2}. Next, we are going to realize the diagram \eqref{main-diagram-eq}
on the level of complexes. 
More precisely, we will represent each $R\pi_*(L_j)$ by
a complex $K_j=[A_j\to B_j]$ of vector bundles on $S$, concentrated
in degrees $[0,1]$, in such a way that the map $Z_j:R\pi_*(L_j)\to\OO_S^{\Sigma_j}$ is realized by 
a surjective chain map of complexes 
$Z_j:K_j\to \OO_S^{\Sigma_j}$. Then the subcomplex
$K'_j=\ker(Z_j)=[A'_j\to B_j]$ will represent $R\pi_*(L_j(-\Si_j))$ and the map 
$\nat:M(K'_\bullet)\hra M(K_\bullet)$ will be the natural inclusion.
For each monomial $M$ appearing in $\w$ 
the map $Z_M:M(K_\bullet)\to \OO_S^{\Sigma_M}$ 
will be realized by the composition of 
$$M(Z_\bullet):M(K_\bullet)\to M(\OO_S^{\Si_1},\ldots,\OO_S^{\Si_n})$$ 
with the natural epimorphism $M(\OO_S^{\Si_1},\ldots,\OO_S^{\Si_n})\to \OO_S^{\Sigma_M}$.
Also, the complex $K_M=\Cone(Z_M)[-1]$ 
will represent the object $E_M\in D^b(S)$ and the maps $\tau_M$ and $\iota_M$ 
will be realized by chain map of complexes, so that the diagram
\eqref{main-diagram-eq} will be commutative in the category of complexes.

To construct appropriate complexes representing
$R\pi_*(L_j(-\Sigma_j))$ and $R\pi_*(L_j)$ we will need the following result.

\begin{lem}\label{generation-lem} Let $T$ be a proper smooth 
DM-stack over $\C$ with projective coarse moduli $\ov{T}$.
Let $p:T\to\ov{T}$ be the projection, 
$\OO(1)$ an ample line bundle on $\ov{T}$ and $\OO_T(1):=p^*\OO(1)$.

(i) There exists a vector bundle $\VV$ on $T$ such that for any coherent sheaf $\FF$ on $T$ the
natural map
\begin{equation}\label{gen-mor-eq}
H^0(T,\VV^\vee\ot \FF(n))\ot\VV(-n)\to\FF
\end{equation}
is surjective for $n\gg 0$. Also, 
$$H^{>0}(T,\FF(n))=0$$
for $n\gg 0$.

(ii) Let $\pi:C\to T$, $p_1,\ldots,p_r:T\to C$ be a family of stable curves with marked points over $T$.
Then for every vector bundle $E$ on $C$ there exists an embedding $E\to F$ of vector bundles
on $C$ such that $R^1\pi_*(F)=0$.
\end{lem}

\Pf . (i) It is well known that $T$ is a quotient stack (see \cite[Thm.\ 4.4]{Kresch}). Hence, by
\cite[Thm.\ 1]{KV}, there exists a scheme $Z$ and a finite flat surjective morphism $q:Z\to T$.
Since the projection $p:T\to \ov{T}$ is proper and quasi-finite, it follows that the map
$p\circ q:Z\to \ov{T}$ is finite, so $Z$ is projective and $\OO_Z(1)=(pq)^*\OO(1)$ is ample.
Therefore, for any coherent sheaf $\FF$ on $T$ the natural morphism of sheaves on $Z$
\begin{equation}\label{H-0-Z-q-F-eq}
H^0(Z,q^*\FF(n))\otimes\OO_Z(-n)\to q^*\FF
\end{equation}
is surjective and $H^{>0}(Z,q^*\FF(n))=0$ for $n\gg 0$.
Note that by projection formula
$$H^i(Z,q^*\FF(n))\simeq H^i(T,q_*(\OO_Z)\ot \FF(n)).$$
Since $\OO_T$ is a direct summand of $q_*(\OO_Z)$, we deduce the vanishing 
of $H^{>0}(T,\FF(n))$ for $n\gg 0$. On the other hand, applying the push-forward by $q$
to \eqref{H-0-Z-q-F-eq} we obtain a surjective morphism
$$H^0(Z,q^*\FF(n))\otimes q_*(\OO_Z)(-n)\to q_*(\OO_Z)\ot\FF,$$
which implies the surjectivity of \eqref{gen-mor-eq} for $n\gg 0$
with $\VV=q_*(\OO_Z)^\vee$.

\noindent
(ii) Let $L=\om_{C/T}(p_1+\ldots+p_r)$. Since the line bundle $L$
is ample on fibers of $\pi$, for sufficiently large $n$ we have
$R^1\pi_*(E^\vee\ot L^n)=R^1\pi_*(L^n)=0$ and the morphism
$$\pi^*\pi_*(E^\vee\ot L^n)\ot L^{-n}\to E^\vee$$
is surjective. Thus, setting $F=\pi^*(\pi_*(E^\vee\ot L^n))^{\vee}\ot L^n$ we obtain
an embedding of vector bundles $E\to F$ and 
$$R^1\pi_*(F)\simeq (\pi_*(E^\vee\ot L^n))^\vee\ot R^1\pi_*(L^n)=0.$$
\ed

Without loss of generality we can assume that $S$ is connected.
Our family of $\Ga$-spin curves $(\CC/S,p_1,\ldots,p_r; P,\vareps)$ induces a map
$S\to \SS$ to a connected component of the moduli stack of $\Ga$-spin structures.
Note that the data $(C/S, p_1,\ldots,p_r; L_1,\ldots,L_n,\ov{\phi}_M)$ is obtained by the base
change from the corresponding universal data over $\SS$. By Proposition \ref{spin-stack-prop},
$\SS$ is a proper smooth DM-stack with projective coarse moduli, so
Lemma \ref{generation-lem} can be applied over $\SS$. 
Hence, similar assertions hold for sheaves over $S$ 
that are obtained by the pull-back from sheaves over $\SS$.

By Lemma \ref{generation-lem}(ii),
for each $j$ we can choose an embedding of vector bundles
$L_j\to P_j$, where $R^1\pi_*(P_j)=0$. Let us consider the induced embedding
of $L_j(-\Si_j)$ into $P_j$ and set 
$$Q_j=P_j/L_j(-\Si_j)\simeq\coker(L_j\to P_j\oplus L_j|_{\Si_j}).$$
Then $Q_j$ has finite support over $\SS$, so it is also $\pi$-acyclic.
Thus, we get $\pi$-acyclic resolutions
$$L_j(-\Si_j)\to [P_j\to Q_j] \text{ and }$$
$$L_j\to [P_j\oplus L_j|_{\Si_j}\to Q_j].$$
Note that the exact sequence
$$0\to L(-\Si_j)\to L_j\to L_j|_{\Si_j}\to 0$$
is realized by the exact sequence of $\pi$-acyclic resolutions
$$0\to [P_j\to Q_j]\to [P_j\oplus L_j|_{\Si_j}\to Q_j]\to L_j|_{\Si_j}\to 0.$$
Now consider the sheaves 
$A_j=\pi_*(P_j\oplus L_j|_{\Si_j})$, $B_j=\pi_*(Q_j)$ and 
$A'_j=\pi_*(P_j)$ on $S$. Since $R^{>0}\pi_*(P_j)=0$, it follows that
$A_j$ and $A'_j$ are vector bundles. 
We have an exact sequence of sheaves on $C$
$$0\to L_j|_{\Si_j}\to Q_j\to Q'_j\to 0,$$
where $Q'_j=P_j/L_j$.
Since $Q'_j$ is a vector bundle on $C$ with $R^{>0}\pi_*(Q'_j)=0$ it follows
that $\pi_*(Q'_j)$ is a vector bundle, hence $B_j=\pi_*(Q_j)$ is also
a vector bundle.

Thus, we get the complexes $K'_j=[A'_j\to B_j]$ and $K_j=[A_j\to B_j]$
representing $R\pi_*(L_j(-\Si_j))$ and $R\pi_*(L_j)$, and the map
$Z_j:A_j=\pi_*(P_j\oplus L_j|_{\Si_j})\to \pi_*(L_j|_{\Si_j})\simeq\OO_S^{\Si_j}$
induced by the projection, such that $K'_j=\ker(Z_j)$.

At this point we can realize on the level of complexes the part of the diagram \eqref{main-diagram-eq}
not involving the maps $t_M$ and $\tau_M$.
First, by taking the external tensor products we get for each monomial $M$ a
$\pi^M$-acyclic resolution
$$L^M_\bullet\to R^M \ \ \text{ and } L_\bullet(-\Si_\bullet)^M\to \bar{R}^M,$$
where $\bar{R}^M$ is a subcomplex in $R^M$. Recall that the object $E_M\in D^b(S)$ is given
by $\Sym(M)$-invariants of the push-forward $R\pi^M_*(\FF_M)$, so we also need a $\pi^M$-acyclic resolution of $\FF_M$.
Note that the map $L^M_\bullet\to \De_*\OO_{\Si_M}$ is realized by a surjective
chain map of complexes $r(M):R^M\to \De_*\OO_{\Si_M}$ vanishing on $\bar{R}^M$. 
This implies that
the subcomplex $\ker(r(M))\sub R^M$ is a $\pi^M$-acyclic resolution of $\FF_M$, and
the embedding $L_\bullet(-\Si_\bullet)^M\hra\FF_M$ is realized by an embedding of resolutions
$\bar{R}^M\hra\ker(r(M))$.

Now the $\Sym(M)$-invariants of the push-forwards with respect to $\pi^M$
of $R^M$, $\bar{R}^M$ and $\ker(r(M))$ will represent $M(R\pi_*(L_\bullet))$, 
$M(R\pi_*(L_\bullet(-\Si_\bullet)))$ and $E_M$, respectively. Now all the maps in
the diagram \eqref{main-diagram-eq} except for $t_M$ and $\tau_M$, are realized on the level
of complexes. Furthermore, we have natural isomorphisms
$$\pi^M_*(R^M)^{\Sym(M)}\simeq M(K_\bullet), \ \ 
\pi^M_*(\bar{R}^M)^{\Sym(M)}\simeq M(K'_\bullet)$$
and an exact sequence
$$0\to \pi^M_*(\ker(r(M))^{\Sym(M)}\to M(K_\bullet)\rTo{Z_M}\OO^{\Si_M}\to 0.$$
It follows that if instead of $\pi^M_*(\ker(r(M))^{\Sym(M)}$ we use
$K_M=\Cone(Z_M)[-1]$, we will still have a representation of 
the part of the diagram \eqref{main-diagram-eq} not involving
the maps $t_M$ and $\tau_M$, on the level of complexes 
(note that $\iota_M$ becomes represented by the natural inclusion
$M(K'_\bullet)\to K_M$).

Note that the above realization depends only on complexes $[A_j\to B_j]$ and surjective maps
$Z_j:A_j\to\OO_S^{\Si_j}$ (recall that $A'_j=\ker(Z_j)$). 
To realize the maps $\tau_M$ by chain maps we will modify this realization by replacing
the complexes $[A_j\to B_j]$ with new complexes $[\bar{A}_j\to\bar{B}_j]$ equipped with
surjective quasi-isomorphisms $[\bar{A}_j\to\bar{B}_j]\to [A_j\to B_j]$.
To do this we will need the following technical assertion
similar to Proposition 4.7 of \cite{PV}. 

\begin{lem}\label{change-com-van-lem} 
Let a stack $T$ and a vector bundle $\VV$ be such that the assertions of Lemma \ref{generation-lem}(i)
hold.
Let $[C_0\to C_1]$ be a complex of vector bundles on $T$.
For each integer $d>0$ there exists
$m_0>0$ such that for any $m_1\ge m_0$ and any surjection 
$$\ov{C}_1=\VV^\vee(-m_1)^{\oplus N}\rTo{\si} C_1$$ 
one has
\begin{equation}\label{ample-van-eq}
H^{>0}(T,(\ov{C}_0^{\ot q_1})^{\vee}\ot\VV^{\ot q_2}(m))=0 \text{ for }m\ge m_0  \text{ and } q_1+q_2\le d,
\end{equation}
where the bundle $\ov{C}_0$ is the fiber product of $C_0$ and $\ov{C}_1$
over $C_1$, so that we have a quasi-isomorphism of complexes 
$$[\ov{C}_0\to\ov{C}_1]\to [C_0\to C_1].$$
\end{lem}

\Pf . Set $K=\ker(\si)$, so that we have exact sequences of vector bundles
\begin{equation}\label{K-0-seq}
0\to K\to \ov{C}_0\to C_0\to 0,
\end{equation}
\begin{equation}\label{K-1-seq}
0\to K\to \ov{C}_1\to C_1\to 0.
\end{equation}
From the sequence \eqref{K-0-seq} we see that \eqref{ample-van-eq}
would follow from the vanishing of 
$$H^{>0}(T,(C_0^{\ot q_1})^{\vee}\ot (K^{\ot q_2})^{\vee}\ot \VV^{\ot q_3}(m))$$
for $m\ge m_0$ and $q_1+q_2+q_3\le d$. Taking tensor powers of
the sequence dual to \eqref{K-1-seq} we get a resolution of the bundle
$(K^{\ot q_2})^{\vee}$ with terms that are direct sums of vector bundles of the form 
$(C_1^{\ot (q_2-s)})^{\vee}\otimes \VV^{\ot s}(m_1s)$.
Thus, it would be enough to find $m_0$ such that
$$H^{>0}(T,(C_0^{\ot q_1})^{\vee}\ot (C_1^{\ot q_2})^{\vee}\ot \VV^{\ot q_3}(m))=0$$
for $m\ge m_0$ and $q_1+q_2+q_3\le d$. But this is possible by 
Lemma \ref{generation-lem}(i).
\ed

Let us apply Lemma \ref{change-com-van-lem} to the complex $[\oplus_j A_j\to \oplus_j B_j]$
and $d$ equal to the maximum of the degrees $|M|$ of all the monomials $M$ occurring in $\w$.
Then we can choose large enough $m_0$ and surjections 
$\ov{B}_j=\VV^\vee(-m_0)^{\oplus N_j}\to B_j$ and replace each $[A_j\to B_j]$ by a
quasi-isomorphic complex $[\ov{A}_j\to \ov{B}_j]$ such that
$$\Ext^{>0}_S((\oplus_j \ov{A}_j)^{\ot q_1}\ot (\VV^{\ot q_2})^\vee(-m),\OO_S)=0$$
for $m\ge m_0$ and $q_1+q_2\le d$. 
This implies that
$$\Ext^{>0}_S((\oplus_j \ov{A}_j)^{\ot q_1}\ot(\oplus_j \ov{B}_j)^{\ot q_2},\OO_S)=0$$
for $q_1+q_2\le d$ and $q_2\ge 1$. Hence, 
for every monomial $M$ appearing in $\w$ 
the terms of the complex 
$$\ov{K}_M=\Cone(M([\ov{A}_\bullet\to \ov{B}_\bullet])\rTo{Z_M}\OO_S^{\Si_M})[-1]$$
representing $E_M$, satisfy $\Ext^{>0}(\ov{E}_M^i,\OO_S)=0$ for $i\ge 2$. 
This easily implies (using the standard
spectral sequence) that the space of morphisms
$$\Hom_{D(S)}(\ov{K}_M,\OO_S[-1])$$
in the derived category is the same as in the homotopy category of complexes.
Thus, replacing $[A_j\to B_j]$ with $[\ov{A}_j\to\ov{B}_j]$
we can realize the map $\tau_M$ by a chain map 
$K_M\to \OO_S[-1]$. 

\

\noindent
{\bf Step 3}. Let $X$ be the total space of the bundle $A_1\oplus\ldots\oplus A_n$ over $S$
and let $p:X\to S$ be the projection.
Note that there is a natural map $Z: X\to \A^{\ov{\ga}}$ induced by the morphisms
$Z_j:A_j\to\OO_S^{\Si_j}$ constructed in Step 2. 
Recall (see Section \ref{w-moduli-sec}) that the group $G(\ov{\ga})=\prod_{i=1}^rG_{I(\ga_i)}$
acts on the set of restricted rigidifications of a given $\Ga$-spin structure in such a way
that an element $\la=(\la(i))\in G(\ov{\ga})$ 
changes the trivializations $(e_j(i))$ to $(\la_j(i)\cdot e_j(i))$, where 
$\la_j(i)$, for $(p_i,j)\in\Si$, are the components of $\la(i)\in G_{I(\ga_i)}$.
Suppose that a subgroup 
$G_S\sub G(\ov{\ga})$ acts on $S$, so that the map $S\to\SS$
is $G_S$-invariant and the action on fibers is induced by the above rescaling action
of $G(\ov{\ga})$.
Then all complexes $[A_j\to B_j]$ are $G_S$-equivariant,
so we have an action of $G_S$ on $X$. 
Since the maps $Z_j:R\pi_*(L_j)\to\OO_S^{\Si_j}$ were defined
using the trivializations
of $L_j|_{\Sigma_j}$ given by $(e_j(i))$, we have
$$\la^*Z_j=\la_j^{-1}\cdot Z_j$$
for $\la\in G_S$, where $\la_j=(\la_j(i))_{i\in\Si_j}$ acts diagonally on $\OO_S^{\Si_j}$.
Hence, the map $Z:X\to\A^{\ov{\ga}}$ satisfies
$$Z\circ\la=\la^{-1}\cdot Z.$$
Thus, the map $(p,Z):X\to S\times\A^{\ov{\ga}}$ is $G_S$-equivariant, where
$\la\in G_S\subset G(\ov{\ga})$ acts on $S\times\A^{\ov{\ga}}$ by 
$$\la\cdot (s,z)=(\la\cdot s,\la^{-1}\cdot z).$$
The natural action of $\G_m^n$ on the fibers of $p:X\to S$ induces an action
of the group $\Ga_\w$ (see Section \ref{sym-sec}) on $X$ such that
the map $Z$ is $\Ga_\w$-equivariant. 
Now we will construct a $G_S\times\Ga_\w$-equivariant Koszul matrix factorization of
the potential $-Z^*\w_{\ov{\ga}}$ on $X$.

The complex $K_M=\Cone(M([A_\bullet\to B_\bullet])\rTo{Z_M}\OO_S^{\Si_M})[-1]$ has
the following form
\begin{align*}
&M(A_\bullet)\rTo{(Z_M,-\de)} 
\OO_S^{\Sigma_M}\oplus \bigoplus_{j=1}^n\pa_jM(A_\bullet)\ot B_j\to\\
&\bigoplus_{j<j'}\bigl(\pa_j\pa_{j'}M(A_\bullet)\ot B_j\ot B_{j'}\bigr)\oplus\bigoplus_j
\bigl(\pa_j^2M(A_\bullet)\ot {\bigwedge}^2B_j\bigr)\to\ldots
\end{align*}
where the first term is in degree $0$ (here $\de$ is the differential on the complex
$M([A_\bullet\to B_\bullet])$).
Since the chain map $\tau_M$ is equal to $\Tr$ on $\OO_S^{\Si_M}[-1]\sub K_M$
(see diagram \eqref{main-diagram-eq}),
it corresponds to a map
$$\a_M=(\a_{M,j}):\oplus_{j=1}^n B_j\ot\pa_jM(A_\bullet)\to\OO_X$$
such that the following diagram is commutative
\begin{equation}\label{Z-a-diag}
\begin{diagram}
M(A_\bullet)&\rTo{\de}& \bigoplus_{j=1}^n\pa_jM(A_\bullet)\ot B_j\to\ldots\\
\dTo{Z_M}&&\dTo{\a_M}\\
\OO_S^{\Sigma_M}&\rTo{\ \ \ \ \ \ \ \ \ \Tr}&\OO_S
\end{diagram}
\end{equation}
Set 
$$\a'_M:=(k_j\cdot \a_{M,j})_{j=1,\ldots,n},$$
where $M=x_1^{k_1}\cdot\ldots\cdot x_n^{k_n}$.
Let us view the differential 
$$\b=\oplus_j\b_j:\oplus_j A_j\to\oplus_j B_j$$ 
as a section of 
the bundle $p^*(\oplus_j B_j)$ on $X$ (linear along fibers). Similarly,
we can view $\a'_M$ as a section of the bundle $p^*(\oplus_j B_j^{\vee})$ on $X$
(polynomial along fibers). Let $\lan\cdot,\cdot\ran$ denote the natural pairing between
$p^*(\oplus_j B_j^{\vee})$ and $p^*(\oplus_j B_j)$.

\begin{lem} One has $\lan \a'_M,\b\ran=Z^*M^{\oplus \Si_M}$.
\end{lem}

\Pf . The components of the differential $\de$ have form
$$\de(M(a_\bullet))_j=\b_j(a_j)\ot\pa_jM(a_1,\ldots,a_n)\in B_j\ot \pa_jM(A_\bullet),$$
where
$a_j\in A_j$, $j=1,\ldots,n$. Hence, the composition
$\a_M\circ\de(M(a_\bullet)):M(A_\bullet)\to\OO_S$ corresponds
to the function $\lan \a'_M,\b\ran$ on $X$. On the other hand, 
the components of the map $Z_M$ correspond to the functions
$Z^*M(x_1(i),\ldots,x_n(i))$ (where $i\in\Sigma_M$). Thus, the assertion follows from the
commutativity of the diagram \eqref{Z-a-diag}.
\ed

Now we set 
$$\a_\w=\sum_M c_M\a'_M, \text{  so that}$$
$$\lan \a_\w,\b\ran=\sum_M c_M Z^*M^{\oplus \Si_M}=Z^*\w_{\ov{\ga}}.$$
Hence, we have a Koszul
matrix factorization $\{-\a_\w,\b\}$ of $-Z^*\w_{\ov{\ga}}$. 
We claim that it can be equipped with a $G_S\times\Ga_\w$-equivariant
structure with respect to the character $\chi_\w:\Ga_\w\to\G_m$ (and trivial on $G_S$).
The bundles $\bigoplus A_j$ and $\bigoplus B_j$ are equipped with a
$G_S\times\Ga_\w$-equivariant structure, in such a way that the action of
$\Ga_\w$ is induced by the embedding $\Ga_\w\sub\G_m^n$.
Then $\b$ can be viewed as a $G_S\times\Ga_\w$-invariant section of $\bigoplus_j B_j$.
On the other hand, $\a_\w$ gives 
a $G_S\times\Ga_\w$-invariant section of $\chi_\w\ot\bigoplus_j B_j^{\vee}$.
Thus, we obtain a $G_S\times\Ga_\w$-equivariant structure on the \mf\ $\{-\a_w,\b\}$.
More explicitly, we have
$$\{-\a_w,\b\}_0=\bigoplus_i {\bigwedge}^{2i}(p^*(\bigoplus_j B_j^{\vee}))\ot\chi_\w^i,\ \
\{-\a_w,\b\}_1=\bigoplus_i {\bigwedge}^{2i+1}(p^*(\bigoplus_j B_j^{\vee}))\ot\chi_\w^i,$$
where the differential is given by
$$\de=\iota(\b)-\a_\w\we.$$ 

\

\noindent 
{\bf Step 4}. Now we will show that the 
matrix factorization $\{-\a_\w,\b\}$ constructed in Step 3 
is supported on the zero section in $X$.
This will allow us to apply the push-forward functor (see Example \ref{G-push-forward-constr})
for the projection $(p,Z):X\to S\times \A^{\ov{\ga}}$.
For the universal family of $\Ga$-spin curves over the moduli space $\SS_g^{\rig,0}$ we get an object
\begin{equation}\label{fund-mf-res-eq}
\begin{array}{l}
\bP_g^{\rig,0}(\ov{\ga}):=(p,Z)_*\{-\a_\w,\b\}
\in 
\ov{\DMF}_{G(\ov{\ga})\times\Ga_\w}(\SS_g^{\rig,0}(\ov{\ga})\times 
\A^{\ov{\ga}}, -\w_{\ov{\ga}}).
\end{array}
\end{equation}

Recall that by Step 2 we have $A'_j=\ker(A_j\to \OO_S^{\Si_j})$ and 
the complex $[A'_j\stackrel{\b'_j}{\to} B_j]$ represents
$R\pi_*(L_j(-\Sigma_j))$, where $\b'_j=\b_j|_{A'_j}$.
Thus,  the total space $X_0$ of the vector bundle 
$A'_1\oplus\ldots\oplus A'_n$ over $S$ coincides with the preimage of the origin $Z^{-1}(0)\sub X$.
Since the critical locus of $\w_{\ov{\ga}}$ on $\A^{\ov{\ga}}$ is the origin $0\in\A^{\ov{\ga}}$, it is enough to
consider the zero locus of the restriction of $\{-\a_\w,\b\}$ to $X_0=Z^{-1}(0)$ (by \cite[Cor.\ 5.3]{PV-stacks}).
Let $\a_0$ be the restriction of $\a_\w$ to $X_0$.
The equation \eqref{tau-t-M-eq} implies that $\a_0$ and $\b'=(\b'_j)$ are exactly the sections
of the pull-backs of $\bigoplus_j B_j^{\vee}$ and $\bigoplus_j B_j$
obtained by the construction of section \ref{simple-constr-sec} applied to the collection of morphisms
\eqref{M-L-Si-om-eq}. It follows from Proposition \ref{zero-locus-prop}
that the zero locus $Z(s_0)$ is exactly the zero section $S\sub X_0\sub X$, as claimed.
Hence, by Lemma \ref{Koszul-mf-support-lem}, $\{-\a_w,\b\}|_{X_0}$ is supported on the zero section.

Recall that we have a natural morphism
$$\SS_g^{\rig}(\ov{\ga})\to\SS_g^{\rig,0}(\ov{\ga})$$
compatible with the homomorphism $\prod_{i=1}^r G/\lan\ga_i\ran\to G(\ov{\ga})$.
Hence, by taking the pull-back of $\bP_g^{\rig,0}(\ov{\ga})$ we obtain a \mf\
\begin{equation}\label{fund-mf-rig-eq}
\bP^{\rig}_g(\ov{\ga})\in 
\ov{\DMF}_{\Ga_\w}(\SS_g^{\rig}(\ov{\ga})\times 
\A^{\ov{\ga}}, -\w_{\ov{\ga}})
\end{equation}
which is equivariant with respect to the action of $\prod_{i=1}^r G/\lan\ga_i\ran$.

This finishes the construction of the fundamental \mf s. In 
Section \ref{independence-sec} we will show 
that it does not depend on the choices made (up to an isomorphism).

\begin{ex} In the case when $\A^{\ov{\ga}}=0$ (and $\w_{\ov{\ga}}=0$) we have
$\SS_g^{\rig,0}(\ov{\ga})=\SS_g(\ov{\ga})$ and
the category $\ov{\DMF}_{\Ga_\w}(\SS_g(\ov{\ga}),0)$ is (non-canonically) equivalent to
the bounded derived category of $G_\w$-equivariant
coherent sheaves on $\SS_g(\ov{\ga})$. For example, for $\w=x^n$, $G_\w=\Z/n$ 
(and $\Ga_\w=\G_m$) this will be the case whenever all $\ga_i\in\Z/n$ are nontrivial.
In this case the Chern class of $\bP_g(\ov{\ga})$ 
is closely related to the {\it Witten's virtual top Chern class} on the moduli spaces of higher spin curves
(see \cite{JKV}, \cite{PV} and \cite{Chiodo}). To get the Witten's virtual top Chern class one has
to twist it with a certain Todd class (see \eqref{lambda-1-eq} below). 
\end{ex}

\subsection{Independence of choices}\label{independence-sec}

Here we will show that the isomorphism class of the fundamental \mf\  $\bP_g^{\rig,0}(\ov{\ga})$
does not depend on the choices made in Step 2 
when realizing the diagram \eqref{main-diagram-eq} on
the level of complexes. 

First, since $\Hom_{D^b(S)}(K_M,\OO_S[-1])$ can be computed in the homotopy category,
all chain maps $K_M\to\OO_S[-1]$ representing $\tau_M$ are homotopic.
A homotopy between two such maps is given by a $G(\ov{\ga})\times \Ga_\w$-equivariant map
$$h:\bigoplus_{j<j'}\bigl(\pa_j\pa_{j'}M(A_\bullet)\ot B_j\ot B_{j'}\bigr)\oplus
\bigoplus_j\bigl(\pa_j^2M(A_\bullet)\ot {\bigwedge}^2B_j\bigr)\to\OO_S\ot\chi_\w.$$
After dualization $h$ can be viewed as a section of $p^*\bigwedge^2(\bigoplus_j B_j^{\vee})\ot\chi_\w$.
Now the operator $\exp(-h)\we ?$ 
induces a $G(\ov{\ga})\times \Ga_\w$-equivariant  
isomorphism between the matrix factorizations associated with two homotopic choices of $\tau_M$
(cf.\ \cite[Prop. 4.2]{PV} or the proof of \cite[Lem.2.5.5]{PV-mf}).

To prove independence of the choice of presentations $R\pi_*(L_j)=[A_j\to B_j]$
we will use the following property of Koszul \mf s (which is
analogous to the results of \cite[sec. 3.2]{PV}).

\begin{prop}\label{independence-prop}
Let $V$ be a vector bundle on a smooth FCDRP-stack $X$, $W\in H^0(X,L)$ a potential,
and let $\{\a,\b\}$ be the Koszul \mf\
associated with sections $\a\in H^0(X,V^{\vee}\ot L)$ and $\b\in H^0(X,V)$ 
such that $\lan \a, \b\ran=W$.
Let $V_1\sub V$ be a subbundle such that $\b\mod V_1$ is a regular section
of $V/V_1$. Assume that the zero locus $X'=Z(\b\mod V_1)$ is smooth and consider
the induced sections 
$$\b'=\b|_{X'}\in H^0(X',V')$$ 
of the bundle $V'=V_1|_{X'}$ 
(note that the restriction $\b|_{X'}$ belongs to $V_1|_{X'}\sub V|_{X'}$) and 
$$\a'=\a\mod (V_1^{\perp}|_{X'})$$ of the bundle
$L\ot (V^{\vee}/V_1^{\perp})|_{X'}\simeq L\ot (V')^{\vee}$.
Assume also that either $W|_{X'}$ is a non-zero-divisor 
or $W=0$ and the zero loci $Z(\a,\b)$ and $Z(\a',\b')$ are proper.
Then one has an isomorphism
\begin{equation}\label{independ-mf-isom}
\{\a,\b\}\simeq i_*\{\a', \b'\}
\end{equation}
in $\DMF(X,W)$, where $i:X'\hra X$ is the natural embedding.
\end{prop}

\Pf . We have a natural morphism 
\begin{equation}\label{independ-mf-map-eq}
{\bigwedge}^{\bullet}(V^{\vee}\ot L^{1/2})(L^{-1/2})\to 
i_*i^*{\bigwedge}^{\bullet}(V^{\vee}\ot L^{1/2})(L^{-1/2})\to 
i_*{\bigwedge}^\bullet((V')^{\vee}\ot L^{1/2})(L^{-1/2}), 
\end{equation}
compatible with the differentials in $\{\a,\b\}$ and $\{\a',\b'\}$.

Assume first that $W|_{X'}$ is a non-zero-divisor.
To show that the map \eqref{independ-mf-map-eq}
is an isomorphism
in $\DMF(X,W)$ we can argue locally.
Thus, we can assume that $V=V_1\oplus V_2$, so we can write 
$\a=(\a_1,\a_2)$ 
$\b=(\b_1,\b_2)$ with $\a_i\in H^0(X,V_i^\vee\ot L)$ and $\b_i\in H^0(X,V_i)$, $i=1,2$.
Our assumptions mean that $\b_2$ is a regular section of $V_2$ and $X'$ is its zero locus.
We have a natural isomorphism of \mf s
$$\{\a,\b\}\simeq\{\a_1,\b_1\}\ot\{\a_2,\b_2\},$$
where $\{\a_i,\b_i\}$ is a \mf\ of $W_i:=\lan\a_i,\b_i\ran$. Note that $W=W_1+W_2$ and
$W_2|_{X'}=0$. Thus, we only need to show that the natural morphism
$$\{\a_1,\b_1\}\ot\{\a_2,\b_2\}\to i_*\{i^*\a_1,i^*\b_1\}$$
is an isomorphism in $\DMF(X,W)$. But this follows from Proposition \ref{koszul-reg-prop}(i).

Next, consider the case $W=0$. We use a deformation argument.
Namely, let us consider a family of chain maps
$$f_t:\{t\a,\b\}\to i_*\{t\a',\b'\},$$
where $t\in\A^1$. We can view the complex $\KK=\Cone(f_t)$ as a complex of sheaves
on $X\times\A^1$, flat over $\A^1$. Note that $f_0$ is quasi-isomorphism, since
$\{0,\b\}\simeq \{0,\b_1\}\ot \{0,\b_2\}$ is the tensor product of the usual Koszul complex and $\b_2$ is regular, hence
$\KK|_{X\times\{0\}}$ is acyclic. On the other hand, 
by Lemma \ref{Koszul-mf-support-lem}, the cohomology of
$\KK_{X\times (\A^1\setminus \{0\})}$ is supported on $Z\times (\A^1\setminus\{0\})$,
where $Z=Z(\a,\b)\cup Z(\a',\b')$. By our assumption, $Z$ is proper, 
hence there exists an open neighborhood $U$ of
$0\in\A^1$ such that $\KK|_{X\times U}$ is acyclic. In particular, there exists $t\neq 0$
such that $f_t$ is quasi-isomorphism. Using an isomorphism
between $f_1$ and $f_t$ for $t\neq 0$ (given by multiplying by $t^i$ on $\bigwedge^i$), we derive
that $f_1$ is also a quasi-isomorphism.
\ed

Now we are ready to show that different choices of 
complexes of bundles $K_j=[A_j\to B_j]$ realizing $R\pi_*(L_j)$ and surjective
maps $K_j\to \OO_S^{\Si_j}$ realizing $Z_j$, which we made in Step 2 realizing the
diagram \eqref{main-diagram-eq} at the level of complexes,
lead to isomorphic \mf s. It is enough to check this in the case when one of the complexes $K_j$ is replaced by a complex $\wt{K}_j$ in degrees $[0,1]$ such that we have
 a quasi-isomorphism $\wt{K}_j\to K_j$ for which
the composition $\wt{K}_j\to K_j\to \OO_S^{\Si_j}$ is still surjective.
By Lemma 4.4 of \cite{PV}, we can represent the map $\wt{K}_j\to K_j$ in the homotopy category 
as a composition of an embedding of complexes 
$\wt{K}_j\to K_j\oplus U$ followed by the projection $K_j\oplus U\to K_j$,
where $U$ is of the form $U=[F\rTo{\id}F]$. Thus, it suffices to consider separately two cases:
the case when $\wt{K}_j\to K_j$ is an embedding and the case when it is a projection
onto a direct summand.

First, suppose we have a quasi-isomorphism $[\wt{A}_1\to \wt{B}_1]\to [A_1\to B_1]$, such that
$\wt{A}_1\to A_1$ and $\wt{B}_1\to B_1$ are
embeddings of bundles. Let us denote $C=A_1/\wt{A}_1=B_1/\wt{B}_1$ and
let $\wt{Z}_1:\wt{A}_1\to\OO_S^{\Sigma_1}$ be the restriction of $Z_1$ (recall that
$\wt{Z}_1$ is still surjective).
Then for each monomial $M$, the components
$$\wt{\a}_j:\partial_jM(\wt{A}_1,A_2,\ldots,A_n)\to B_j^{\vee},$$
for $j\neq 1$, are obtained from $\a_j$ by restriction to $\wt{A}_1$, while to get
$\wt{\a}_1$ we also have to use the map $B_1^{\vee}\to (\wt{B}_1)^{\vee}$.
In the new realization the space $X$ is replaced by its subspace
$\wt{X}\sub X$ which is the total space of $\wt{A}_1\oplus \bigoplus_{j>1} A_j$. Note that $\wt{X}$
is the zero locus of the regular section of $p^*B_1/p^*\wt{B}_1=p^*C$ on $X$,
induced by $\b_1$. 
Thus, we can apply Proposition \ref{independence-prop} to the subbundle
$$p^*\wt{B}_1\oplus\bigoplus_{j>1} B_j\sub \bigoplus_j B_j$$
to obtain an isomorphism
$$\{\a,\b\}\simeq i_*\{\wt{\a},\wt{\b}\},$$
where $i:\wt{X}\hra X$ is the embedding.
Hence, the push-forwards of $\{\a,\b\}$ and $\{\wt{\a},\wt{\b}\}$ to 
$S\times \A^{\ov{\ga}}$ are isomorphic.

The case of the quasi-isomorphism of the type
$K_j\oplus U\to K_j$ is similar because we can apply the above
argument to the embedding $K_j\to K_j\oplus U$ which is also a quasi-isomorphism.

This finishes the proof of independence of the isomorphism class of the fundamental
\mf\ in the corresponding derived category of the choices made in Step 2 of the construction.

\section{Cohomological field theories associated with a quasihomogeneous isolated singularity}
\label{axioms-sec}

\subsection{Construction of \CFTs}

Let $R$ be a commutative $\C$-algebra and
$\HH$ a finitely generated $\Z/2$-graded projective $R$-module
equipped with a perfect symmetric $R$-bilinear pairing
$$b: \HH\ot_R \HH\to R$$
(i.e., $b$ induces an isomorphism $\HH\simeq\Hom_R(\HH,R)$). 
Let $\De^R\in\HH\ot _{R}\HH$ be the Casimir element corresponding to $b$.
Recall (see \cite[III.4]{Manin}) that a {\it complete Cohomological Field Theory} (\CFT) 
on the state space $(\HH,b)$ with coefficients in $R$ 
is a collection of even $R$-linear maps
\begin{equation}\label{R-linear-correlators-eq}
\La_{g,r}:\HH^{\ot_R n}\to H^*(\ov{\MM}_{g,r})\ot R
\end{equation}
and a fixed element $\unit\in\HH$ (called {\it flat unit})
satisfying certain properties. Here the Casimir $\De^R$ is used to formulate the factorization properties
and the insertion of the identity $\unit$ corresponds to
forgetting a marked point.

An example of \CFT\ with coefficients is provided by 
the $G$-equivariant Gromov-Witten theory considered in \cite{Giv-eqGW},
where $R=H^*_G(pt)$. 
Note that if we have a homomorphism $R\to R'$ and a \CFT\ with coefficients in $R$
then by extending scalars we can obtain
a \CFT\ with coefficients in $R'$ and the state space $\HH\ot_R R'$.

Let us fix a  quasihomogeneous polynomial $\w(x_1,\ldots,x_n)$ with an isolated singularity.
We assume that the degrees $d_j=\deg(x_j)$ are positive.
Let $G\sub G_\w$ be a finite subgroup in the group of diagonal symmetries of $\w$, such that
$G$ contains 
the exponential grading element $J$ (defined by \eqref{J-eq}).
Let $\Ga\sub\Ga_\w$ be the subgroup associated with $G$ by Lemma \ref{J-lem},
so that there is an exact sequence
\begin{equation}\label{group-ext-seq}
1\to G\to \Ga \rTo{\chi}\G_m\to 1.
\end{equation}
For each $\ga\in G$ let us set
$$\HH_{\ga}:=HH_*(\MF_{\Ga}(\w_{\ga})),$$
where $\w_{\ga}$ is the restriction of $\w$ to the fixed point locus $(\A^n)^{\ga}\sub\A^n$
(we consider the induced action of $\Ga$ on $(\A^n)^{\ga}$).
Recall from Section \ref{Hoch-sec} that $\HH_{\ga}$ has a natural $\widehat{G}$-action, so we can view it
as a module over the ring $R=\C[\widehat{G}]$.

We will construct a \CFT\ with coefficients in $R$ and
the state space
\begin{equation}\label{state-space-eq}
\HH=\HH(\w,G):=\bigoplus_{\ga\in G}\HH_{\ga}.
\end{equation}

To define the corresponding $R$-linear maps
\be\label{lambda-g-eq}
\La_{g,\Ga,\chi}=\La_g^R(\ov{\ga}):\HH_{\ga_1}\ot_R\ldots\ot_R \HH_{\ga_n}\to H^*(\ov{\MM}_{g,r},R)=
H^*(\ov{\MM}_{g,r})\ot R
\end{equation}
for $\ov{\ga}=(\ga_1,\ldots,\ga_r)\in G^r$ (the components of the \CFT\ maps
\eqref{R-linear-correlators-eq}) we will use the moduli spaces 
$$\SS_g^{\rig}(\ov{\ga})=\SS_{g,r,\Ga,\chi}^{\rig}(\ov{\ga})$$
introduced in Section \ref{w-moduli-sec}
and the fundamental matrix factorizations $\bP^{\rig}_g(\ov{\ga})$ (see \eqref{fund-mf-rig-eq}).

Retaining
only the
action of the subgroup $\Ga\sub\Ga_\w$ on $\bP^{\rig}_g(\ov{\ga})$
we obtain an object
$$\bP^{\rig}_{g,\Ga}(\ov{\ga})\in
\ov{\DMF}_{\Ga}(\SS^{\rig}_{g}(\ov{\ga})\times\A^{\ov{\ga}},-\w_{\ov{\ga}}),$$
which is also invariant under the action of $G^r$ (see Step 4 of the construction of
Section \ref{constr-sec}).
This object gives a functor
\begin{equation}\label{Phi-functor-eq}
\Phi_g(\ov{\ga}):\ov{\DMF}_{\Ga}(\A^{\ov{\ga}},\w_{\ov{\ga}})\to D_{G}(\SS^{\rig}_{g}(\ov{\ga})):
\bar{E}\mapsto (p_1)_*(p_2^*\bar{E}\ot \bP_{g,\Ga}^{\rig}(\ov{\ga})),
\end{equation}
where $p_1$ and $p_2$ are the projections from the product 
$\SS^{\rig}_{g}(\ov{\ga})\times\A^{\ov{\ga}}$ onto its factors.
Here we use the push-forward functor \eqref{G-mf-0-push-forward}. Note that the functor $\Phi_g(\ov{\ga})$ is compatible
with the action of $G^r$ on both categories.

For a stack $\XX$ let us denote by $HH_*(\XX)$ the Hochschild homology of the dg-version
of the derived category of coherent sheaves on $\XX$. Note that for a morphism
$f:\XX\to\XX'$ we have pull-back maps $f^*:HH_*(\XX')\to HH_*(\XX)$ induced
by the pull-back functor between the derived categories of coherent sheaves.

For a smooth proper DM-stack $\XX$ we denote by $H^*(\XX,\C)$ and $H_*(\XX,\C)$
the usual cohomology and homology of its coarse moduli space. Note that since this moduli space
is a rational homology manifold, we have Poincar\'e duality isomorphism
$H^*(\XX,\C)\rTo{\sim} H_*(\XX,\C)$ given by the cap product with the fundamental class $[\XX]$
(note that the definition of $[\XX]$ takes into account the order of generic automorphism groups
over connected components of $\XX$; see \cite[Sec.\ 2.2]{AGV}).
For a morphism $f:\XX\to \YY$ of smooth proper DM-stacks we can use
the Poincar\'e duality to view the push-forward map $f_*:H_*(\XX,\C)\to H_*(\YY,\C)$
as a map between {\it cohomology} spaces.

Recall that for a smooth projective variety $X$ one has the Hochschild-Kostant-Rosenberg isomorphism
$$I_{HKR}:HH_*(X)\simeq H^*(X,\C),$$
sending $HH_i(X)$ to $\bigoplus_{q-p=i}H^{p,q}(X)$.
Now let $\XX$ be a smooth proper connected DM-stack with projective coarse moduli space. 
Since we work over $\C$, by \cite[Thm.\ 4.4]{Kresch}, $\XX$ is a quotient stack. Hence, 
by \cite[Thm.\ 1]{KV},
there exists a finite flat surjective morphism $\pi:X\to\XX$ with $X$
a (connected) smooth projective variety. 
Using the pull-back
$\pi^*:HH_*(\XX)\to HH_*(X)$ and the push-forward
$\pi_*:H^*(X,\C)\to H^*(\XX,\C)$ we obtain a map
\begin{equation}\label{stack-map}
\a_{\XX}:HH_*(\XX)\rTo{\pi^*} HH_*(X)\rTo{I_{HKR}} H^*(X,\C)\rTo{\frac{1}{\deg\pi}\pi_*} H^*(\XX,\C)
\end{equation}
which
we will call the {\it HKR map}.
It is easy to see that this map does not depend on a choice of a morphism $X\to \XX$.
Indeed, for any other map $X'\to\XX$ we can find $\wt{X}\to\XX$ dominating both $X$ and $X'$.
Thus, we just need to use the fact that for a flat surjective morphism $f:\wt{X}\to X$ the map
$f^*:HH_*(X)\to HH_*(\wt{X})$ on Hochschild homology is compatible with the usual pull-back
map $f^*:H^*(X,\C)\to H^*(\wt{X},\C)$ via the Hochschild-Kostant-Rosenberg isomorphisms.
Similarly one can check that the maps \eqref{stack-map} are compatible with pull-backs and push-forwards
with respect to finite \'etale morphisms.
Also, the construction of the map \eqref{stack-map} easily extends to the case when $\XX$ is not connected.

\begin{lem} Let $\XX$ be a smooth proper DM-stack with projective coarse moduli space.
For an object $F\in D(\XX)$ 
let $\ch^{HH}(F)\in HH_*(\XX)$ be the categorical Chern character 
(see Section \ref{dg-Hoch-sec}) and let 
$\ch^{top}(F)\in H^*(\XX,\C)$ be the usual topological Chern character.
Then $\a_{\XX}(\ch^{HH}(F))=\ch^{top}(F)$.
\end{lem}

\Pf . We can assume $\XX$ is connected. Let $\pi:X\to \XX$ be a finite flat surjective morphism as above.
First, note that by \cite[Thm.\ 4.5]{Cal}, for any $G\in D(X)$ we have
$$I_{HKR}(\ch^{HH}(G))=\ch^{top}(G)\in H^*(X,\C).$$
Applying this to $G=\pi^*F$ we obtain
$$I_{HKR}(\pi^*\ch^{HH}(F))=\ch^{top}(\pi^*G)=\pi^*\ch^{top}(F).$$
Hence,
$$\a_{\XX}(\ch^{HH}(F))=\frac{1}{\deg\pi}\pi_*\pi^*\ch^{top}(F)=\ch^{top}(F).$$
\ed

Recall that by Corollary \ref{Kunneth-cor}, we have an isomorphism
\begin{equation}\label{iota-isom}
\iota(\ov{\ga}):\HH_{\ga_1}\ot_R\ldots\ot_R\HH_{\ga_r}\wt{\to} 
HH_*(\MF_{\Ga}(\A^{\ov{\ga}},\w_{\ov{\ga}}))^{G^r}\sub
HH_*(\MF_{\Ga}(\A^{\ov{\ga}},\w_{\ov{\ga}})).
\end{equation}
Let 
\begin{equation}\label{hoch-phi-map}
\begin{array}{l}
\phi_g(\ov{\ga}):\HH_{\ga_1}\ot_R\ldots\ot_R\HH_{\ga_r}\rTo{\iota(\ov{\ga})} HH_*(\MF_{\Ga}(\A^{\ov{\ga}},\w_{\ov{\ga}}))\rTo{\Phi_g(\ov{\ga})_*}
HH_*(\SS^{\rig}_{g}(\ov{\ga}))\ot R\to \\
H^*(\SS^{\rig}_{g}(\ov{\ga}),\C)\ot R
\end{array}
\end{equation}
be the composition of this embedding, the map $\Phi_g(\ov{\ga})_*$
induced on the Hochschild homology spaces by the functor $\Phi_g(\ov{\ga})$ 
and the HKR map
\eqref{stack-map}  for $\XX=\SS^{\rig}_{g}(\ov{\ga})$.

Let $\st_g:\SS_{g}^{\rig}(\ov{\ga})\to\ov{\MM}_{g,r}$ be the projection.
Let us consider the corresponding push-forward map 
$$\st_{g*}:H^*(\SS_{g}^{\rig}(\ov{\ga}),\C)\ot R\to H^*(\ov{\MM}_{g,r},\C)\ot R.$$

Now we are ready to define the \CFT\ on the state space
$\HH=\bigoplus_{\ga\in G}\HH_\ga$.
The maps \eqref{lambda-g-eq} are given by 
\begin{equation}\label{lambda-st-g-eq}
\La^R_g(\ov{\ga})=\frac{1}{\deg(\st_g)}\cdot \st_{g*}\circ\phi_g(\ov{\ga}).
\end{equation}
Note that these maps are even with respect to the natural $\Z_2$-gradings on \eqref{state-space-eq}.
For the exponential grading element $J\in G$ we have $\HH_J=R$ since 
$(\A^n)^J=0$.
We take the element $1\in\HH_J\sub\HH$ to be the flat unit $\unit$ of our \CFT.
To define a metric on $\HH$ we consider the element $\zeta_\bullet\in(\C^*)^n$
with components
\begin{equation}\label{zeta-def-eq}
\zeta_j=\exp(\pi iq_j) \text{ for } j=1,\ldots,n,
\end{equation}
where $q_j=d_j/d$ (see Section \ref{quasihom-mf-sec}). Note that $(\zeta_\bullet)^2=J$ and $\chi(\zeta_\bullet)=-1$,
so
\begin{equation}\label{zeta-eq}
\w(\zeta_\bullet x)=-\w(x).
\end{equation}
We set 
\begin{equation}\label{metric-eq}
(x,y)=\sum_{\ga}((\zeta_\bullet)_*x_{\ga},y_{\ga^{-1}})^R_{\w_{\ga}}
\end{equation}
where 
$(\cdot,\cdot)^R_{\w_{\ga}}$ 
is the canonical $R$-valued bilinear form on $\HH_{\ga}$
(see Section \ref{bilinear-sec}).

\begin{thm}\label{CohFT-thm}
Let $\w(x_1,\ldots,x_n)$ be a quasihomogeneous polynomial with isolated singularity and
$G\sub G_\w$ a finite subgroup containing $J$. 
The state space $\HH=\HH(\w,G)$, the metric \eqref{metric-eq}, the flat unit $1\in\HH_J\sub\HH$
and the collection of maps $\La^R_g(\ov{\ga})$ define the \CFT\ with coefficients in
$R=\C[\widehat{G}]$.
\end{thm}

In Sections \ref{factor-1-sec}, \ref{factor-2-sec} and \ref{tails-sec} we will check the factorization and
the flat unit axioms for this \CFT. 
In Section \ref{genus-0-3-sec} we will finish the proof of the theorem by verifying the remaining axiom relating the metric on $\HH$
and the maps $\La_0^R(\ga,\ga^{-1},J)$ (for $\ga_1\neq\ga_2^{-1}$ the moduli space 
$\SS_0(\ga_1,\ga_2,J)$ is empty).
To do this we will show that the fundamental matrix factorization on 
each component of $\SS_{0}^{\rig}(\ga,\ga^{-1},J)$
is obtained from the diagonal \mf\ $\De^{\st}_{\w_\ga,\zeta}$
by the action of an element in $G\times G$.

Since the algebra $R$ is isomorphic to the direct sum of algebras
$R\simeq\bigoplus_{\ga'\in G} \C\cdot e_{\ga'}$,
where $e_{\ga'}$ are idempotents \eqref{e-idemp}, our \CFT\ decomposes into
a direct sum of \CFTs\ indexed by elements of $G$. 
By Theorem \ref{hoch-prop}(ii),
$$e_{\ga'}\cdot\HH_{\ga}=H(\w_{\ga,\ga'})^{G},$$
where $\w_{\ga,\ga'}$ is the restriction of $\w$ to the subspace of $\{\ga,\ga'\}$-invariants.
Thus, the sub-\CFT\ corresponding to an element $\ga'\in G$ has the state space
\begin{equation}\label{ga'-state-space-eq}
\HH(\w,G,\ga'):=\bigoplus_{\ga\in G} H(\w_{\ga,\ga'})^{G}.
\end{equation}
and the components of the maps \eqref{R-linear-correlators-eq} are given by 
\begin{equation}\label{lambda-spec-eq}
\La^{\ga'}_g(\ov{\ga})=\frac{1}{\deg(\st_g)}\cdot \left.(\st_g)_*
\pi_{\ga'}\phi_g(\ov{\ga})\right|_{e_{\ga'}\HH_{\ga_1}\ot\ldots\ot e_{\ga'}\HH_{\ga_r}}~,
\end{equation}
where $\pi_{\ga'}:R\to\C$ is the specialization homomorphism.

The \CFT\ corresponding to $\ga'=1$ can be twisted 
to produce a theory satisfying an analog of the {\it concavity}
axiom from \cite{FJR} (see section \ref{concavity-sec}).
Namely, we define twisted maps
$$\phi^{tw}_g(\ov{\ga}):e_1\HH_{\ga_1}\ot\ldots\ot e_1\HH_{\ga_r}
\to H^*(\SS^{\rig}_{g}(\ov{\ga}))$$
for $\ov{\ga}\in G^r$ by
\begin{equation}\label{phi-twisted-eq}
\phi^{tw}_g(\ov{\ga})=\Td\bigl(R\pi_*(\bigoplus_{j=1}^n\LL_j)\bigr)^{-1}\cdot
\left.\pi_{1}\phi_g(\ov{\ga})\right|_{e_1\HH_{\ga_1}\ot\ldots\ot e_1\HH_{\ga_r}}~,
\end{equation}
where $(\LL_\bullet)$ is the universal $\w$-structure and $\Td$ is the
Todd class, and we set
\begin{equation}\label{lambda-1-eq}
\la_g(\ov{\ga})=\exp(\pi i \wt{D}_g(\ov{\ga}))\cdot
\frac{1}{\deg(\st_g)}\cdot (\st_g)_*\phi^{tw}_g(\ov{\ga}),
\end{equation}
where 
$\wt{D}_g(\ov{\ga})$ is given by \eqref{mod-index-eq}.

\begin{thm}\label{twisted-CohFT-thm}
Let $\w$ and $G$ be as in Theorem \ref{CohFT-thm}.
The collection of maps $\la_g(\ov{\ga})$ defines a \CFT\ on the
state space $\HH(\w,G,1)$ with the metric obtained by restricting the metric \eqref{metric-eq}
and the flat unit element $1\in H(\w_{J,1})^{G}$.
\end{thm} 

We will call the \CFT\ of this Theorem the {\it reduced \CFT\ } associated with $(\w,G)$.

The proof of this theorem will be given in Sections \ref{factor-1-sec}, \ref{factor-2-sec}, \ref{tails-sec} and \ref{genus-0-3-sec}
simultaneously with the proof of Theorem \ref{CohFT-thm}.

Note that the state space of the reduced theory
$$\HH(\w,G,1)=\bigoplus_{\ga\in G} H(\w_{\ga})^G$$
can be identified with the state space of the \CFT\ constructed in \cite[sec. 5.1]{FJR}.
Hypothetically, the theories themselves also match.

\medskip

\noindent
{\bf Conjecture.} {\it The reduced \CFT\ associated with $(\w,G)$
is isomorphic to the FJR-theory for the same pair constructed in \cite{FJR}.} 

\medskip

One of the obstacles on the way to proving this conjecture is that it is not clear 
how to verify the homogeneity property of the maps \eqref{lambda-1-eq} in our setup in general.
For some particular cases this homogeneity is proved in Section \ref{dim-sec}.
In particular, it holds for all simple singularities. Combining this with the reconstruction theorem
\cite[Thm.\ 6.2.10]{FJR} we will show in Section \ref{comparison-sec} that
 the above conjecture holds for them. 
In the case of the simple singularity of type $A$, $\w=x^n$, the conjecture
also follows from the results
of Chiodo \cite[sec. 5]{Chiodo} (which imply that in this case
our definition of the reduced \CFT\ is compatible with the construction of \cite{PV}) together
with the reconstruction results of \cite{JKV} and \cite{FSZ} (see also \cite{FJR-An}).

\begin{rems}
1. If $\w'$ is obtained from $\w$ by rescaling of variables then the above constructions for $\w$ and $\w'$ are
naturally identified. In particular, the maps \eqref{hoch-phi-map}, \eqref{lambda-st-g-eq},
\eqref{phi-twisted-eq} and \eqref{lambda-1-eq} are $G_\w$-invariant, where 
$G_\w$ acts naturally on the state spaces (and trivially on the cohomology of the relevant moduli stacks).

\noindent
2. It would be interesting to find quantum K-theory versions of \CFT\ from Theorems \ref{CohFT-thm} and
\ref{twisted-CohFT-thm}
(see \cite{Giv-WDVV}, \cite{Lee}). 
\end{rems}

\subsection{Behavior of fundamental matrix factorizations under gluing}
\label{factor-1-sec}

Let $(\wt{\pi}:\wt{\CC}\to S, p_1,\ldots,p_r,p,q)$ 
be a family of stable orbicurves with $r+2$ marked orbipoints 
over a connected base $S$. 
Also, let
$f:\wt{\CC}\to \CC$ be a morphism of families of stable orbicurves 
over $S$, where the family $\CC$ has a node $\si$ such that 
$f$ is an isomorphism over $\CC\setminus\si$ and $f^{-1}(\si)$ is the union of $p$ and $q$.
We will say that $\CC$ is obtained from $\wt{\CC}$ by gluing points $p$ and $q$.
The remaining marked points $p_1,\ldots,p_r$ can be viewed as marked points on $\CC$. 
Let $(\wt{P},\wt{\vareps})$ be a $\Ga$-spin structure of type
$\wt{\ga}=(\ga_1,\ldots,\ga_r,\ga_p,\ga_q)$
on
$(\wt{\CC},p_1,\ldots,p_r,p,q)$ equipped with a rigidification
$$
(e(1),\ldots,e(r),e(p),e(q)), \text{ where }
e(i)\in \wt{P}|_{p_i}, i=1,\ldots, r, \ e(p)\in \wt{P}|_p \text{ and } e(q)\in \wt{P}|_q.
$$
Note that $\ga_q=\ga_p^{-1}$ since the orbicurve $\CC$ is balanced 
at the node $\si$.

Now we will define a glued $\Ga$-spin structure $(P,\vareps)$ on
$(\CC,p_1,\ldots,p_r)$ such that
$\wt{P}\simeq f^*P$.
Recall (see e.g., \cite[Ch.\ 10.2]{ACG2}) 
that there is a canonical isomorphism $f^*\om_{\CC/S}\simeq \om_{\wt{\CC}/S}(p+q)$, such that
the following diagram is commutative
\begin{equation}\label{opp-res-diagram}
\begin{diagram}
(f^*\om_{\CC/S})|_p&\rTo{}&\om_{\wt{\CC}/S}(p+q)|_p&\rTo{\Res_p}&\OO_S\\
\dTo{}&&&&\dTo{-\id}\\
(f^*\om_{\CC/S})|_q&\rTo{}&\om_{\wt{\CC}/S}(p+q)|_q&\rTo{\Res_q}&\OO_S\\
\end{diagram}
\end{equation}
Hence, we obtain an isomorphism
$f^*\om^{\log}_{\CC/S}\simeq \om^{\log}_{\wt{\CC}/S}$
with a similar property (the residues are opposite).
To get a $\Ga$-spin structure on $\CC$ we will use the element $\zeta_\bullet\in \Ga\sub(\C^*)^n$ 
defined by \eqref{zeta-def-eq} which satisfies
$$\chi(\zeta_\bullet)=-1.$$
Consider the isomorphism
\begin{equation}\label{u-P-eq}
u:\wt{P}|_p\rTo{\sim}\wt{P}|_q
\end{equation}
defined by
\be\label{e-j-q-p-eq}
e(q)=\zeta_\bullet\cdot u(e(p)),
\end{equation}
and let $P$ be the $\Ga$-bundle on $\CC$ obtained from 
$\wt{P}$ by gluing with respect to $u$.
Then the induced trivializations $\chi_*(e(p))$ and $\chi_*(e(q))$ of $\G_m$-torsors
$\chi_*(\wt{P})|_p$ and $\chi_*(\wt{P}_\bullet)|_q$, respectively, satisfy
$$\chi_*(e(q))=-\chi_*(u)(\chi_*(e(p))).$$
Hence, by commutativity of the diagram \eqref{opp-res-diagram}, the isomorphism
$\wt{\vareps}:\chi_*(\wt{P})\to P(\om^{\log}_{\wt{\CC}/S})$ descends to an isomorphism
$\vareps:\chi_*(P)\to P(\om^{\log}_{\CC/S})$.
Thus, we obtain a $\Ga$-spin structure on $(\CC,p_1,\ldots,p_r)$ of type
$\ov{\ga}=(\ga_1,\ldots,\ga_r)$. Furthermore, the 
trivializations $e(i)$, $i=1,\ldots,r$ define a rigidification structure of $(P,\vareps)$.

Now we are going to compare 
the matrix factorizations
$$
\bP\in \ov{\DMF}_{\Ga_\w}(S\times (\A^n)^{\ov{\ga}},-\w_{\ov{\ga}}) \ \text{ and }
\wt{\bP}\in \ov{\DMF}_{\Ga_\w}(S\times (\A^n)^{\wt{\ga}},-\w_{\wt{\ga}}),
$$
associated with the $\Ga$-spin structures $(P,\vareps)$ and $(\wt{P},\wt{\vareps})$, respectively,
by the construction of Section \ref{constr-sec}.

\begin{thm}\label{factor-thm} 
One has isomorphisms
\begin{equation}\label{factor-isom}
\bP\simeq (\id_S\times \pr_{1,\ldots,r}^p)_*(\id_S\times \id\times\De^{\zeta})^*\wt{\bP} 
\end{equation}
and
\begin{equation}\label{main-factor-isom}
\bP\simeq (\id_S\times\pr_{1,\ldots,r})_*\bigl(\wt{\bP}\otimes (\pr_p,\pr_q)^*\De^{\st}_{\w_{\ga},\zeta}\bigr)
\end{equation}
in $\ov{\DMF}_{\Ga_\w}(S\times (\A^n)^{\ov{\ga}},-\w_{\ov{\ga}}),$
where 
$\pr_{1,\ldots,r}^p:(\A^n)^{\ov{\ga}}\times(\A^n)^{\ga_p}\to (\A^n)^{\ov{\ga}}$ and
$\pr_{1,\ldots,r}:(\A^n)^{\wt{\ga}}\to (\A^n)^{\ov{\ga}}$ are
the coordinate projections,
$$\De^{\zeta}:(\A^n)^{\ga_p}\to (\A^n)^{\ga_p}\times (\A^n)^{\ga_q}: x\mapsto (\zeta_\bullet\cdot x,x)$$
is the shifted diagonal, and
$$\De^{\st}_{\w_{\ga},\zeta}:=(\zeta_\bullet,\id)^*\De^{\st}_{\w_{\ga}}$$
is the shifted diagonal \mf.
\end{thm}

We are going to prove Theorem \ref{factor-thm} by going through the steps of the construction 
of Section \ref{constr-sec} for both $\Ga$-spin structures $(P,\vareps)$ and $(\wt{P},\wt{\vareps})$.

For $j=1,\ldots, n$, let
$$\Sigma_j\sub p_1+\ldots+ p_r\sub \CC \ \text{ and }
\wt{\Sigma}_j\sub p_1+\ldots+ p_r+ p+ q\sub \wt{\CC}$$
be the subdivisors defined by \eqref{Sigma-j-eq} for
  the collections $\ov{\ga}$ and 
$\wt{\ga}=(\ov{\ga},\ga_p,\ga_q)$, respectively.

For each monomial $M$ in $\w$ we set $\Sigma_M=\cap_{\deg_{x_j}M>0}\Sigma_j$ and
$\wt{\Sigma}_M=\cap_{\deg_{x_j}M>0}\wt{\Sigma}_j$. Note that 
$$\Sigma_j=\wt{\Sigma}_j\cap\{p_1+\ldots+p_r\} \text{ and }
\Sigma_M=\wt{\Sigma}_M\cap\{p_1+\ldots+p_r\}.$$

Let $C$ and $\wt{C}$ be the curves obtained from $\CC$ and $\wt{\CC}$ by forgetting the orbifold structure
at $p_1,\ldots,p_r$, let $\rho:\CC\to C$, $\wt{\rho}:\wt{\CC}\to\wt{C}$ be the natural projections. 
We still denote by $f:\wt{C}\to C$ the morphism induced by $f:\wt{\CC}\to\CC$ and by $\si:S\to C$
the node obtained by gluing $p$ and $q$. Also, 
let $(\LL_{\bullet})$ (resp., $(\wt{\LL}_{\bullet})$) be the line bundles associated with the $\Ga$-bundle
$P$ (resp., $\wt{P}$), and let
$L_j=\rho_*\LL_j$ and $\wt{L}_j=\wt{\rho}_*\wt{\LL}_j$.

We have the natural exact sequence
\begin{equation}\label{f*-seq}
0\to \OO_{C}\to f_*\OO_{\wt{C}}\rTo{\delta} \si_*\OO_S\to 0.
\end{equation}
where the map $\delta$ is the difference of evaluations along $p$ and $q$.
Tensoring this sequence with $L_j$ and using the projection formula we get an exact triangle
\begin{equation}\label{L-f-triangle-0}
R\pi_*(L_j)\to R\wt{\pi}_*(\wt{L}_j)\to \pi_*(\si^*L_j)\to\ldots
\end{equation}
The third term depends on the action of the local group at $\si$ on $\si^*L_j$.
We have two cases:

\noindent 
{\bf Case 1}.
$\wt{\Sigma}_j=\Sigma_j$ (i.e., $\ga_p$ has a nontrivial $j$th component). Then the third
term of the triangle \eqref{L-f-triangle-0} vanishes, so we have an isomorphism
\be\label{push-forward-equality}
R\pi_*(L_j)\simeq R\wt{\pi}_*(\wt{L}_j).
\end{equation}

\noindent
{\bf Case 2}.
$\wt{\Sigma}_j=\Sigma_j+p+q$ (i.e., $\ga_p$ has trivial $j$th component).
Then we get an exact triangle
\be\label{L-f-triangle}
R\pi_*(L_j)\to R\wt{\pi}_*(\wt{L}_j)\rTo{r_p-u_j^{-1}r_q} (\wt{L}_j)|_p\to\ldots,
\end{equation}
where $r_p:R\wt{\pi}_*(\wt{L}_j)\to \wt{L}_j|_p$ and $r_q:R\wt{\pi}_*(\wt{L}_j)\to \wt{L}_j|_q$
are the restriction maps and $u_j:\wt{L}_j|_p\to \wt{L}_j|_q$ is the isomorphism induced by
\eqref{u-P-eq}.

Let $Z_j:R\pi_*(L_j)\to\OO_S^{\Si_j}$ and $\wt{Z}_j:R\wt{\pi}_*\wt{L}_j\to\OO_S^{\wt{\Si}_j}$
be the maps \eqref{Z-j-eq} induced by the rigidifications of the $\Ga$-spin structures $(P,\vareps)$ and 
$(\wt{P},\wt{\vareps})$.
In the second case we will also have to consider
the maps 
\begin{equation}\label{Z'-j-eq}
Z'_j:R\pi_*(L_j)\to R\wt{\pi}_*\wt{L}_j\rTo{\wt{Z}_j}\OO_S^{\wt{\Si}_j}.
\end{equation}
Note that the components of this map corresponding to the points $p$ and $q$ satisfy
\begin{equation}\label{Z'-p-q-eq}
Z'_j(p)=\zeta_j\cdot Z'_j(q)
\end{equation}
(this follows from \eqref{e-j-q-p-eq}).

Now let $\FF_M$ (resp., $\wt{\FF}_M$) be the
coherent sheaf on $C_M$ (resp., $\wt{C}_M$) defined by \eqref{F-M-eq}
for our $\Ga$-spin structures over $\CC$ and $\wt{\CC}$. By Step 1 of the construction of \ref{constr-sec},
we have canonical maps 
$$\tau_M:E_M=R\pi^M_*(\FF_M)^{\Sym(M)}\to\OO_S[-1] \ \text{ and } \ 
\wt{\tau}_M:\wt{E}_M=R\wt{\pi}^M_*(\FF_M)^{\Sym(M)}\to\OO_S[-1].$$
We are going to establish a certain compatibility
between these maps (see Lemma \ref{factor-lem} below).

First, we need some preparations.
Define a subbundle $\OO_S^{\wt{\Si}_M}(p,q)\sub\OO_S^{\wt{\Si}_M}$ by
$$\OO_S^{\wt{\Si}_M}(p,q)=\begin{cases}\OO_S^{\wt{\Si}_M}, & \text{if }\wt{\Si}_M=\Si_M,\\
\ker(\pi_p+\pi_q:\OO_S^{\wt{\Si}_M}\to\OO_S), & \text{if }\wt{\Si}_M=\Si_M+p+q,\end{cases}$$
where $\pi_p,\pi_q:\OO_S^{\wt{\Si}_M}\to \OO_S$ are the projections corresponding to 
the marked points $p$ and $q$, respectively. 

Suppose first that $\wt{\Si}_M=\Si_M+p+q$. This happens exactly
when for each $j$ with $\deg_{x_j}M>0$ the element $\ga_p$ has trivial
$j$th component. In this case the action of the local group at the node
$\sigma$ on $\si^*M(L_1,\ldots,L_n)\simeq\OO_{\si}$ is trivial.
By \eqref{Z'-p-q-eq}, for every monomial $M$ in $\w$ we have
$$M(Z'_\bullet)_q=-M(Z'_\bullet)_p,$$
so $M(Z'_\bullet)$ induces a morphism
$$Z'_M:M(R\pi_*(L_\bullet))\to \OO_S^{\wt{\Si}_M}(p,q).$$
Let us define a coherent sheaf $\FF'_M$ on $C^M$ from the exact sequence
$$0\to \FF'_M\to L_\bullet^M\to (\De_M)_*(\OO_{\Si_M}\oplus\OO_\si)\to 0$$
and set 
$$E'_M=R\pi^M_*(\FF'_M)^{\Sym(M)},$$
where $\Sym(M)$ is the product of symmetric groups \eqref{Sym-M-eq}.
Then the above exact sequence gives rise to an exact triangle
$$E'_M\to M(R\pi_*(L_\bullet))\rTo{Z'_M}\OO_S^{\wt{\Si}_M}(p,q)\to E'_M[1].$$
Now let $f^M:\wt{C}^M\to C^M$ be the map induced by the morphism $f:\wt{C}\to C$.
We have a commutative diagram with exact rows
\begin{equation}\label{FFF-diagram-eq}
\begin{diagram}
0&\rTo{}&\FF_M&\rTo{\eps}&L_\bullet^M&\rTo{}&(\De_M)_*\OO_{\Sigma_M}&\rTo{} &0\\
&&\uTo{}&&\uTo{\id}&&\uTo{}\\
0&\rTo{}&\FF'_M&\rTo{}& L_\bullet^M&\rTo{}&(\De_M)_*(\OO_{\Sigma_M}\oplus\OO_{\si})&\rTo{} &0\\
&&\dTo{}&&\dTo{}&&\dTo{(\id,k)}\\
0&\rTo{} &f^M_*\wt{\FF}_M&\rTo{}&(f_*\wt{L}_\bullet)^M&\rTo{}&
\De_*(\OO_{\Sigma_M}\oplus\OO_{\si}^{\oplus 2})&\rTo{}& 0
\end{diagram}
\end{equation}
where $k:\OO_{\si}\to\OO_{\si}^{\oplus 2}: x\mapsto (x,-x)$.
It induces a commutative diagram in $D^b(S)$ whose rows are exact triangles
\begin{equation}\label{triangles-diagram-eq}
\begin{diagram}
E_M&\rTo{\eps}&M(R\pi_*(L_\bullet))&\rTo{Z_M}&\OO_S^{\Sigma_M}&
\rTo{i}&E_M[1]\\
\uTo{f}&&\uTo{\id}&&\uTo{\pr_{\Si_M}}&&\uTo{f[1]}\\
E'_M&\rTo{\eps'}&M(R\pi_*(L_\bullet))&\rTo{Z'_M}&\OO_S^{\wt{\Sigma}_M}(p,q)&
\rTo{i}&E'_M[1]\\
\dTo{g}&&\dTo{}&&\dTo{}&&\dTo{g[1]}\\
\wt{E}_M&\rTo{\wt{\eps}}&M(R\wt{\pi}_*(\wt{L}_\bullet))&\rTo{\wt{Z}_M}&
\OO_S^{\wt{\Sigma}_M}&
\rTo{i}&\wt{E}_M[1]
\end{diagram}
\end{equation}
where the map $\pr_{\Sigma_M}$ is induced by the natural
projection $\OO_S^{\wt{\Si}_M}\to\OO_S^{\Si_M}$ and $\wt{Z}_M$ is defined similarly to
$Z_M$ for the data associated with $\wt{C}$.

In the case $\wt{\Si}_M=\Si_M$ we still have the commutative diagram \eqref{triangles-diagram-eq}
with $E'_M=E_M$ and the morphism from the second row to the first being identity.

\begin{lem}\label{factor-lem} 
The diagram
\begin{equation}\label{factor-diagram}
\begin{diagram}
E'_M&\rTo{f}& E_M\\
\dTo{g}&&\dTo{\tau_M}\\
\wt{E}_M&\rTo{\wt{\tau}_M}&\OO_S[-1]
\end{diagram}
\end{equation}
is commutative. 
\end{lem}

\Pf . Recall  (see \cite{Hart-RD}) that there is a natural map 
$$\Tr_f:f_*\om_{\wt{C}/S}\to\om_{C/S}$$
such that the composition
$$f_*\om_{\wt{C}/S}\rTo{\Tr_f}\om_{C/S}\to f_*f^*\om_{C/S}\simeq f_*(\om_{\wt{C}/S}(p+q))$$
is the natural map induced by the embedding $\om_{\wt{C}/S}\to\om_{\wt{C}/S}(p+q)$, while
the composition
$$R\pi_*(f_*\om_{\wt{C}/S})\rTo{\Tr_f}R\pi_*(\om_{C/S})\rTo{\Tr_{C/S}}\OO_S[-1]$$
can be identified with $\Tr_{\wt{C}/S}$.
By definition of the maps $\tau_M$ and $\wt{\tau}_M$, 
to show that diagram \eqref{factor-diagram} is commutative
it is sufficient to check commutativity of the diagram
\begin{diagram}
\FF'_M&\rTo{}& \FF_M\\
\dTo{}&&&\rdTo{\kappa_M}\\
f^M_*\wt{\FF}_M&\rTo{\wt{\kappa}_M}&\De_*f_*\om_{\wt{C}/S}&\rTo{\Tr_f}&\De_*(\om_{C/S})
\end{diagram}
This follows from the commutativity of the diagram
\begin{diagram}
\FF'_M&\rTo{}& L_\bullet^M&\rTo{}&\De_*(\om_{C/S}(\Si_M))\\
\dTo{}&&\dTo{}&&\dTo{}\\
f^M_*\wt{\FF}_M&\rTo{}&(f_*\wt{L_\bullet})^M&\rTo{}&\De_*f_*(\om_{\wt{C}/S}(\wt{\Si}_M))
\end{diagram}
(see \eqref{FFF-diagram-eq}) since the natural morphism of sheaves
$$\De_*(\om_{C/S})\to\De_*(\om_{C/S}(\Si_M))\to\De_*f_*(\om_{\wt{C}/S}(\wt{\Si}_M)).$$
is injective.
\ed

As in Step 2 of the construction of the fundamental \mf\ in Section \ref{constr-sec}, we 
choose for each $j$ a $\pi$-acyclic resolution 
\begin{equation}\label{L-D-Q-res}
L_j\to [P_j\oplus \left.L_j\right|_{\Si_j\cup\si}~\to Q_j].
\end{equation}
However, we will modify these resolutions using the map $f:\wt{C}\to C$.

In Case 1, when $\wt{\Si}_j=\Si_j$, the pull-back of \eqref{L-D-Q-res} by $f$ gives a $\wt{\pi}$-acyclic
resolution
$$\wt{L}_j\to [f^*P_j\oplus\left.\wt{L}_j\right|_{\wt{\Si}_j+p+q}~\to f^*Q_j]$$
and a compatible resolution of $\wt{L}_j(-\wt{\Si}_j)$.
In Case 2, when $\wt{\Si}_j=\Si_j+p+q$, 
the pull-back of \eqref{L-D-Q-res} gives a $\wt{\pi}$-acyclic resolution 
$$\wt{L}_j\to [f^*P_j\oplus\left.\wt{L}_j\right|_{\wt{\Si}_j}~\to f^*Q_j]$$
and a compatible resolution of $\wt{L}_j(-\wt{\Si}_j)$. 
In both cases pushing forward by $f$ we get a two-term resolution for
$f_*\wt{L}_j$ such that both the difference map $f_*\wt{L}_j\to \left.\wt{L}_j\right|_{\si}~$ and
the map $f_*\wt{L}_j\to L_j|_{\Si_j}$ factor through this resolution. This leads
to a new $\pi$-acyclic resolution for $L_j$ of the form
\begin{equation}\label{L-j-new-res-eq}
L_j\to \bigl[\ker\bigl(f_*\OO\ot (P_j\oplus L_j|_{\Si_j\cup\si})\to L_j|_{\si}\bigr)\to f_*\OO\ot Q_j\bigr]
\end{equation}
and a compatible resolution for $L_j(-\Si_j)$,
such that the map $L_j\to L_j|_{\Si_j\cup\si}$ is induced by a surjective map from the resolution
in \eqref{L-j-new-res-eq}.
Pushing forward these resolutions to $S$ we obtain complexes 
\begin{equation}\label{4-res-eq}
[\wt{A}'_j\to B_j]\sub [A'_j\to B_j]\sub [A_j\to B_j]\sub [\wt{A}_j\to B_j]
\end{equation}
representing $R\wt{\pi}_*(\wt{L}_j(-\wt{\Si_j}))$, $R\pi_*(L_j(-\Si_j))$, $R\pi_*(L_j)$ and 
$R\wt{\pi}_*\wt{L}_j$, respectively.
Note that we have an exact triple
$$0\to \wt{A}'_j\to\wt{A}_j\rTo{\wt{Z}_j}\OO_{S}^{\wt{\Si}_j}\to 0$$
and a similar exact triple relating $A'_j$ and $A_j$.
In the case $\wt{\Si}_j=\Si_j$ we have an equality $A_j=\wt{A_j}$, while in
the case $\wt{\Si}_j=\Si_j+p+q$ we have
an exact triple of complexes
$$0\to [A_j\to B_j]\hra [\wt{A}_j\to B_j]\rTo{\de_j(p,q)}\OO_S\to 0$$
representing the exact triangle \eqref{L-f-triangle}, where the morphism $\de_j(p,q):\wt{A}_j\to \OO_{S}$
is given by
$$\de_j(p,q)=\wt{Z}_j(p)-\zeta_j\cdot\wt{Z}_j(q).$$
(here $\wt{Z}_j(p)$ and $\wt{Z}_j(q)$ are the components of $\wt{Z}_j$).
The maps $Z'_j:A_j\to\OO_S^{\wt{\Si}_j}$ realizing \eqref{Z'-j-eq} 
are defined as compositions 
$$Z'_j:A_j\to\wt{A}_j\rTo{\wt{Z}_j}\OO_S^{\wt{\Si}_j}$$
while $Z_j:A_j\to\OO_S^{\Si_j}$ are obtained from them by projecting to $\OO_S^{\Si_j}$.

Note that all four resolutions \eqref{4-res-eq} have the same second term $B_j$. Thus,
as in Step 2 
we can assume that $B_j=\VV^{\vee}(-m_0)^{\oplus N_j}$ for large enough $m_0$, so that
$$
\Ext^{>0}_S((\oplus_j\wt{A}_j)^{\ot q_1}\ot(\VV^{\vee})^{\ot q_2}(-m),\OO_S)=
\Ext^{>0}_S((\oplus_j A_j)^{\ot q_1}\ot(\VV^{\vee})^{\ot q_2}(-m),\OO_S)=
0 $$
for $m\ge m_0$ and $q_1+q_2\le d$.
Hence, we can assume that 
$$\Ext^{>0}_S((\oplus_j\wt{A}_j)^{\ot q_1}\ot(\oplus_j B_j)^{\ot q_2},\OO_S)=
\Ext^{>0}_S((\oplus_j A_j)^{\ot q_1}\ot(\oplus_j B_j)^{\ot q_2},\OO_S)=0$$
for $q_1+q_2\le d$ and $q_2\ge 1$.

As in Step 2 we realize $E_M$, $\wt{E}_M$ and $E'_M$ by the complexes
\begin{align*}
&K_M=\Cone(M([A_\bullet\to B_\bullet])\rTo{Z_M}\OO_S^{\Si_M})[-1], \ 
\wt{K}_M=\Cone(M([\wt{A}_\bullet\to B_\bullet])\rTo{\wt{Z}_M}\OO_S^{\wt{\Si}_M})[-1],
\ \text{ and }\\
&K'_M=\Cone(M([A_\bullet\to B_\bullet])\rTo{Z'_M}\OO_S^{\wt{\Si}_M}(p,q))[-1],
\end{align*}
respectively.
Then for all $M$ appearing in $\w$ we will have
\begin{equation}\label{ext-van-factor2-eq}
\Ext^{>0}_S(\wt{K}^i_M,\OO_S)=
\Ext^{>0}_S(K^i_M,\OO_S)=
\Ext^{>0}_S(K^{\prime i}_M,\OO_S)=0 \ \text{ for } i\ge 2.
\end{equation}

As before we consider the spaces
$$X:=\tot(A_1\oplus\ldots\oplus A_n)\sub \wt{X}:=\tot(\wt{A}_1\oplus\ldots\oplus\wt{A}_n)$$
over $S$.
We have the following diagram with a cartesian square
\begin{equation}\label{factor-diagram2}
\begin{diagram}
X&\rTo{q}&\prod_{i=1}^r(\A^n)^{\ga_i}\times (\A^n)^{\ga_p}&\rTo{\pr_{1,\ldots,r}^p}&
\prod_{i=1}^r(\A^n)^{\ga_i}\\
\dInto{}& & \dTo{\id\times\De^{\zeta}}\\
\wt{X} &\rTo{\wt{Z}}& \prod_{i=1}^r(\A^n)^{\ga_i}\times(\A^n)^{\ga_p}\times(\A^n)^{\ga_q}
\end{diagram}
\end{equation}
where the maps $\wt{Z}$ and $q$ are given by $\wt{Z}_1,\ldots,\wt{Z}_n$ and
$Z'_1,\ldots,Z'_n$, respectively, and 
the composition of two horizontal arrows in the first row is equal to $Z$.

\bigskip

\noindent
{\it Proof of Theorem \ref{factor-thm}.}
Isomorphism \eqref{main-factor-isom} follows from
\eqref{factor-isom}  by Proposition \ref{diag-mf-lem}(ii), so we only need to prove 
\eqref{factor-isom}.

As in Step 2 of the construction,
the vanishing \eqref{ext-van-factor2-eq} can be used to choose for each monomial $M$ 
chain maps representing $\tau_M$ and $\wt{\tau}_M$. Also, the morphisms in derived category
from $K'_M$ to $\OO_S[-1]$ can be calculated in the homotopy category.
Hence, we obtain a realization of the commutative diagram \eqref{factor-diagram} in 
the homotopy category.
Since the complex $K'_M$ differs from $K_M$ only in terms of degree $0$ and $1$,
we can replace $\tau_M$ by a homotopic chain map, so that the diagram \eqref{factor-diagram}
will be commutative in the category of complexes.

Let $\wt{p}:\wt{X}\to S$ and $p:X\to S$ be the projections.
Recall that we have a section $\wt{\b}$ (resp., $\b$) 
of $\wt{p}^*(\bigoplus_j B_j)$ (resp., $p^*(\bigoplus_j B_j)$), corresponding to the differential
$\bigoplus_j \wt{A}_j\to \bigoplus_j B_j$ (resp., $\bigoplus_j A_j\to\bigoplus_j B_j$), and
$\b$ is equal to the restriction of $\wt{\b}$ to $X$ by construction.
As in Step 3 of the construction, 
from the chain map $\wt{\tau}_M$ (resp., $\tau_M$) we get a section $\wt{\a}_M$ 
(resp., $\a_M$) of the bundle $\wt{p}^*(\bigoplus_j B_j^\vee)$ 
(resp., $p^*(\bigoplus_j B_j^\vee)$)
satisfying
$$\lan\wt{\a}_M,\wt{\b}\ran=\wt{Z}^*M^{\oplus \wt{\Si}_M}
\text{ (resp., } \lan \a_M,\b\ran= Z^*M^{\oplus \Si_M}).$$ 
The
components of degree $1$ of the diagram \eqref{factor-diagram} form the following 
commutative square
\begin{diagram}
\OO_S^{\wt{\Si}_M}(p,q)\oplus\bigoplus_j\pa_jM(A_\bullet)\ot B_j &\rTo{}&
\OO_S^{\Si_M}\oplus\bigoplus_j\pa_jM(A_\bullet)\ot B_j \\
\dTo{}&&\dTo_{(\Tr,\a_M)}\\
\OO_S^{\wt{\Si}_M}\oplus\bigoplus_j\pa_jM(\wt{A}_\bullet)\ot B_j 
&\rTo{(\Tr,\wt{\a}_M)}&\OO_S
\end{diagram}
This implies
that the restriction of $\wt{\a}_M$ to $X$ is equal to $\a_M$. 
Therefore, the matrix factorization $\{-\a_\w,\b\}$ on $X$ is isomorphic to the restriction 
of $\{-\wt{\a}_\w,\wt{\b}\}$.
As in Step 4, we have
$$\wt{\bP}=(\wt{p},\wt{Z})_*\{-\wt{\a}_\w,\wt{\b}\} \ \text{ and }\ \bP=(p,Z)_*\{-\a_\w,\b\}.$$
Since $Z$ is the composition of arrows in the first row of the diagram \eqref{factor-diagram2},
we have
$$\bP\simeq (\id_S\times\pr_{1,\ldots,r}^p)_*\bigl((p,q)_*\{-\a_\w,\b\}\bigr).$$
On the other hand, by the base change formula (see Proposition \ref{base-change-prop}) in the cartesian square 
\begin{diagram}
X&\rTo{(p,q)}& S\times \prod_{i=1}^r(\A^n)^{\ga_i}\times (\A^n)^{\ga_p}\\
\dInto{}& & \dTo{\id_S\times\id\times\De^{\zeta}}\\
\wt{X} &\rTo{(\wt{p},\wt{Z})}& S\times\prod_{i=1}^r(\A^n)^{\ga_i}\times(\A^n)^{\ga_p}\times(\A^n)^{\ga_q}
\end{diagram}
we have
$$(p,q)_*\{-\a_\w,\b\}\simeq (\id_S\times\id\times\De^{\zeta})^*\wt{\bP},$$
which implies the required isomorphism \eqref{factor-isom}.
\ed

\subsection{Verification of the factorization axiom of \CFT}\label{factor-2-sec}

The main axiom of \CFT\ describes the factorization property of the maps
$\La_g(\ov{\ga})$ under the gluing morphisms
$$\rho_{\tree}:\ov{\MM}_{g_1,r_1+1}\times\ov{\MM}_{g_2,r_2+1}\to\ov{\MM}_{g,r} \ \ \text{ and } \ \ 
\rho_{\petlya}:\ov{\MM}_{g-1,r+2}\to\ov{\MM}_{g,r},$$
where $g_1+g_2=g$ and $r_1+r_2=r$ (see \cite[2.2.6, 2.2.7]{KonM}). Here we will verify this axiom for the maps \eqref{lambda-st-g-eq}
of our \CFT\ on the state space $\HH=\HH(\w,G)$ (see \eqref{state-space-eq}) and also for the twisted maps
\eqref{lambda-1-eq} on the state space $\HH(\w,G,1)$.
Specifically, in this section we prove the equalities 
\begin{equation}\label{factor-axiom-eq}
\begin{array}{l}
(\rho_{\tree})^*\circ\La_g(\ga_1,\ldots,\ga_{r_1};\ga'_1,\ldots,\ga'_{r_2})=\\
\sum_{\ga\in G}(\La_{g_1}(\ga_1,\ldots,\ga_{r_1},\ga)\ot\La_{g_2}(\ga'_1,\ldots,\ga'_{r_2},\ga^{-1}))\circ
(\id^{\ot r}\ot T_{\w_{\ga},\zeta}).
\end{array}
\end{equation}
\begin{equation}\label{factor-axiom-loop-eq}
(\rho_{\petlya})^*\circ\La_g(\ga_1,\ldots,\ga_r)=
\sum_{\ga\in G}(\La_{g-1}(\ga_1,\ldots,\ga_r,\ga,\ga^{-1})\circ
(\id^{\ot r}\ot T_{\w_{\ga},\zeta}),
\end{equation}
where $\ga_i$ are elements of $G$ and
\begin{align}\label{T-w-zeta-eq}
&T_{\w,\zeta}=\frac{1}{|G|}\cdot\sum_{h\in G}(\id\times h)^*\ch(\De^{\st}_{\w,\zeta})\in &
HH_*(\MF_{\Ga}(\A^n\times\A^n,\w\oplus\w))^{G\times G}\simeq \nonumber \\ 
&&HH_*(\MF_{\Ga}(\w))\ot_R
HH_*(\MF_{\Ga}(\w))
\end{align}
(the last identification follows from Corollary \ref{Kunneth-cor}).
The proof of the axiom will be finished once we check that $(T_{\w_{\ga},\zeta})_{\ga\in G}$ 
are exactly the components of the Casimir element for the metric on $\HH$. 
This will be done in Lemma \ref{metric-lem} below.

To verify \eqref{factor-axiom-eq} we first observe that due to condition \eqref{gamma-J-eq} 
the only potentially nontrivial summand in the right-hand side of \eqref{factor-axiom-eq} corresponds to
$$\ga=(\ga_1\cdot\ldots\cdot\ga_{r_1})^{-1}J^{2g_1-1+r_1}=
\ga'_1\cdot\ldots\cdot\ga'_{r_2}J^{-2g_2+1-r_2}.$$
We have the following commutative diagram
involving the moduli spaces of $\ov{\w}$-structures with rigidifications:
\begin{equation}\label{rig-factor-diagram}
\begin{diagram}
\SS_{g_1}^{\rig}(\ga_1,\ldots,\ga_{r_1},\ga)\times
\SS_{g_2}^{\rig}(\ga'_1,\ldots,\ga'_{r_2},\ga^{-1}) &\rTo{\rho_{\tree,\ga}^{\rig}}&
\SS_{g}^{\rig}(\ga_1,\ldots,\ga_{r_1};\ga'_1,\ldots,\ga'_{r_2})\\
\dTo{\st_{g_1}\times\st_{g_2}}&&\dTo{\st_g}\\
\ov{\MM}_{g_1,r_1+1}\times\ov{\MM}_{g_2,r_2+1}&\rTo{\rho_{\tree}}&\ov{\MM}_{g,r}
\end{diagram}
\end{equation}
where the maps
$\st_g$, $\st_{g_1}$ and $\st_{g_2}$ are the natural projections.
The map $\rho^{\rig}_{\tree,\ga}$ is given by the gluing construction described in the beginning of
Section \ref{factor-1-sec}. 
Theorem \ref{factor-thm} applied to 
$$
S=\SS_{g_1}^{\rig}(\ga_1,\ldots,\ga_{r_1},\ga)\times
\SS_{g_2}^{\rig}(\ga'_1,\ldots,\ga'_{r_2},\ga^{-1})$$
gives the following relation between the fundamental matrix factorizations:
\begin{equation}\label{factor-rig-isom}
\begin{array}{l}
(\rho^{rig}_{\tree,\ga}\times\id)^*\bP^{\rig}_{g,\Ga}(\ga_1,\ldots,\ga_{r_1};\ga'_1,\ldots,\ga'_{r_2})\simeq \\
(\id_S\times\pr_{1,\ldots,r})_*
[\pi_1^*\bP^{\rig}_{g_1,\Ga}(\ga_1,\ldots,\ga_{r_1},\ga)\ot\pi_2^*\bP^{\rig}_{g_2,\Ga}(\ga'_1,\ldots,\ga'_{r_2},\ga^{-1})\ot
\pr_{p,q}^*\De^{\st}_{\w_{\ga},\zeta}],
\end{array}
\end{equation}
where $\pi_1$ and $\pi_2$ are the projections of
$S\times\prod_{i=1}^{r_1}(\A^n)^{\ga_i}\times\prod_{i=1}^{r_2}(\A^n)^{\ga'_i}\times
(\A^n)^{\ga}\times (\A^n)^{\ga^{-1}}$ to
$S\times\prod_{i=1}^{r_1}(\A^n)^{\ga_i}\times (\A^n)^{\ga}$ and
$S\times\prod_{i=1}^{r_2}(\A^n)^{\ga'_i}\times (\A^n)^{\ga^{-1}}$, respectively, and
$\pr_{p,q}$ is the projection to $(\A^n)^{\ga}\times(\A^n)^{\ga^{-1}}$.

The functor 
$$\Phi: 
\ov{\DMF}_{\Ga}(\A^{\ov{\ga}}\times\A^{\ov{\ga'}},\w_{\ov{\ga}}\oplus \w_{\ov{\ga'}})\to D_{G}(S)$$
given by the kernel in the left-hand side of \eqref{factor-rig-isom},
where $\ov{\ga}=(\ga_1,\ldots,\ga_{r_1})$ and $\ov{\ga}'=(\ga'_1,\ldots,\ga'_{r_2})$,
is isomorphic to the composition 
$$\ov{\DMF}_{\Ga}(\A^{\ov{\ga}}\times\A^{\ov{\ga'}},\w_{\ov{\ga}}\oplus \w_{\ov{\ga'}})
\rTo{\Phi_g(\ov{\ga},\ov{\ga}')}D_{G}(\SS_{g}^{\rig}(\ov{\ga},\ov{\ga}'))
\rTo{(\rho^{\rig}_{\tree,\ga})^*}D_{G}(S).
$$
After passing to Hochschild homology and using the HKR map \eqref{stack-map} we obtain that
the map induced by $\Phi$
$$\phi:\bigotimes_{k=1}^{r_1}\HH_{\ga_k}\ot
\bigotimes_{l=1}^{r_2}\HH_{\ga'_l}\rTo{\iota} HH_*(\MF_{\Ga}(\A^{\ov{\ga}}\times\A^{\ov{\ga'}},\w_{\ov{\ga}}\oplus \w_{\ov{\ga'}}))
\rTo{\Phi_*} H^*(S,\C)\ot R$$
(where all tensor products are taken over $R$ and $\iota=\iota(\ov{\ga},\ov{\ga'})$ is the map \eqref{iota-isom})
is equal to the composition 
$(\rho^{\rig}_{\tree,\ga})^*\circ\phi_g(\ov{\ga},\ov{\ga}')$, where $\phi_g(\ov{\ga},\ov{\ga'})$ is the map \eqref{hoch-phi-map}.
On the other hand, let 
$$\phi':\bigotimes_{k=1}^{r_1}\HH_{\ga_k}\ot
\bigotimes_{l=1}^{r_2}\HH_{\ga'_l}\rTo{\iota} 
HH_*(\MF_{\Ga}(\A^{\ov{\ga}}\times\A^{\ov{\ga'}},\w_{\ov{\ga}}\oplus \w_{\ov{\ga'}}))\rTo{\Phi'_*}
H^*(S,\C)\ot R$$ 
be the composition of the HKR map
with the map induced on Hochschild homology 
by the functor $\Phi'$ given by the kernel in
the right-hand side of \eqref{factor-rig-isom}.
Using the projection formula for $\id_S\times \pr_{1,\ldots,r}$ (see Proposition \ref{proj-formula-prop})
we see that $\Phi'_*$ 
is equal to the composition 
\begin{align*}
&HH_*(\MF_{\Ga}(\A^{\ov{\ga}}\times\A^{\ov{\ga'}},\w_{\ov{\ga}}\oplus \w_{\ov{\ga'}}))\rTo{\ot\pr_{p,q}^*\ch(\De^{\st}_{\w_{\ga},\zeta})}
\\
&HH_*(\MF_{\Ga}(\A^{\ov{\ga}}\times\A^{\ov{\ga'}}\times(\A^n)^{\ga}\times(\A^n)^{\ga^{-1}},
\w_{\ov{\ga}}\oplus \w_{\ov{\ga'}}\oplus\w_{\ga}\oplus\w_{\ga^{-1}}))\rTo{\wt{\Phi}_*}H^*(S,\C)\ot R,
\end{align*}
where $\wt{\Phi}_*$ is induced by the functor $\wt{\Phi}$ associated with the kernel
$\pi_1^*\bP^{\rig}_{g_1,\Ga}(\ga_1,\ldots,\ga_{r_1},\ga)\ot
\pi_2^*\bP^{\rig}_{g_2,\Ga}(\ga'_1,\ldots,\ga'_{r_2},\ga^{-1})$.
Since $\wt{\Phi}_*$ is invariant with respect to the action of $G\times G$ on the factors $(\A^n)^{\ga}\times (\A^n)^{\ga^{-1}}$,
we can replace $\ch(\De^{\st}_{\w_{\ga},\zeta})$ by its $G\times G$-averaging $T_{\w_{\ga},\zeta}$
in the above formula for $\Phi'_*$. Since $\iota$ is exactly the embedding of $G^{r_1+r_2}$-invariants (see Corollary \ref{Kunneth-cor}), 
we obtain
that $\phi'$ is equal to the composition
$$\bigotimes_{k=1}^{r_1}\HH_{\ga_k}\ot
\bigotimes_{l=1}^{r_2}\HH_{\ga'_l}\rTo{\id^{\ot r}\ot T_{\w_{\ga},\zeta}} 
\bigotimes_{k=1}^{r_1}\HH_{\ga_k}\ot
\bigotimes_{l=1}^{r_2}\HH_{\ga'_l}\ot \HH_{\ga}\ot_R\HH_{\ga^{-1}}
\rTo{\phi_{g_1}\ot\phi_{g_2}} H^*(S,\C)\ot R,$$
where $\phi_{g_1}=\phi_{g_1}(\ga_1,\ldots,\ga_{r_1},\ga)$ and
$\phi_{g_2}=\phi_{g_2}(\ga'_1,\ldots,\ga'_{r_2},\ga^{-1})$.
Since the functors $\Phi$ and $\Phi'$ are isomorphic, we have 
\begin{equation}\label{phi-phi'-eq}
\phi=\phi'.
\end{equation}
The desired formula \eqref{factor-axiom-eq} is obtained from this by applying $(\st_{g_1}\times\st_{g_2})_*$
taking into account the relation
\begin{equation}\label{fact-st-eq}
\frac{1}{\deg(\st_g)}\rho_{\tree}^*(\st_g)_*=
\sum_{\ga\in G}\frac{1}{\deg(\st_{g_1})\deg(\st_{g_2})}(\st_{g_1}\times\st_{g_2})_*
(\rho^{\rig}_{\tree,\ga})^*,
\end{equation}
which holds because the space in the left upper corner of diagram
\eqref{rig-factor-diagram} is an \'etale covering of the fibered product of $\rho_{\tree}$ and
$\st_g$.

The proof of \eqref{factor-axiom-loop-eq} is analogous: one has to apply Theorem
\ref{factor-thm} to compare \mf s over
$$S=\SS_{g-1}^{\rig}(\ga_1,\ldots,\ga_r,\ga,\ga^{-1}).$$
Also, one has to replace \eqref{rig-factor-diagram}
with the commutative diagram 
\begin{diagram}
\sqcup_{\ga\in G}\SS_{g-1}^{\rig}(\ga_1,\ldots,\ga_{r},\ga,\ga^{-1}) &\rTo{\rho_{\petlya,\ga}^{\rig}}&
\SS_{g}^{\rig}(\ga_1,\ldots,\ga_r)\\
\dTo{\st_{g-1,\ga}}&&\dTo{\st_g}\\
\ov{\MM}_{g-1,r+2}&\rTo{\rho_{\petlya}}&\ov{\MM}_{g,r}
\end{diagram}
and use the corresponding equation
\begin{equation}\label{fact-st-loop-eq}
\frac{1}{\deg(\st_g)}\rho_{\petlya}^*(\st_g)_*=
\sum_{\ga\in G}\frac{1}{\deg(\st_{g-1,\ga})}(\st_{g-1,\ga})_*
(\rho^{\rig}_{\petlya,\ga})^*.
\end{equation}

To check the factorization axiom for the twisted maps \eqref{lambda-1-eq} 
we note that in the situation of Section \ref{factor-1-sec} we have
$$\Td(R\pi_*(\LL_j))=\Td(R\wt{\pi}_*(\wt{\LL}_j)).$$ 
Indeed, in Case 1 this follows from the isomorphism
\eqref{push-forward-equality}, and in Case 2 --- from the
exact triangle \eqref{L-f-triangle} using the fact that $\wt{L}_j|_p$ is trivial and so has the trivial Todd class.
Now it remains to apply the equalities \eqref{phi-phi'-eq} and \eqref{fact-st-eq} (resp., \eqref{fact-st-loop-eq}),
taking into account Lemma \ref{mod-index-factor-lem}.

\subsection{Forgetting tails}\label{tails-sec}

Here we will check the forgetting tails axiom of \CFT\ (see \cite[4.2]{FJR})
which corresponds to the projection 
$$\th:\ov{\MM}_{g,r}\to\ov{\MM}_{g,r-1}.$$
To do this we have to compare the fundamental \mf s in the following situation.
Let $S=\SS_g^{\rig,0}(\ov{\ga},J)$ be the moduli space of $\Ga$-spin curves with a restricted rigidification structure
$\psi$ (see Section \ref{w-moduli-sec}) of type $(\ov{\ga},J)$, where $\ov{\ga}=(\ga_1,\ldots,\ga_{r-1})$ and $J$ is the exponential
grading element (see \eqref{J-eq}).
Let $(\CC\to S, p_1,\ldots, p_r; P,\vareps)$ be the universal $\Ga$-spin curve.
Let $\CC_r\to S$
be the corresponding family obtained by forgetting the orbifold structure at $p_r$.
Let $\rho_r:\CC\to\CC_r$ and $\rho':\CC_r\to C$ 
be the corresponding projections, so that $\rho'\circ\rho_r=\rho:\CC\to C$
is the morphism of forgetting the orbifold structure at all the marked points 
(see Section \ref{w-moduli-sec}). 

Let $(\LL_1,\ldots,\LL_n)$ be the collection of line bundles on $\CC$ associated
with the $\Ga$-bundle $P$.

\begin{lem} 
On the family $(\CC_r\to S, p_1,\ldots, p_{r-1})$ 
there is
a natural $\Ga$-spin structure $(\ov{P}, \bar{\vareps})$ with a restricted rigidification $\bar{\psi}$
such that the collection line bundles on $\CC_r$ associated with $\ov{P}$
is isomorphic to $(\rho_{r*}\LL_1,\ldots, (\rho_r)_*\LL_n)$. 
Hence, we obtain a morphism 
$$\wt{\th}:S=\SS_g^{\rig,0}(\ov{\ga},J)\to\SS_g^{\rig,0}(\ov{\ga})$$
covering the forgetting tail map $\th$.
\end{lem}

\Pf . By Corollary \ref{quasihom-monom-cor}, a $\Ga$-spin structure
$(P,\vareps)$ can be described as an additional structure for the line bundles $(\LL_1,\ldots,\LL_n)$ given
by isomorphisms 
\begin{equation}\label{w-str-isom-bis-eq}
M_i(\LL_1,\ldots,\LL_n)\simeq\om_{\CC/S}^{\log} \text{ for }\ i=0,\ldots,n-1,
\end{equation}
where $M_0,\ldots,M_{n-1}$ 
are the Laurent monomials \eqref{M-i-def-eq}.
In the notation of Proposition \ref{deg-prop}(i), we have $l_{M_i}(\bq)=1$ for $i=0,\ldots,n-1$.
Therefore, the isomorphisms \eqref{w-str-isom-bis-eq} induce isomorphisms
$$M_i(\rho_{r*}\LL_1,\ldots,\rho_{r*}\LL_n)\simeq \rho_{r*}(\om_{\CC/S}^{\log})(-p_r)
\simeq \om_{\CC_r/S}^{\log},  \text{ for }\ i=0,\ldots,n-1,$$
where the log-structure on $\CC_r$ is given by the markings $p_1,\ldots,p_{r-1}$. 
This is established by the argument similar to the proof of \eqref{M-L-C-main-isom} applied only at the marked
point $p_r$.
Thus, we get a $\Ga$-spin structure
$(\ov{P},\bar{\vareps})$ on $(\CC_r\to S, p_1,\ldots, p_{r-1})$, which inherits the rigidification at the marked
points $p_1,\ldots,p_{r-1}$.
\ed

Since $(\A^n)^{J}=0$ and $\w_J=0$, the fundamental \mf s associated with the data
$(\CC\to S, p_1,\ldots,p_r; P,\vareps, \psi)$ and 
$(\CC_r\to S, p_1,\ldots,p_{r-1}; \ov{P},\bar{\vareps}, \bar{\psi})$ (see Section \ref{constr-sec}) 
belong to the same category,
so we can compare them.

\begin{prop} One has an isomorphism
$$\bP^{\rig,0}_g(\ov{\ga},J)
\simeq (\wt{\th}\times\id)^*\bP^{\rig,0}_g(\ov{\ga}).$$
\end{prop}

\Pf . The natural isomorphisms $\rho'_*(\rho_r)_*\LL_j\simeq\rho_*\LL_j$ are compatible with the maps
\eqref{over-phi-eq} constructed using the $\w$-structures on $\CC$ and on $\CC_r$ and with
trivializations \eqref{e-j-Si-j-eq}. Since the rest of the construction 
in Section \ref{constr-sec} depends only on these data, the assertion follows.
\ed

The above proposition immediately implies that
$$\La_g^R(\ov{\ga})=\La_g^R(\ov{\ga},J)\circ (\id^{\ot r-1}\ot\unit) \text{  and }$$
$$\la_g(\ov{\ga})=\la_g(\ov{\ga},J)\circ (\id^{\ot r-1}\ot\unit)$$
which is the forgetting tails axiom for the \CFTs\ of Theorems \ref{CohFT-thm} 
and \ref{twisted-CohFT-thm}.

\subsection{Concavity}\label{concavity-sec}

Here we consider a special class of families of $\Ga$-spin structures, called {\it concave},
and derive a formula connecting the class of the fundamental \mf\ 
in cohomology of the base with the Chern character of a certain vector bundle.
This is a generalization of the concavity property of \cite[Thm.\ 4.1.5]{FJR}.

Let $S$ be a DM-stack admitting a finite flat covering by a smooth projective scheme,
and let $G$ be a finite group.
For a $G$-equivariant vector bundle $V$ on $S$ (where the action of $G$ on $S$ is trivial)
we have a canonical decomposition
$$V=\bigoplus_{\eta\in \Irr(G)} V_{\eta}\ot \eta$$
compatible with the action of $G$ (where $\Irr(G)$ is the set
of isomorphism classes of irreducible representations of $G$, and 
the bundles $V_{\eta}$ have the trivial $G$-action).
Hence, the abstract Chern character (see Section \ref{dg-Hoch-sec})
$\ch_G^{HH}(V)\in HH_*(D_G(S))$ decomposes as
$$\ch_G^{HH}(V)=\sum_{\eta\in \Irr(G)}\ch^{HH}(V_{\eta})[\eta].$$
Let $R(G)$ be the representation ring of $G$ over $\C$.
Applying the map $\a_{S}\ot\id:HH_*(S)\ot R(G)\to H^*(S,\C)\ot R(G)$
(see \eqref{stack-map})
we obtain an element with values in $H^*(S,\C)\ot R(G)$:
\begin{equation}\label{ch-G-def}
\ch_G(V)=(\a_{S}\ot\id)\ch_G^{HH}(V)=\sum_{\eta\in \Irr(G)}\ch(V_{\eta})[\eta]\in 
H^*(S,\C)\ot R(G).
\end{equation}

\begin{rem} The above notion is different from
the usual $G$-equivariant Chern character of $V$ 
with values in $H^*_G(S,\C)$, which for a finite group
acting trivially on $S$ is equal to the non-equivariant
Chern character, because in this case
$H^*_G(S,\C)=H^*(S,\C)\ot H^*(BG,\C)=H^*(S,\C)$.
\end{rem}

If $G$ is commutative then $R(G)=R=\C[\widehat{G}]$ and for every $\ga\in G$ we have the evaluation homomorphism
$\pi_\ga:R\to\C$. Thus, we can consider the components 
$$\ch_G(V)_\ga:=\pi_\ga(\ch_G(V))=\sum_{\eta\in \widehat{G}}\eta(\ga)\ch(V_{\eta})\in H^*(S,\C).$$
Note that for $\ga=1$ the component $\ch_G(V)_1$ is the usual (non-equivariant) Chern character.

Now assume that we have a family
$(\pi:\CC\to S, p_1,\ldots,p_r)$ of orbicurves of genus $g$ with $r$ marked points and a 
rigidified $\Ga$-spin structure $(P,\vareps)$, induced by a
map
$$f:S\to\SS^{\rig}_{g}(\ov{\ga}),$$
where $\ov{\ga}=(\ga_1,\ldots,\ga_r)\in G^r$.
Such a family is called {\it concave} if $\pi_*(\bigoplus_{j=1}^n \LL_j)=0$ for 
$j=1,\ldots,n$, where $(\LL_1,\ldots,\LL_n)$ are the line bundles on $\CC$ associated with 
the $\Ga$-bundle $P$.
In this case
$$\VV=R^1\pi_*(\bigoplus_{j=1}^n \LL_j)$$
is a vector bundle on $S$, equipped
with a $G$-equivariant structure via the embedding $G\sub\G_m^n$.
Denote by
$\bP$ the pull-back of the fundamental \mf\ \eqref{fund-mf-rig-eq} to $S\times\A^{\ov{\ga}}$.
For each $\ga\in G$ let us denote by 
$$\kappa_\ga:\HH_\ga=HH_*(\MF_\Ga((\A^n)^\ga,\w_\gamma))\to R$$
the map induced by the restriction to the origin functor.
Note that $\kappa_\ga$ is given by the canonical pairing with the Chern character of 
the stabilization of the residue field $\C^{\st}$ (see Example \ref{k-st-ex}).

\begin{prop}\label{concavity-prop} For a concave family of $\Ga$-spin curves over $S$ induced by
the morphism $f:S\to \SS^{\rig}_{g}(\ov{\ga})$, the map
\begin{equation}\label{f-phi-concavity-eq}
f^*\phi_g(\ov{\ga}):\HH_{\ga_1}\ot_R\ldots\ot_R\HH_{\ga_r}\to H^*(S,\C)\ot R
\end{equation}
is given by
\begin{equation}\label{concavity-R-gen-eq}
f^*\phi_g(\ga_1,\ldots,\ga_r)=\ch_{G}\bigl({\bigwedge}^\bullet_\chi \VV^\vee\bigr)\ot
\kappa_{\ga_1}\ot\ldots\ot\kappa_{\ga_r},
\end{equation}
where $\ch_{G}\in H^*(S,\C)\ot R$ is given by \eqref{ch-G-def} and
${\bigwedge}^\bullet_\chi(\VV^\vee):=\bigoplus_i {\bigwedge}^i(\VV^\vee)\ot \chi^{\lfloor i/2 \rfloor}.$
\end{prop}

\Pf . By the base change formula, the \mf\ $\bP$ is obtained by applying the construction of Section \ref{constr-sec}
directly to the family of $\Ga$-spin curves over $S$.
Note that in our case for each $j=1,\ldots,n$, the complex $[A_j\to B_j]$ representing $R\pi_*(\LL_j)$ has the property that $\b_j:A_j\to B_j$ is the embedding of a subbundle. 
Hence, the map $\b=\oplus\b_j$ viewed as a section of $p^*(\bigoplus_j B_j)$ on $X=\tot(\bigoplus_j A_j)$,
is a regular section of the subbundle $p^*(\bigoplus_j A_j)\sub p^*(\bigoplus_j B_j)$.
Let $i:S\to X$ be the zero section.
Assume first that $\A^{\ov{\ga}}\neq 0$ and let $H\sub \A^{\ov{\ga}}$ be the hypersurface
$\w_{\ov{\ga}}=0$. Let  $i':S\to Z^{-1}(H)$ denote the natural embedding.
Since $B_j/A_j\simeq R^1\pi_*(\LL_j)$,
Proposition \ref{koszul-deformed-prop}(i) implies that
$$[\fC(\{-\a_\w,\b\})]=[i'_*{\bigwedge}_\chi^\bullet(\VV^\vee)]$$ 
in the Grothendieck group of $D_{\Sg}(Z^{-1}(H)/\Ga)$. 
Hence, the class of $\fC(\bP)$
in the Grothendieck group of $D_{\Sg}(S\times H/\Ga)$
is given by
$$[\fC(\bP)]= 
[(p,Z)_*\fC(\{-\a_\w,\b\})]=[(\id_S\times k)_*{\bigwedge}_\chi^\bullet\VV^\vee],$$
where $k:\{0\}\to H$ is the embedding.
This implies the result using Lemma \ref{stable-functor-lem} and Remark \ref{dg-fun-rem}.2.
In the case $\A^{\ov{\ga}}=0$, by Proposition \ref{koszul-deformed-prop}(ii), we have
$$[\{-\a_\w,\b\}]=[i_*\maf {\bigwedge}_\chi^\bullet(\VV^\vee)]$$ 
and so by Remark \ref{com-G-mf-rem},
$$[\com_G(\bP)]=[{\bigwedge}_\chi^\bullet(\VV^\vee)].$$ 
\ed

\begin{cor}\label{concavity-cor} 
In the situation of Proposition \ref{concavity-prop}
assume in addition that
$(\A^n)^{\ga_i}=0$ for all $i=1,\ldots,r$. 
In this case $\HH_{\ga_i}=R$ for every $i$, so we can view
the map \eqref{f-phi-concavity-eq} as an element of $H^*(S,\C)\ot R$. 
We have
\begin{equation}\label{concavity-R-eq}
f^*\phi_g(\ga_1,\ldots,\ga_r)=\ch_{G}({\bigwedge}_\chi^\bullet\VV^\vee),
\end{equation}
where $\ch_{G}$ is given by \eqref{ch-G-def}.
The twisted element $\phi^{tw}_g(\ov{\ga})\in H^*(S,\C)$ (see \eqref{phi-twisted-eq}) is equal to the top Chern class of
$\VV^{\vee}$:
\begin{equation}\label{concavity-1-eq}
\phi^{tw}_g(\ov{\ga})=(-1)^{D}c_{D}(\VV)=c_{D}(\VV^{\vee}).
\end{equation}
where $D$ is the rank of $\VV$.
\end{cor}

\Pf .  The formula \eqref{concavity-R-eq} is a direct consequence \eqref{concavity-R-gen-eq}.
The formula \eqref{concavity-1-eq}
follows from the fact that the component $\ch_{G}({\bigwedge}^*\VV^*)_1$
is the usual Chern character and from the standard relation
$$c_{top}(\VV)=\Td(\VV)\cdot\ch({\bigwedge}^\bullet\VV^\vee).$$
\ed

\subsection{Homogeneity conjecture}
\label{dim-sec}

Our conjecture that the reduced \CFT\ of Theorem
\ref{twisted-CohFT-thm} is isomorphic to the one constructed in \cite{FJR} implies a certain homogeneity
property of the maps $\la_g(\ov{\ga})$. This suggests the following analog of the Dimension axiom
of \cite[Thm.\ 4.1.5]{FJR} for the maps 
$$\phi^{tw}_g(\ov{\ga}): e_1\HH_{\ga_1}\otimes\ldots\otimes e_1\HH_{\ga_r}\to 
H^*(\SS^{\rig}_{g}(\ov{\ga})),\C)$$ 
(see \eqref{phi-twisted-eq}).

\medskip

\noindent
{\bf Homogeneity Conjecture}. {\it The image of the map
$\phi^{tw}_g(\ov{\ga})$ is contained in $H^{2\wt{D}_g(\ov{\ga})}(\SS^{\rig}_{g}(\ov{\ga})),\C)$,
where 
$$\wt{D}_g(\ov{\ga})=D_g(\ov{\ga})+\frac{1}{2}\sum_{i=1}^r N_{\ga_i} \text{ and}$$
$$D_g(\ga_1,\ldots,\ga_r)=(g-1)\hat{c}_{\w}+\iota_{\ga_1}+\ldots+\iota_{\ga_r}$$
(see Section \ref{degrees-sec}).}

\medskip

Note that for $\ov{\ga}$ such that $\SS_g(\ov{\ga})$ is nonempty, we have
$2\wt{D}_g(\ov{\ga})\in\Z$ and
$$2\wt{D}_g(\ov{\ga})\equiv N_{\ga_1}+\ldots+N_{\ga_r}\ \mod 2,$$
where $N_\ga=\dim (\A^n)^{\ga}$.
Since the map $\phi^{tw}_g(\ov{\ga})$ is even with respect to the natural 
$\Z/2$-grading,
we know that the above conjecture holds modulo $2$, i.e., the image of
$\phi^{tw}_g(\ov{\ga})$ is contained in 
$$\bigoplus_{m\in\Z} H^{2m+N_{\ga_1}+\ldots+N_{\ga_r}}(\SS^{\rig}_{g}(\ov{\ga}),\C).$$

Now we are going to prove a certain homogeneity property related to the above conjecture.
Let us say that a Koszul \mf\ $\{\a,\b\}$ of $\w$ has {\it rank }$k$ if $\a$ and $\b$ are sections of dual vector bundles of rank $k$.

\begin{prop}\label{homog-prop} Let $\bar{E}_i$, for $i=1,\ldots,r$, be
a Koszul \mf\ of the potential $\w_{\ga_i}$ of rank $k_i$. 
Then
$$\phi^{tw}_g(\ov{\ga})(\ch(\bar{E}_1),\ldots,\ch(\bar{E}_r))\in H^{2(D_g(\ov{\ga})+k_1+\ldots+k_r)}(\SS^{\rig}_{g}(\ov{\ga}),\C).$$ 
\end{prop}

\Pf . Let us write for brevity $\SS=\SS^{\rig}_{g}(\ov{\ga})$.
Recall that the fundamental \mf\ \eqref{fund-mf-rig-eq} is a $\Ga$-equivariant \mf\ of $-\w_{\ov{\ga}}$
on $\SS\times\prod_{i=1}^r (\A^n)^{\ga_i}$ of the form
$$\bP=(p,Z)_*\{\a,\b\},$$
where $p:X=\tot(A)\to \SS$ is the projection from the total space of a vector bundle
$A=\bigoplus_j A_j$, and $\{\a,\b\}$ is the Koszul \mf\ of $0$,  
$$\{\a,\b\}=\bigl({\bigwedge}^\bullet_\chi(p^*B^\vee),\de\bigr),$$ 
where $B=\bigoplus_j B_j$.
Here the complexes $[A_j\to B_j]$ represent the objects $R\pi_*(\LL_j)\in D^b(\SS)$.
Since $\phi^{tw}_g$ is defined using the specialization by $\pi_1:R\to\C$ corresponding to the element $1\in G$,
we can disregard the $G$-equivariant structure on $\bP$. Therefore, we have
\begin{equation}\label{phi-Td-Ch-eq}
\phi^{tw}_g(\ov{\ga})(\ch(\bar{E}_1),\ldots,\ch(\bar{E}_r))=\exp(\pi i\wt{D}_g(\ov{\ga}))\cdot
\frac{1}{\deg(\st_g)}\st_{g*}(\Td(R\pi_*(\bigoplus_{j=1}^n\LL_j))^{-1}\cdot \ch(K)),
\end{equation}
where $K=(p_{\SS})_*(\bP\ot\bar{E}_1\ot\ldots\ot\bar{E}_r)$. Here 
$p_{\SS}: \SS\times\A^{\ov{\ga}}\to\SS$ is the projection. We view $K$ as a $\Z/2$-graded complex of quasicoherent sheaves on $\SS$ with
coherent cohomology. The Chern character of such a complex $K$ can be calculated as
$\ch(K)=\ch(H^{even}K)-\ch(H^{odd}K)$.
Thus, we can replace the expression in the right-hand side of \eqref{phi-Td-Ch-eq} with a similar expression of characteristic classes in
the Chow group $A^*(\SS)\ot\Q$.
It is enough to show that
$$\Td(R\pi_*(\bigoplus_{j=1}^n\LL_j))^{-1}\cdot \ch(K)=\\
\Td(B)\cdot\Td(A)^{-1}\cdot\ch(K)\in A^{D_g(\ov{\ga})+k_1+\ldots+k_r}(\SS)\ot\Q,$$
where $D_g(\ov{\ga})=\rk B-\rk A$ by \eqref{index-eq}. 
Note that 
$$K=(p_{\SS})_*(p,Z)_*(\{\a,\b\}\ot Z^*(\bar{E}_1\ot\ldots\ot\bar{E}_r))\simeq p_*\{\a',\b'\},$$
where $\{\a',\b'\}$ is a Koszul \mf\ of zero of rank $\rk B+k_1+\ldots+k_r$ on $X$ (supported at the zero section $\SS\sub X$). 
By \cite[Lemma 5.3.8]{Chiodo}, we have
$$\Td(A)^{-1}\cdot\ch(K)=\ch_{\SS}^X(\{\a',\b'\})\cdot [p],$$
where $\ch_{\SS}^X(\{\a',\b'\})\in A^*(\SS\to X)$ is
the localized Chern character of the $\Z/2$-graded complex $\{\a',\b'\}$ 
(see \cite[sec. 2.2]{PV}), and $[p]\in A^{-\rk A}(X\to\SS)$ is the orientation class of $p$.
Now by \cite[Thm.\ 3.2]{PV}, the class
$$\Td(B)\cdot \ch_{\SS}^X(\{\a',\b'\})\in A^*(\SS\to X)$$
is concentrated in degree $\rk B+k_1+\ldots+k_r$. Hence, the class
$$\Td(R\pi_*(\bigoplus_{j=1}^n\LL_j))^{-1}\cdot \ch(\bP)=\Td(B)\cdot \ch_{\SS}^X(\{\a',\b'\})\cdot [p]\in A^*(\SS)$$
lives in degree $\rk B-\rk A+k_1+\ldots+k_r$ as claimed.
\ed

\begin{cor}\label{homog-NS-cor} The Homogeneity Conjecture holds in the case when
$(\A^n)^{\ga_i}=0$ for every $i=1,\ldots,r$.
\end{cor}

Proposition \ref{homog-prop} is compatible with the above Homogeneity Conjecture due to the following vanishing property
of Koszul \mf s.

\begin{thm}\label{vanishing-thm}
Let $W\in \C[[x_1,\ldots,x_n]]$ be an isolated singularity and let $\{\a,\b\}$ be a Koszul \mf\ of $W$ of rank 
$r\neq n/2$. Then $\ch(\{\a,\b\})=0$.
\end{thm}

Note that we do not assume here that $W$ is quasi-homogeneous. 
We start by computing
the maps on the Hochschild homology of categories of matrix factorizations induced by homomorphisms of rings.

\begin{lem}\label{mf-homomorphism-lem}
Let 
$\phi:A=\C[[y_1,\ldots,y_m]]\to B=\C[[x_1,\ldots,x_n]]$ 
be a homomorphism of $\C$-algebras 
sending an isolated singularity $q\in A$ to an isolated singularity $W\in B$,
and let 
\begin{equation}\label{Phi-MF-functor}
\Phi:\MF(q)\to \MF(W): \ov{E}\mapsto \ov{E}\ot_A B
\end{equation}
be the functor induced by $\phi$.
Then, under the identifications $HH_*(\MF(q))\simeq H^*(\Om^\bullet_{A/k}, \wedge dq)$ and
$HH_*(\MF(W))\simeq H^*(\Om^\bullet_{B/k}, \wedge dW)$
(see Section \ref{Hoch-sec}), 
the homomorphism induced by $\Phi$ on the Hochschild homology coincides with
 the map
\begin{equation}\label{coh-MF-map}
H^*(\Om^\bullet_{A/k}, \wedge dq)\rTo H^*(\Om^\bullet_{B/k}, \wedge dW)
\end{equation}
induced by $\phi$.
\end{lem}

\Pf .  It is useful to think about the functor $\Phi$ geometrically as the pull-back with respect to
the morphism $f:Y=\Spec(B)\to X=\Spec(A)$ given by $\phi$.
We denote by 
$$f^*:\Om^\bullet_{A/k}\to \Om^\bullet_{B/k}$$ 
the corresponding pull-back map on forms, so that
\eqref{coh-MF-map} is induced by $f^*$. First, we claim that if $\ov{E}=(E,\de)$ is a \mf\ of $q$
then 
\begin{equation}\label{ch-f-E-eq}
\ch(f^*\ov{E})=f^* \ch(\ov{E}) \mod \Om^\bullet\wedge dW
\end{equation}
in $H^*(\Om^\bullet_{B/k})$. Indeed, the formula \cite[(0.2)]{PV-mf} for $\ch(\ov{E})$ can be 
rewritten in a coordinate-free way as
$$\ch(\ov{E})=\str(\exp[\nabla,\de]) \mod \Om^\bullet_{A/k}\wedge dq,$$
where $\nabla:E\to E\ot\Om^1_{A/k}$ is a connection associated with a 
choice of an $A$-basis in $E$
(see \cite{Platt} for a more general formula). 
Since  the right-hand side of the above equation 
is compatible with 
pull-backs, our claim follows.
Applying \eqref{ch-f-E-eq} to the morphism 
$$\id\times f:X\widehat{\times} Y\to X\widehat{\times} X$$
and the diagonal \mf\ of $(-q)\oplus q$ on $X\widehat{\times} X$ (here we take completed products, so that
the corresponding rings are power series in all variables), 
we get that
\begin{equation}\label{ch-f-De-eq}
\ch((\id\times f)^*\De^{\st}_q)=(\id\times f)^* \ch(\De^{\st}_q)
\end{equation}
in the Hochschild homology of $\MF((-q)\oplus W)$.
Note that for $\ov{E}\in\MF(q)$ we have
$$p_X^*\ov{E}\ot (\id\times f)^*\De^{\st}_q\simeq (\id\times f)^*(p_X^*\ov{E}\ot \De^{\st}_q),$$
where $p_X$ is the projection to $X$. Using the base change formula and the fact
that $\De^{\st}_q$ corresponds to the identity functor (see Proposition \ref{diag-mf-lem}), we get that
the functor $\MF(q)\to\MF(W)$ associated with the kernel $(\id\times f) ^*\De^{\st}_q$ is
exactly the pull-back functor $f^*$. Hence, by Lemma \ref{kernel-pairing-lem}, 
the map on the Hochschild homology induced by $f^*$ is given by
$$a\mapsto \tr_{12}(a\ot \ch((\id\times f)^*\De^{\st}_q)),$$
where $\tr_{12}(a\ot b\ot c)=(a, b) c$ 
and $(\cdot,\cdot)$ is the canonical pairing between the Hochschild homology of
 the categories $\MF(q)$ and
$\MF(q)^{op}\simeq \MF(-q)$. Let $(e_i)$ be a basis of 
 $HH_*(\MF(q))$, 
and let $(e^i)$ be the dual basis of 
$HH_*(\MF(-q))$, 
so that $(e_i,e^j)=\de_{ij}$. Then
$\ch(\De^{\st}_q)=\sum_i e^i\ot e_i$ (see \cite[(1.19)]{PV-mf}).
Hence, using \eqref{ch-f-De-eq} we get
$$\tr_{12}(a\ot \ch((\id\times f)^*\De^{\st}_q))=
\tr_{12}(a, (\id\otimes f^*)\Bigl(\sum_i e^i\ot e_i)\Bigr)=
\sum_i (a, e^i) f^*(e_i)= f^*(a)$$
as claimed.
\ed
 
\noindent
{\it Proof of Theorem \ref{vanishing-thm}}.
Let  $\{\a,\b\}=\{a_1,\ldots,a_r;b_1,\ldots,b_r\}$ be a Koszul \mf\  of $W$. 
If 
$a_i$ or $b_i$ is invertible for some $i$,
then $\{\a,\b\}$ is contractible and $\ch(\{\a,\b\})=0$. Otherwise,
we have a homomorphism of $\C$-algebras
$$\phi:A=\C[[u_1,\ldots,u_r;v_1,\ldots,v_r]]\to B=\C[[x_1,\ldots.x_n]],$$ 
given by $\phi(u_i)=a_i$, $\phi(v_i)=b_i$, 
such that
$$W=\phi(q), \ \text{ where } q=u_1v_1+\ldots +u_rv_r\in A.$$
Hence, 
for the  functor \eqref{Phi-MF-functor} corresponding to $\phi$ we have
$$\{a_1,\ldots,a_r;b_1,\ldots,b_r\}\simeq \Phi(\{u_1,\ldots,u_r;v_1,\ldots,v_r\}).$$
Therefore,
$$\ch(\{a_1,\ldots,a_r;b_1,\ldots,b_r\})=\Phi_*(\{u_1,\ldots,u_r;v_1,\ldots,v_r\},$$
where $\Phi_*$ is the map on the Hochschild homology induced by $\Phi$.
By Lemma \ref{mf-homomorphism-lem}, $\Phi_*$ 
coincides with 
the natural map
\eqref{coh-MF-map} induced by $\phi$.
But the source of this map is concentrated in degree $2r$ and the target---in degree $n\neq 2r$, so 
$\Phi_*=0$.
\ed

\begin{cor}\label{simple-homog-cor} 
Assume that for each $\ga\in G$, the space $HH_*(\MF(\w_{\ga}))^G$ is generated by the Chern
characters of Koszul \mf s. 
Then the Homogeneity Conjecture holds for the \CFT\ associated with $\w$ and $G$.

In particular, it holds for all simple singularities.
\end{cor}

\Pf . The first assertion follows from Proposition \ref{homog-prop} and Theorem \ref{vanishing-thm}.
In Section \ref{simple-sing-sec} we will verify that this criterion can be applied to all simple singularities.
\ed

\subsection{Index zero}\label{index-zero-sec}

In the case when $D_g(\ov{\ga})=0$ and $(\A^n)^{\ga_i}=0$ for every $i=1,\ldots,r$, the Homogeneity Conjecture predicts that 
$\phi^{tw}_g(\ov{\ga})$ belongs to $H^0(\ov{\SS}_{g}(\ov{\ga}),\C)$, so it should be a multiple of the fundamental class on each connected component. 
Following \cite{FJR}, we will identify this multiple in terms of the degree of a certain map between affine spaces.

Since we are computing a class in degree $0$ cohomology, it is enough to consider a single
$\Ga$-spin curve $(\CC, p_1,\ldots,p_r; P,\vareps)$ corresponding to a point in the moduli space $\SS=\SS_g(\ov{\ga})$, such that passing to the coarse moduli space $\rho:\CC\to C$ we get
a smooth curve $C$. Assume that 
$D_g(\ov{\ga})=0$ and $(\A^n)^{\ga_i}=0$ for every $i=1,\ldots,r$.
Let $(\LL_1,\ldots,\LL_n)$ be the line bundles on $\CC$ associated with
the $\Ga$-bundle $P$ and let $L_j=\rho_*(\LL_j)$ be
the corresponding line bundles on $C$. 
Then for every monomial $M$ in $\w$ and every $j=1,\ldots,n$ we have a morphism
$$\a_j(M):\pa_jM(H^0(L_\bullet))\to H^0(\pa_jM(L_\bullet))\to H^0(\om_C\ot L_j^{-1})\simeq H^1(L_j)^*,$$
which can be viewed as a section of the vector bundle $H^1(L_j)^*\ot\OO_X$ on the affine space $X=\bigoplus_{j=1}^n H^0(L_j)$ (see Section \ref{canonical-sec}).
As in Section \ref{canonical-sec} we take linear combinations of these sections
$$\a_j^{\w}=\sum_k c_k m_{kj}\a_j(M_k),$$
where $\w=\sum_{k=1}^N c_k M_k$ and $M_k=x_1^{m_{k1}}\ldots x_n^{m_{kn}}$.
Then by Proposition \ref{zero-locus-prop}, 
the section $\a=(\a_1^\w,\ldots,\a_n^\w)$ of $\bigoplus_{j=1}^n H^1(L_j)^*\ot\OO_X$ has zero locus supported only at the origin.
We can view $\a$ as a $\Ga$-equivariant morphism between affine spaces
\begin{equation}\label{Witten-map-eq}
\a:X=\bigoplus_{j=1}^n H^0(L_j)\to Y=\bigoplus_{j=1}^n H^1(L_j)^*\ot \chi.
\end{equation}
This morphism has the property that the subscheme $\a^{-1}(0)$ is concentrated at the origin in $X$.
Note that our assumption $D_g(\ov{\ga})=0$ implies that $X$ and $Y$ are affine spaces of the same dimension.

Let $X(\Ga)$ denote the group of algebraic characters of $\Ga$. 
The natural embedding $\iota:G\to\Ga$ induces a surjective homomorphism of the group rings 
$$\iota^*:\C[X(\Ga)]\to\C[\widehat{G}]=R.$$
Also, consider the homomorphism 
$$\varphi:\G_m\to\Ga:\la\mapsto (\la^{d_1},\ldots,\la^{d_n})$$ 
and the induced homomorphism
$$\varphi^*:\C[X(\Ga)]\to\C[X(\G_m)]=\C[t,t^{-1}],$$
where we choose a generating character $t\in X(\G_m)$ to be $t(\la)=\la^{-1}$.
Note that $\varphi^*(\chi)=t^{-d}$.
Define the $\Z$-grading on $\C[X(\Ga)]$ by
$$\varphi^*(\xi)=t^{\deg(\xi)} \text{ for }\xi\in X(\Ga),$$
so that
$\varphi^*$ becomes a homomorphism of $\Z$-graded algebras.
Let $\C[X(\Ga)]^{\compl}$ be the completion of $\C[X(\Ga)]$ with respect to the degree filtration
$(\C[X(\Ga)])_{\ge q}$. The homomorphism $\varphi^*$ extends to a homomorphism
from $\C[X(\Ga)]^{\compl}$ to the ring of Laurent series in $t$. 

We have the following analog of the Index zero Axiom of \cite[Thm.\ 4.1.5]{FJR}.

\begin{prop}\label{index-zero-prop} 
Assume that $D_g(\ov{\ga})=0$ and $(\A^n)^{\ga_i}=0$ for every $i=1,\ldots,r$. Let 
$x=(\CC,\LL_\bullet)\in\SS$ be a point with smooth $\CC$. Let
$$a_j=h^0(C,L_j),\ \ b_j=h^1(C,L_j) \ \text{ and } \ h=\sum_{j=1}^n a_j=\sum_{j=1}^n b_j=\dim X,$$
where we use the notation above)

\noindent
(i) The restriction of the class $\phi_g(\ov{\ga})\in H^*(\SS,\C)\ot R$ to the point $x\in\SS$ is given by
$$\phi_g(\ov{\ga})|_x=(-1)^h\cdot [H^0(\a^{-1}(0),\OO)]_G\in R,$$
where $[V]_G$ denotes the class of a $G$-module in the representation ring $R$ of $G$.
In particular, the restriction of the specialized 
class $\phi^{tw}_g(\ov{\ga})\in H^*(\SS,\C)$ to $x\in\SS$ is equal to
\begin{equation}\label{index-zero-axiom-eq}
\phi^{tw}_g(\ov{\ga})|_x=(-1)^{h}\cdot\ell(\a^{-1}(0))\in \C,
\end{equation}
where $\ell(\a^{-1}(0))$ is the length of the zero-dimensional scheme $\a^{-1}(0)$.

\noindent
(ii) Let $t_j\in X(\Ga)$ denote the inverse of the character of $\Ga$ induced by the $j$th projection 
$\G_m^n\to\G_m$.
Consider the element
$$P=\prod_{j=1}^n\frac{(1-\chi t_j)^{b_j}}{(1-t_j)^{a_j}}\in\C[X(\Ga)]^{\compl}.$$
Then $P$ belongs to $\C[X(\Ga)]$ and
$$\phi_g(\ov{\ga})|_x=(-1)^h\cdot [H^0(\a^{-1}(0),\OO)]_G=\iota^*P.$$

\noindent
(iii) Consider the subgroup $\lan J\ran=G\cap\varphi(\G_m)\sub\Ga$.
Since $J$ has order $d$ we can identify the representation ring of $\lan J\ran$ with $\C[u]/(u^d-1)$
using the character $u=t|_{\lan J\ran}$. 
We have the specialization homomorphism $\C[t,t^{-1}]\to\C[u]/(u^d-1)$ sending $t$ to $u$.
Consider the Laurent series
$$Q(t)=\prod_{j=1}^n \frac{(1-t^{-d+d_j})^{b_j}}{(1-t^{d_j})^{a_j}}.$$
Then $Q(t)$ is a  Laurent polynomial
and viewing $H^0(\a^{-1}(0),\OO)$ as a representation of $\lan J\ran$
we obtain
$$(-1)^h\cdot [H^0(\a^{-1}(0),\OO)]_{\lan J\ran}=Q|_{t=u}.$$

\noindent
(iv) One has
$$\ell(\a^{-1}(0))=\prod_{j=1}^n\frac{(1-q_j)^{b_j}}{q_j^{a_j}}.$$
\end{prop}

\Pf . (i) Let $\bP$ be the restriction of the fundamental \mf\ over $\SS$ to $x$.
By our construction, $\bP$ is the $\Ga$-equivariant \mf\ of $0$ over a point given by 
\begin{equation}\label{fund-mf-index-zero-eq}
\bP=p_*\{-\a,0\}\simeq p_*\maf(K^\bullet(-\a)),
\end{equation}
where $p:X\to pt$ is the projection (see Lemma \ref{Koszul-complex-lem}). Since
$$H^0(\maf(K^\bullet(-\a)))\simeq H^{even}(K^\bullet(-\a)) \text{ and }
H^1(\maf(K^\bullet(-\a)))\simeq H^{odd}(K^\bullet(-\a))$$
as $G$-equivariant sheaves (see \eqref{coh-maf-eq}), 
the assertion follows from the fact that the only nonzero cohomology of the
Koszul complex $K^\bullet(-\a)$ is the sheaf
$\OO_{\a^{-1}(0)}$ in degree $h(s)=\dim X$.

\noindent
(ii) Let $\CC(\Ga)$ be the category of representations $V$ of $\Ga$ of the form 
$$V=\oplus_{\xi\in X(\Ga), \deg(\xi)\ge -N}V_{\xi}\ot\xi,$$
where all multiplicities $V_{\xi}$ are finite-dimensional. The assignment
$$V\mapsto [V]_{\Ga}=\sum_{\xi\in X(\Ga)}\dim V_{\chi}\cdot \xi$$ 
gives an additive function 
$K_0(\CC(\Ga))\to \CC[X(\Ga)]^{\compl}$, compatible with tensor products.
It is easy to see that 
$$P=[p_*K(-\a)]_{\Ga}.$$
The isomorphism \eqref{fund-mf-index-zero-eq} implies that
the class of $\bP$ in $R$ is equal to $\iota^*P$.
Also, the cohomology of $p_*K(-\a)$ is finite-dimensional, so $P\in\CC[X(\Ga)]$.

\noindent
(iii) This follows from (ii) applying the specialization with respect to the homomorphism $\varphi^*:X(\Ga)\to X(\G_m)$ because $\varphi^*(\chi)=t^{-d}$
and $\varphi^*(t_j)=t^{d_j}$. 

\noindent
(iv) This follows from (iii) by specializing to $t=1$.
\ed

\begin{rem} The {\it Witten map} $\DD$ considered in \cite[Thm.\ 4.1.5]{FJR} (restricted to a point) 
is equal to the complex conjugate of our map $\a$ (see \eqref{Witten-map-eq}).
Thus, the degree of $\DD$ differs from the algebraic degree of $\a$ by the factor $(-1)^{h}$, so our formula \eqref{index-zero-axiom-eq} agrees with
Index Zero Axiom of \cite[Thm.\ 4.1.5]{FJR}.
\end{rem}

\begin{ex} Consider the case when $G=\lan J\ran$ and 
$\w(x_1,\ldots,x_n)$ is homogeneous of degree $d$, so that $d_j=1$ and
$q_j=1/d$. Then the degrees of all the line bundles $L_j$ are the same, so the index zero condition means
that $\deg(L_j)=g-1$ for every $j$. In this case the formula of Proposition \ref{index-zero-prop}(iii) gives
$$\phi_g(\ov{\ga})|_x=\Bigl(\frac{1-t^{-d+1}}{1-t}\Bigr)^{\sum_j a_j}|_{t=u}=
\bigl(-u(1+u+\ldots+u^{d-2})\bigr)^{\sum_j a_j}.$$
Specialization at $u=1$ gives in this case
$$\phi_g^{tw}(\ov{\ga})=(-d+1)^{\sum_j a_j}.$$
\end{ex}

\subsection{Sums of singularities}\label{sum-sing-sec}

Assume that we have a decomposition $\w=\w'(x_1,\ldots,x_{n'})\oplus \w''(y_1,\ldots,y_{n''})$ with
$n'>0$ ad $n''>0$.
Then both polynomials $\w'$ and $\w''$ are also quasi-homogeneous with respect to
the restrictions of the degree vector. Assume also that $G=G'\times G''$, where 
$G'\sub G_{\w'}$ (resp., $G''\sub G_{\w''}$) is a finite subgroup 
containing the exponential grading element.
Let $\Ga'\sub\Ga_{\w'}$ (resp., $\Ga''\sub\Ga_{\w''}$) be the subgroup associated with $G'$ (resp., $G''$).
Then the character $\chi:\Ga\to\G_m$ factors through each of the natural projections $\Ga\to \Ga'$ and $\Ga\to \Ga''$ and we have a cartesian square of commutative algebraic groups
\begin{diagram}
\Ga &\rTo{}&\Ga'\\
\dTo{}&&\dTo{\chi'}\\
\Ga''&\rTo{\chi''}&\G_m
\end{diagram}
Hence, 
for $\ov{\ga'}\in (G')^r$ and $\ov{\ga''}\in (G'')^r$
we have natural isomorphisms of the moduli spaces
\begin{equation}\label{spin-sum-moduli-eq}
\SS^{\rig}_{g,\Ga}(\ov{\ga})\simeq \SS^{\rig}_{g,\Ga'}(\ov{\ga'})\times_{\ov{\MM}_{g,r}}
\SS^{\rig}_{g,\Ga''}(\ov{\ga''}),
\end{equation}
where $\ov{\ga}=(\ov{\ga'},\ov{\ga''})\in G^r$.
Also, we have a natural decomposition of the state space $\HH(\w,G)$ into a tensor product:
\begin{equation}\label{state-space-tensor-prod-eq}
\HH(\w,G)\simeq\HH(\w',G')\otimes_{\C} \HH(\w'',G''),
\end{equation}
compatible with the decomposition $R=R'\ot R''$, where $R'=\C[\widehat{G'}]$, $R''=\C[\widehat{G''}]$.

The following result is an analog of \cite[Thm.\ 4.1.5(8)]{FJR} and \cite[Thm.4.2.2]{FJR}.

\begin{thm}\label{sum-sing-thm} 
(i) For $\ov{\ga}=(\ov{\ga'},\ov{\ga''})\in G^r$ one has an isomorphism
\begin{equation}\label{sum-sing-isom}
\bP^{\rig}_{g,\Ga}(\ov{\ga})\simeq p_{\Ga'}^*
\bP^{\rig}_{g,\Ga'}(\ov{\ga'})\ot p_{\Ga''}^*\bP^{\rig}_{g,\Ga''}(\ov{\ga''})
\end{equation}
in $\ov{\DMF}_{\Ga}(\SS^{\rig}_{g,\Ga}(\ov{\ga})\times\A^{\ov{\ga'}}\times\A^{\ov{\ga''}},
-\w_{\ov{\ga}})$, where
$p_{\Ga'}$ and $p_{\Ga''}$ are the projections to $\SS^{\rig}_{g,\Ga'}(\ov{\ga'})\times\A^{\ov{\ga'}}$
and $\SS^{\rig}_{g,\Ga''}(\ov{\ga''})\times\A^{\ov{\ga''}}$, respectively.

\noindent
(ii) Under the decomposition \eqref{state-space-tensor-prod-eq}, the
map $\La^R_{g,\Ga}(\ov{\ga})$ becomes the tensor product (over $\C$)
of the maps $\La^{R'}_{g,\Ga'}(\ov{\ga'})$ and $\La^{R''}_{g,\Ga''}(\ov{\ga''})$.
The similar result holds for the twisted map $\la_g(\ov{\ga})$ (see \eqref{lambda-1-eq}).
\end{thm}                       

\Pf . (i) Let 
$S=\SS^{\rig}_{g,\Ga}(\ov{\ga})$ be the moduli space of $\Ga$-spin curves of genus $g$ and type
$\ov{\ga}$. 
Recall (see Section \ref{constr-sec}) that
$$\bP=(p,Z)_*\bar{E},$$
where $\bar{E}$ is a Koszul \mf\ on the
total space of the vector bundle 
$p:X\to S$ and $Z:X\to \A^{\ov{\ga}}$ is a linear map. We denote by $X',p',\bar{E}'$, etc.
(resp., $X'',p'',\bar{E}''$, etc.) the similar data constructed for $(\w',\Ga')$ (resp., $(\w'',\Ga'')$)
using the induced $\Ga'$-spin structure (resp., $\Ga''$-spin structure) 
on the universal family of $\Ga$-spin curves over $S$.
Going through the steps of the construction
of these data we see that 
$X=X'\times_S X''$, $Z=(Z',Z'')$ and 
$$\bar{E}\simeq p_{X'}^*(\bar{E}')\ot p_{X''}^*(\bar{E}''),$$
where $p_{X'}:X\to X'$ and $p_{X''}:X\to X''$ are the projections.
Applying Proposition \ref{rel-Kun-lem} to the maps $(p',X'):X'\to S\times\A^{\ov{\ga'}}$ and
$(p'',X''):X''\to S\times\A^{\ov{\ga''}}$
we obtain an isomorphism
$$\bP\simeq p_{12}^*(\bP')\ot p_{13}^*(\bP''),$$
which leads to \eqref{sum-sing-isom}, since $\bP'$ and $\bP''$ are the pull-backs of the fundamental
\mf s associated with the universal families over 
$\SS^{\rig}_{g,\Ga'}(\ov{\ga'})$ and $\SS^{\rig}_{g,\Ga''}(\ov{\ga''})$.

\noindent
(ii) Denote by 
$S'=\SS^{\rig}_{g,\Ga'}(\ov{\ga'})$ (resp., $S''=\SS^{\rig}_{g,\Ga''}(\ov{\ga''})$ the corresponding
moduli spaces, and let $p_{S'}:S\to S'$ and $p_{S''}:S\to S''$ be the projections.
Using (i) we obtain that
for $\bar{A}'\in \MF_{\Ga'}(\A^{\ov{\ga'}})$ and $\bar{A}''\in \MF_{\Ga''}(\A^{\ov{\ga''}})$ one has
$$\Phi_g(\ov{\ga})(\bar{A}'\ot\bar{A}'')\simeq p_{S'}^*(\Phi'_g(\ov{\ga'})(\bar{A}'))\ot
p_{S''}^*(\Phi''_g(\ov{\ga''})(\bar{A}'')),$$
where $\Phi'_g(\ov{\ga})$ (resp., $\Phi''_g(\ov{\ga''})$) is the functor given by the kernel 
$\bP^{\rig}_{g,\Ga'}(\ov{\ga'})$ (resp., $\bP^{\rig}_{g,\Ga''}(\ov{\ga''})$). 
Hence, the induced map 
$$\phi_g(\ov{\ga}):\HH_{\ga_1}\ot_R\ot\ldots\ot_R\HH_{\ga_r}\to H^*(S,\C)\ot R$$
(see \eqref{hoch-phi-map}) is given by 
$$(h'_1\ot h''_1)\ot\ldots\ot(h'_r\ot h''_r)\mapsto p_{S'}^*(\phi'_g(h'_1\ot\ldots\ot h'_r))\cdot
p_{S''}^*(\phi''_g(h''_1\ot\ldots\ot h''_r)),$$
where $\phi'_g$ and $\phi''_g$ are the corresponding maps for $(\w',\Ga')$ and $(\w'',\Ga'')$.
It remains to take the push-forward with respect to
the projection $\st_g:S\to\ov{\MM}_{g,r}$ and use the isomorphism \eqref{spin-sum-moduli-eq}.
In the case of the twisted map $\la_g(\ov{\ga})$
 one has to use in addition the multiplicativity of the Todd class and the additivity
of the numbers $\wt{D}_g(\ov{\ga})$ under the product of symmetry groups.
\ed

\section{Calculations for genus zero and three points}\label{calc-sec}

In this section we will compute the fundamental \mf\ and the maps of the \CFT\ of Theorem 
\ref{CohFT-thm}
corresponding to $\P^1$ with three marked points under some technical assumptions on the 
type $(\ga_1,\ga_2,\ga_3)$. In particular,
we will verify the metric axiom of the \CFT.
The results of Section \ref{3-pt-corr-sec} will be used in Section \ref{simple-sing-sec} to determine
the Frobenius algebras associated with our \CFT\ for
all simple singularities.

\subsection{Metric axiom}\label{genus-0-3-sec}

In this section  
we will check that the components of the metric on the state space of 
the \CFT\ of Theorem \ref{CohFT-thm} are equal
to the maps $\La_{0}^R(\ga,\ga^{-1},J)$, as required by the metric axiom of \CFT.

Since $\ov{\MM}_{0,3}$ is a point, the moduli space $\SS_{0,3}$ a finite stack, where a point of
$\SS_{0,3}$ corresponds to a $\Ga$-spin curve $(\CC,p_1,p_2,p_3; P,\vareps)$
with the projection $\rho:\CC\to C=\P^1$ 
obtained by forgetting the orbi-structure at $p_1$, $p_2$ and $p_3$. 
By Proposition \ref{deg-prop}, the type 
$(\ga_1,\ga_2,\ga_3)\in G^3$ of such a $\Ga$-spin curve should satisfy $\ga_1\ga_2\ga_3=J$,
and there exists a $\Ga$-spin curve of every such type, unique up to an isomorphism.

First, let us consider the situation where $\ga_3=J$ and therefore $\ga_1\cdot\ga_2=1$.
Thus, we will denote $\ga_1=\ga$ and $\ga_2=\ga^{-1}$.
As we have seen in Section \ref{tails-sec}, in this case we can drop the point $p_3$ when 
calculating the fundamental \mf. 

Let $L_j=\rho_*(\LL_j)$ be the corresponding line bundles on $\P^1$ (where $(\LL_1,\ldots,\LL_n)$
are the line bundles on $\CC$ associated with $P$), and
let $S\sub\{1,\ldots,n\}$ be the set of $j$ such that the $j$th component of $\ga$ is trivial.
Formula \eqref{deg-eq} for $g=0$ and $r=2$ implies that
$$\deg L_j=
\begin{cases} 0, & j\in S,\\ -1, & j\not\in S.\end{cases}$$
Note that in the notation of Section \ref{constr-sec}
$$\Si_j=\begin{cases} \{p_1,p_2\}, &j\in S \\ \emptyset, &j\not\in S.\end{cases}$$
Hence, $\Si_M=\{p_1,p_2\}$ when $j\in S$ for all $j$ with $\deg_{x_j}M>0$, and $\Si_M$ is empty otherwise. 
Thus, for each $M$ we have one of the two cases depending on whether $l_M(\bdeg)$ is zero
or not: 

\noindent
(i) $\Si_M=\{p_1,p_2\}$ and $l_M(w_1)=l_M(w_2)=0$;

\noindent
(ii) $\Si_M=\emptyset$ and both $l_M(w_1)$ and $l_M(w_2)$ are positive.

Permuting the indices $\{1,\ldots,n\}$ corresponding to the variables $x_1,\ldots,x_n$ we can assume that 
$S=\{1,\ldots,k\}$, so that
$L_j=\OO_{\P^1}$ for $j=1,\ldots,k$ and $L_j=\OO_{\P^1}(-1)$ for $j=k+1,\ldots,n$.

The moduli space $\SS_{0}^{\rig}(\ga,\ga^{-1},J)$ is a collection of points $\psi$ corresponding to different choices of rigidification. We are going to compute the restriction $\bP(\psi)$ of 
the fundamental \mf\ $\bP^{\rig}_{0,\Ga}(\ga,\ga^{-1},J)$ to $\{\psi\}\times\A^{\ov{\ga}}$.
We claim that there exists a rigidification $\psi_0$ of $(P,\vareps)$ such that for $j=1,\ldots,k$
the induced trivialization of $L_j|_{p_1}$ 
 comes from a global trivialization of $L_j=\OO$ and the induced trivialization of
$L_j|_{p_2}$ differs  from the restriction of a global trivialization by the factor
$\zeta_j=\exp(\pi i q_j)$. Indeed, by Corollary \ref{quasihom-monom-cor},
a choice of a restricted rigidification with this property is equivalent to choosing trivializations of the line
bundles $L_1,\ldots,L_k$ at $p_1$ and $p_2$, such that the induced trivializations
of $M_i(L_1,\ldots,L_k)\simeq\om(p_1+p_2)$ at these points coincide with the canonical trivialization,
where $M_0,\ldots,M_{k-1}$ are Laurent monomials in the subset of variables $x_1,\ldots,x_k$, 
chosen as in \eqref{M-i-def-eq}. If we choose arbitrary global trivializations 
$e_j:\OO\to L_j$ for $j=1,\ldots,k$, then we will get a collection of nonzero residues
$$r_i=\Res_{p_1}(M_i(e_1,\ldots,e_k)), \ i=0,\ldots,k-1.$$
Since the map $(M_0,\ldots,M_{k-1}):(\C^*)^k\to (\C^*)^k$ is surjective, we can rescale $e_j$'s,
so that $r_i=1$ for each $i=0,\ldots,k-1$. Since the residues of global sections of $\om(p_1+p_2)$ 
at $p_1$ and $p_2$ are opposite, the trivializations $e_j|_{p_1}\in L_j|_{p_1}$ and
$\zeta_j e_j|_{p_2}\in L_j|_{p_2}$ define a restricted rigidification of our $\Ga$-spin structure. 
Using surjectivity of the map \eqref{sur-rig-mor-eq} we extend it to a rigidification $\psi_0$.
We are going to show that $\bP(\psi_0)$ is essentially the diagonal matrix factorization on
$(\A^n)^{\ga}\times(\A^n)^{\ga^{-1}}=\A^k\times\A^k$.  

To apply the construction of Section \ref{constr-sec} we need to choose resolutions
$[A_j\to B_j]$ for $R\Ga(L_j)$ such that the restriction maps $Z_j:R\Ga(L_j)\to L_j|_{\Si_j}$ are
realized by surjective maps $A_j\to L_j|_{\Si_j}$. When $j=1,\ldots,k$, we have $L_j\simeq \OO$ and
we take the resolution
\begin{equation}\label{O-resolutions-eq}
R\Ga(L_j)\to [L_j|_{p_1}\oplus L_j|_{p_2}\rTo{\de}L|_{p_1}],
\end{equation}
where $\de$ is the difference map that uses the natural identification $L|_{p_1}\simeq L|_{p_2}$.
Using the rigidification $\psi_0$ we can identify $L_j|_{p_1}\oplus L_j|_{p_2}$ with $\C^2$.
Then the above resolution becomes $[\C^2\rTo{\b_j}\C]$,
where
$\b_j(x,y)=y-\zeta_j x$ by the choice of $\psi_0$. Note that that the restriction map
$Z_j: R\Ga(L_j)\to L_j|_{p_1+p_2}\simeq\C^2$ is realized by the identity map $\C^2\to \C^2$.

When $j=k+1,\ldots,n$, there exists an isomorphism $L_j=\OO(-1)$.
Since for such $j$ we have $\Si_j=\emptyset$, we can simply set $A_j=B_j=0$ in this case. 

With these choices of resolutions 
the space $X$ can be identified with $\A^k\times\A^k$, so that the map
$Z$ becomes the identity.
The bundle $p^*(\bigoplus B_j)$ on $\A^k\times\A^k$ is the trivial bundle of rank $k$
with the basis $e_1,\ldots,e_k$, and the differential $\b:\bigoplus A_j\to\bigoplus B_j$ corresponds
to the section 
$$\b=\sum_{j=1}^k(y_j-\zeta_jx_j)e_j\in H^0(X,p^*(\bigoplus B_j)).$$
The $\w$-structure also induces a section $-\a_\w\in H^0(X,p^*(\bigoplus B_j^\vee))$ such that
$\{-\a_w,\b\}$ is a \mf\ of $-\w_{\ga}(x)-\w_{\ga^{-1}}(y)$ on $\A^k\times\A^k$ (see 
Section \ref{constr-sec}).
Since $\b$ is regular, the Koszul \mf\ $\{-\a_w,\b\}$ 
is a stabilization of the structure sheaf of the shifted diagonal $y=\zeta_\bullet x$ in $\A^k\times\A^k$
(see Theorem \ref{factor-thm}). Therefore,
$$\bP(\psi_0)=\{-\a_w,\b\}\simeq \De^{\st}_{-\w_{\ga},\zeta},$$
where $\De^{\st}_{-\w_{\ga},\zeta}\simeq (\zeta_\bullet,\id)^*\De^{\st}_{-\w_\ga}$.

Any other rigidification $\psi\in \SS_{0}^{\rig}(\ga,\ga^{-1},J)$ 
is obtained from $\psi_0$ by the action of an element $(g_1,g_2)\in G\times G$.
Hence,
$$\bP(\psi)\simeq (g_1\times g_2)^*\De^{\st}_{-\w_{\ga},\zeta}\simeq (\id\times g)^*
\De^{\st}_{-\w_{\ga},\zeta},$$
where $g=g_2g_1^{-1}\in G$.
Now we are ready to calculate the maps $\La^R_0(\ga,\ga^{-1},J)(\cdot,\cdot,\unit)$.

\begin{lem}\label{metric-lem}
For $h\in\HH(\w_{\ga})$ and $h'\in\HH(\w_{\ga^{-1}})$ one has 
\begin{equation}\label{genus-0-formula}
\La^R_0(\ga,\ga^{-1},J)(h,h',\unit)=
((\zeta_\bullet)_*h,h')^R_{\w_{\ga}},
\end{equation}
where 
$$(\cdot,\cdot)^R_{\w_\ga}:\HH(-\w_\ga)\ot\HH(\w_\ga)\to R$$ 
is the pairing \eqref{R-pairing-formula}, and
$$(\zeta_\bullet)_*:\HH(\w_{\ga^{-1}})=\HH(\w_{\ga})\to \HH(-\w_\ga)=\HH(-\w_{\ga^{-1}})$$
is the isomorphism induced by 
the automorphism $x_j\mapsto \zeta_jx_j$ of $\C[x_1,\ldots,x_n]$.

The Casimir element corresponding to this metric on the state space
$\HH=\bigoplus\HH(\w_{\ga})$ has
components $T_{\w_{\ga},\zeta}\in \HH(\w_{\ga})\ot\HH(\w_{\ga^{-1}})$ given by
\eqref{T-w-zeta-eq}.
\end{lem}

\Pf . Recall that $\La^R_0(\ga,\ga^{-1},J)$ is defined using the functor 
$$\Phi_0(\ga,\ga^{-1},J):\DMF_{\Ga}(\A^k\times \A^k,\w_{\ga}\oplus\w_{\ga^{-1}})\to
D_G(\SS_{0}^{\rig}(\ga,\ga^{-1},J))$$
associated with the kernel $\bP^{\rig}_{0,\Ga}(\ga,\ga^{-1},J)$ (see \eqref{Phi-functor-eq}).
Let us consider the component of this functor
$$\Phi(\psi): \MF_{\Ga}(\A^k\times \A^k,\w_{\ga}\oplus\w_{\ga^{-1}})\to \Com_f(G-\mod)$$
corresponding to the point $\psi\in \SS_{0}^{\rig}(\ga,\ga^{-1},J)$.
As we saw above, $\Phi(\psi)$ is given by tensoring with
$$\bP(\psi)\simeq (\id\times g)^*\De^{\st}_{-\w_{\ga},\zeta}$$
 for some
$g\in G$. Hence, by Proposition \ref{diag-mf-lem}(ii), we have
$$\Phi(\psi)=\pi_*\circ\De^*\circ(\zeta_\bullet\times g^{-1})^*,$$
where $\pi:\A^k\to pt$ is the projection.
Recall that the canonical pairing $(\cdot,\cdot)^R_{\w_{\ga}}$ is induced by the functor $\pi_*\circ \De^*$,
whereas the functor $(g^{-1})^*$ induces the
identity map on the Hochschild homology $HH_*(\MF_{\Ga}(\w_{\ga}))$.
Formula \eqref{genus-0-formula} follows from this and from the
fact that under the isomorphisms \eqref{graded-G-hoch-isom}
the map
$$(\zeta_\bullet)_*: HH_*(\MF_{\Ga}(\w_{\ga}))\to HH_*(\MF_{\Ga}(-\w_{\ga}))$$
decomposes into the direct sum of morphisms
$$(\zeta_\bullet)_*:H(\w_{\ga,\ga'})\to H(-\w_{\ga,\ga'})=H(\w_{\ga,\ga'})$$
induced by the automorphism of $\C[x]$ sending $x_i$ to $\zeta_ix_i$,
where $\w_{\ga,\ga'}=\left.\w\right|_{(\A^n)^{\{\ga,\ga'\}}}~$.

The second assertion follows from Proposition \ref{Casimir-prop}.
\ed

\subsection{Three-point correlators}\label{3-pt-corr-sec}

Now we are going to consider the maps $\La^R_0(\ov{\ga})$ corresponding to 
arbitrary markings $\ov{\ga}=(\ga_1,\ga_2,\ga_3)$. 
As in the beginning of the section, we work with a $\w$-curve 
 $\CC$ with three marked points, where we no longer assume that $\ga_3=J$.
Note that in this case \eqref{deg-eq} takes form
\begin{equation}\label{deg-0-3-eq}
\bdeg=\bq-\bth_1-\bth_2-\bth_3.
\end{equation}
Recall that for every $j=1,\ldots,n$, the subset $\Sigma_j\sub\{p_1,p_2,p_3\}$ 
contains the point $p_i$ if and only if $(\th_i)_j=0$. Using the fact that the coordinates of
$\bth_i$ belong the interval $[0,1)$ we obtain the following list of possible cases.

\begin{lem} For every $j=1,\ldots,n$, one of the following possibilities is realized:

\noindent
(i) $|\Sigma_j|=2$ and $L_j=\OO$;

\noindent
(ii) $|\Sigma_j|=1$ and $L_j$ is either $\OO$ or $\OO(-1)$;

\noindent
(iii) $\Sigma_j=\emptyset$ and $L_j$ is either $\OO$, or $\OO(-1)$, or $\OO(-2)$.
\end{lem}

\ed

As before, to compute $\La^R_0(\ov{\ga})$ 
we need to choose resolutions $[A_j\to B_j]$ for $R\Ga(L_j)$.

In case (i) we will use the resolution \eqref{O-resolutions-eq}.
In case (ii) when $L_j=\OO$ we can use the resolution
$$A_j=H^0(\OO)\to B_j=0$$ 
with $Z_j: A_j\to\C$ equal to the identity map.
In case (ii) when $L_j=\OO(-1)$ we will use the resolution
$$A_j=H^0(\OO)\rTo{\ev_{\infty}} B_j=\C,$$
where $\ev_{\infty}$ is the evaluation at the point $\infty\in \P^1$ (that we assume to be distinct from
$p_1$, $p_2$, $p_3$) with $Z_j$ still equal to the identity.
Finally, in case (iii) we will take the resolution with $A_j=H^0(L_j)$, $B_j=H^1(L_j)$ and zero differential.
 
 Assume that for any $j$ such that $L_j=\OO$ one has $|\Sigma_j|\ge 1$.
Then the corresponding space $X=\bigoplus_j A_j$ gets identified with 
$\A^{\ov{\ga}}=\prod_j \A^{\Sigma_j}$.
We will use the coordinates $x_j(i)$ on this affine space indexed by $(i,j)$ such that $p_i\in \Si_j$.
For $i=0,1,2$ let us set
$$S_i=\{j\ |\ |\Sigma_j|=i \text{ and } L_j\simeq \OO(i-2)\}\subset\{1,\ldots,n\}$$
The bundle $\bigoplus_j B_j$ has generators $e_j$ labeled by
$j\in S_0\sqcup S_1\sqcup S_2$.
For $j\in S_2$ the coefficient of $e_j$ in $\b$
has the form $a_jx_j(i_1)+b_jx_j(i_2)$,
where $\Sigma_j=\{p_{i_1},p_{i_2}\}$, $i_1<i_2$, and $a_j$ and $b_j$ are nonzero constants
(depending on a rigidification $\psi$).
For $j\in S_1$ the coefficient of $e_j$ in $\b$ 
is $x_j(i)$, where $\Sigma_j=\{p_i\}$.
Finally, for $j\in S_0$ the coefficient of $e_j$ in $\b$ is zero.

\begin{prop}\label{3-pt-prop}
(i) Assume that one has $|\Sigma_j|\ge 1$ for each $j$ such that $L_j=\OO$.
Let $\bP(\psi)$ be  
the restriction of the fundamental \mf\ to a point $\psi\in\SS^{\rig}_0(\ov{\ga})$.
Then $\bP(\psi)$ is isomorphic to a Koszul \mf\ 
$$\bP(\psi)=\{\a,\b\}$$
of the potential $\w_{\ga_1}(x_\bullet(1))+\w_{\ga_2}(x_\bullet(2))+\w_{\ga_3}(x_\bullet(3))$,
where $\b$ is a section of the trivial bundle with generators $e_j$ numbered by $j\in S_0\sqcup S_1\sqcup S_2$, 
of the form
$$\b=\sum_{j\in S_1}x_j(i)e_j+\sum_{j\in S_2}
(a_jx_j(i_1)+b_jx_j(i_2))e_j$$
for some $a_j,b_j\in \C^*$,
where in the first (resp., second) sum we assume that $\Si_j=\{p_i\}$ (resp., 
$\Si_j=\{p_{i_1},p_{i_2}\}$). 

Furthermore, the map 
$$\La^R_0(\ov{\ga}):\HH_{\ga_1}\ot_R\HH_{\ga_2}\ot_R\HH_{\ga_3}\to R$$
is equal to $i_{\psi}^*\phi_0(\ov{\ga})$ (see \eqref{hoch-phi-map}) for any choice of a rigidification $\psi$,
where
$i_{\psi}:\{\psi\}\hra \SS^{\rig}_0(\ov{\ga})$ is the embedding.

\noindent
(ii) Now assume that for any $j$ such that $L_j=\OO$ one has $|\Sigma_j|=2$.
Consider the composition
$$\phi(S_1,S_2):\HH_{\ga_1}\ot_R\HH_{\ga_2}\ot_R\ot\HH_{\ga_3}\rTo{\iota(\ov{\ga})} 
HH_*(\MF_{\Ga}(\A^{\ov{\ga}},\w_{\ov{\ga}}))\to R,$$
where $\iota(\ov{\ga})$ is the map \eqref{iota-isom} and the second arrow is 
induced by the functor
$$\MF_{\Ga}(\A^{\ov{\ga}},\w_{\ov{\ga}})\rTo{\com_G\circ \pi_*i^*}\Com_f(G-\mod),$$
where $i:\A^{S_2}\to\A^{\ov{\ga}}$ is the embedding of the subspace defined by the linear equations
\begin{align*}
&a_jx_j(i_1)+b_jx_j(i_2)=0 & \text{ for } j\in S_2, \Sigma_j=\{p_{i_1},p_{i_2}\},\\
&x_j(i)=0 & \text{ for } j\in S_1, \Sigma_j=\{p_i\}
\end{align*}
and $\pi:\A^{S_2}\to pt$ is the projection.
Then one has 
\begin{equation}\label{3-pt-formula}
\La^R_0(\ov{\ga})=\prod_{j\in S_0}(1-t_j)\cdot\phi(S_1,S_2),
\end{equation}
where $(t_1,\ldots,t_n)$ are the inverses of the characters of $G$ corresponding to the coordinates
of the natural map $G\to\G_m^n$.
\end{prop}

\Pf . (i) The first assertion follows immediately from the discussion preceding the proposition.
For the second, we observe that the group $G^3$ acts transitively on $\SS^{\rig}_0(\ov{\ga})$
and for $\bar{g}\in G^3$ one has
$$\bP(\bar{g}\psi)=\bar{g}^*\bP(\psi).$$
This implies the statement since the functors $g^*$ on 
$\MF_{\Ga}((\A^n)^{\ga_i}, \w_{\ga_i})$ for $g\in G$ induce the identity
maps on the Hochschild homology and $\La^R_0(\ov{\ga})$ is equal the average of the maps
$i_{\psi}^*\phi_0(\ov{\ga})$ over $\psi\in\SS^{\rig}_0(\ov{\ga})$.

\noindent
(ii) Under our assumptions we have an identification 
$$\A^{S_1}\times \A^{S_2}\times\A^{S_2}\simeq \A^{\ov{\ga}}$$
so that the section 
$$\b\in\bigoplus_j B_j=\bigoplus_{j\in S_0\cup S_1\cup S_2}\OO\cdot e_j$$ 
is of the form
$$\sum_{j\in S_1} x_je_j +\sum_{j\in S_2} (a_jx_j+b_jy_j)e_j,$$
where $(x_j)_{j\in S_1}$ are coordinates on $\A^{S_1}$,
$(x_j,y_j)_{j\in S_2}$ are coordinates on $\A^{S_2}\times \A^{S_2}$, and $e_j$ is the generator of the one-dimensional
representation $\eta_j$ of $G$.
By the second assertion in part (i), it is enough to make calculations for one rigidification 
$\psi$.
By part (i), we have
$$\bP(\psi)\simeq \{\a,\b\}=\{\a_0,0\}\ot\{\a_{12},\b_{12}\}\simeq K^\bullet(\a_0)\ot\{a_{12},\b_{12}\},$$
where we 
use the decomposition 
$$\bigoplus_j B_j=(\bigoplus_{j\in S_0} \OO_X\cdot e_j)\oplus(\bigoplus_{j\in S_1\sqcup S_2}
\OO_X\cdot e_j)$$
and write in components $\a=(\a_0,\a_{12})$, $\b=(0,\b_{12})$. 
Applying Lemma \ref{stable-functor-lem} we obtain an isomorphism
$$\Phi_{\bP(\psi)}(\bar{E}_1,\bar{E}_2,\bar{E}_3)\simeq \pi_*(i^*(\bar{E}_1\ot \bar{E}_2\ot 
\bar{E}_3)\ot K^\bullet(\a_0))$$
for the functor associated with $\bP(\psi)$.
Since $i^*(\bar{E}_1\ot \bar{E}_2\ot \bar{E}_3)$ is supported at the origin, on the level of Hochschild homology we can replace the complex $K^\bullet(\a)$ by the alternating sum of its terms.
This gives rise to the factor $\prod_{j\in S_0}(1-t_j)$ in the formula \eqref{3-pt-formula}.
\ed

Note that we have a natural identification
$$\A^{S_2}=(\A^n)^{\ga_1,\ga_2}\times (\A^n)^{\ga_1,\ga_3}\times(\A^n)^{\ga_2,\ga_3},$$
so that the embedding $i$ is the product of three diagonal embeddings complemented
by zero in the remaining coordinates.

\begin{rem} 
The assumptions of Proposition \ref{3-pt-prop}(ii) are automatically satisfied if $d_j=\deg(x_j)=1$
for $j=1,\ldots,n$ and $\ga_i^d=1$ for $i=1,2,3$, where $d=\deg(\w)$.
\end{rem}

\subsection{Case of homogeneous polynomials and the scalar group action}\label{hom-pol-sec}

Let us consider the case when $\w(x_1,\ldots,x_n)$ is homogeneous (i.e., $\deg(x_j)=1$) and 
$G=\Z/d$, where $d=\deg(\w)$, such that $m\in\Z/d$ acts by the scalar
multiplication with $\exp(2\pi i m/d)$. The group $\Ga$ in this case is
$\G_m$ acting on $\A^n$ via scalar multiplications, the character $\chi:\G_m\to\G_m$
sends $\la$ to $\la^d$. The element $\zeta_\bullet\in(\C^*)^n$ has all the components equal to
$\exp(\pi i/d)$.
Our \CFT\  has the state space 
$$\HH=\bigoplus_{m\in\Z/d}\HH_m,$$
where $\HH_m=R$ for $m\not\equiv 0$ and
$$\HH_0=HH_*(\MF_{\G_m}(\A^n,\w)).$$
Let us set for $m\in\Z/d$, $m\not\equiv 0$
$$e(m):=1\in\HH_m.$$
Recall that $\HH_0$ is an $R$-module equipped with the canonical $R$-valued metric 
$(\cdot,\cdot)^R$.
We also have a special element
$$e(0):=\ch(\C^{\st})\in\HH_0.$$
Using Proposition \ref{3-pt-prop}(ii) we can calculate all the maps
$\La^R_0(m_1,m_2,m_3)$, or equivalently the $R$-algebra structure on $\HH$ in terms of these data.
Recall that the $R$-valued metric on $\HH$ restricts to the canonical metric on $\HH_0\times\HH_0$ 
(twisted by $(\zeta_\bullet)_*$ in the first factor)
and also restricts to the natural pairing between $\HH_m$ and $\HH_{-m}$, so that
$$(e(m),e(-m))^R=1 \text{ for } m\not\equiv 0.$$
The flat unit element in this case is $e(1)\in\HH_1$.
Let us pick as a generator $t$ of the group $\widehat{G}$ the character
\begin{equation}\label{t-eq}
t(m):=\exp(-2\pi i m/d)
\end{equation}
and use the corresponding identification $R=\C[t]/(t^d-1)$. 
Below we will often omit $(m_1,m_2,m_3)$ from the notation $\La^R_0(m_1,m_2,m_3)(x,y,z)$,
where $x\in\HH_{m_1}$, $y\in\HH_{m_2}$ and $z\in\HH_{m_3}$.

\begin{thm}\label{homog-thm}
(a) All the nonzero maps $\La^R_0(\cdot,\cdot,\cdot)$ are given by
\begin{equation}\label{La-e-e-e-eq}
\La^R_0(e(m_1),e(m_2),e(m_3))=\begin{cases} 1, & m_1+m_2+m_3=d+1, \\
(1-t)^n, & m_1+m_2+m_3=2d+1,
\end{cases}
\end{equation}
\begin{equation}\label{La-x-e-e-eq}
\La^R_0(x,e(m),e(d+1-m))=(x,e(0))^R,
\end{equation}
\begin{equation}\label{La-x-y-e-eq}
\La^R_0(x,y,e(1))=((\zeta_\bullet)_*x,y)^R,
\end{equation}
where $m_1,m_2,m_3\in [1,d-1]$, $m\in [2,d-1]$ and $x,y\in\HH_0$.

\noindent
(b) The maps $1\mapsto e(1), u\mapsto e(2)$ induce an isomorphism of $R$-algebras
$$R[\HH_0,u]/I\simeq \HH,$$
where $I$ is the ideal generated by the relations
$$u^{d-1}=e(0), \ \ xu=(x,e(0))^R, \ \ xy=((\zeta_\bullet)_*x,y)^Ru^{d-2}, {\text where }\ x,y\in\HH_0.$$
\end{thm}

\Pf . (a) Recall that the map $\La^R_0(m_1,m_2,m_3)$
is zero unless
$$m_1+m_2+m_3\equiv 1\mod(d).$$

Assume first that $m_1=m_2=0$. Then $m_3=1$, and the corresponding 
equality \eqref{La-x-y-e-eq} follows from Lemma \ref{metric-lem}.

Next, assume that $m_1=0$ while $m_2,m_3\in [1,d-1]$ are such that $m_2+m_3=d+1$.
Then \eqref{deg-eq} gives $L_j=\OO(-1)$ for all $j$. Hence,
Proposition \ref{3-pt-prop}(ii) applies with $S_0=S_2=\emptyset$ and $S_1=\{1,\ldots,n\}$.
Thus, $\La^R_0(0,m_2,m_3)(\cdot,1,1):\HH_0\to R$ is the map induced on Hochschild
homology by the functor of restriction to $0\in\A^n$. Thus, \eqref{La-x-e-e-eq} follows from
Example \ref{k-st-ex}.

Finally, in the case when $m_1,m_2,m_3\in [1,d-1]$ we have
either $L_j=\OO(-1)$ for all $j$ (when $m_1+m_2+m_3=d+1$) or
$L_j=\OO(-2)$ for all $j$ (when $m_1+m_2+m_3=2d+1$). In the former case
$S_0=S_1=S_2=\emptyset$, while in the latter case $S_1=S_2=\emptyset$ and
$S_0=\{1,\ldots,n\}$. Now \eqref{La-e-e-e-eq} follows from Proposition \ref{3-pt-prop}(ii).

\noindent
(b) From (a) and the definition of the metric on $\HH$
we obtain the following multiplication rules in $\HH$. 
First, we get that $e(1)$ is a unit.
For $i,j\in[1,d-1]$ we have
$$e(i)\cdot e(j)=\begin{cases} e(i+j-1), & i+j\le d+1,\\ (1-t)^ne(i+j-1), & i+j>d+1\end{cases}$$
(for $i+j=d+1$ we use the fact that $(\zeta_\bullet)^*e(0)=e(0)$).
Also, for $i\in [2,d-1]$ and $x\in \HH_0$ one has
$$x\cdot e(i)=(x,e(0))^R\cdot e(i-1).$$
Finally, for $x,y\in\HH_0$ one has
$$x\cdot y=((\zeta_\bullet)_*x,y)^R\cdot e(d-1).$$
This is equivalent to our assertion.
\ed

Note that the associativity of the product rules obtained in the proof of Theorem \ref{homog-thm}(b)
amounts to the identity
\begin{equation}\label{metric-identity}
(x,e(0))^R\cdot (y,e(0))^R=((\zeta_\bullet)_*x,y)^R\cdot (1-t)^n=(x,y)^R\cdot (1-t)^n.
\end{equation}
This can be checked independently using Example \ref{k-st-ex}. Namely, the left-hand side corresponds
to the functor of restriction to the origin in $\A^n\times\A^n$, and the right-hand-side is obtained
by first restricting to the shifted diagonal and then to the origin.

\begin{rem}
The specialization of $\HH$ at $t=1$ is the algebra
$$\C[\HH_{0,0},u]/(u^{d-1}, xu, xy-((\zeta_\bullet)_*x,y)u^{d-2}),$$
where $x,y\in\HH_{0,0}=H(\w)^{\Z/d}$. Indeed, this
follows easily from the fact that $e(0)$ has zero component in
$\HH_{0,0}$ (see \cite{PV-mf}).
In the case when $d=n$ the category of $\G_m$-equivariant
\mf s of $\w$ is equivalent to the derived category of coherent
sheaves on the corresponding Calabi-Yau hypersurface $X\sub\P^n$
(see \cite{Orlov-graded})
and the above algebra is isomorphic to $H^*(X,\C)$. More precisely, the subalgebra generated
by $u$ corresponds to classes in $H^*(X,\C)$ restricted from $\P^n$, while the
subspace $\HH_{0,0}$ corresponds to the primitive part of the middle cohomology $H^n(X,\C)$.
\end{rem}

\section{Simple singularities}\label{simple-sing-sec}

In this section we will calculate the Frobenius algebra structure on the $R$-module
$\HH=\bigoplus_{\ga\in G}\HH_{\ga}$
given by the \CFT\ of Theorem
\ref{CohFT-thm} for all simple singularities, i.e., for singularities of type $A_n$, $D_n$, $E_6$,
$E_7$ and $E_8$.
 We will use the notation and the results of Section \ref{calc-sec}.

\subsection{Singularity of type $A$}\label{A-sec} 

Consider the $A_{d-1}$ singularity 
$\w=x^d$ with the symmetry group $G=G_\w=\lan J\ran$, where $J=\exp(2\pi i m/d)$
(this is the only possible symmetry group allowed by our construction). 
Since $\w$ a homogeneous polynomial, we can apply Theorem \ref{homog-thm}. We keep the notation
of Section \ref{hom-pol-sec}.
Thus, for $m\in\Z/d$ we have
$$\HH_m=\begin{cases} R, & m\not\equiv 0,\\ 
HH_*(\MF_{\G_m}(\A^1,x^d)), &m\equiv 0.\end{cases}$$
From the description of the Hochschild homology of $\MF_{\G_m}(\A^1,x^d)$ 
(see Theorem \ref{hoch-prop}) and the calculation of
the Chern character $\ch(\C^{\st})$ in \cite[Ex. 4.2.2]{PV-mf} 
we have an identification of $R$-modules
$$R/(1+t+\ldots+t^{d-1})\rTo{\sim}\HH_0: 1\mapsto e(0)=\ch(\C^{\st}).$$

The metric on $\HH_0$ is given by
\begin{equation}\label{A-n-metric-eq}
(e(0),e(0))=(e(0),\zeta^*e(0))^R=(e(0),e(0))^R=1-t,
\end{equation}
where $\zeta=\exp(\pi i/d)$ (see \eqref{metric-eq}).
Indeed, since $e(0)=\ch(\C^{\st})$, by Lemma \ref{stable-functor-lem}, 
$(e(0),e(0))^R$ is equal to the class of the $\Z/d$-representation
$(\C^{\st})|_0$ in $R$. Also, $\zeta^*$ acts trivially on the regular Koszul \mf\ $\C^{\st}$.

Thus, we have an isomorphism of $R$-modules
$$\HH=R/(1+t+\ldots+t^{d-1})\cdot e(0)\oplus\bigoplus_{m\in\Z/d, m\not\equiv 0} R\cdot e(m).$$
The metric on $\HH$ is given by
$$(e(m),e(-m))=\begin{cases} 1, &m\not\equiv 0,\\ 1-t, &m\equiv 0.\end{cases}$$
By Theorem \ref{homog-thm}, setting $u=e(2)$ we obtain an isomorphism of $R$-algebras
$$\HH\simeq R[u]/(u^d-1+t,(1+t+\ldots+t^{d-1})u^{d-1}).$$
If we specialize with respect to $\pi_1:t\mapsto 1$ we get the Frobenius ring
$\C[u]/(u^{d-1})$, i.e., the Milnor ring of the same singularity. 
On the other hand, the specialization with respect to
$\pi_{\om}:t\mapsto \om$, where $\om$ is a nontrivial $d$th root of unity,
gives the semisimple ring $\C[u]/(u^d-1+\om)$.

\subsection{Singularity of type $D$ with the maximal group of diagonal symmetries}
\label{D-max-sec}

Consider the $D_{d+1}$ singularity $\w=x^d+xy^2$. 
The group $G=G_\w$ of diagonal symmetries is isomorphic to $\Z/2d$, where $m\in\Z/2d$ acts on $\A^2$
by $(\exp(-2\pi im/d),\exp(\pi im/d))$. 
In this case we denote the line bundles on $\P^1$ associated with a $\Ga$-spin structure
as $L_x$ and $L_y$.
We will use the identification $R=\C[t]/(t^{2d}-1)$, where
$t:\Z/2d\to\C^*$ is given by \eqref{t-eq} with $d$ replaced by $2d$.

Let us calculate the $R$-modules $\HH_m$ in the decomposition
of the state space $\HH=\bigoplus_{m\in \Z/2d}\HH_m$.
For $m=0$ we claim that there is an isomorphism of $R$-modules 
$$R/(1-t+t^2-t^3+\ldots-t^{2d-1})\to\HH_0:1\mapsto e(0),$$
where $e(0)=\ch(\bar{E})$ with $\bar{E}=\{x^{d-1}+y^2,x\}$.
Indeed, consider the decomposition
$$\HH_0=\bigoplus_{m\in\Z/2d}H(\w_m)^{\Z/2d},$$
where $\w_m$ is the restriction of $\w$ to the subspace of $m$-invariants.
We have
$$H(\w_0)^{\Z/2d}=\C\cdot y\, dx\we dy, \ \ H(\w_d)^{\Z/2d}=0, \text{ and }\ H(\w_m)^{\Z/2d}=\C
\text{ for }\ m\not\equiv 0,d.
$$
Now our claim follows from the fact that $e(0)$ has nonzero components in
$H(\w_m)^{\Z/2d}$ for all $m\not\equiv d$ (this computation is analogous to
the case $m=0$ and $d=3$ considered in \cite[Ex. 4.1.8]{PV-mf}).

For $m\not\equiv 0,d$ the components $\HH_m$ can be identified with $R\cdot e(m)$, and for $m\equiv d$
we have an isomorphism
$$\HH_d=HH_*(\MF_{\Z/2d}(\A^1,x^d))\simeq R/(1+t^2+t^4+\ldots+t^{2d-2})\cdot e(d),$$
where $e(d)=\ch(\C^{\st})$ (here $m\in \Z/2d$ acts on $\A^1$ by $\exp(2\pi i m/d)$).
Note that $(\zeta_\bullet)^*e(m)=e(m)$ since $e(m)$ is the Chern class of an object invariant under
$(\zeta_\bullet)^*$.

The metric on $\HH$ is given by
$$(e(m),e(-m))=1 \text{ for } m\not\equiv 0,d,$$
$$(e(d),e(d))=1-t^{-2},$$
$$(e(0),e(0))=-(1+t)t^{-2}.$$
The second equality follows from the case of singularity $A_{d-1}$.
To get the last equality we observe that by Lemma \ref{stable-functor-lem}, $(e(0),e(0))^R$ is equal to
the class of the $\Z/2$-graded complex of $G$-modules (with finite-dimensional cohomology)
\begin{equation}\label{e-0-restriction-eq}
\pi_*(\bar{E}|_{x=0})=\pi_*(\{y^2,0\})\simeq \C[y]/(y^2)\ot t^{-2}\, [1],
\end{equation}
and the class of $\C[y]/(y^2)$ is equal to $1+t$.

To calculate $\La^R_0(e(m_1),e(m_2),e(m_3))\in R$ we can apply
Proposition \ref{3-pt-prop}(ii) in most cases.
From now on we assume that $0\le m_i<2d$ for $i=1,2,3$.
Note that $\bq=(\frac{1}{d},\frac{d-1}{2d})$ while
$$\bth_i=\bth_{m_i}=\begin{cases} (0,0), &m_i=0,\\ (1-\frac{m_i}{d},\frac{m_i}{2d}), &0<m_i\le d, \\ (2-\frac{m_i}{d},\frac{m_i}{2d}), &d<m_i<2d. \end{cases}$$
Thus, the condition \eqref{deg-0-3-eq} implies that $L_x=\OO$ only when two of the $m_i$'s are zero.
On the other hand, $L_y=\OO$ if and only if $m_1+m_2+m_3=d-1$, in which case it is possible that only one
or none of $m_i$'s is trivial.
Thus, we have two cases not covered by Proposition \ref{3-pt-prop}(ii): 

\noindent
{\bf Case 1a}. $m_1+m_2+m_3=d-1$ with $m_i>0$ for $i=1,2,3$;

\noindent
{\bf Case 1b}. $m_1=0$ and $m_2+m_3=d-1$ with $m_2>0$ and $m_3>0$.

Consider first Case 1a (which can occur only for $d\ge 4$). 
We have $\bth_i=(1-m_i/d,m_i/(2d))$, so $\bth_1+\bth_2+\bth_3=(2+1/d,(d-1)/(2d))$.
Thus, we have $L_x=\OO(-2)$ and $L_y=\OO$. Hence, this is the case of index zero considered in
Section \ref{index-zero-sec} (since $\chi(\OO(-2))+\chi(\OO)=0$).
Also, since $0<m_i<d-1$, all these group elements have trivial invariants in $\A^2$. Thus, applying
Proposition \ref{index-zero-prop}(ii), we obtain  
$$\La^R_0(e(m_1),e(m_2),e(m_3))=\left.\frac{1-t_1^{-d}\cdot t_1}{1-t_2}\right|_{t_1=t^{-2}, t_2=t}~,
$$
where $t_1$ and $t_2$ are characters of the group $\Ga$ satisfying $t_1^{d-1}=t_2^2$. Thus, we have
$$\La^R_0(e(m_1),e(m_2),e(m_3))=\left.\frac{1-t_2^{-2}}{1-t_2}\right|_{t_2=t}~=-(1+t)t^{-2}.
$$

In Case 1b we have 
$\bth_1+\bth_2+\bth_3=(1+1/d,(d-1)/2d)$, so $L_x=\OO(-1)$ and $L_y=\OO$. Also, 
we have $\Si_x=\Si_y=\{p_1\}$. 
By Proposition \ref{3-pt-prop}(i), 
in this case the fundamental \mf\ $\bP(\psi)$ is the Koszul \mf\ $\{-x^{d-1}-y^2,x\}$ on $\A^2$. 
Hence, by Lemma \ref{stable-functor-lem},
the corresponding functor $\MF_{\Ga}(\w)\to\Com_f(G-\mod)$ associates with $\bar{E}$
the restriction $\com_G(E|_{x=0})$, which we computed in \eqref{e-0-restriction-eq}. This gives
$$\La^R_0(e(0),e(m),e(d-1-m))=-(1+t)t^{-2} \text{ for } 0<m<d-1.$$

In all other cases we can apply Proposition \ref{3-pt-prop}(ii). Note that $\La^R_0(e(m_1),e(m_2),e(m_3))$
is only nonzero when $m_1+m_2+m_3\equiv d-1 (2d)$ (since $J$ corresponds to $d-1\in\Z/2d$). 
Note also that by \eqref{deg-0-3-eq}, 
$$L_y=\OO(-a),$$
where $m_1+m_2+m_3=d-1+2da$. 

\noindent
{\bf Case 2a}. $m_1+m_2+m_3=3d-1$ with $0<m_1<d$, $0<m_2<d$ and $d<m_3<2d$.

\noindent
In this case $L_x=L_y=\OO(-1)$, so $S_0=S_1=S_2=\emptyset$ and we get
$$\La^R_0(e(m_1),e(m_2),e(m_3))=1.$$

\noindent
{\bf Case 2b} (occurs only for $d\ge 4$). $m_1+m_2+m_3=3d-1$ with $0<m_1<d$, $d<m_2<2d$ and $d<m_3<2d$.

\noindent
In this case $L_x=\OO(-2)$ and $L_y=\OO(-1)$, so $S_1=S_2=\emptyset$ and $S_0=\{x\}$. Thus, we get
$$\La^R_0(e(m_1),e(m_2),e(m_3))=1-t^{-2}.$$

\noindent
{\bf Case 2c}. $m_1+m_2+m_3=5d-1$ with $d<m_i<2d$ for $i=1,2,3$.

\noindent
In this case $L_x=\OO(-1)$ and $L_y=\OO(-2)$, so
$S_1=S_2=\emptyset$ and $S_0=\{y\}$. Thus, 
$$\La^R_0(e(m_1),e(m_2),e(m_3))=1-t.$$

\noindent
{\bf Case 3}. $m_1=d$, $0<m_2<d$, $d<m_3<2d$ with $m_1+m_2=2d-1$.
In this case $L_x=L_y=\OO(-1)$
$S_1=\{x\}$ and $S_0=S_2=\emptyset$, so
we have to calculate the class of the restriction of $\C^{\st}\in\MF_{\G_m}(\A^1,x^d)$ to $0\in\A^1$.
Thus,
$$\La^R_0(e(d),e(m),e(2d-1-m))=1-t^{-2} \text{ for } 0<m<d-1.$$

\noindent
{\bf Case 4}. $m_1=0$, $m_2+m_3=3d-1$, where $d<m_2<2d$ and $d<m_3<2d$.
In this case $L_x=L_y=\OO(-1)$, so
$S_1=\{x,y\}$ and $S_0=S_2=\emptyset$. Hence, we are reduced to 
computing the class of the restriction of $\bar{E}$ to $x=y=0$. Therefore, we get
$$\La^R_0(e(0),e(m),e(3d-1-m))=1-t^{-2} \text{ for } d<m<2d-1.$$

\noindent
{\bf Case 5}. $m_1=0$, $m_2=d$, $m_3=2d-1$.
In this case $L_x=\OO$, $L_y=\OO(-1)$, so we have
$S_1=\{y\}$, $S_2=\{x\}$ and $S_0=\emptyset$. Thus, we have to calculate the class
of the restriction 
$$\bar{E}\boxtimes \C^{\st}\in\MF(\A^2\times\A^1,\w(x_1,y)\oplus x_2^d)$$
to the linear subspace $x_1=x_2$, $y=0$.
Since tensoring with $\C^{\st}$ has the same effect as the restriction to the origin,
this class is equal to the class of the restriction of $\bar{E}$ to $x=y=0$, so
we obtain as in the previous case
$$\La^R_0(e(0),e(d),e(2d-1))=1-t^{-2}.$$

The remaining two cases $(m_1,m_2,m_3)=(0,0,d-1)$ or $(d,d,d-1)$
can be computed using \eqref{genus-0-formula}, since $e(d-1)$ is the flat unit $\unit$ for our theory.

Now we can determine the ring structure on $\HH$.
The element $e(d-1)$ is a unit. We always have 
$$e(m_1)e(m_2)=r(m_1,m_2)e(m_1+m_2-d+1)$$
with some $r(m_1,m_2)\in R$.
From the formulas for $\La^R_0(?,?,?)$ we get the relations
$$e(d-2)e(m)=e(m-1) \text{ for } 1\le m\le d-2,$$
$$e(d-2)^d=e(d-2)e(0)=-(1+t)t^{-2} e(-1),$$
$$e(d-2)e(-m)=e(-m-1) \text{ for } 1\le m\le d-1,$$
$$e(-1)^2=(1-t)e(d-1).$$
Thus, $\HH$ is generated as an $R$-algebra by elements
$u=e(d-2)$ and $v=e(-1)$ subject to the relations
$$v^2=1-t, \ \ u^d=-(1+t)t^{-2}v, \ \  (1-t)(1+t^2+\ldots+t^{2d-2})u^{d-1}=(1+t^2+\ldots+t^{2d-2})u^{d-1}v=0.$$

The specialization $t=1$ gives the Frobenius algebra $\C[u]/(u^{2d-1})$ 
with the pairing given by
$$(u^{2d-2},1)=-2(u^{d-2}e(-1),e(d-1))=-2(e(-d+1),e(d-1))=-2.$$

\subsection{Singularity of type $D$ with the non-maximal symmetry group}
\label{D-J-sec}

The group of diagonal symmetries $G_\w$ of the $D_{d+1}$ singularity $\w=x^d+xy^2$
is not generated by the exponential grading element $J=(\exp(2\pi i/d),\exp(2\pi i k/d))$
precisely when $d=2k+1$ is odd.   In this case $J$ has order $d$, so
the subgroup $G=\lan J\ran$ has index two in $G_\w$.
In this section we will calculate the Frobenius algebra corresponding to this subgroup.
We will use the identification $R=\C[t]/(t^d-1)$, where the generating character $t$ is defined by
$t(J)=\exp(-2\pi i/d)$.
As before we denote the line bundles of a $\Ga$-spin structure
by $L_x$ and $L_y$.
For $m\in\Z/d$ we set $\HH_m=\HH_{J^m}$. Since all nontrivial powers of $J$ have
no invariants on $\A^2$, we have $\HH_m=R$ for $m\not\equiv 0$. For such $m$ 
we denote by $e(m)$ the element $1\in\HH_m$. To compute $\HH_0$ as an $R$-module
consider $\lan J\ran$-equivariant \mf s
$$\bar{E}_{\pm}=\{x^{k+1}\pm i xy, x^k\mp iy\}.$$
and denote their Chern characters by $e_\pm(0)=\ch(\bar{E}_{\pm})$.
We claim that the map
$$R\oplus R\to\HH_0: (r_1,r_2)\mapsto r_1e_+(0)+r_2(e_+(0)-e_{-}(0))$$
induces an isomorphism of $R$-modules
$$R\oplus R/(t-1)\to\HH_0.$$
Indeed, the components of the decomposition
$$\HH_0=\bigoplus_{m\in\Z/d} H(\w_{J^m})^{J}$$
are 
$$H(\w)^J=(\C\cdot x^{k}+\C\cdot y)\cdot dx\we dy, \ \ H(\w_{J^m})^{J}=\C \text{ for }m\not\equiv 0.$$
Now using \cite[Thm.\ 3.3.3]{PV-mf} 
we obtain that $\ch(\bar{E}_{\pm})$ have nonzero $m$-components for
$m\not\equiv 0$, and their $0$-component are given by
$$\ch(\bar{E}_{\pm})_0=(\mp idx^{k}+y)\cdot dx\we dy.$$
This immediately implies our claim.

The metric on $\HH$ is given by
$$(e(m),e(-m))=1 \text{ for }\ m\not\equiv 0,$$
$$(e_\pm(0),e_\pm(0))=-(1+t+\ldots+t^k)t^{k}$$
$$(e_+(0),e_-(0))=1+t+\ldots+t^{k-1}.$$
Indeed, the last two equalities follow from the quasi-isomorphisms
$$\pi_*\bar{E}_{\pm}|_{x^k\mp iy=0}=\pi_*\{2x^{k+1},0\}\simeq \C[x]/(x^{k+1})\ot t^{-k-1}[1],$$
$$\pi_*\bar{E}_+|_{x^k+ iy=0}=\pi_*\{0,2x^{k}\}\simeq \C[x]/(x^k).$$

The calculation of the three-point correlators
$\La_0^R(m_1,m_2,m_3)$ is done similarly to the case of the maximal symmetry group.
We have $\bq=(\frac{1}{d},\frac{k}{d})$ and
$$\bth_{m_i}=\begin{cases}(0,0), &m_i=0, \\ (\frac{2l}{d},1-\frac{l}{d}), &m_i=2l, 0<l\le k, \\
(1-\frac{2l}{d}, \frac{l}{d}), &m_i=-2l, 0<l\le k.\end{cases}$$
The condition \eqref{deg-0-3-eq} implies that $L_x=\OO$ only when two of the $m_i$'s are trivial.

We have the following cases with $L_y=\OO$ and $|\Si_y|\le 1$.

\noindent
{\bf Case 1a.} $m_i=-2l_i$, $i=1,2,3$, where $l_1+l_2+l_3=k$, $l_i>0$.
In this case $L_x=\OO(-2)$ and $\Si_x=\Si_y=\emptyset$. Hence, by Proposition \ref{index-zero-prop}(ii),
we obtain
$$\La_0^R(e(-2l_1),e(-2l_2),e(-2l_3))=\left.\frac{1-t^{-d}\cdot t}{1-t^k}\right|_{t^d=1}~=-(1+t^k)t.$$

\noindent
{\bf Case 1b.} $m_1=0$, $m_2=-2l_2$, $m_3=-2l_3$, where $l_2+l_3=k$, $l_i>0$. In this case
$L_x=\OO(-1)$ and $\Si_x=\Si_y=\{p_1\}$. By Proposition \ref{3-pt-prop}(i), we have
$$\bP(\psi)=\{-x^{d-1}-y^2,x\}.$$
Therefore, by Lemma \ref{stable-functor-lem}, the corresponding functor $\MF_{\Ga}(\w)\to\Com_f(G-\mod)$ is given by
the restriction to $x=0$. We have
$$\pi_*\bar{E}_\pm|_{x=0}=\pi_*\{0,\mp iy\}\simeq\C.$$
Hence,
$$\La_0^R(e_\pm(0),e(-2l_2),e(-2l_3))=1.$$

In the cases when $L_y\neq\OO$ or $|\Si_y|\ge 2$
 we can apply Proposition \ref{3-pt-prop}(ii). Note that
$\La_0^R(m_1,m_2,m_3)$ is nonzero only when $m_1+m_2+m_3\equiv 1 \mod d$.

\noindent
{\bf Case 2a.} $m_1=2l_1$, $m_2=-2l_2$, $m_3=-2l_3$, where $1\le l_i\le k$ and 
$l_2+l_3-l_1=k$. In this case $L_x=L_y=\OO(-1)$, so $S_0=S_1=S_2=\emptyset$, and we get
$$\La_0^R(e(2l_1),e(-2l_2),e(-2l_3))=1.$$

\noindent
{\bf Case 2b.} $m_1=2l_1$, $m_2=2l_2$, $m_3=-2l_3$, where $1\le l_i\le k$ and
$-l_3+l_1+l_2=k+1$. In this case $L_x=\OO(-2)$ and $L_y=\OO(-1)$, so $S_1=S_2=\emptyset$ and
$S_0=\{x\}$. Hence, 
$$\La_0^R(e(2l_1),e(2l_2),e(-2l_3))=1-t.$$

\noindent
{\bf Case 2c.} $m_i=2l_i$, where $1\le l_i\le k$ and $l_1+l_2+l_3=k+1$. In this case $L_x=\OO(-1)$ and 
$L_y=\OO(-2)$, so $S_1=S_2=\emptyset$ and
$S_0=\{y\}$. Hence, 
$$\La_0^R(e(2l_1),e(2l_2),e(2l_3))=1-t^k.$$

\noindent
{\bf Case 3.} $m_1=0$, $m_2=2l$, $m_3=2(k+1-l)$, where $1\le l\le k$.
In this case $L_x=L_y=\OO(-1)$, so $S_0=S_2=\emptyset$ and $S_1=\{x,y\}$. Therefore, 
$\La_0^R(0,2l,2(k+1-l))$ sends the class of a \mf\ of $\w$ on $\A^2$ to the class of its restriction to the origin.
Thus,
$$\La_0^R(e_\pm(0),e(2l),e(2(k+1-l)))=1-t^k.$$

The remaining case $m_1=m_2=0$, $m_3=1$ follows from Lemma
\ref{metric-lem} since $e(1)$ is the flat unit.

Now let us determine the ring structure on $\HH$. 
The element $e(1)$ is a unit.
Note that the product of elements in $\HH_{m_1}$ and $\HH_{m_2}$ lies in $\HH_{m_1+m_2-1}$.
Using the above calculations we obtain
$$e_\pm(0)^2=-(1+t+\ldots+t^k)t^k e(2k), \ \ e_+(0)e_-(0)=(1+t+\ldots+t^{k-1})e(2k),$$
$$e_\pm(0)e(2l)=(1-t^k)e(-2(k+1-l)), $$
$$e_\pm(0)e(-2l)=e(2(k-l)) \text{ for } l<k,$$ 
$$e(2l_1)e(2l_2)=\begin{cases} (1-t^k) e(-2(k+1-l_1-l_2)), &  l_1+l_2<k+1, \\
(1-t)e(2(l_1+l_2-k-1)), & l_1+l_2>k+1, \\
\frac{1-t}{2}\cdot (e_+(0)+e_-(0)), & l_1+l_2=k+1, \end{cases}$$
$$e(-2l_1)e(-2l_2)=\begin{cases} -t(1+t^k)e(2(k-l_1-l_2)), & l_1+l_2<k, \\
e(-2(l_1+l_2-k)), & l_1+l_2>k, \\
-\frac{t(1+t^k)}{2}\cdot(e_+(0)+e_-(0)), & l_1+l_2=k,
\end{cases}$$
$$e(2l_1)e(-2l_2)=\begin{cases} e(2(k+l_1-l_2)), & l_1\le l_2,\\
(1-t)e(-2(k+1-l_1+l_2)), & l_1>l_2, \end{cases}$$
where $1\le l\le k$, $1\le l_i\le k$.

Assume first that $k>1$. Then setting $u=e(3)=e(-2(k-1))$ we obtain the following relations
$$e(-2l)=u^{k-l} \text{ for } 1\le l\le k-1; \ \ e(2l)=e_\pm(0)u^l \text{ for } 1\le l\le k,$$ 
$$e_\pm(0)^2u=1-t^k, \ \ e_\pm(0)u^{k+1}=1-t,$$
$$u^k=-\frac{t(1+t^k)}{2}\cdot(e_+(0)+e_-(0)),$$
$$e_-(0)(e_+(0)-e_-(0))=-e_+(0)(e_+(0)-e_-(0))=(1+t+\ldots+t^{2k})e_\pm(0)u^k.$$
Thus, $\HH$ is generated as an $R$-algebra by the elements $u$ and $v=e_+(0)-e_-(0)$, subject
to the relations
$$(t-1)v=uv=0,$$
$$u^{2k+1}=-t(1+t^k)(1-t),$$
$$t(1+t^k)v^2=2(1+t+\ldots+t^{2k})u^{2k}.$$
The specialization $t=1$ gives the Frobenius
algebra $\C[u,v]/(uv,v^2-du^{d-1})$ with the pairing 
$$(u^{d-1},1)=-((e_+(0)+e_-(0))u^k,e(1))=-2(e(2k),e(1))=-2.$$
Note that this Frobenius algebra is isomorphic to the Milnor ring of the singularity $D_{d+1}$.

In the $D_4$ case, i.e., when $k=1$, we obtain that $\HH$ is generated as an $R$-algebra by
the elements $e_\pm(0)$, subject to the relations
$$(1-t)(e_+(0)-e_-(0))=0,$$
$$e_+(0)^2=e_-(0)^2=-t(1+t)e_+(0)e_-(0),$$
$$e_+(0)^3=e_-(0)^3=1-t.$$
If we take generators $u=(e_+(0)+e_-(0))/2$ and $v=(e_+(0)-e_-(0))/2$, the relations
become
$$uv=(1+t+t^2)u^2+(1-t-t^2)v^2=0,$$
$$u^3=1, \ \ v^3=-t.$$
The specialization $t=1$ gives the Frobenius algebra
$\C[u,v]/(uv, 3u^2-v^2)$
with the pairing
$$(u^2,1)=-\frac{1}{2}.$$
This Frobenius algebra is isomorphic to the Milnor ring of the $D_4$ singularity.

\subsection{$E_7$ singularity}\label{E7-sec}

In the case of the $E_7$ singularity $\w=x^3+xy^3$
the maximal symmetry group $G=G_\w$ is generated by $J=(\exp(2\pi i/3),\exp(4\pi i/9))$.
We will use the identification of $G$ with $\Z/9$ where $m\in\Z/9$ acts on $\A^2$ by 
$(\exp(-2\pi im/3), \exp(2\pi im/9))$. 
We will denote by $t$ the character of $\Z/9$ given by \eqref{t-eq} (with $d=9$) and use the 
identification $R=\C[t]/(t^9-1)$. 

First, let us determine the $R$-module $\HH=\bigoplus_{m\in\Z/9}\HH_m$.
We claim that there is an isomorphism of $R$-modules
$$R/((t-1)(1+t^3+t^6))\to \HH_0: 1\mapsto e(0),$$
where $e(0)=\ch(\bar{E})$ with $\bar{E}=\{x^2+y^3,x\}$. Indeed, the summands of the 
decomposition \eqref{graded-G-hoch-isom} of 
$\HH_0$ are
$$H(\w_0)^G=\C\cdot y^2\cdot dx\we dy, \ H(\w_{\pm 3})^G=0,\
H(\w_m)^G=\C \text{ for }m\not\equiv 0,\pm 3,
$$
where $\w_m$ is the restriction of $\w$ to the subspace 
of $m$-invariants in $\A^2$. Our claim follows from the fact that $e(0)$ has nonzero components in
$H(\w_m)$ for all $m\not\equiv\pm 3$. 

The components $\HH_m$ for $m$ not divisible by $3$ can be identified
with $R\cdot e(m)$, while
$$\HH_{\pm 3}=HH_*(\MF_{\Z/9}(\A^1,x^3))\simeq R/(1+t^3+t^6)\cdot e(\pm 3),$$
where $e(\pm 3)=\ch(\C^{\st})$.

The metric on $\HH$ is given by
$$(e(m),e(-m))=1 \text{ for }\ m\not\equiv 0,\pm 3,$$
$$(e(3),e(-3))=1-t^{-3},$$
$$(e(0),e(0))=-(1+t+t^2)t^{-3},$$
where the last equality follows from the quasi-isomorphism
$$\pi_*\bar{E}|_{x=0}=\pi_*\{y^3,0\}\simeq \C[y]/(y^3)\ot t^{-3}[1].$$

Now let us compute the three-point correlators $\La_0^R(e(m_1),e(m_2),e(m_3))$.
We have $\bq=(\frac{1}{3},\frac{2}{9})$ and
$$\bth_{m_i}=\begin{cases} (0,\frac{l}{3}), &m_i=3l, 0\le l\le 2,\\ 
(\frac{2}{3},\frac{3l+1}{9}), &m_i=3l+1, 0\le l\le 2, \\ (\frac{1}{3},\frac{3l+2}{9}), 
&m_i=3l+2, 0\le l\le 2. \end{cases}$$
Note that $L_x=\OO$ implies that $|\Si_x|=2$. Hence, we have only one case not covered by
Proposition \ref{3-pt-prop}(ii).

\noindent
{\bf Case 1.} $m_1=0$, $m_2=1$, $m_3=1$.
In this case
$L_x=\OO(-1)$, $L_y=\OO$, $\Si_x=\Si_y=\{p_1\}$.
Thus, we have $S_0=S_2=\emptyset$ and $S_1=\{x\}$.
Then by Proposition \ref{3-pt-prop}(i), 
$\bP(\psi)$ is the Koszul \mf\ $\{-x^2-y^3,x\}$. 
Thus, $\La_0^R(e(0),e(1),e(1))$ is given by
the class of $\pi_*(\bar{E}|_{x=0})$ computed above.
Hence,
$$\La_0^R(e(0),e(1),e(1))=-(1+t+t^2)t^{-3}.
$$

In all remaining cases we can use Proposition \ref{3-pt-prop}(ii). 
Note that $\La_0^R(e(m_1),e(m_2),e(m_3))$
is nonzero only when $m_1+m_2+m_3\equiv 2 \mod 9$.

\noindent
{\bf Case 2a.} $m_1=3l_1+1$, $m_2=3l_2+2$, $m_3=3l_3+2$, where $0\le l_i\le 2$, $l_1+l_2+l_3=2$. 
We have $L_x=L_y=\OO(-1)$, $S_0=S_1=S_2=\emptyset$. Hence, in this case
$$\La_0^R(e(3l_1+1),e(3l_2+2),e(3l_3+2))=1.$$

\noindent
{\bf Case 2b.} $m_1=3l_1+1$, $m_2=3l_2+2$, $m_3=3l_3+2$, where $0\le l_i\le 2$, $l_1+l_2+l_3=5$. 
We have $L_x=\OO(-1)$, $L_y=\OO(-2)$, $S_1=S_2=\emptyset$, $S_0=\{y\}$. 
Hence, in this case
$$\La_0^R(e(3l_1+1),e(3l_2+2),e(3l_3+2))=1-t.$$

\noindent
{\bf Case 3a.} $m_1=3l_1$, $m_2=3l_2+1$, $m_3=3l_3+1$, where $0\le l_i\le 2$, $l_1>0$, $l_1+l_2+l_3=3$. 
We have $L_x=L_y=\OO(-1)$, $S_0=S_2=\emptyset$, $S_1=\{x\}$.
Hence, in this case we have to compute the class of the restriction of $e(3l_1)$ to the origin, which gives
$$\La_0^R(e(3l_1),e(3l_2+1),e(3l_3+1))=1-t^{-3}.$$

\noindent
{\bf Case 3b.} $m_1=6$, $m_2=m_3=7$. 
We have $L_x=\OO(-1)$, $L_y=\OO(-2)$, $S_2=\emptyset$, $S_0=\{y\}$, $S_1=\{x\}$.
Hence, 
$$\La_0^R(e(6),e(7),e(7))=(1-t)(1-t^{-3}).$$

\noindent
{\bf Case 3c.} $m_1=0$, $m_2=3l_2+1$, $m_3=3l_3+1$, where $0\le l_i\le 2$, $l_2+l_3=3$. 
We have $L_x=L_y=\OO(-1)$, $S_0=S_2=\emptyset$, $S_1=\{x,y\}$.
Hence, we have to compute the restriction of $\bar{E}$ to the origin, which gives
$$\La_0^R(e(0),e(3l_2+1),e(3l_3+1))=1-t^{-3}.$$

\noindent
{\bf Case 4a.} $m_1=3l_1$, $m_2=3l_2$, $m_3=3l_3+2$, where $0\le l_i\le 2$, $l_1>0$, $l_2>0$, $l_1+l_2+l_3=3$. 
We have $L_x=\OO$, $L_y=\OO(-1)$, $S_0=S_1=\emptyset$, $S_2=\{x\}$.
Hence, in this case we have to compute the class of 
$\pi_*(\C^{\st}\otimes\C^{\st})$, which is equivalent to computing the pairing $((e(3),e(-3))$. Therefore, we get
$$\La_0^R(e(3l_1),e(3l_2),e(3l_3+2))=1-t^{-3}.$$

\noindent
{\bf Case 4b.} $m_1=6$, $m_2=6$, $m_3=8$. 
We have $L_x=\OO$, $L_y=\OO(-2)$, $S_1=\emptyset$, $S_0=\{y\}$, $S_2=\{x\}$.
Hence, 
$$\La_0^R(e(6),e(6),e(8))=(1-t)(1-t^{-3}).$$

\noindent
{\bf Case 4c.} $m_1=0$, $m_2=3l_2$, $m_3=3l_3+2$, where $l_2>0$, $l_3>0$, $l_2+l_3=3$. 
We have $L_x=\OO$, $L_y=\OO(-1)$, $S_0=\emptyset$, $S_1=\{y\}$, $S_2=\{x\}$.
Hence, we have to compute the restriction of the \mf\ $\bar{E}\boxtimes\C^{\st}$ on the space
$\A^2\times\A^1$ with coordinates $(x,y,x')$ to
the subspace $ax+bx'=y=0$. Restricting to $y=0$ we get $\C^{\st}\boxtimes\C^{\st}$, so the answer
is the same as in Case 4a:
$$\La_0^R(e(0),e(3l_2),e(3l_3+2))=1-t^{-3}.$$

The remaining case $m_1=m_2=0$, $m_3=2$ follows from the metric axiom.

Now we can determine the ring structure on $\HH$. Note that $e(m_1)\cdot e(m_2)$ is always
proportional to $e(m_1+m_2-2)$ and $e(2)$ is a unit of $\HH$. From the above computation
of $\La_0^R(e(m_1),e(m_2),e(m_3))$ and the formulas for the metric
we obtain the following multiplication table:
$$e(1)e(5)=e(4),\ \ e(1)e(8)=e(4)e(5)=e(7),$$
$$e(4)e(8)=e(7)e(5)=(1-t)e(1), \ \ e(7)e(8)=(1-t)e(4),$$
$$e(5)^2=e(8), \ \ e(5)e(8)=(1-t), \ \ e(8)^2=(1-t)e(5),$$
$$e(1)^2=e(0), \ \ e(1)e(5)=e(3), \ \ e(1)e(7)=e(4)^2=e(6), \ \ e(4)e(7)=(1-t)e(0),$$
$$e(0)e(1)=-(1+t+t^2)t^{-3}e(-1), \ \ e(6)e(7)=(1-t)(1-t^{-3}),$$
$$e(0)e(4)=e(3)e(1)=(1-t^{-3}), \ \ e(0)e(7)=e(3)e(4)=e(6)e(1)=(1-t^{-3})e(5),$$
$$e(3)e(7)=e(6)e(4)=(1-t^{-3})e(8),$$
$$e(3)e(5)=e(6), \ \ e(3)e(8)=e(6)e(5)=(1-t)e(0), \ \ e(6)e(8)=(1-t)e(3),$$
$$e(3)^2=e(0)e(6)=(1-t^{-3})e(4), \ \ e(3)e(6)=(1-t^{-3})e(7), \ \ e(6)^2=(1-t)(1-t^{-3})e(1),
$$
$$e(0)e(3)=(1-t^{-3})e(1), \ \ e(0)^2=-(1+t+t^2)t^{-3}e(7).$$

It follows that $\HH$ is generated as an $R$-algebra by the elements $u=e(5)$ and $v=e(1)$ subject to
the relations
$$u^3=1-t, \ \ v^3=-(1+t+t^2)t^{-3}u^2, \ \ (1+t^3+t^6)uv^2=0.$$

The specialization $t=1$ gives the Frobenius algebra $\C[u,v]/(v^3+3u^2, uv^2)$ with the pairing
$$(u^2v,1)=(e(7),e(2))=1.$$
Thus, we again obtain the Milnor ring of the same singularity.

\subsection{$E_6$ and $E_8$ singularities}\label{E68-sec}

In the case of the $E_6$ and $E_8$ singularities $\w=x^3+y^4$ and $\w=x^3+y^5$ the maximal group $G=G_\w$ of diagonal
symmetries coincides with $\lan J\ran$. Thus, in both cases we have $\w=\w_1\oplus\w_2$ and
$G=G_1\times G_2$, where $(\w_i,G_i)$ is an $A_n$-singularity for some $n$. 
By Theorem \ref{sum-sing-thm}, the corresponding Frobenius algebras over $R$ are tensor products (over $\C$)
of the Frobenius algebras corresponding to $(\w_1,G_1)$ and $(\w_2,G_2)$.
Thus, for $E_6$ we have $R=\C[t]/(t^{12}-1)$, and $\HH_{E_6}$ is the $R$-algebra generated by 
$u$ and $v$ subject to the relations
$$u^3=1-t^4, \ \ (1+t^4+t^8)u^2=0, \ \ v^4=1-t^3, \ \ (1+t^3+t^6+t^9)v^3=0.$$
For $E_8$ we have $R=\C[t]/(t^{15}-1)$, and $\HH_{E_8}$ is generated over $R$ by $u$ and $v$
with the relations
$$u^3=1-t^5, \ \ (1+t^5+t^{10})u^2=0, \ \ v^5=1-t^3, \ \ (1+t^3+t^6+t^9+t^{12})v^4=0.$$

In both cases the specialization $t=1$ gives the Milnor ring of the corresponding singularity.

\subsection{Comparison with the Fan-Jarvis-Ruan theory}\label{comparison-sec}

Let $\w$ be a simple singularity, and let $G\sub G_\w$ be a subgroup
containing $J$. Here we will show that our reduced \CFT\ for the pair $(\w,G)$ is
isomorphic to the FJR-theory for the same data constructed in \cite{FJR}.
Recall that the state space 
$$\HH^{FJR}=\bigoplus_{\ga\in G}\HH^{FJR}_{\ga}$$ 
of the FJR theory coincides with our state space
$$\HH^{red}=\HH(\w,G,1)=\bigoplus_{\ga\in\Ga}\HH^{red}_\ga$$
(see \eqref{ga'-state-space-eq}). However, the obvious identification of the state spaces is not compatible
with the operations of \CFT. 

\begin{rem} 
The state space $\HH^{FJR}$ is defined in \cite{FJR}
in terms of some relative cohomology groups, whereas
our state space $\HH^{red}$ is defined using Milnor rings of singularities.
These spaces can be identified by the so-called
Wall's isomorphism (see \cite[Eq.\ (74)]{FJR}).
This identification however should not be confused with the mirror symmetry
isomorphisms  between two different types of Landau-Ginzburg models,
the FJR-theory of the pair $(\w,G)$ (the A-model) and the B-model
associated with the dual pair $(\hat{\w},\hat{G})$ (see \cite{FJR},
\cite{Krawitz} and \cite{PAKWR} for details). Even though the state
space of the B-model is also constructed from the Milnor rings of singularities,
its ring structure is completely different from the one induced on
$\HH^{red}$ by our \CFT.
\end{rem}

\begin{thm}\label{comparison-thm} 
The reduced \CFT\ associated with the pair $(\w,G)$ 
is isomorphic to the FJR-theory for the same pair.
\end{thm}

\Pf . In Sections \ref{A-sec}--\ref{E68-sec} we showed that for simple singularities
all the components $\HH^{red}_{\ga}$ are generated by the Chern characters of Koszul
\mf s of rank $1$. Therefore, by Corollary \ref{simple-homog-cor},
the Homogeneity Conjecture holds for ($\w,G)$. 
Together with the results of Sections \ref{concavity-sec}--\ref{sum-sing-sec} this implies that
our reduced \CFT\ has all the properties established in \cite[Sec.\ 4]{FJR} for the
FJR-theory. Therefore, the Reconstruction Theorem \cite[Thm.\ 6.2.10]{FJR}, proved 
for the FJR-theory, is valid for our theory as well. 

We claim that in order to construct an isomorphism of the theories
it is enough to find an isomorphism $\psi:\HH^{red}\rTo{\sim}\HH^{FJR}$
of Frobenius algebras respecting metrics such that for every $\ga\in G$ with $(\A^n)^\ga=0$,
the restriction of $\psi$ to $\HH^{red}_\ga$
is the identity map $\C=\HH^{red}_\ga\to\HH^{FJR}_\ga=\C$. 
Indeed, the latter condition guarantees that $\psi$ respects the \CFT\ maps in the concave case
(see Corollary \ref{concavity-cor}).
For all $(\w,G)$ except the $D_4$ singularity with the group $G=\lan J\ran$, the Reconstruction Theorem 
implies that both theories are determined by the Frobenius algebra structure on the state 
space along with a certain four-point correlator which can computed using the Concavity property.
In the remaining case $(D_4, \lan J\ran)$ a similar statement is true as follows from
\cite[Thm.\ 4.5]{FJMR}. (The analog of this theorem holds also for our reduced \CFT.)
Now let us construct the required isomorphism $\HH^{red}\simeq\HH^{FJR}$
for singularities of each type. 

\noindent
{\bf 1.  Types $A$, $E_6$ and $E_8$.}
In these cases for every $\ga\in G$ we have either
$(\A^n)^{\ga}=0$ or $\HH^{red}_\ga=\HH^{FJR}_\ga=0$. Thus, we should take $\psi$ to be 
the natural identification $\HH^{red}=\HH^{FJR}$. 
All genus zero correlators
in both theories are computed from the Concavity property  (see \cite{JKV}, \cite{FJR-An}), 
so they coincide under this identification of the state spaces.

\noindent
{\bf 2. Type $D_{d+1}$, $G=G_\w$.} In this case we still have $(\A^2)^{\ga}=0$ for all $\ga\neq 1$,
but the component corresponding to $\ga=1$ is nonzero.
As shown in Section \ref{D-max-sec},
the algebra $\HH^{red}$ is generated by the element $e(d-2)$ and we have 
\begin{align*}
&e(d-2)^i=e(d-1-i) &\text{ for } 1\le i\le d-1 \text{ and }\\
&e(d-2)^{d-1+i}=-2e(-i) &\text{ for } 1\le i\le d-1. 
\end{align*}
Note that these relations follow formally from the formulas for the metric, for the correlators
of Case 2a in Section \ref{D-max-sec} and for one of the correlators of Case 1b,
\begin{equation}\label{D-our-cor-eq}
\la(e(0),e(d-2),e(1))=-2.
\end{equation} 
Now we define the map $\psi:\HH^{red}\to\HH^{FJR}$ by
$$\psi(e(i))= \begin{cases} \bfe_{-i}, &i\not\equiv 0,\\ \eps\cdot 2y\bfe_0, &i\equiv 0\end{cases}$$
with $\eps=\pm 1$, where in the right-hand side we use the notation of 
\cite[Sec.\ 5.3.1]{FJR} with $n=d$.
Our calculations in Section \ref{D-max-sec} together with calculations of  \cite[Sec.\ 5.3.1]{FJR}
imply that this map is compatible with the metrics and with the correlators of Case 2a. Furthermore,
as shown in \cite[Sec.\ 5.3.1]{FJR}, 
\begin{equation}\label{D-FJR-cor-eq}
\lan y\bfe_0, \bfe_{d+2}, \bfe_{-1}\ran=\pm 1.
\end{equation}
Therefore, we can choose $\eps$ such that $\psi$ is compatible with the correlators
\eqref{D-our-cor-eq} and \eqref{D-FJR-cor-eq}.
Such $\psi$ sends powers of $e(d-2)$ to the corresponding powers of $\bfe_{d+2}$, and so it is
a ring isomorphism satisfying our requirements.

\noindent
{\bf 3. Type $D_{d+1}$, $G=\lan J\ran$, $d=2k+1$.} Assume first that $k>1$.
Then the algebra $\HH^{red}$ is generated by the elements $u=e(3)$ and 
$v=e_+(0)-e_-(0)$ subject to the relations $uv=0$, $v^2=du^{d-1}$ (see Section \ref{D-J-sec}).
Therefore, using the notation and the calculations of \cite[Sec.\ 5.2.4]{FJR}
we see that there is an algebra isomorphism
$$\psi:\HH^{red}\to\HH^{FJR}: u\mapsto\bfe_3, v\mapsto i\a x^k\bfe_0+i\b y\bfe_0.$$
Note that we have 
$$u^l=\begin{cases} e(2l+1) & 0\le l\le k-1,\\ -2e(2l-2k) & k+1\le l\le 2k, \end{cases}$$
and similar relations, expressing $\bfe_j$ in terms of $\bfe_3$, hold in $\HH^{FJR}$.
Therefore, $\psi(e(j))=\bfe_j$ for $j\not\equiv 0$.
Finally, the metric on $\HH^{red}$ is determined by $(u^{2k},1)=-2$, and we have
the similar relation for the metric on $\HH^{FJR}$. Hence, $\psi$ respects the metrics.

In the case $k=1$ the algebra $\HH^{red}$ is generated by the elements $u$ and $v$ such that
$uv=0$ and $u^2=v^2/3=-e(2)/2$. Calculations of \cite[Sec.\ 5.2.4]{FJR} show that $\HH^{FJR}$ is generated by the elements
$X=x\bfe_3$ and $Y=y\bfe_3$ with the relations 
$$XY=0, \ \ X^2=\frac{1}{6}\bfe_2, \ \text{ and } Y^2=-\frac{1}{2}\bfe_2.$$
Therefore, the map $\psi$ defined by
$$\psi(u)=i\sqrt{3}\cdot X, \psi(v)=\sqrt{3}\cdot Y$$
gives an algebra isomorphism sending
$e(2)=\bfe_2$. This $\psi$ also respects the metrics.

\noindent
{\bf 4. Type $E_7$.} In the notation of \cite[Sec.\ 5.2.2]{FJR}, let us define the algebra
isomorphism $\psi:\HH^{red}\to\HH^{FJR}$ by sending $u=e(5)$ to $\bfe_7$ and
$v=e(1)$ to $\bfe_5$. Using the relations of Section \ref{E7-sec} and of \cite[Sec.\ 5.2.2]{FJR} one immediately
checks that $\psi$ sends $e(2j)$ to $\bfe_j$ for $j\not\equiv 0 \mod 3$. 
Since the metrics are determined by the relations $(e(-2),1)=(\bfe_8,1)=1$,
the map $\psi$ also preserves the metrics.

\section{APPENDIX. Functoriality of Hochschild homology}

We will use the notation of Section \ref{dg-Hoch-sec}.
Let $\CC$ and $\DD$ be small dg-categories which are dg Morita equivalent to smooth
and proper dg-algebras, and let $F:\Per_{dg}(\CC)\to\Per_{dg}(\DD)$ be a dg-functor.
In this appendix we recall the construction of the map on Hochschild homology
$$F_*:HH_*(\CC)\to HH_*(\DD)$$ 
given  in \cite[Sec.\ 1.2]{PV-mf} and will show
that it agrees with the similar map constructed using the standard Hochschild complexes
(see \cite[Sec.\ 2.3]{Shk}). 

Recall that our construction in \cite[Sec.\ 1.2]{PV-mf} uses the fact that every dg-functor $F$ 
can be realized as the tensor
product functor with a perfect $\CC-\DD$-bimodule $X$: 
$$F(M)=M\ot_{\CC} X.$$
Let us consider the $\DD-\CC$-bimodule $X^T$ given by
$$X^T(D,C^\vee)=\Hom_{\DD^{op}-\mod}(X(C,?),h_D),$$
where $h_D$ is the representable right $\DD$-module associated with $D\in\DD$.
In \cite[Sec.\ 1.2]{PV-mf} we constructed canonical morphisms 
$$u:\De_{\CC}\to X\ot_\DD X^T \text{ and } c: X^T\ot_\CC X\to\De_{\DD}$$
in the derived categories of $\CC-\CC$ and $\DD-\DD$-bimodules.
The map $F_*$ is defined as the composition
$$\Tr_{\CC}(\De_{\CC})\rTo{\Tr_{\CC}(u)}\Tr_{\CC}(X\ot_{\DD} X^T)\simeq\Tr_{\DD}(X^T\ot_{\CC} X)
\rTo{\Tr_{\DD}(c)}\Tr_{\DD}(\De_{\DD}),$$
where the isomorphism in the middle is the canonical isomorphism constructed in 
\cite[Lem.\ 1.1.3]{PV-mf}.

Let us describe a modification of this construction, convenient  for our purposes.
Consider the functor
$$F^{(2)}:\Per(\CC^{op}\ot\CC)\to\Per(\DD^{op}\ot\DD)$$
defined using the tensor products with $X$ and $X^T$:
$$F^{(2)}(M)=X^T\ot^{\dL}_{\CC} M\ot^{\dL}_{\CC} X\simeq M\ot^{\dL}_{\CC^{op}\ot\CC}(X\ot X^T)
$$
for $M\in\Per(\CC^{op}\ot\CC)$.
Note that $F^{(2)}$
sends a representable $\CC-\CC$-bimodule $h_{C_1^\vee\ot C_2}$
to $h_{F(C_1)^\vee\ot F(C_2)}$. 

Consider the canonical morphism of functors from $\Per(\CC^{op}\ot\CC)$ to $\Per(k)$
$$t_F:\Tr_{\CC}\to\Tr_{\DD}\circ F^{(2)}$$
induced by the map $u$ and the isomorphism
$$X\ot_\DD X^T\simeq (X\ot X^T)\ot_{\DD\ot\DD^{op}}\De_{\DD}$$
(see \cite[eq.\ (1.8)]{PV-mf}).
Note that on a representable bimodule $h_{C_1^\vee\ot C_2}$
the morphism $t_F$ is given by the map
$$\Hom_{\CC}(C_1,C_2)\rTo{F}\Hom_{\DD}(F(C_1),F(C_2)).$$
Also, consider the canonical morphism in $\Per(\DD^{op}\ot\DD)$
$$c_F:F^{(2)}(\De_{\CC})\to\De_{\DD}.$$
given by the composition
$$F^{(2)}(\De_{\CC})\simeq\De_\CC\ot^{\dL}_{\CC^{op}\ot\CC}(X\ot X^T)\simeq X^T\ot_{\CC} X
\rTo{c}\De_{\DD}.$$

\begin{prop}\label{hoch-funct-prop}
The map $F_*$ is equal to the composition
$$\Tr_{\CC}(\De_{\CC})\rTo{t_F(\De_{\CC})} \Tr_{\DD}F^{(2)}(\De_{\CC})\rTo{\Tr_{\DD}(c_F)} 
\Tr_{\DD}(\De_{\DD}).$$
\end{prop}

\Pf . By definition, $F_*$ is the composition of the following four morphisms
\begin{align*}
&\Tr_{\CC}(\De_{\CC})\rTo{\Tr_{\CC}(u)} (X\ot_{\DD} X^T)\ot^{\dL}_{\CC^{op}\ot\CC}\De_{\CC} 
\rTo{\sim} (X\ot X^T)\ot^{\dL}_{\DD^{op}\ot\CC^{op}\ot\CC\ot\DD}(\De_{\DD}\ot\De_\CC)\\
&\rTo{\sim} (X^T\ot_{\CC} X)\ot^{\dL}_{\DD^{op}\ot\DD}\De_{\DD}
\rTo{\Tr_{\DD}(c)}\Tr_{\DD}(\De_{\DD}).
\end{align*}
It remains to notice that the composition of the first two arrows is $t_F(\De_{\CC})$ while the composition
of the last two arrows is $\Tr_\DD(c_F)$.
\ed

Now we will compare our map $F_*$ with the map on Hochschild homology constructed
using the standard complexes.
Recall that the diagonal bimodule $\De_{\CC}$ has the bar-resolution by representable 
$\CC-\CC$-bimodules (see \cite[Sec.\ 6.6]{Keller-derived}):
\begin{equation}
\ldots\to \Barc_1(\CC)(P^\vee,Q)\to\Barc_0(\CC)(P^\vee,Q)\to\De_{\CC}(P^\vee,Q)=\CC(P,Q),
\end{equation}
where we use the notation $\CC(?,?)=\Hom_{\CC}(?,?)$, and for $P,Q\in\CC$,
$$\Barc_n(\CC)(P^\vee,Q)=\bigoplus_{C_0,\ldots,C_n\in\CC}
\CC(C_n,Q)\ot\CC(C_{n-1},C_n)\ot\ldots\ot\CC(C_0,C_1)\ot\CC(P,C_0).$$
Computing $\Tr_{\CC}(\De_{\CC})$ with the help of this resolution
and taking into account the identification
$$\Barc_n(\CC)\ot_{\CC^{op}\ot\CC}\De_{\CC}=
\bigoplus_{C_0,\ldots,C_n\in\CC} \CC(C_n,C_0)\ot\CC(C_{n-1},C_n)\ot\ldots\ot\CC(C_0,C_1)$$
leads to the standard {\it Hochschild complex} ${\mathbf C}^H(\CC)$
(see e.g. \cite[Sec.\ 2.3]{Shk}).
Thus, we obtain a canonical isomorphism in $D(k)$
\begin{equation}\label{HC-Tr-isom}
{\mathbf C}^H(\CC)\simeq \Tr_{\CC}(\De_{\CC}).
\end{equation}

Note that if $\CC$ has a compact
generator $G$, then the restriction functor induces an
equivalence of the derived category of $\CC-\CC$-bimodules
with the derived category of bimodules over the dg-algebra $A=\CC(G,G)$ that sends $\De_{\CC}$
to the diagonal bimodule $A$. Hence we obtain an isomorphism in $D(k)$
between $\Tr_{\CC}(\De_{\CC})$ and 
the Hochschild homology of $A$. 
To realize this isomorphism on the chain level we can use
the subcomplex  $\Barc_\bullet(\CC,G)$ in the bar-resolution with
$$\Barc_n(\CC,G)(P^\vee,Q)=\CC(G,Q)\ot \CC(G,G)^{\ot n}\ot\CC(P,G).$$
The corresponding subcomplex $\Tr_{\CC}(\Barc_\bullet(\CC,G))$
in ${\mathbf C}^H(\CC)$ computes the
Hochschild homology of $A$.
Similar subcomplexes can be defined
in the situation when $\CC$ is generated by a finite set of compact objects $G_1,\ldots,G_m$.
The isomorphisms \eqref{HC-Tr-isom}
are compatible with the inclusions 
$\{G_1,\ldots,G_m\}_{\CC}\sub \CC$ 
of the full dg-subcategories with objects $G_1,\ldots,G_m$
(inducing equivalences of derived categories).

The standard complex ${\mathbf C}^H(?)$ is functorial with respect to dg-functors between dg-categories.
Let us show that the induced maps on Hochschild homology coincide with the maps
$F_*$ defined above.

\begin{thm}\label{hoch-funct-thm} 
Let $F:\CC\to\DD$ be a dg-functor. Then the map $F_*$ coincides with the
map on Hochschild homology induced by the chain map of Hochschild complexes 
$${\mathbf C}^H(F):{\mathbf C}^H(\CC)\to {\mathbf C}^H(\DD)$$
given by $F$.
\end{thm}

\Pf . Assume first that both $\CC$ and $\DD$ have finite number of objects.
Then we have the dg-algebras
$$A=\bigoplus_{C_1,C_2\in\CC}\CC(C_1,C_2) \text{ and }
B=\bigoplus_{D_1,D_2\in\DD}\DD(D_1,D_2),$$
so that the categories of modules over $\CC$ and $A$ (resp., over $\DD$ and $B$)
are equivalent. We
can view $F$ as a non-unital homomorphism of dg-algebras $f:A\to B$ and
extend it to a functor between the categories of perfect modules
$$F:\Per(\CC)=\Per(A)\to\Per(B)=\Per(\DD)$$ 
by $F(M)=M\ot_A B$ for $M\in\Per(A)$.
We are going to compute the map $F_*$ in this case using Proposition \ref{hoch-funct-prop}.
In our case the diagonal $\De_{\CC}$ (resp., $\De_{\DD}$) corresponds to the $A-A$-bimodule
$A$ (resp., $B-B$-bimodule $B$), and
$\Tr_{\CC}$ (resp., $\Tr_{\DD}$) is given by the functor $?\ot_{A^e}A$
(resp., $?\ot_{B^e}B$), where $A^e=A^{op}\ot A$. 
The functor $F^{(2)}:\Per(A^e)\to\Per(B^e)$ sends $M\in\Per(A^e)$ to
$M\ot_{A^e} B^e$. 
The natural transformation $t_F:\Tr_{\CC}\to\Tr_{\DD}\circ F^{(2)}$ is given by the morphisms
$$t_F(M):M\ot_{A^e} A\rTo{\id\ot f} M\ot_{A^e} B\simeq (M\ot_{A^e} B^e)\ot_{B^e}B.$$
induced by $f$. Finally, the map $c_F:F^{(2)}(A)\to B$ 
in $D(B^e)$ is the natural map $A\ot_{A^e} B^e\to B$ induced by $f$.
Let $\Barc_\bullet(A)$ (resp., $\Barc_\bullet(B)$)
be the bar-resolution of the bimodule $A$ (resp., $B$). Then $c_F$ 
is realized by the natural morphism of complexes
$$\Barc_\bullet(A) \ot_{A^e} B^e \to \Barc_\bullet(B),$$
given by
$$B\ot (A^{\ot n})\ot B \rTo{\id\ot (f^{\ot n})\ot\id} B\ot (B\ot\ldots\ot B)\ot B.$$
Hence, the morphism $\Tr_{\DD}(c_F)$ is realized by the map
$$\Barc_\bullet(A) \ot_{A^e} B\simeq (\Barc_\bullet(A)\ot_{A^e} B^e) \ot_{B^e} B
\to \Barc_\bullet(B) \ot_{B^e} B$$
given by
$$A^{\ot n}\ot B \rTo{f^{\ot n}\ot\id} B^{\ot n}\ot B.$$
Since the map $t_F(\Barc_\bullet(A))$ is given by
$$A^{\ot n}\ot A \rTo{\id\ot f} A^{\ot n}\ot B,
$$
we see that its composition with 
$\Tr_{\DD}(c_F)$ is equal to the natural map 
$${\mathbf C}^H(f):{\mathbf C}^H(A)\to {\mathbf C}^H(B)$$
of the Hochschild complexes induced by $f$.

Now let us consider the general case. Let $X$ (resp., $Y$)
be a compact generator of $\CC$ (resp., $\DD$). We have the following
commutative diagram of dg-functors
\begin{diagram}
\{X\}_{\CC} &\rTo{}&\CC  \\
\dTo{F'}&&\dTo{F}\\
\{F(X),Y\}_{\DD}&\rTo{}&\DD 
\end{diagram}
where horizontal arrows are inclusions inducing equivalences of derived categories, and $F'$
is the restriction of $F$.
By the first part of the proof, the map $(F')_*$ is represented by the chain map ${\mathbf C}^H(F')$
of the Hochschild complexes. The same is true for both 
horizontal arrows by the discussion preceding the formulation of the
theorem. Since the inclusion $\{X\}_{\CC}\to\CC$ induces an isomorphism
on Hochschild homology, our assertion follows.
\ed

\begin{cor}\label{Groth-cor} 
The map $F_*: HH_*(\CC)\to HH_*(\DD)$ depends only on the class of $F$ in the Grothendieck
group of $\Per(\CC^{op}-\DD)$.
\end{cor}

\Pf . This follows from \cite[Thm.\ 2.4]{Keller-inv}.
\ed

{\sc Department of Mathematics, University of Oregon, Eugene, OR 97405}

{\it Email addresses}: apolish@uoregon.edu, vaintrob@uoregon.edu

\end{document}